\documentclass[10pt]{amsart}
\usepackage[margin=0.95in]{geometry}
\usepackage[T1]{fontenc}
\usepackage[utf8]{inputenc}
\usepackage{lmodern}
\usepackage{microtype}
\usepackage{amsmath,amssymb,amsthm,mathtools,mathrsfs,bm}
\usepackage{enumitem}
\usepackage{booktabs}
\usepackage{array}
\renewcommand{\arraystretch}{1.03}
\usepackage{tikz}
\usetikzlibrary{arrows.meta,positioning}
\usepackage{color}
\usepackage{cite}
\usepackage{hyperref}
\hypersetup{
  colorlinks=true,
  linkcolor=blue,
  citecolor=blue,
  urlcolor=blue,
  pdftitle={Color--Phase Separation for Mixed Random Operators in Two-Speed Stochastic Klein--Gordon Systems},
  pdfauthor={Guangqian Zhao},
  pdfsubject={Stochastic PDEs; paracontrolled distributions; two-speed stochastic Klein--Gordon systems},
  pdfkeywords={stochastic PDEs, paracontrolled distributions, stochastic Klein--Gordon equations, mixed random operators, Gaussian tensors, weak covariance, Galerkin convergence}
}
\numberwithin{equation}{section}

\theoremstyle{plain}
\newtheorem{theorem}{Theorem}[section]
\newtheorem{proposition}[theorem]{Proposition}
\newtheorem{lemma}[theorem]{Lemma}
\newtheorem{corollary}[theorem]{Corollary}
\theoremstyle{definition}
\newtheorem{definition}[theorem]{Definition}

\theoremstyle{remark}
\newtheorem{remark}[theorem]{Remark}

\newcommand{\T}{\mathbb{T}}
\newcommand{\R}{\mathbb{R}}
\newcommand{\Z}{\mathbb{Z}}
\newcommand{\E}{\mathbb{E}}
\newcommand{\Prob}{\mathbb{P}}
\newcommand{\C}{\mathcal{C}}
\newcommand{\calZ}{\mathcal{Z}}
\newcommand{\calK}{\mathcal{K}}
\newcommand{\Btwo}{B_{2,\infty}}
\newcommand{\eps}{\varepsilon}
\newcommand{\la}{\langle}
\newcommand{\ra}{\rangle}
\newcommand{\ii}{\mathrm{i}}
\newcommand{\one}{\mathbf{1}}
\newcommand{\wh}{\widehat}
\DeclareMathOperator{\supp}{supp}
\DeclareMathOperator{\dist}{dist}
\newcommand{\frakm}{\mathfrak{m}}
\newcommand{\HH}{\mathfrak{H}}
\newcommand{\hh}{\mathfrak{h}}
\providecommand{\Sym}{\operatorname{Sym}}

\title[Color--Phase Separation for Mixed Random Operators]{Color--Phase Separation for Mixed Random Operators\\
in Two-Speed Stochastic Klein--Gordon Systems}
\author[G. Zhao]{Guangqian Zhao}
\address{School of Mathematical Sciences, University of Science and Technology of China, Hefei, Anhui 230026, China}
\email{zhaoguangqian@mail.ustc.edu.cn}

\begin{document}

\begin{abstract}
We study a two-component stochastic fractional Klein--Gordon system on the three-dimensional torus with fractional dispersion exponent \(\alpha\), distinct propagation speeds, and a pure cross nonlinearity.  The singular mixed terms combine a low--high paraproduct inside a Duhamel integral with an outer resonant product.  Their structure is governed by two index pairs: the color pair determines the possible Gaussian contraction, whereas the propagator--input pair determines the Duhamel phase.  In every same-color block generated by the interaction graph, the outer propagator has a different speed.  The contracted term is therefore a Fourier-diagonal Volterra multiplier, and integration by parts in time yields a factor of order \(N^{-\alpha}\) at frequency scale \(N\).  For the centered remainder we give a finite Hilbert-space kernel normal form that preserves the Fourier incidence relation and reduces the operator estimate to oriented Gaussian-tensor flattenings.  For independent noises we construct the enhanced stochastic data, prove local well-posedness in the prescribed enhanced paracontrolled/Galerkin class, and establish pathwise convergence of the Fourier--Galerkin equations for \(12/13<\alpha\le1\).  We also prove a logically separate perturbative extension to Fourier-diagonal color covariances with frequency-decaying off-diagonal entries.
\end{abstract}

\maketitle

\section{Introduction}

\paragraph{The equation and the mixed blocks.}

We study the two-component system
\begin{equation}\label{eq:main-system-intro}
        L_i u_i=u_1u_2+\xi_i,\qquad i=1,2,
\end{equation}
with independent real space-time white noises in the independent-color theorem and Fourier-diagonal weak covariance in the covariance-augmented theorem.  The common linear form is
\begin{equation}\label{eq:omega-intro}
        L_i=\partial_t^2+\omega_i(D)^2,
        \qquad
        \omega_i(n)=\bigl(1+c_i^2|n|^{2\alpha}\bigr)^{1/2},
        \qquad c_1\ne c_2.
\end{equation}
The nonlinearity is the pure cross interaction.  The mixed paracontrolled operators generated by this graph are
\[
        T^{i;j,k}(w)=I_i(w<\Psi_j)\circ\Psi_k,
        \qquad
        T^{1;2,2},\ T^{2;1,1},\ T^{1;1,2},\ T^{2;2,1}.
\]
For \(T^{i;j,k}w=I_i(w<\Psi_j)\circ\Psi_k\), the pair \((j,k)\) determines the possible Wick contraction, whereas \((i,j)\) determines the phase difference
\begin{equation}\label{eq:intro-color-phase-denominator}
        \Delta_{i|j}(q,\ell)=\omega_i(q+\ell)-\omega_j(\ell).
\end{equation}
The cross-interaction graph has the property that \(j=k\) implies \(i\ne j\).  Thus every same-color contraction passes through a different Klein--Gordon phase.  Its contraction condition gives $r=-\ell$ and hence $n=q$, so the contracted part is a Fourier-diagonal Volterra multiplier acting on the low input.  In the proportional low--high region,
\begin{equation}\label{eq:intro-speed-gap}
        |\omega_i(\ell+q)-\omega_j(\ell)|\gtrsim N^\alpha,
        \qquad i\ne j,
        \qquad |\ell|\sim N,
        \qquad |q|\le \delta_0N,
\end{equation}
one integration by parts in time gives the Volterra denominator $N^{-\alpha}$.  The high-shell count is then
\begin{equation}\label{eq:intro-diagonal-count}
        N^3\cdot N^{-\alpha}\cdot N^{-2\alpha}\cdot N^{-\alpha}=N^{3-4\alpha},
\end{equation}
corresponding respectively to the high loop, the Duhamel kernel, the covariance factor, and the phase denominator.  After the covariance branch is split as in Lemma~\ref{lem:exact-covariance-branches}, integration by parts in time produces the denominator; Lemma~\ref{lem:volterra-ET} gives the multiplier estimate and Appendix~\ref{app:phase-gap} gives the phase lower bound.  If \(i=j\), then \(j\ne k\), so color independence centers the block.  The generated family is
\begin{center}
\small
\begin{tabular}{@{}ccc@{}}
\toprule
block & color relation & independent-noise decomposition \\
\midrule
$T^{1;2,2}$ & $j=k$, $i\ne j$ & $\mathcal D^{1;2}+\mathcal B^{1;2,2}$ \\
$T^{2;1,1}$ & $j=k$, $i\ne j$ & $\mathcal D^{2;1}+\mathcal B^{2;1,1}$ \\
$T^{1;1,2}$ & $j\ne k$, $i=j$ & $\mathcal B^{1;1,2}$ \\
$T^{2;2,1}$ & $j\ne k$, $i=j$ & $\mathcal B^{2;2,1}$ \\
\bottomrule
\end{tabular}
\end{center}

At the first Picard level the singular source is
\[
        \Theta=\Psi_1\Psi_2.
\]
For independent colors this is already a cross-color Wick product.  For a Fourier-diagonal color covariance
\[
        \E[d\widehat W_a(n,t)d\widehat W_b(m,t)]
        =\mathbf 1_{n+m=0}\mathsf R_{ab}(n)\,dt,
\]
the finite product is split as
\begin{equation}\label{eq:intro-correlated-color-split}
        \widehat\Psi_a(\ell,s)\widehat\Psi_b(r,t)
        =\bigl(:\widehat\Psi_a(\ell,s)\widehat\Psi_b(r,t):\bigr)_{\mathsf R}
        +\mathbf 1_{\ell+r=0}\Sigma_{ab}^{\mathsf R}(\ell;s,t).
\end{equation}
If the off-diagonal covariance satisfies
\[
        |\mathsf R_{12}(n)|+|\mathsf R_{21}(n)|
        \lesssim \langle n\rangle^{-\kappa},
        \qquad \kappa>3-3\alpha+,
\]
then the extra lower-chaos branches are summable.  The covariance extension retains these deterministic branches and estimates the centered remainder by the same second-order random-operator argument as in the independent case.

The stochastic expansion belongs to the Da Prato--Debussche and paracontrolled line of singular SPDE analysis \cite{DPD,Hairer,Paracontrolled,GKO}; related stochastic wave and quantization problems are treated in \cite{GKOT,OhOkamotoTolomeo,OhOkamoto,OhRobertTzvetkov}.  For deterministic Klein--Gordon systems with quadratic or different-speed interactions, see \cite{Shatah,SimonTaflin,Germain,IonescuPausader,Deng}; periodic long-time problems are studied in \cite{DelortSzeftel,FangZhang}.

\paragraph{Contributions.}
The paper makes three model-level contributions.  First, it identifies and quantifies the color--phase mechanism: the pure-cross interaction graph forces every same-color contraction into a different-speed Volterra channel, and the resulting phase denominator yields the summable diagonal bound.  Second, it gives a class-indexed finite Hilbert-kernel normal form for the centered, time-correlated Fourier kernels and proves joint pathwise convergence in time and cutoff in the operator topology needed by the fixed point.  Third, it combines these estimates with the first Picard and cubic symbols to obtain a complete local solution theory, including convergence of the Galerkin fields and their projected nonlinear sources.  The independent-noise theorem is the principal result.  The Fourier-diagonal weak-covariance theorem is a separate perturbative extension with its additional lower-chaos branches displayed explicitly.

\paragraph{The centered operator estimate.}
After the covariance diagonals have been separated, a centered dyadic block has the form
\[
        B_{nq}=\sum_{\nu}\sum_{\ell_0,r_0}\sum_{\mu,\lambda}
        H_\nu(\ell_0,\mu,r_0,\lambda,q,n)
        :G_{j,[\ell_0],\mu}G_{k,[r_0],\lambda}:,
        \qquad
        n=q+\varepsilon_{\nu,1}\ell_0+\varepsilon_{\nu,2}r_0,
\]
where \(\ell_0,r_0\) are fixed representatives of the involution classes \([\ell_0]=\{\ell_0,-\ell_0\}\) and \([r_0]=\{r_0,-r_0\}\), while the signs \(\varepsilon_{\nu,1},\varepsilon_{\nu,2}\) belong to the deterministic pattern \(\nu\).  Thus the Gaussian indices are \(a=([\ell_0],\mu)\) and \(b=([r_0],\lambda)\); in particular, \(\ell_0\) and \(-\ell_0\) never define two independent Gaussian coordinates.  The bounded indices \(\mu,\lambda\) arise from the time-kernel normalization.  Proposition~\ref{prop:abstract-random-tensor} controls the relevant oriented flattenings and yields the operator bound
\[
        \ell_q^2\longrightarrow \ell_n^2.
\]
The normalization is performed separately in each real Fourier involution class.  It is obtained by factoring the finite Gram matrix of the time kernels, without inverting that matrix; a common Gram factorization is used when two time points or two cutoffs are compared.  Thus the auxiliary multiplicity is uniformly bounded, rank loss causes no singular coefficient, and no spatial label is mixed during the reduction.
Only after this operator norm is estimated do we insert the dyadic input bound
\begin{equation}\label{eq:intro-positive-envelope}
        \|P_Qw\|_{L_T^\infty L_x^2}
        \le Q^{-\sigma}\|w\|_{L_T^\infty B_{2,\infty}^\sigma},
\end{equation}
and the dyadic centered profile is
\begin{equation}\label{eq:intro-centered-profile}
        N^{\frac32-3\alpha+}\bigl(M^{\frac32+}+Q^{\frac32+}\bigr)Q^{-\sigma}.
\end{equation}
The two summands correspond to the output and input column counts.  They give the Sobolev source bound and the direct \(L_T^1B_{2,\infty}^{\sigma-\alpha}\) estimate required by the fixed point.

The paracontrolled random-operator formulation, including the centered/contracted decomposition and time integration by parts, follows \cite{GKO}.  The degree-two operator inequality is the corresponding specialization of the random-tensor estimates in \cite{DNY,Kaneshiro}.  What is specific to the present equation is the compatibility between the pure-cross color graph and the two speeds, together with its use throughout the enhanced-data and Galerkin construction.  We use the finite-chaos calculus of \cite{Janson,Nualart,PeccatiTaqqu}, Gaussian decoupling \cite{DeLaPenaGine}, and noncommutative Khintchine inequalities \cite{LustPisier}.

The first Picard object and the cubic source are constructed separately.  Appendix~\ref{app:first-picard-phase} proves
\[
        V_i=I_i(\Psi_1\Psi_2)\in C_T\mathcal C^{\frac72\alpha-3-},
        \qquad
        \partial_tV_i\in C_T\mathcal C^{\frac52\alpha-3-}.
\]
The same construction supplies the Besov bound used when $V_i$ enters the mixed-operator input space.  The H\"older estimate is the separate multiplier estimate used in ordinary deterministic products.  Appendix~\ref{app:cubic-resonance} constructs
\[
        \Gamma_i=V_{3-i}\circ\Psi_i=\Gamma_i^{(3)}+C_i,
\]
where $C_i$ is the first-chaos integrated Volterra-kernel component required by the source topology.

\paragraph{Linear estimates.}
The stochastic construction and the deterministic dispersive estimates are separated in the proof.  In the subunit branch $12/13<\alpha<1$, the deterministic fixed point uses a massive flat-torus version of the compact-manifold fractional-wave Sobolev--Strichartz estimate; Appendix~\ref{app:fracKG-strichartz} includes a self-contained proof, following the exponent range in \cite{Dinh}.  At $\alpha=1$, the homogeneous $(8,8/3)$ and $(4,4)$ estimates follow from the compact Klein--Gordon theorem of \cite{CacciafestaDanesiMeng} after a scalar time rescaling; the dual retarded line is then obtained by $TT^*$ and retarded truncation as in \cite{KeelTao,GKO}.  The color--phase mixed-operator estimates and the covariance split are common to both branches.

\paragraph{Galerkin limit.}
At each cutoff, the enhanced variables are obtained algebraically from the projected Galerkin equation.  The field and its projected nonlinear source are therefore fixed before the limiting argument.

The solution is defined as the limit of the Galerkin equations, and the nonlinear term is the limit of
\[
        N_{i,\Lambda}=\pi_\Lambda(u_{1,\Lambda}u_{2,\Lambda}),
        \qquad i=1,2.
\]
Convergence in the enhanced-data and fixed-point topologies identifies the resulting mild equation with the cutoff limit of these sources.  Compatible coordinate representations and cofinal cutoff sequences consequently give the same field and source on every common interval; the uniqueness assertion is confined to this Galerkin limit.

\begin{figure}[ht]
\centering
\begin{tikzpicture}[
  node distance=7mm and 8mm,
  box/.style={draw,rounded corners,align=center,font=\scriptsize,
              text width=0.25\textwidth,minimum height=8mm,inner sep=3pt},
  arr/.style={-{Latex[length=1.8mm]},line width=0.45pt}
]
\node[box] (diag) {phase gap and Volterra integration\\contracted diagonals};
\node[box,right=of diag] (center) {class-indexed Gram normal form\\centered mixed operators};
\node[box,right=of center] (symbols) {first Picard and cubic estimates\\stochastic symbols};
\node[box,below=of center] (enhanced) {joint enhanced-data convergence\\in time and cutoff};
\node[box,below left=of enhanced] (linear) {fractional/conic Strichartz\\deterministic product bounds};
\node[box,below right=of enhanced] (fixed) {paracontrolled fixed point\\localized stability};
\node[box,below=of fixed] (galerkin) {finite Galerkin bridge\\field and source convergence};
\draw[arr] (diag) -- (enhanced);
\draw[arr] (center) -- (enhanced);
\draw[arr] (symbols) -- (enhanced);
\draw[arr] (enhanced) -- (fixed);
\draw[arr] (linear) -- (fixed);
\draw[arr] (fixed) -- (galerkin);
\end{tikzpicture}
\caption{Dependency structure of the independent-color theorem.  The weak-covariance theorem replaces the stochastic-symbol layer by its covariance-augmented Wick expansion and then uses the same deterministic closure and Galerkin bridge.}
\label{fig:proof-dependency}
\end{figure}

\section{Model, notation, and main results}

\subsection{Equation and cutoff system}
Define the Duhamel operators
\begin{equation}\label{eq:duhamel}
I_i f(t)=\int_0^t \frac{\sin((t-s)\omega_i(D))}{\omega_i(D)}f(s)\,ds.
\end{equation}
The stochastic enhanced objects are constructed from the zero-data stochastic convolutions; deterministic Cauchy data enter through the affine regular component of the fixed point.

\begin{definition}[Admissible cofinal Fourier cutoff family]\label{def:admissible-cutoff-family}
A countable family \((\pi_\Lambda)_{\Lambda\in\mathbb N}\) is called an admissible cofinal Fourier cutoff family if
\[
        \widehat{\pi_\Lambda f}(n)
        =\chi_\Lambda(n)\widehat f(n),
        \qquad n\in\mathbb Z^3,
\]
and the symbols have the following properties:
\begin{enumerate}[label=\textup{(\roman*)},leftmargin=2.2em]
\item each \(\chi_\Lambda\) has finite support and
\(
 \sup_{\Lambda}\|\chi_\Lambda\|_{\ell^\infty(\mathbb Z^3)}
 \le C_{\rm cut}<\infty;
\)
\item \(\chi_\Lambda(-n)=\overline{\chi_\Lambda(n)}\), so that \(\pi_\Lambda\) preserves real-valued fields;
\item \(\chi_\Lambda(n)\to1\) as \(\Lambda\to\infty\) for every fixed \(n\in\mathbb Z^3\).
\end{enumerate}
The family is called standard localized if, in addition, \(|\chi_\Lambda|\le1\), \(\chi_\Lambda(n)=1\) for \(|n|\le\Lambda\), and \(\chi_\Lambda(n)=0\) for \(|n|\ge2\Lambda\).  Sharp spectral projectors and compactly supported smooth spectral cutoffs are the basic examples.  Numerical cutoff rates below are asserted only for the standard localized dyadic family; for a general admissible family we assert convergence without a rate.

An admissible Galerkin convention uses the same multiplier \(\pi_\Lambda\) on the noise and Cauchy data, on every stochastic leg used to construct the finite enhanced datum, and on the outer nonlinear source.  This compatibility, rather than idempotence of \(\pi_\Lambda\), is the property used in the finite-dimensional bridge.
\end{definition}

For the cutoff convergence statement we use the Galerkin cutoff equation
\begin{equation}\label{eq:cutoff-system}
\begin{aligned}
L_1u_{1,\Lambda}&=\pi_\Lambda(u_{1,\Lambda}u_{2,\Lambda})+\pi_\Lambda\xi_1,
& (u_{1,\Lambda},\partial_tu_{1,\Lambda})|_{t=0}&=(\pi_\Lambda y_{1,0},\pi_\Lambda y_{1,1}),\\
L_2u_{2,\Lambda}&=\pi_\Lambda(u_{1,\Lambda}u_{2,\Lambda})+\pi_\Lambda\xi_2,
& (u_{2,\Lambda},\partial_tu_{2,\Lambda})|_{t=0}&=(\pi_\Lambda y_{2,0},\pi_\Lambda y_{2,1}).
\end{aligned}
\end{equation}
The same cutoff acts on the noise, the nonlinear source, and the Cauchy data.  This convention makes the bridge in Proposition~\ref{prop:finite-cutoff-bridge} an exact finite-dimensional identity.  The affine semigroup term is placed in the regular component of the mild map, while the Volterra contractions and centered fluctuation operators are retained in the enhanced datum before the cutoff is removed.

\subsection{Paraproduct and resonant-product notation}\label{subsec:bony-notation}
We use one fixed smooth Littlewood--Paley decomposition throughout the paper, following the standard Bony paraproduct conventions for Besov spaces; see for example \cite{BCD}.  The low--high constant is chosen once and for all as follows.  Let $\delta_0=\delta_0(c_1,c_2)>0$ be small enough for the proportional phase-difference estimate of Lemma~\ref{lem:speed-gap}; then choose $0<c_0<\delta_0/16$.  With $\Delta_N$ denoting a dyadic block and $S_{<c_0N}$ a smooth low-frequency cutoff, the Bony products are written as
\begin{equation}\label{eq:bony-definitions}
 f<g:=\sum_N S_{<c_0N}f\,\Delta_Ng,
 \qquad
 f>g:=g<f,
 \qquad
 f\circ g:=fg-f<g-f>g.
\end{equation}
Thus $<$ is the low--high paraproduct, $>$ is the high--low paraproduct, and $\circ$ is the resonant high--high product.  The symbol $\circ$ is reserved exclusively for this Bony resonant product; ordinary multiplication is written by juxtaposition or by $\cdot$, and composition is written explicitly when needed.  When one factor is a stochastic distribution, such as in $V_{3-i}\circ\Psi_i$ or $I_i(w<\Psi_j)\circ\Psi_k$, the resonant product is interpreted through the finite-cutoff enhanced objects constructed below.  The Fourier multiplier $\chi^{\mathrm{res}}$ appearing later is the smooth cutoff implementing the support condition $N\sim M$ in the resonant product and is part of the fixed Bony convention.
With this definition, the three terms form an exact partition of the product.  The notation $N\sim M$ below refers to the complementary frequency band determined by this fixed partition; its comparison constants may depend on $c_0$ but never on the dyadic scales.

\subsection{Core notation and scale conventions}\label{subsec:notation-guide}
The notation below is fixed once and then used without further relabeling.  We also write $P_N=\Delta_N$ for a dyadic Littlewood--Paley block, and $C_TX$, $C_T^1X$, and $L_T^pX$ mean $C([0,T];X)$, $C^1([0,T];X)$, and $L^p([0,T];X)$.  The symbols $<$, $>$, and $\circ$ always mean the Bony operations in \eqref{eq:bony-definitions}; scalar comparisons and ordinary products are written in their usual contexts.  The bare symbol $\Omega$ in $L^p(\Omega;\cdot)$ denotes the probability space.  The two-frequency phase-difference variables $\Phi_{ij}^{\pm}$ are introduced in Section~\ref{sec:random-operators}.  As usual, $r-$ means $r-\kappa$ for an arbitrarily small fixed loss.  When a named exponent is used repeatedly, the loss is chosen once; in particular one may read
\begin{equation}\label{eq:rhoV-notation-guide}
 \rho_V:=\frac{7\alpha}{2}-3-\kappa_V,\qquad 0<\kappa_V\ll1.
\end{equation}

Fix once and for all a transversal $\mathbb Z^3_{\mathrm{rep}}\subset\mathbb Z^3$ for the involution $m\mapsto-m$, containing $0$.  For every $m\in\mathbb Z^3$, let $m^\circ\in\mathbb Z^3_{\mathrm{rep}}$ be the representative of $[m]=\{m,-m\}$ and write
\begin{equation}\label{eq:fourier-involution-representative}
        m=\varepsilon(m)m^\circ,
        \qquad \varepsilon(m)\in\{-1,1\},
        \qquad \varepsilon(0)=1.
\end{equation}
Gaussian coordinates in the real Fourier convention are indexed by the class $[m]$, whereas the sign $\varepsilon(m)$ remains in the deterministic Fourier coefficient.  Thus $G_{j,[m],\mu}$ and $G_{j,[-m],\mu}$ denote the same coordinate.  This convention is used in every Gram normal form and local covariance-coordinate formula below.

The mixed random operators generated by the pure cross equation are restricted to
\begin{equation}\label{eq:mixed-index-set}
        \mathfrak M
        :=\bigl\{(i;j,k)\in\{1,2\}^3:\ k=3-i,\ j\in\{1,2\}\bigr\}
        =\bigl\{(2;1,1),(2;2,1),(1;1,2),(1;2,2)\bigr\}.
\end{equation}
Thus, whenever \((i;j,k)\in\mathfrak M\) and \(j=k\), one automatically has \(i\ne j\).  The pure cross system excludes the direct same-color same-speed branch \(i=j=k\) from the generated mixed-operator estimates.
\begin{center}
\small
\begin{tabular}{@{}p{0.22\textwidth}p{0.70\textwidth}@{}}
\toprule
symbol & meaning \\
\midrule
$N,Q,M$ & stochastic high scale, input scale, and output scale in Section~\ref{sec:random-operators}; the paraproduct region is $Q\le c_0N$ \\
$\mathsf R(n)$ & Fourier-diagonal two-color covariance matrix in the covariance-augmented theorem \\
$\ell,q,n$ & high stochastic frequency, input frequency, and output frequency; after resonance the second high frequency is determined by these variables \\
$E_T$, $E_T^{2,\sigma}$ & rough input space $C_TL^2\cap C_T^1H^{-\alpha-\eps}$ and its Besov refinement, including $w\in L_T^\infty\Btwo^\sigma$ and $\partial_tw\in L_T^\infty\Btwo^{\sigma-\alpha}$ \\
$T^{i;j,k}$ & mixed operator $I_i(w<\Psi_j)\circ\Psi_k$; covariance contraction depends on the colors $j,k$, while the Duhamel phase-difference gap depends on the pair $(i,j)$ \\
$\mathcal D^{i;j}$, $\mathcal B^{i;j,k}$ & deterministic same-color Volterra diagonal and centered fluctuation; both are controlled in the Besov-input mixed-operator estimates \\
$\Phi_{ij}^{\pm}$ & two-frequency phase-difference variables used for the same-color Volterra diagonal \\
$\sigma$ & dyadic $L^2$ input exponent \\
\bottomrule
\end{tabular}
\end{center}
The phrase ``control norm'' refers to the enhanced-data norm $M_T$ in Definition~\ref{def:enhanced-data-norm}.  The high--high low-output centered sector is included in $\mathcal B^{i;j,k}$ and is controlled in the $E_T^{2,\sigma}$ topology.  The term ``full-input'' refers to the deterministic variable $W_i=V_i+X_i+Y_i$ in the paracontrolled split.

\subsection{Endpoint conventions}\label{subsec:theorem-conventions}
We record \(B_{2,\infty}^s\) endpoint bounds as dyadic \(L_T^\infty\) envelopes and impose time continuity below the endpoint:
\[
        f\in L_T^\infty B_{2,\infty}^{s}
        \quad\text{and}\quad
        f\in \bigcap_{\eta>0}C_TB_{2,\infty}^{s-\eta}.
\]
The first Picard object is used both as a Besov input and as a H\"older multiplier; accordingly, we retain both estimates
\[
        V_i\in L_T^\infty B_{2,\infty}^{\rho_V}
        \cap C_T\mathcal C^{\rho_V}.
\]
For \(12/13<\alpha<1\) the deterministic argument uses the fractional Klein--Gordon estimate of Appendix~\ref{app:fracKG-strichartz}; for \(\alpha=1\) it uses Proposition~\ref{prop:gko-conic-strichartz}.  All stochastic identities are first written at finite Galerkin cutoff and then passed to the limit in the displayed topologies.

\subsection{Main results}\label{sec:main-results}

\begin{theorem}[Independent-color main theorem]\label{thm:main}
\normalfont
Fix distinct speeds \(c_1\ne c_2\), a time horizon \(T_0\le1\), an admissible cofinal Fourier cutoff family \((\pi_\Lambda)_{\Lambda\in\mathbb N}\) as in Definition~\ref{def:admissible-cutoff-family}, and the pure-cross system \eqref{eq:main-system-intro}, with Galerkin cutoff convention \eqref{eq:cutoff-system}.  Assume \(12/13<\alpha\le1\).  For \(\alpha<1\), choose a strict admissible tuple as in Definition~\ref{def:admissible}; for \(\alpha=1\), use the endpoint window in \eqref{eq:conic-window}.  Let the two driving Gaussian noises be independent.
For deterministic data
\[
        y_i\in (H^{s_2}\cap B_{2,\infty}^{\sigma})
        \times (H^{s_2-\alpha}\cap B_{2,\infty}^{\sigma-\alpha}),
        \qquad i=1,2,
\]
there exists a full-probability event \(\Omega_{\mathrm{ind},T_0}\), chosen after the parameters and cutoff convention are fixed, with the following properties.

For every \(\omega\in\Omega_{\mathrm{ind},T_0}\) there is a random time \(T_\ast=T_\ast(\omega,y)\in(0,T_0]\) such that, for every \(0<T\le T_\ast\), the finite Galerkin solutions converge to a Galerkin limit.  More precisely:
\begin{enumerate}[label=\textup{(\roman*)},leftmargin=2.2em]
\item The finite stochastic enhanced data converge in the corresponding enhanced distance,
\[
        d_T(\Xi_\Lambda(\omega),\Xi(\omega))\longrightarrow0.
\]
The mixed source operators have cutoff limits as bounded maps
\[
        E_T^{2,\sigma}\longrightarrow
        C_TH^{s_2-\alpha}\cap L_T^1B_{2,\infty}^{\sigma-\alpha}.
\]
In the independent case the mixed components have the diagonal/centered form
\[
        T^{i;j,k}={\bf 1}_{j=k}\mathcal D^{i;j}+\mathcal B^{i;j,k}.
\]

\item The Galerkin fields converge in the reconstruction topology:
\[
        d_T^{\rm rec}
        \bigl((\Xi_\Lambda,Z_\Lambda),(\Xi,Z)\bigr)\longrightarrow0.
\]
The limiting field is reconstructed as
\[
        u_i=\Psi_i+V_i+X_i+Y_i,
        \qquad i=1,2,
\]
where \(Z=(X_1,X_2,Y_1,Y_2)\) is the deterministic fixed point driven by \((\Xi,y)\).

\item The projected nonlinear sources
\[
        N_{i,\Lambda}
        :=\pi_\Lambda(u_{1,\Lambda}u_{2,\Lambda}),
        \qquad i=1,2,
\]
converge to limits \(N_i^{\rm Gal}\) in the component source topology of the paracontrolled decomposition.  The limiting fields satisfy the mild equations with limiting source
\[
        u_i=\Psi_i+S_i(t)y_i+I_iN_i^{\rm Gal},
        \qquad i=1,2.
\]

\item The reconstructed field and nonlinear source are independent of the compatible enhanced-coordinate representation and of the cofinal cutoff subsequence.  Thus the prescribed Galerkin scheme has a unique limit on every common interval.  The solution map is locally Lipschitz in the deterministic data on localized-smallness classes.
\end{enumerate}
The constants in the Volterra estimates deteriorate as the speed gap closes, because the phase lower bound depends on \(|c_1-c_2|\).
\end{theorem}

\begin{proof}
The stochastic lift and mixed operators are given by Theorems~\ref{thm:stochastic-lift}, \ref{thm:first-picard-main}, and~\ref{thm:operator-main}, together with Proposition~\ref{prop:cubic-resonance}.  Proposition~\ref{prop:common-independent-event} places all cutoff limits on one full-probability event.  The deterministic fixed point is Theorem~\ref{thm:deterministic-closure}.  The Galerkin comparison is Proposition~\ref{prop:finite-cutoff-bridge} and Corollary~\ref{cor:cutoff-stability}.
\end{proof}

\begin{remark}[Scope of uniqueness]\label{rem:uniqueness-scope}
The theorem asserts uniqueness of the deterministic fixed point reconstructed from the prescribed enhanced datum and, equivalently, uniqueness of the corresponding Galerkin field and source limit on common existence intervals.  It does not assert unconditional pathwise uniqueness among arbitrary distributional solutions that are not known to induce this enhanced datum.
\end{remark}

\begin{theorem}[Stochastic lift]\label{thm:stochastic-lift}
Assume $\alpha>3/4$.  With
\[
        \Psi_i=I_i(\xi_i),
        \qquad
        \Theta=\Psi_1\Psi_2,
        \qquad
        V_i=I_i(\Theta),
\]
the cutoff approximations converge almost surely in
\begin{equation}\label{eq:stochastic-lift-paths}
        \Psi_i\in C_T\C^{\alpha-3/2-},
        \qquad
        \Theta\in C_T\C^{2\alpha-3-},
        \qquad
        V_i\in C_T\C^{3\alpha-3-},
        \qquad
        \partial_tV_i\in C_T\C^{2\alpha-3-}.
\end{equation}
The quadratic lift $\Theta$ is the cross-Wick product determined by the finite cutoff algebra.
\end{theorem}

\begin{theorem}[First Picard smoothing]\label{thm:first-picard-main}
The phase analysis in Appendix~\ref{app:first-picard-phase} yields the following first Picard bounds:
\begin{equation}\label{eq:first-picard-gain-intro}
        V_i\in C_T\C^{\frac{7\alpha}{2}-3-},
        \qquad
        \partial_tV_i\in C_T\C^{\frac{5\alpha}{2}-3-}.
\end{equation}
Consequently, for every admissible $\sigma<\rho_V:=7\alpha/2-3-$,
\begin{equation}\label{eq:first-picard-besov-intro}
        V_i\in L_T^\infty B_{2,\infty}^{\rho_V},
        \qquad
        \partial_tV_i\in L_T^\infty B_{2,\infty}^{\rho_V-\alpha}.
\end{equation}
The H\"older estimate is used for deterministic multiplication, and the Besov estimate supplies the dyadic Besov input bound.
\end{theorem}

\begin{proposition}[Cubic resonant symbols]\label{prop:cubic-resonance}
The cubic symbols
\[
        \Gamma_1=V_2\circ\Psi_1,
        \qquad
        \Gamma_2=V_1\circ\Psi_2
\]
have finite colored Wick decompositions
\begin{equation}\label{eq:cubic-main-decomp-intro}
        \Gamma_i=\Gamma_i^{(3)}+C_i.
\end{equation}
Appendix~\ref{app:cubic-resonance} proves
\begin{equation}\label{eq:cubic-components}
        \Gamma_i^{(3)}\in C_T\C^{\frac{9\alpha}{2}-\frac92-},
        \qquad
        C_i\in C_T\C^{5\alpha-\frac92-},
\end{equation}
and hence
\begin{equation}\label{eq:cubic-total-regularity}
        \Gamma_i\in C_T\C^{\frac{9\alpha}{2}-\frac92-}.
\end{equation}
The first-chaos term $C_i$ is defined as an integrated Volterra-kernel limit after the remaining stochastic convolution is inserted.
\end{proposition}

\begin{theorem}[Mixed operators]\label{thm:operator-main}
Let
\begin{equation}\label{eq:ET-def-intro}
        E_T^{2,\sigma}:=C_TL^2\cap C_T^1H^{-\alpha-\eps}
        \cap L_T^\infty B_{2,\infty}^{\sigma}
        \cap\{\partial_t w\in L_T^\infty B_{2,\infty}^{\sigma-\alpha}\}.
\end{equation}
The $B_{2,\infty}$ component is an $L_T^\infty$ shell envelope, with lower-index time continuity.  Assume
\begin{equation}\label{eq:operator-cond-main}
        s_2<4\alpha-3-10\eps,
        \qquad
        3-3\alpha+10\eps<\sigma<4\alpha-3-10\eps.
\end{equation}
Then each generated mixed operator $T^{i;j,k}$, $(i;j,k)\in\mathfrak M$, has a cutoff limit as a bounded operator
\begin{equation}\label{eq:operator-bound-main}
        E_T^{2,\sigma}\to C_TH^{s_2-\alpha}\cap L_T^1B_{2,\infty}^{\sigma-\alpha}.
\end{equation}
For $j=k$ the limit is the sum of a phase-difference Volterra diagonal and a centered fluctuation; for $j\ne k$ it is centered by color independence.  The centered proof estimates finite dyadic kernels in the operator norm $\ell_q^2\to\ell_n^2$ and preserves the incidence $n=q+\ell+r$.
\end{theorem}

\begin{proposition}[Joint cutoff convergence for the independent enhanced datum]\label{prop:common-independent-event}
Fix $T_0\le1$, the admissible parameters, and one admissible cofinal cutoff family in the sense of Definition~\ref{def:admissible-cutoff-family}.  All independent-color enhanced objects may be constructed on a single full-probability event $\Omega_{\mathrm{ind},T_0}$.  On this event their cutoff approximations converge simultaneously in the topologies listed below.  In the third column, the displayed quantity is the high-shell envelope after the target regularity weight has been inserted; every exponent is understood with the fixed strict loss.
\begin{center}
\scriptsize
\renewcommand{\arraystretch}{1.16}
\begin{tabular}{@{}>{\raggedright\arraybackslash}p{0.16\textwidth}>{\raggedright\arraybackslash}p{0.32\textwidth}>{\raggedright\arraybackslash}p{0.43\textwidth}@{}}
\toprule
object & cutoff-limit topology & summable shell envelope \\
\midrule
$\Psi_i$ & $C_{T_0}\mathcal C^{s_\Psi}$, $s_\Psi<\alpha-3/2$ & $N^{s_\Psi+3/2-\alpha+}$ \\
$\Theta$ & $C_{T_0}\mathcal C^{s_\Theta}$, $s_\Theta<2\alpha-3$ & $N^{s_\Theta+3-2\alpha+}$ \\
$V_i,\partial_tV_i$ & $C_{T_0}\mathcal C^{s_V}\times C_{T_0}\mathcal C^{s_{\dot V}}$, $s_V<7\alpha/2-3$, $s_{\dot V}<5\alpha/2-3$ & $N^{s_V+3-7\alpha/2+}$ and $N^{s_{\dot V}+3-5\alpha/2+}$ \\
$\Gamma_i^{(3)},C_i$ & $C_{T_0}\mathcal C^{s_3}\times C_{T_0}\mathcal C^{s_1}$, $s_3<9\alpha/2-9/2$, $s_1<5\alpha-9/2$ & $N^{s_3+9/2-9\alpha/2+}$ and $N^{s_1+9/2-5\alpha+}$ \\
$\mathcal D^{i;j}$ & $\mathcal L(E_{T_0}^{2,\sigma},C_{T_0}H^{s_2-\alpha}\cap L_{T_0}^1B_{2,\infty}^{\sigma-\alpha})$ & $N^{3-4\alpha+}$ and $N^{s_2+3-4\alpha+}$ \\
$\mathcal B^{i;j,k}$ & the same two operator targets & $N^{s_2+3-4\alpha+}$, $N^{3-3\alpha-\sigma+}$, $M^{\sigma+3-4\alpha+}$, and $M^{3-4\alpha+}$ \\
\bottomrule
\end{tabular}
\end{center}

More precisely, reserve the small losses so that every exponent in the last column is at most $-\delta$ for some $\delta>0$.  For every coefficient-valued block and every finite $p\ge2$, the fixed-time estimate and a time-increment estimate hold with the same shell envelope, the latter multiplied by $|h|^\eta N^{\alpha\eta}$ for a sufficiently small $\eta>0$.  For centered mixed kernels the same assertion holds in
\[
 C_{t,s}\mathcal L(\ell_q^2,\ell_n^2)
\]
after the two-parameter time lift.  A difference of two standard cutoffs is supported where a stochastic high frequency or the output frequency is comparable to the smaller cutoff.  Consequently, along dyadic cutoffs $\Lambda=2^m$, every row of the table has a tail bound $C(\omega)\Lambda^{-\theta}$ in its stated topology for some $\theta>0$ determined by the strict margins.  For a general admissible cofinal cutoff family the same estimates give the Cauchy property without asserting a numerical rate.

The event $\Omega_{\mathrm{ind},T_0}$ is independent of the deterministic input inserted later into the mixed operators.  It is obtained before the fixed point and works simultaneously for all dyadic blocks, all members of the chosen cutoff family, all rational time endpoints in $[0,T_0]$, and all strict target exponents by taking a countable intersection.
\end{proposition}

\begin{proof}
Fix first a rational strict choice of all target exponents.  Since the list of
objects and summation sectors in the table is finite, the losses may be
reserved so that, after decreasing $\delta>0$ if necessary, every weighted
shell envelope is bounded by $2^{-4\delta m}$ when the largest active dyadic
scale is $2^m$.  Hypercontractivity converts each coefficient-level
$L^2(\Omega)$ estimate into an $L^p(\Omega)$ estimate for every finite $p$,
without changing this dyadic power.  The coefficient criteria in
Lemma~\ref{lem:coeff-chaos-criterion}, Appendix~\ref{app:first-picard-phase},
and Appendix~\ref{app:cubic-resonance} give these estimates for the
distribution-valued rows and their cutoff differences.  The diagonal row is
Proposition~\ref{prop:det-contraction}.  For the centered row,
Proposition~\ref{prop:fixed-time-centered-kernel} estimates the matrix norm
before an input is inserted, Lemmas~\ref{lem:centered-increment-coefficients}
and~\ref{lem:operator-valued-time-lift} provide the two time increments, and
Lemma~\ref{lem:centered-deterministic-summation-table} gives the four
summable envelopes displayed in the table.

We spell out the passage to a common pathwise event.  At maximal scale $2^m$
there are at most $C(1+m)^3$ relevant dyadic triples.  Let $d_t=1$ for a
coefficient-valued path and $d_t=2$ for a centered operator kernel depending
on $(t,s)$.  Choose $A>\alpha$ and use the grid of mesh $2^{-Am}$ in each time
variable.  The number of grid points is at most $C2^{d_tAm}$.  At a grid
point, Chebyshev's inequality at threshold $2^{-2\delta m}$ gives a bound
$C_p2^{-2\delta pm}$ after the $2^{-4\delta m}$ moment envelope is inserted.
Choose $p$ so large that
\[
        2\delta p>d_tA+2.
\]
The sum over the grid points, the dyadic triples, and $m$ is then finite.

For neighboring grid points the increment estimates have the additional
factor
\[
        |h|^\eta 2^{\alpha\eta m}
        \le 2^{-\eta(A-\alpha)m}.
\]
Increasing $A$ first and then $p$ makes the corresponding probability sum
finite as well.  Borel--Cantelli and the standard dyadic chaining argument
therefore yield the asserted one-parameter continuous paths and the
two-parameter paths in
$C_{t,s}\mathcal L(\ell_q^2,\ell_n^2)$.  This argument also shows explicitly
that the same event controls all rational time endpoints.

For a difference of two standard cutoffs, the finite-cutoff identities place
the difference on a stochastic leg or on the output.  If the smaller cutoff
is $2^k$, the removed region belongs to the shell tail $m\ge k-O(1)$; summing
the preceding $2^{-4\delta m}$ envelopes gives $O(2^{-3\delta k})$.  Applying
the same grid and Borel--Cantelli argument gives an almost sure
$O(2^{-\theta k})$ bound for some $\theta>0$.  On any fixed finite set of
dyadic blocks, every admissible cutoff multiplier converges uniformly to one
because only finitely many lattice frequencies occur.  Finite-block convergence plus
the uniform high-shell bound proves the Cauchy property for every member of
the chosen countable cofinal family.

Finally, there are only finitely many object types and countably many cutoff
levels, rational time points, and rational strict exponent tuples.  Their
full-probability events may therefore be intersected.  An arbitrary strict
target exponent is dominated by a nearby rational stronger exponent, so the
same event gives all stated strict regularities.  Since the mixed-operator
norm is estimated before the deterministic input is inserted, this event is
common to every $w\in E_{T_0}^{2,\sigma}$.
\end{proof}

\begin{definition}[Fourier-diagonal weak covariance]\label{def:wc-covariance}
A Fourier-diagonal weak covariance is a real Gaussian forcing whose Fourier Brownian motions satisfy
\begin{equation}\label{eq:wc-covariance-def-main}
        \E[d\widehat W_a(n,t)d\widehat W_b(m,t)]
        =\mathbf 1_{n+m=0}\mathsf R_{ab}(n)\,dt,
        \qquad a,b\in\{1,2\}.
\end{equation}
For every real Fourier involution class \(\{n,-n\}\), the corresponding finite real sine--cosine covariance matrix is assumed to be positive semidefinite.  In complex notation this is encoded by a Hermitian two-color matrix \(\mathsf R(n)=(\mathsf R_{ab}(n))_{a,b=1}^2\), with \(\mathsf R_{aa}(n)=1\),
\[
        \mathsf R_{ba}(n)=\overline{\mathsf R_{ab}(n)},
        \qquad
        \mathsf R_{ab}(-n)=\overline{\mathsf R_{ab}(n)}.
\]
The first identity is covariance symmetry and the second is the Fourier reality condition.  They are equivalent to a real symmetric covariance matrix on the sine--cosine coordinates of the class; its explicit $4\times4$ form is given in \eqref{eq:real-fourier-covariance-matrix}.  Thus the complex off-diagonal entry \(\mathsf R_{12}(n)\) may be non-real and need not be even.  The finite Gaussian space is the quotient of this positive semidefinite local Hilbert space by its zero-norm directions.  The off-diagonal block is assumed to satisfy
\begin{equation}\label{eq:wc-kappa-main}
        |\mathsf R_{12}(n)|+|\mathsf R_{21}(n)|
        \le C_{\mathsf R}\langle n\rangle^{-\kappa}.
\end{equation}
The Fourier-diagonal support relation \(n+m=0\) preserves the spatial incidence in the centered tensor estimates.  The diagonal independent model is recovered when \(\mathsf R_{12}=\mathsf R_{21}=0\).

When two covariance fields are compared, we use a common Gaussian realization.  For each real Fourier involution class, fix a standard real Brownian vector $B_{[n]}$ and set
\begin{equation}\label{eq:canonical-covariance-coupling}
        dW^{\mathsf R}_{[n]}(t)
        =\mathsf C_{\mathsf R}(n)^{1/2}\,dB_{[n]}(t),
\end{equation}
where \(\mathsf C_{\mathsf R}(n)\) is the real sine--cosine covariance matrix and the positive semidefinite square root is used.  All covariance-difference estimates below refer to this coupling.
\end{definition}

\begin{definition}[Weak covariance gain and covariance norm]\label{def:wc-gain}
For a covariance satisfying Definition~\ref{def:wc-covariance}, set
\begin{equation}\label{eq:wc-covariance-norm}
        \|\mathsf R\|_{\kappa}:=1+
        \sup_{n\in\mathbb Z^3}\langle n\rangle^{\kappa}
        \bigl(|\mathsf R_{12}(n)|+|\mathsf R_{21}(n)|\bigr).
\end{equation}
When the theorem is used with the fixed loss \(\eps\), write
\begin{equation}\label{eq:wc-gain-margin}
        \delta_{\mathsf R}:=
        \kappa-(3-3\alpha)-10\eps.
\end{equation}
The weak-covariance regime is the open condition \(\delta_{\mathsf R}>0\).  This margin controls the zero-mode quadratic covariance branch; the constants in the covariance-augmented estimates may depend polynomially on \(\|\mathsf R\|_\kappa\) and may deteriorate as \(\delta_{\mathsf R}\downarrow0\).
\end{definition}

\begin{definition}[Covariance-augmented enhanced datum]\label{def:wc-enhanced-datum-main}
A covariance-augmented enhanced datum is the tuple
\begin{equation}\label{eq:wc-enhanced-datum-main-tuple}
\Xi^{\mathsf R}:=
\bigl(
\Psi_i^{\mathsf R},\Theta^{\mathsf R},C_\Theta^{\mathsf R},V_i^{\mathsf R},\partial_tV_i^{\mathsf R},
\Gamma_{i,\mathsf R}^{(3)},C_{i,\mathrm{diag}}^{\mathsf R},C_{i,\mathrm{off}}^{\mathsf R},R_{i,\Theta}^{\mathsf R},
\mathcal D_{\mathsf R}^{i;j,k},\mathcal B_{\mathsf R}^{i;j,k}
\bigr).
\end{equation}
Here \(i\in\{1,2\}\) and \((i;j,k)\in\mathfrak M\).  The superscript on \(\Psi_i^{\mathsf R}\) records the joint Gaussian realization; the marginal linear convolution has the same one-color law as in the independent case because \(\mathsf R_{ii}=1\).  The datum is admissible on \([0,T]\) if every component is the cutoff limit of its finite \(\mathsf R\)-Wick counterpart in the topology entering the covariance-augmented distance \(d_T^{\rm wc}\), and if the bounded-size and localized-smallness controls of Definition~\ref{def:wc-bounded-size} are finite.
\end{definition}

\begin{theorem}[Weak-covariance extension]\label{thm:weak-covariance-main}
Assume Definition~\ref{def:wc-covariance} and the weak-covariance gain condition
\begin{equation}\label{eq:wc-main-threshold}
        \delta_{\mathsf R}=\kappa-(3-3\alpha)-10\eps>0.
\end{equation}
Keep the fixed distinct speeds, the pure cross interaction, the admissible cutoff family of Definition~\ref{def:admissible-cutoff-family}, and the admissible parameters of Definition~\ref{def:admissible}.  Consider the finite Galerkin equation \eqref{eq:cutoff-system}, but driven by the two-color Gaussian forcing with covariance \eqref{eq:wc-covariance-def-main}.  Then the deterministic solution theorem and Galerkin convergence remain valid after replacing the independent enhanced datum by the covariance-augmented datum \(\Xi^{\mathsf R}\).

At finite cutoff every stochastic product is first expanded by the finite \(\mathsf R\)-Wick identity.  The additional lower-chaos components are: the deterministic zero-mode branch \(C_\Theta^{\mathsf R}\) in the quadratic symbol, the covariance-weighted Volterra diagonals \(\mathcal D_{\mathsf R}^{i;j,k}\) in the mixed operators, and the off-diagonal cubic first-chaos branches \(C_{i,\mathrm{off}}^{\mathsf R}\) together with the regular branch \(R_{i,\Theta}^{\mathsf R}\).  The centered remainders \(\mathcal B_{\mathsf R}^{i;j,k}\) are written in covariance coordinates chosen separately in each Fourier involution class.  This local Gram factorization preserves the incidence \(n=q+\ell+r\), so the same row/column tensor theorem applies.  Only covariance square roots are used.

The new branches and the estimates used to remove the cutoff are summarized here.  The table is an inventory of the finite \(\mathsf R\)-Wick decomposition, rather than an additional hypothesis.
\begin{center}
\scriptsize
\renewcommand{\arraystretch}{1.16}
\begin{tabular}{@{}>{\raggedright\arraybackslash}p{0.18\textwidth}>{\raggedright\arraybackslash}p{0.18\textwidth}>{\raggedright\arraybackslash}p{0.23\textwidth}>{\raggedright\arraybackslash}p{0.32\textwidth}@{}}
\toprule
branch & chaos and support & high-shell envelope & cutoff-limit topology \\
\midrule
$C_\Theta^{\mathsf R}$ & deterministic, $n=0$ & $N^{3-3\alpha-\kappa+}$ & $C_T C^r$ for every $r$, by Proposition~\ref{prop:wc-theta-branch} \\
$\mathcal D_{\mathsf R}^{i;j,k}$ & deterministic, $n=q$ & $N^{3-4\alpha+}$ if $j=k$; $N^{3-4\alpha-\kappa+}$ if $j\ne k$ & $\mathcal L(E_T^{2,\sigma},C_TH^{s_2-\alpha}\cap L_T^1B_{2,\infty}^{\sigma-\alpha})$ \\
$\mathcal B_{\mathsf R}^{i;j,k}$ & centered second chaos, $n=q+\ell+r$ & $N^{3/2-3\alpha+}(M^{3/2+}+Q^{3/2+})$ & the same two operator targets, by Proposition~\ref{prop:wc-centered-component} \\
$\Gamma_{i,\mathsf R}^{(3)}$ & centered third chaos & coefficient shell $M^{6-9\alpha+}$ & $C_T\mathcal C^{9\alpha/2-9/2-}$ and both source spaces \\
$C_{i,\mathrm{diag}}^{\mathsf R}$ & first chaos, same-color pairing & independent integrated-kernel profile & $C_T\mathcal C^{5\alpha-9/2-}$ and both source spaces \\
$C_{i,\mathrm{off}}^{\mathsf R}$ & first chaos, off-diagonal pairing & integrated kernel $\langle n\rangle^{-\alpha}N^{3-4\alpha-\kappa+}$ & $C_T\mathcal C^{5\alpha+\kappa-9/2-}$ and both source spaces \\
$R_{i,\Theta}^{\mathsf R}$ & zero-mode multiplier times first chaos & finitely many spatial blocks & every topology required of the cubic source \\
\bottomrule
\end{tabular}
\end{center}

The cutoff and difference estimates behind this inventory are explicit rather than formal identifications.  The quadratic rows are closed by Propositions~\ref{prop:wc-theta-branch} and~\ref{prop:wc-quadratic-lift}; the deterministic and centered mixed rows by Propositions~\ref{prop:wc-deterministic-diagonals} and~\ref{prop:wc-centered-component}; and all four cubic rows by Proposition~\ref{prop:wc-cubic-symbols}, whose proof uses the integrated off-diagonal kernel tail of Lemma~\ref{lem:wc-offdiag-cubic-kernel-tails}.  Proposition~\ref{prop:wc-enhanced-distance} then sums these branchwise cutoff differences in the single augmented metric, while Proposition~\ref{prop:wc-covariance-approximation} treats covariance perturbations and rank-deficient limits.  Thus every row entering the fixed point has both a construction bound and a cutoff-difference bound in its stated topology.

For every fixed \(T_0\le1\) there is a full-probability event \(\Omega^{\rm wc}_{*,T_0}\) on which
\begin{equation}\label{eq:wc-main-pathwise-dT}
        d_{T_0}^{\rm wc}(\Xi_\Lambda^{\mathsf R}(\omega),\Xi^{\mathsf R}(\omega))\to0.
\end{equation}
For every deterministic datum \(y\) in the Cauchy-data class of Theorem~\ref{thm:main}, there is a random time \(T_*(\omega,y,\mathsf R)>0\) such that the correlated Galerkin solutions converge in the reconstruction topology to
\begin{equation}\label{eq:wc-main-solution-formula}
        u_i^{\mathsf R}(\omega,y)=\Psi_i^{\mathsf R}(\omega)+V_i^{\mathsf R}(\omega)+X_i^{\mathsf R}(\omega,y)+Y_i^{\mathsf R}(\omega,y).
\end{equation}
Here \(X_i^{\mathsf R},Y_i^{\mathsf R}\) are the fixed-point variables associated with the covariance-augmented datum \(\Xi^{\mathsf R}\), and \(\Psi_i^{\mathsf R}\) denotes the corresponding linear stochastic convolution.
Compatible covariance-augmented coordinates and cofinal cutoff sequences produce the same Galerkin field and source limit on every common interval.  If \(\mathsf R_{12}=\mathsf R_{21}=0\), all off-diagonal covariance branches vanish and the theorem reduces to the independent-color theorem.
\end{theorem}

\begin{remark}[Logical status of the covariance extension]
The independent-color theorem and its full-probability event do not use Theorem~\ref{thm:weak-covariance-main}.  The latter is a separate perturbative result whose additional finite Wick branches, local covariance coordinates, and cutoff limits are proved in Section~\ref{subsec:weak-covariance-extension} and Appendix~\ref{app:weak-covariance-symbols}.
\end{remark}

\begin{proposition}[Galerkin cutoff bridge]\label{prop:finite-cutoff-bridge}
At finite cutoff, the Galerkin equation \eqref{eq:cutoff-system} with the same cutoff on the noise, deterministic data, and each nonlinear source is algebraically equivalent to the finite enhanced $X$--$Y$ system obtained from \eqref{eq:X-eq}--\eqref{eq:Y-eq}.  The reconstruction is
\[
        u_{i,\Lambda}=\Psi_{i,\Lambda}+V_{i,\Lambda}+X_{i,\Lambda}+Y_{i,\Lambda}.
\]
All Wick products and contractions in this statement are finite-dimensional identities.
\end{proposition}

\begin{proof}
Substitute the decomposition into the finite Galerkin mild equation.  The equations for $\Psi_{i,\Lambda}$ and $V_{i,\Lambda}$ remove the cutoff noise and the cutoff cross stochastic product.  The fixed Bony decomposition of the remaining finite Fourier polynomial gives the finite $X$--$Y$ system with the same outer cutoff.  Summing the component equations reconstructs the Galerkin PDE.
\end{proof}

\begin{lemma}[Deterministic Galerkin data convergence]\label{lem:deterministic-galerkin-data-convergence}
If $\sigma<s_2$, then $\pi_\Lambda f\to f$ in $H^{s_2}\cap B_{2,\infty}^{\sigma}$ for every $f\in H^{s_2}$.  Similarly, $\pi_\Lambda g\to g$ in $H^{s_2-\alpha}\cap B_{2,\infty}^{\sigma-\alpha}$ for every $g\in H^{s_2-\alpha}$.
\end{lemma}

\begin{proof}
The uniform symbol bound and pointwise convergence imply Sobolev convergence by dominated convergence on Fourier coefficients.  We give the corresponding finite-block plus tail argument for the Besov norm.  For every dyadic $K\ge1$, uniform multiplier boundedness gives
\[
 \sup_\Lambda\sup_{M>K}
 M^\sigma\|P_M(\pi_\Lambda-1)f\|_{L^2}
 \lesssim (1+C_{\rm cut})
 K^{\sigma-s_2}\|f\|_{H^{s_2}},
\]
which tends to zero as $K\to\infty$ because $\sigma<s_2$.  For the finitely many blocks $M\le K$, the lattice frequency set is finite, so \(\chi_\Lambda(n)\to1\) implies uniform convergence of the multiplier and hence convergence in every block.  Taking first \(\Lambda\to\infty\) and then \(K\to\infty\) proves the Besov convergence.  The velocity estimate is identical.
\end{proof}

\begin{proposition}[Cutoff convention equivalence]\label{prop:cutoff-convention-equivalence}
Let $\mathbb E_\Lambda^{\rm stoch}$ be the data built with stochastic-leg cutoffs and let $\mathbb E_\Lambda^{\rm Gal}$ be the Galerkin-compatible data with the outer cutoff in Proposition~\ref{prop:finite-cutoff-bridge}.  Then
\begin{equation}\label{eq:cutoff-convention-equivalence}
        d_T(\mathbb E_\Lambda^{\rm stoch},\mathbb E_\Lambda^{\rm Gal})\to0
        \quad\hbox{a.s.}
\end{equation}
The representative first Picard differences are
\begin{align}
        \pi_\Lambda I_i(\Psi_{1,\Lambda}\Psi_{2,\Lambda})-I_i(\Psi_{1,\Lambda}\Psi_{2,\Lambda})
        &=(\pi_\Lambda-1)I_i(\Psi_{1,\Lambda}\Psi_{2,\Lambda}),\label{eq:cutoff-diff-V-output}\\
        I_i(\Psi_{1,\Lambda}\Psi_{2,\Lambda})-I_i(\Psi_1\Psi_2)
        &=I_i((\Psi_{1,\Lambda}-\Psi_1)\Psi_{2,\Lambda})+I_i(\Psi_1(\Psi_{2,\Lambda}-\Psi_2)).\label{eq:cutoff-diff-V-leg}
\end{align}
For the cubic first-chaos contraction, \(C_{i,\Lambda}^{\rm Gal}-C_{i,\Lambda}^{\rm stoch}\) is the sum of the output and stochastic-leg tails of the integrated Volterra kernel.
The object-by-object source of the cutoff discrepancy is as follows; the convergence topologies and summable powers are those in Proposition~\ref{prop:common-independent-event}.
\begin{center}
\scriptsize
\renewcommand{\arraystretch}{1.14}
\begin{tabular}{@{}p{0.19\textwidth}p{0.31\textwidth}p{0.40\textwidth}@{}}
\toprule
object & possible removed tail & controlling estimate \\
\midrule
$\Psi_i,\Theta$ & stochastic legs and output & Lemma~\ref{lem:coeff-chaos-criterion} and the quadratic coefficient tails \\
$V_i,\partial_tV_i$ & two stochastic legs and output & Lemma~\ref{lem:B-tail-sums} and the first-Picard coefficient criterion \\
$\Gamma_i^{(3)}$ & three stochastic legs and output & Lemma~\ref{lem:C-third-tail} \\
$C_i$ & internal loop, remaining Gaussian leg, and output & Lemmas~\ref{lem:C-K21-conv} and~\ref{lem:C-first-tail} \\
$\mathcal D^{i;j}$ & contracted high loop and output diagonal & Proposition~\ref{prop:det-contraction} \\
$\mathcal B^{i;j,k}$ & either Gaussian high leg or the output & Lemma~\ref{lem:time-sup-centered-kernel} and Theorem~\ref{thm:pathwise-fluct} \\
\bottomrule
\end{tabular}
\end{center}
For the centered mixed operators, the identity is first read after finite dyadic truncation,
\begin{equation}\label{eq:cutoff-diff-B-operator}
        \sum_{N,Q,M\le 2^K}P_M
        (\mathfrak B^{\rm Gal}_{\Lambda,N,Q,M}-\mathfrak B^{\rm stoch}_{\Lambda,N,Q,M})P_Q,
\end{equation}
and then $K\to\infty$ in the operator topology \eqref{eq:operator-bound-main}.
\end{proposition}

\begin{proof}
For fixed finitely many dyadic blocks, all expressions are finite Fourier sums and the cutoff symbols converge entrywise.  The output tails for $\Theta,V,\Gamma^{(3)}$ use the coefficient tails of Appendices~\ref{app:first-picard-phase} and~\ref{app:cubic-resonance}; the first-chaos $C_i$ uses the integrated-kernel tail; $D^{i;j}$ uses the Volterra diagonal majorant; and $B^{i;j,k}$ uses the centered operator-tail majorant of Section~\ref{sec:random-operators}.  The finite-set/tail argument with these summable envelopes gives \eqref{eq:cutoff-convention-equivalence}.
\end{proof}

\subsection{Galerkin reconstruction}\label{subsec:galerkin-source-limit}
For enhanced data \(\Xi\) and fixed-point variables \(Z=(X_1,X_2,Y_1,Y_2)\), set
\[
        R_i(\Xi,Z)=\Psi_i+V_i+X_i+Y_i,
        \qquad
        d_T^{\rm rec}\bigl((\Xi,Z),(\widetilde\Xi,\widetilde Z)\bigr)
        =d_T(\Xi,\widetilde\Xi)+\|Z-\widetilde Z\|_{Z_T}.
\]
At finite cutoff define
\[
        N_{i,\Lambda}=\pi_\Lambda(u_{1,\Lambda}u_{2,\Lambda}).
\]
The Bony expansion of this source is a finite algebraic identity.  Hence convergence in the enhanced and fixed-point norms implies componentwise convergence of \(N_{i,\Lambda}\) in the source topology.

\begin{corollary}[Galerkin convergence]\label{cor:cutoff-stability}
Suppose the cutoff enhanced data converge in the augmented distance and form a common localized-smallness class on \([0,T]\).  Then the fixed points converge in \(\mathcal Z_T\), the Galerkin fields converge in reconstruction distance, and the projected nonlinear sources converge in their component source topologies.  Compatible enhanced coordinates and cofinal cutoff sequences yield the same field and source limits on common intervals.  The statement also holds in the Fourier-diagonal weak-covariance regime with \((\Xi,d_T)\) replaced by \((\Xi^{\mathsf R},d_T^{\rm wc})\).
\end{corollary}

\begin{proof}
Corollary~\ref{cor:deterministic-cutoff-passage} gives convergence of the fixed points.  Propositions~\ref{prop:finite-cutoff-bridge} and~\ref{prop:cutoff-convention-equivalence} identify their reconstruction with the Galerkin fields.  Each term in the finite Bony expansion of \(N_{i,\Lambda}\) converges in its assigned source topology.  In the weak-covariance regime the finite \(\mathsf R\)-Wick expansion gives the same algebraic identification, with convergence measured by \(d_T^{\rm wc}\).
\end{proof}

\section{Parameter constraints and critical exponents}\label{sec:parameters}

The parameter constraints have three independent layers.  First, the mixed-operator estimates give the color--phase constraints
\begin{equation}\label{eq:operator-cond-constraints}
        s_2<4\alpha-3-,\qquad 3-3\alpha+<\sigma<4\alpha-3-.
\end{equation}
The upper bound on $s_2$ is the output-count condition for both the deterministic Volterra diagonal and the centered operator.  The lower bound on $\sigma$ pays for the low-output column contribution of the centered row/column estimate.  The upper bound on $\sigma$ is the direct $B_{2,\infty}^{\sigma-\alpha}$ source summability condition.

Second, the full-input formulation uses the first Picard regularity
\begin{equation}\label{eq:Vi-regularity-constraints}
        \rho_V=\frac72\alpha-3- ,
        \qquad
        V_i\in C_T\mathcal C^{\rho_V}\cap L_T^\infty B_{2,\infty}^{\rho_V}.
\end{equation}
The H\"older entry is used in deterministic products involving $V_i$; the Besov entry supplies the input Besov estimate.  The admissible $\sigma$ is therefore chosen below $\rho_V$ and below the Sobolev indices used in the deterministic closure:
\begin{equation}\label{eq:combined-sigma-table}
        3-3\alpha+<\sigma<\min\{\rho_V,4\alpha-3,2\alpha-3/2,s_1,s_2\}-.
\end{equation}

Third, the cubic source uses
\begin{equation}\label{eq:Gamma-regularity-constraints}
        \Gamma_i\in C_T\mathcal C^{\frac92\alpha-\frac92-},
        \qquad
        s_2<\frac{11}{2}\alpha-\frac92-.
\end{equation}
Together with the deterministic lower bound $s_2>3/2-\alpha+$, this gives
\begin{equation}\label{eq:alpha-threshold}
        \frac32-\alpha < \frac{11}{2}\alpha-\frac92
        \quad\Longleftrightarrow\quad
        \alpha>\frac{12}{13}.
\end{equation}
Thus \(12/13\) is the threshold produced by the present fixed-point scheme; the equation-level optimal threshold remains open.  The same boundary is also reached by the Besov-input interval
\[
        3-3\alpha<\sigma<\rho_V=\frac72\alpha-3-
\]
and, in the near-$(4,3)$ deterministic choice, by the quadratic parameter window.  The mixed-operator sector alone closes under the weaker conditions \(s_2<4\alpha-3\) and \(3-3\alpha<\sigma<4\alpha-3\); the full threshold is imposed jointly by the cubic, first-Picard input, and deterministic quadratic constraints.

At the endpoint \(\alpha=1\), choose the conic window
\begin{equation}\label{eq:conic-window}
        \frac14+<s_1<\frac12-,\qquad
        \frac12+<s_2<\min\{1,s_1+\tfrac14\}-,
        \qquad
        0<\sigma<\min\{s_1,s_2,1/4\}-.
\end{equation}
This branch is separated from Appendix~\ref{app:fracKG-strichartz}; it uses the conic estimates recorded in Section~\ref{sec:deterministic-closure}.

\section{Paracontrolled decomposition}\label{sec:paracontrolled}
This algebraic section displays the random-operator structure.  On $[0,T]$, deterministic Cauchy data are propagated by the homogeneous semigroup inside the $Y$ component, so the estimates below apply to general deterministic data in the stated class.

\subsection{Exact full-input split}
Define
\begin{equation}\label{eq:Theta-def}
\Theta=\Psi_1\Psi_2,
\qquad
V_i=I_i(\Theta).
\end{equation}
The fixed point uses
\begin{equation}\label{eq:three-layer}
u_i=\Psi_i+W_i,
\qquad
W_i=V_i+U_i,
\qquad
U_i=X_i+Y_i.
\end{equation}
Equivalently,
\begin{equation}\label{eq:W-full-def}
W_i=V_i+X_i+Y_i.
\end{equation}
Since $L_i\Psi_i=\xi_i$ and $L_iV_i=\Theta$, the variable $W_i$ satisfies
\begin{equation}\label{eq:Wi-eq}
L_iW_i=\Theta+\Psi_1W_2+\Psi_2W_1+W_1W_2.
\end{equation}
Applying the Bony decomposition from Subsection~\ref{subsec:bony-notation} and subtracting $L_iV_i=\Theta$ gives the $X$--$Y$ system below. The low--high input is the full $W_i$, including the first Picard component and the two deterministic remainders.

The algebraic dependency is
\[
 W<\Psi \longrightarrow X,\qquad
 X\circ\Psi \longrightarrow T(W),\qquad
 V\circ\Psi \longrightarrow \Gamma,\qquad
 Y\circ\Psi \ \hbox{classical}.
\]
Thus \(V_i\) is part of the regular input in the full low--high source \(W_i<\Psi_j\).

\subsection{The \texorpdfstring{$X$}{X} and \texorpdfstring{$Y$}{Y} equations}
The low--high products between the full first remainders $W_i=V_i+X_i+Y_i$ and the noises are placed into $X_i$:
\begin{equation}\label{eq:X-eq}
L_iX_i=W_2<\Psi_1+W_1<\Psi_2.
\end{equation}
The remainder is then
\begin{align}\label{eq:Y-eq}
L_iY_i
&=W_2>\Psi_1+W_1>\Psi_2+W_1W_2
\notag\\
&\quad+W_2\circ\Psi_1+W_1\circ\Psi_2.
\end{align}
The resonant terms are decomposed as
\begin{align}\label{eq:res-split-W}
W_2\circ\Psi_1+W_1\circ\Psi_2
&=\Gamma_1+\Gamma_2+Y_2\circ\Psi_1+Y_1\circ\Psi_2
\notag\\
&\quad+X_2\circ\Psi_1+X_1\circ\Psi_2,
\end{align}
where
\[
\Gamma_1=V_2\circ\Psi_1,
\qquad
\Gamma_2=V_1\circ\Psi_2
\]
are stochastic symbols constructed in Appendix~\ref{app:cubic-resonance}. Substituting the Duhamel formula for $X_i$ into the last two terms gives
\begin{align}\label{eq:X-res-random-operators}
X_2\circ\Psi_1&=T^{2;1,1}(W_2)+T^{2;2,1}(W_1),
\notag\\
X_1\circ\Psi_2&=T^{1;1,2}(W_2)+T^{1;2,2}(W_1).
\end{align}
Here $<$, $>$, and $\circ$ are the Bony products fixed in Subsection~\ref{subsec:bony-notation}.  In particular, $W_i\circ\Psi_j$ denotes the Bony resonant product; Wick contractions enter through the finite enhanced objects.  The singular resonant pieces are interpreted by the enhanced symbols $\Gamma_i$ and the mixed random operators $T^{i;j,k}$.  Thus the mixed random operators act on the full input $W_i=V_i+X_i+Y_i$.  In the $[0,T]$ formulation, the affine homogeneous part appears in the regular component of the mild fixed-point map.

\section{Stochastic lift for independent colors}\label{sec:stoch-lift}

\begin{lemma}[Coefficient-level dyadic chaos criterion]\label{lem:coeff-chaos-criterion}
Let $Z$ be a spatially homogeneous random distribution on $\T^3$ whose dyadic block $P_NZ(t)$ belongs to a fixed Wiener chaos of order at most $m$.  Suppose that, uniformly in $t\in[0,T]$ and $|n|\sim N$,
\begin{equation}\label{eq:coeff-criterion-fixed}
  \E |\widehat Z(n,t)|^2\lesssim_T N^p,
\end{equation}
and that for some $0<\eta\le1$,
\begin{equation}\label{eq:coeff-criterion-inc}
  \E |\widehat Z(n,t)-\widehat Z(n,t')|^2
  \lesssim_T |t-t'|^{2\eta}N^{p+2\alpha\eta}.
\end{equation}
Then for every $\kappa>0$ and every $\gamma<\eta$ there is an almost surely finite random constant $C_\omega$ such that
\begin{equation}\label{eq:coeff-criterion-block}
  \|P_NZ\|_{C^\gamma([0,T];L^\infty_x)}
  \le C_\omega N^{\frac32+\frac p2+\alpha\eta+\kappa}.
\end{equation}
In particular, after choosing $\eta>0$ sufficiently small, $Z\in C_T\C^s$ for every
\begin{equation}\label{eq:coeff-criterion-regularity}
  s< -\frac32-\frac p2.
\end{equation}
The same statement applies to dyadic cutoff tails.
\end{lemma}

\begin{proof}
For fixed $x$, the variance of $P_NZ(t,x)$ is bounded by the sum of the coefficient variances on the shell, hence by $N^{3+p}$.  Nelson hypercontractivity on a fixed Wiener chaos upgrades this to all finite $L^r(\Omega)$ norms with the same power of $N$, up to constants depending on $r$ and the chaos order.  Applying the same estimate to a finite spatial net and to finitely many derivatives gives the dyadic $L^\infty_x$ bound with an arbitrary $N^\kappa$ loss.  The time estimate follows from \eqref{eq:coeff-criterion-inc} in the same way, with the usual fractional KG cost $N^{\alpha\eta}$, and Kolmogorov plus Borel--Cantelli over dyadic $N$ yields the almost sure estimate.  Summability of $N^s\|P_NZ\|_{C_TL^\infty}$ gives \eqref{eq:coeff-criterion-regularity}.  Tail estimates are identical because at least one Fourier variable is restricted to the removed dyadic tail.
\end{proof}

\subsection{Linear stochastic convolutions}
The Fourier coefficients of $\Psi_i$ are
\begin{equation}\label{eq:Psi-fourier}
\widehat\Psi_i(n,t)=\int_0^t \frac{\sin((t-s)\omega_i(n))}{\omega_i(n)}\,d\beta_n^i(s),
\end{equation}
where the Brownian families for $i=1,2$ are independent.  We use the real-valued Fourier convention
\[
        \beta_{-n}^i(t)=\overline{\beta_n^i(t)}.
\]
The non-conjugated covariance is
\[
        \E[d\beta_n^i(t)\,d\beta_m^j(t)]
        =\mathbf 1_{i=j}\mathbf 1_{n+m=0}\,dt,
\]
whereas the conjugated covariance is
\[
        \E[d\beta_n^i(t)\,d\overline{\beta_m^j(t)}]
        =\mathbf 1_{i=j}\mathbf 1_{n=m}\,dt.
\]
The former is used for Wick contractions inside products, and the latter is used in second-moment and operator-norm estimates.  By It\^o isometry,
\begin{equation}\label{eq:Psi-var}
\E|\widehat\Psi_i(n,t)|^2
\lesssim_T \omega_i(n)^{-2}
\lesssim \langle n\rangle^{-2\alpha}.
\end{equation}
A standard dyadic Bernstein--hypercontractive argument then yields the first statement in Theorem \ref{thm:stochastic-lift}.

\begin{proof}[Proof of Theorem \ref{thm:stochastic-lift} for the stochastic convolution]
For $|n|\sim N$, \eqref{eq:Psi-var} is the coefficient estimate \eqref{eq:coeff-criterion-fixed} with $p=-2\alpha$.  The standard kernel increment bound
\begin{equation*}
\E|\widehat\Psi_i(n,t+h)-\widehat\Psi_i(n,t)|^2\lesssim |h|^{2\eta}\langle n\rangle^{-2\alpha+2\eta\alpha}
\end{equation*}
is \eqref{eq:coeff-criterion-inc}.  Lemma~\ref{lem:coeff-chaos-criterion} therefore gives
\begin{equation*}
 \|P_N\Psi_i\|_{C_TL^\infty_x}\lesssim_\omega N^{3/2-\alpha+},
\end{equation*}
and the Besov summability condition is $s<\alpha-3/2$.  This coefficient-level Gaussian block estimate is the input for the Besov summation.
\end{proof}

\subsection{The cross quadratic symbol \texorpdfstring{$\Theta=\Psi_1\Psi_2$}{Theta = Psi1 Psi2}}
Because the two colors are independent, the product $\Psi_1\Psi_2$ is already centered:
\begin{equation}\label{eq:cross-centered}
\E[\Psi_1(x,t)\Psi_2(x,t)]=0.
\end{equation}
For its Fourier coefficient we have
\begin{equation}\label{eq:Theta-fourier}
\widehat\Theta(n,t)=\sum_{k+\ell=n}\widehat\Psi_1(k,t)\widehat\Psi_2(\ell,t),
\end{equation}
and therefore
\begin{equation}\label{eq:Theta-var}
\E|\widehat\Theta(n,t)|^2
\lesssim_T \sum_{k+\ell=n}\langle k\rangle^{-2\alpha}\langle \ell\rangle^{-2\alpha}
\lesssim \langle n\rangle^{3-4\alpha+}.
\end{equation}
For sufficiently small $\eta>0$, expand the time difference in the two factors and use independence together with the stochastic-convolution increment bound.  The same convolution estimate gives
\begin{equation}\label{eq:Theta-increment}
 \E|\widehat\Theta(n,t)-\widehat\Theta(n,t')|^2
 \lesssim_{T,\eta}
 |t-t'|^{2\eta}\langle n\rangle^{3-4\alpha+2\alpha\eta+}.
\end{equation}
Lemma~\ref{lem:coeff-chaos-criterion} therefore yields
$\Theta\in C_T\C^{2\alpha-3-}$ once $\alpha>3/4$.

The same estimates give convergence of the cutoff products in this path space.  Indeed, a cutoff difference places at least one of $k,\ell$ above the lower cutoff $L$.  For every sufficiently small $\vartheta>0$,
\[
 \mathbf 1_{\{|k|>L\}}+\mathbf 1_{\{|\ell|>L\}}
 \lesssim L^{-2\vartheta}
 \bigl(\langle k\rangle^{2\vartheta}+\langle\ell\rangle^{2\vartheta}\bigr).
\]
Thus the fixed-time and increment bounds for the cutoff difference gain $L^{-2\vartheta}$ after an arbitrarily small loss in the spatial exponent.  The strict inequality $s<2\alpha-3$ absorbs this loss.  Hypercontractivity and Borel--Cantelli over the countable Galerkin cutoffs then show that $\Theta_\Lambda$ is almost surely Cauchy in $C_T\C^s$ for every such $s$.

\subsection{The first Picard objects \texorpdfstring{$V_i=I_i(\Theta)$}{Vi = Ii(Theta)}}
Using one Duhamel factor, we obtain
\begin{equation}\label{eq:Vi-var}
\E|\widehat V_i(n,t)|^2
\lesssim_T \langle n\rangle^{-2\alpha}\langle n\rangle^{3-4\alpha+}
=\langle n\rangle^{3-6\alpha+},
\end{equation}
which implies
\begin{equation}\label{eq:Vi-reg}
V_i\in C_T\C^{3\alpha-3-},
\qquad
\partial_tV_i\in C_T\C^{2\alpha-3-}.
\end{equation}
This completes the proof of Theorem \ref{thm:stochastic-lift}.

\subsection{Cubic symbols \texorpdfstring{$\Gamma_i=V_{3-i}\circ\Psi_i$}{Gamma i = V 3-i o Psi i}}
We display $\Gamma_1=V_2\circ\Psi_1$; the other case is symmetric. In Fourier variables,
\begin{equation}\label{eq:Gamma-fourier}
\widehat\Gamma_1(n,t)=\sum_{m+r=n}\chi^{\mathrm{res}}(m,r)\widehat V_2(m,t)\widehat\Psi_1(r,t).
\end{equation}
After expanding $V_2=I_2(\Psi_1\Psi_2)$, this becomes
\begin{equation}\label{eq:Gamma-expanded-body}
 \widehat\Gamma_1(n,t)
 =\sum_{a+b+r=n}\chi^{\mathrm{res}}(a+b,r)
 \int_0^t K_2(t-s,a+b)
 \widehat\Psi_1(a,s)\widehat\Psi_2(b,s)\widehat\Psi_1(r,t)\,ds.
\end{equation}
The direct same-color contraction pairs the inner color-$1$ field $\widehat\Psi_1(a,s)$ with the outer color-$1$ field $\widehat\Psi_1(r,t)$; independence prevents pairings with the color-$2$ leg. Therefore
\begin{equation}\label{eq:Gamma-wick-body}
 \Gamma_1=\Gamma_1^{(3)}+C_1,
\end{equation}
where $\Gamma_1^{(3)}$ is the centered third homogeneous chaos and $C_1$ is a first-chaos contraction. The contraction condition is $a+r=0$, hence $b=n$, and
\begin{equation}\label{eq:Gamma-C1-body}
 \widehat C_1(n,t)
 =\int_0^t
 \left[\sum_{r\in\mathbb Z^3}\chi^{\mathrm{res}}(n-r,r)K_2(t-s,n-r)\sigma_1(r;s,t)\right]
 \widehat\Psi_2(n,s)\,ds.
\end{equation}
Thus the contracted component is a deterministic Volterra kernel acting on the remaining color-$2$ stochastic convolution, and it has mean zero.  The infinite contraction sum is defined in Appendix~\ref{app:cubic-resonance} after inserting this remaining stochastic convolution; this integrated-kernel formulation is the convergent object.  This is the simplest instance of the colored Wick algebra used later in Appendix~\ref{app:cubic-resonance}: only the two color-$1$ variables pair, reducing the chaos order by two and leaving the color-$2$ stochastic convolution as a first-chaos factor.

The direct product estimate using $V_2\in C_T\C^{3\alpha-3-}$ gives $\Gamma_i\in C_T\C^{4\alpha-9/2-}$.  Incorporating the first Picard smoothing estimate \eqref{eq:Vi-regularity-constraints}, Appendix~\ref{app:cubic-resonance} proves
\begin{equation}\label{eq:Gamma-regularity-body}
 \Gamma_i^{(3)}\in C_T\C^{\frac{9\alpha}{2}-\frac92-},
 \qquad
 C_i\in C_T\C^{5\alpha-\frac92-},
 \qquad
 \Gamma_i\in C_T\C^{\frac{9\alpha}{2}-\frac92-}.
\end{equation}
The proof has two separate components. The third-chaos estimate uses the first Picard coefficient bound
\[
 \sup_{t\le T}\E|\widehat V_j(m,t)|^2\lesssim \langle m\rangle^{3-7\alpha+}
\]
inside the resonant convolution with $\Psi_i$. The first-chaos contraction is treated by its exact kernel \eqref{eq:Gamma-C1-body}: low-output high--high interactions have a different-phase gap, while balanced-output interactions are controlled by a one-frequency phase-layer count for the factor $\min\{1,|\Omega|^{-1}\}$.

Consequently the resonant cubic source is compatible with the source norm once
\begin{equation}\label{eq:Gamma-source-body}
 s_2<\frac{11\alpha}{2}-\frac92-
\end{equation}

\subsection{The product \texorpdfstring{$V_1V_2$}{of the first Picard terms} in the admissible range}\label{subsec:first-picard-product}
The product $V_1V_2$ is estimated directly from the first Picard smoothing.  Appendix~\ref{app:first-picard-phase} proves
\[
  V_i\in C_T\C^{\rho_V},\qquad \rho_V=\frac{7\alpha}{2}-3-\kappa_V,
\]
\begin{equation}\label{eq:first-picard-product-regularity}
  V_1V_2\in C_T\C^{\rho_V}.
\end{equation}
Since the admissible parameters below always satisfy $s_2<\alpha$, the source exponent $s_2-\alpha$ is negative and
\begin{equation}\label{eq:V1V2-source-condition}
  C_T\C^{\rho_V}\hookrightarrow C_TH^{s_2-\alpha}.
\end{equation}
Hence $V_1V_2\in C_TH^{s_2-\alpha}$ throughout the admissible range.

\section{Mixed random operators}\label{sec:random-operators}

Here $i$ is the Duhamel index, $j,k$ are the stochastic colors, and $N,Q,M$ denote the high, input, and output scales, respectively.

\subsection{Generated mixed blocks}
From \eqref{eq:X-eq},
\begin{equation}\label{eq:X-formulas}
X_1=I_1(W_2<\Psi_1)+I_1(W_1<\Psi_2),
\qquad
X_2=I_2(W_2<\Psi_1)+I_2(W_1<\Psi_2).
\end{equation}
Substitution into the resonant part of \eqref{eq:Y-eq} gives
\begin{equation}\label{eq:generated-mixed-blocks}
X_2\circ\Psi_1=T^{2;1,1}(W_2)+T^{2;2,1}(W_1),
\qquad
X_1\circ\Psi_2=T^{1;1,2}(W_2)+T^{1;2,2}(W_1),
\end{equation}
where
\begin{equation}\label{eq:T-def}
T^{i;j,k}(w):=I_i(w<\Psi_j)\circ\Psi_k.
\end{equation}
The branch structure is as follows.  The blocks $T^{1;2,2}$ and $T^{2;1,1}$ have same-color stochastic legs and different Duhamel-source phases, so their contracted parts form Volterra diagonals.  The blocks $T^{1;1,2}$ and $T^{2;2,1}$ are cross-color blocks and are centered in the independent model.  These alternatives exhaust the covariance and phase channels generated by the cross interaction, equivalently the index set $\mathfrak M$ in \eqref{eq:mixed-index-set}.  In the Fourier-diagonal weak-covariance extension, the cross-color blocks acquire additional $\mathsf R$-weighted diagonal branches before the centered remainder is estimated.

\subsection{Finite-cutoff kernels and the deterministic/centered split}\label{subsec:finite-kernels}
The scale notation in this section is fixed throughout: $N$ is the high stochastic scale, $Q$ is the low input scale, and $M$ is the output scale.  The frequency variables are also fixed: $q$ denotes the input mode of $w$, $\ell$ is the high mode carried by the stochastic factor inside the paraproduct, and $n$ is the output mode after the exterior resonant product.  The second high stochastic mode is then determined by the convolution relation $r=n-q-\ell$.  The paraproduct cutoff in $w<\Psi_j$ imposes $Q\le c_0N$ for a fixed small constant $c_0$.

We use the following dyadic kernels.
\begin{center}
\small
\begin{tabular}{@{}p{0.18\textwidth}p{0.24\textwidth}p{0.45\textwidth}@{}}
\toprule
symbol & type & meaning \\
\midrule
$A_{N,Q,M}^{i;j,k}$ & raw dyadic kernel & finite-cutoff kernel before separating contractions \\
$D_{N,Q}^{i;j}$ & deterministic diagonal & same-color covariance branch; after contraction $n=q$ \\
$B_{N,Q,M}^{i;j,k}$ & centered kernel & raw kernel after subtracting $D$ in same-color blocks; raw $A$ in cross-color blocks \\
$\mathcal D^{i;j}$ & operator & Volterra Fourier multiplier acting diagonally on the low input \\
$\mathcal B^{i;j,k}$ & operator & centered second-chaos operator estimated in $\ell^2_q\to\ell^2_n$ norm \\
\bottomrule
\end{tabular}
\end{center}
With these conventions, set
\begin{equation}\label{eq:A-kernel}
A_{N,Q,M}^{i;j,k}(n,q;t,s)=\chi_M(n)\chi_Q(q)
\sum_{\ell+r=n-q}\rho_N(\ell)\rho_N(r)\chi^{\mathrm{res}}(q+\ell,r)
K_i(t-s,q+\ell)\widehat\Psi_j(\ell,s)\widehat\Psi_k(r,t),
\end{equation}
where $K_i(\tau,m)=\sin(\tau\omega_i(m))/\omega_i(m)$.  If $j=k$, the deterministic diagonal contraction is
\begin{equation}\label{eq:D-kernel}
D_{N,Q}^{i;j}(q;t,s)=\chi_Q(q)
\sum_{\ell}\rho_N(\ell)\rho_N(-\ell)\chi^{\mathrm{res}}(q+\ell,-\ell)
K_i(t-s,q+\ell)\sigma_j(\ell;s,t),
\end{equation}
where
\begin{equation}\label{eq:sigma-def}
\sigma_j(\ell;s,t):=\E\bigl[\widehat\Psi_j(\ell,s)\widehat\Psi_j(-\ell,t)\bigr].
\end{equation}
The deterministic operator is diagonal:
\begin{equation}\label{eq:D-operator-diagonal}
\widehat{\mathcal D_{N,Q}^{i;j}w}(q,t)=\int_0^tD_{N,Q}^{i;j}(q;t,s)\widehat w(q,s)\,ds.
\end{equation}
The centered kernel is
\begin{equation}\label{eq:B-kernel}
B_{N,Q,M}^{i;j,j}(n,q;t,s)=A_{N,Q,M}^{i;j,j}(n,q;t,s)-\one_{\{n=q\}}\chi_M(n)D_{N,Q}^{i;j}(q;t,s),
\end{equation}
and for $j\ne k$ we set $B_{N,Q,M}^{i;j,k}=A_{N,Q,M}^{i;j,k}$.  The associated centered operator is denoted by $\mathcal B_{N,Q,M}^{i;j,k}$.

\begin{remark}[Finite-cutoff convention for the dyadic kernels]\label{rem:finite-cutoff-dyadic-kernel-convention}
All identities in this subsection are at finite cutoff.  Stochastic-leg cutoff factors are suppressed until the cutoff-tail estimate.
\end{remark}

\begin{lemma}[Exact finite-cutoff Wick split for one dyadic block]\label{lem:finite-cutoff-wick-block}
Fix a finite Fourier cutoff $\Lambda$ and a dyadic triple $(N,Q,M)$.  In a same-color block, before any limiting operation is taken, the product of the two cutoff stochastic legs satisfies the finite Wick identity
\begin{equation}\label{eq:finite-block-wick-identity}
        c_\Lambda(\ell)c_\Lambda(r)\widehat\Psi_j(\ell,s)\widehat\Psi_j(r,t)
        =c_\Lambda(\ell)c_\Lambda(r):\widehat\Psi_j(\ell,s)\widehat\Psi_j(r,t):\;+
        c_\Lambda(\ell)c_\Lambda(r)\mathbf 1_{\ell+r=0}\,\sigma_j(\ell;s,t).
\end{equation}
in the real Fourier convention, up to the fixed finite sign copies described in Lemma~\ref{lem:same-color-hilbert-kernel-reduction}.  Consequently the finite kernels obey
\begin{equation}\label{eq:finite-block-A-B-D-identity}
        A_{\Lambda,N,Q,M}^{i;j,j}
        =B_{\Lambda,N,Q,M}^{i;j,j}
         +\mathbf 1_{\{n=q\}}\chi_M(n)D_{\Lambda,N,Q}^{i;j}
\end{equation}
as an identity of finite matrices in the variables $(n,q)$.  For $j\ne k$, the covariance vanishes by color independence and the finite raw block is already centered.
\end{lemma}

\begin{proof}
The product formula for first Gaussian integrals gives \eqref{eq:finite-block-wick-identity}.  The only non-conjugated same-color pairing is the involution $r=-\ell$, so the convolution relation $n=q+\ell+r$ reduces to $n=q$ on the contracted term.  Substituting this contracted contribution into the finite raw kernel gives precisely the finite diagonal kernel described in Remark~\ref{rem:finite-cutoff-dyadic-kernel-convention}.  The output cutoff then reads $\chi_M(q)$ because $n=q$; for fixed low input scale $Q$ only finitely many output scales $M\sim Q$ contribute, with a constant depending only on the Littlewood--Paley partition.  Thus the diagonal is counted as a scalar multiplier on the low mode.  In a cross-color block the two stochastic legs belong to orthogonal Gaussian color spaces, so the block is centered in the independent model.
\end{proof}

\begin{remark}[Finite-cutoff algebra for Fourier-diagonal covariance]\label{rem:weakly-correlated-colors}
The independent theorem corresponds to the diagonal covariance case, where a cross-color block with \(j\ne k\) is already centered.  Under Fourier-diagonal weak covariance the same finite raw kernel is split by the \(\mathsf R\)-Wick identity.  If
\[
  \E[d\widehat W_j(n,t)d\widehat W_k(m,t)]
  =\mathbf 1_{n+m=0}\mathsf R_{jk}(n)\,dt,
  \qquad j\ne k,
\]
then at finite cutoff
\[
  \widehat\Psi_j(\ell,s)\widehat\Psi_k(r,t)
  =\bigl(:\widehat\Psi_j(\ell,s)\widehat\Psi_k(r,t):\bigr)_{\mathsf R}
  +\mathbf 1_{\ell+r=0}\Gamma_{jk}^{\mathsf R}(\ell;s,t).
\]
The contracted support imposes \(n=q\), producing a Fourier-diagonal Volterra branch with the covariance weight \(\mathsf R_{jk}\).  If the off-diagonal block satisfies
\[
  |\mathsf R_{jk}(n)|\lesssim \langle n\rangle^{-\kappa},
\]
then this branch carries the additional dyadic gain used in Subsection~\ref{subsec:weak-covariance-extension} and Appendix~\ref{app:weak-covariance-symbols}.  The Wick-centered remainder has the same second-chaos tensor form as the independent centered block after a local covariance-coordinate factorization inside each real Fourier involution class.
\end{remark}

\subsection{Same-color phase-difference deterministic contractions}\label{subsec:same-color-det}

\begin{lemma}[Proportional low--high phase-difference gap]\label{lem:speed-gap}
There exist $\delta_0>0$ and $N_0\ge1$, depending only on the fixed annular cutoffs, $\alpha$, and the speed gap, such that if $i\ne j$, $|\ell|\sim N$, $|q|\le\delta_0N$, and $N\ge N_0$, then
\begin{equation}\label{eq:speed-gap}
 |\omega_i(q+\ell)-\omega_j(\ell)|\gtrsim N^\alpha,
 \qquad
 |\omega_i(q+\ell)+\omega_j(\ell)|\gtrsim N^\alpha.
\end{equation}
For one fixed speed $\omega$, the corresponding same-speed difference satisfies
\begin{equation}\label{eq:same-speed-comparison}
        |\omega(q+\ell)-\omega(\ell)|\lesssim |q|N^{\alpha-1}
\end{equation}
on the same proportional low--high window.  Finitely many lower shells are absorbed into the low-frequency constants.
\end{lemma}

\begin{proof}
The different-speed estimates are Proposition~\ref{prop:phase-gap-annulus}; take $\delta_0$ to be the proportional aperture supplied there.  The same-speed comparison is Proposition~\ref{prop:same-speed-contrast}.  The lower shells below $N_0$ are finite-dimensional.
\end{proof}

\begin{lemma}[Exact same-color covariance branch decomposition]\label{lem:exact-covariance-branches}
Let $0\le s\le t\le T$, $|\ell|\sim N$, and let
\[
        \sigma_j(\ell;s,t)=\E\bigl[\widehat\Psi_j(\ell,s)\widehat\Psi_j(-\ell,t)\bigr].
\]
Then
\begin{equation}\label{eq:exact-covariance-branches}
\sigma_j(\ell;s,t)=
\frac{s}{2\omega_j(\ell)^2}\cos((t-s)\omega_j(\ell))
+\frac{\sin((t-s)\omega_j(\ell))}{4\omega_j(\ell)^3}
-\frac{\sin((t+s)\omega_j(\ell))}{4\omega_j(\ell)^3}.
\end{equation}
We denote the three summands by $\sigma_j^{(0)}$, $\sigma_j^{(1)}$, and $\sigma_j^{(2)}$.  The main branch $\sigma_j^{(0)}$ is the only branch for which we use the Volterra phase-difference denominator.  The two endpoint branches satisfy the absolute bound
\begin{equation}\label{eq:endpoint-covariance-absolute}
        |\sigma_j^{(1)}(\ell;s,t)|+|\sigma_j^{(2)}(\ell;s,t)|\lesssim N^{-3\alpha}.
\end{equation}
For each sign branch in the product of $K_i(t-s,q+\ell)$ with $\sigma_j^{(0)}$, the integrand can be written as
\begin{equation}\label{eq:main-covariance-sign-branch}
        e^{\ii(t-s)\Phi_{ij}^{\pm}(q,\ell)}a_{q,\ell}(s)\widehat w(q,s),
        \qquad
        \Phi_{ij}^{\pm}(q,\ell)=\omega_i(q+\ell)\pm\omega_j(\ell),
\end{equation}
where
\begin{equation}\label{eq:main-covariance-amplitude-bounds}
        |a_{q,\ell}(s)|\lesssim_T N^{-3\alpha},
        \qquad
        |\partial_sa_{q,\ell}(s)|\lesssim N^{-3\alpha}.
\end{equation}
The constants may depend on $T$ and the fixed cutoff profile and are uniform in $q,\ell,N$.
\end{lemma}

\begin{proof}
By It\^o isometry,
\[
\sigma_j(\ell;s,t)=\int_0^s
\frac{\sin((s-r)\omega_j(\ell))\sin((t-r)\omega_j(\ell))}{\omega_j(\ell)^2}\,dr.
\]
The identity \eqref{eq:exact-covariance-branches} follows from
$2\sin A\sin B=\cos(A-B)-\cos(A+B)$ and an elementary integration in $r$.
The endpoint estimate is immediate from $\omega_j(\ell)\sim N^\alpha$.  For the main branch, expand both the Duhamel sine and the covariance cosine into exponentials.  The factor $K_i$ contributes $\omega_i(q+\ell)^{-1}\lesssim N^{-\alpha}$ and $\sigma_j^{(0)}$ contributes $s\omega_j(\ell)^{-2}\lesssim_T N^{-2\alpha}$.  If the oscillation is placed in $e^{\ii(t-s)\Phi_{ij}^{\pm}}$, the derivative in $s$ falls only on the smooth scalar factor $s$ and on fixed cutoffs, giving the same $N^{-3\alpha}$ amplitude bound.
\end{proof}

\begin{lemma}[Volterra multiplier on the proportional low--high window]\label{lem:volterra-ET}
Assume $j=k$ and $i\ne j$.  For every small $\eps>0$,
\begin{align}\label{eq:volterra-ET}
\|\mathcal D_{N,Q}^{i;j}P_Qw\|_{C_TH^{s_2-\alpha}}
&\lesssim
N^{3-4\alpha+\eps}\|P_Qw\|_{L_T^\infty L^2}\notag\\
&\quad+N^{3-4\alpha+\eps}Q^{s_2+\eps}\|\partial_tw\|_{C_TH^{-\alpha-\eps}}.
\end{align}
whenever $Q\le\delta_0N$.
\end{lemma}

\begin{proof}
We use the exact covariance branch decomposition of Lemma~\ref{lem:exact-covariance-branches}.  The endpoint branches $\sigma_j^{(1)}$ and $\sigma_j^{(2)}$ are estimated absolutely.  They contribute
\[
        \underbrace{N^3}_{\text{high loop}}
        \underbrace{N^{-\alpha}}_{K_i}
        \underbrace{N^{-3\alpha}}_{\sigma_j^{(1)},\sigma_j^{(2)}}
        =N^{3-4\alpha},
\]
which is already at the desired shell size.

For the main branch $\sigma_j^{(0)}$, expand the Duhamel sine and the covariance cosine into sign branches as in \eqref{eq:main-covariance-sign-branch}.  Lemma~\ref{lem:speed-gap} gives $|\Phi_{ij}^{\pm}(q,\ell)|^{-1}\lesssim N^{-\alpha}$ on $Q\le\delta_0N$.  Hence
\begin{align*}
\int_0^t e^{i(t-s)\Phi}a(s)\widehat w(q,s)\,ds
&=\frac{a(t)\widehat w(q,t)-e^{it\Phi}a(0)\widehat w(q,0)}{i\Phi}\\
&\quad -\int_0^t\frac{e^{i(t-s)\Phi}}{i\Phi}
       \bigl(a'(s)\widehat w(q,s)+a(s)\partial_s\widehat w(q,s)\bigr)\,ds.
\end{align*}
We record the three contributions separately.  The boundary terms have the high-loop count $N^3$, the Duhamel factor $N^{-\alpha}$, the covariance factor $N^{-2\alpha}$, and the phase-difference denominator $N^{-\alpha}$.  Thus their multiplier size is $N^{3-4\alpha}$.  The amplitude-derivative term has the same shell size, because \eqref{eq:main-covariance-amplitude-bounds} gives $|a'(s)|\lesssim N^{-3\alpha}$ after the Duhamel factor has already been included in the amplitude.  The endpoint covariance branches were estimated absolutely above and have the same shell size.

For the input-derivative term, the only additional quantity is $\partial_tw$ at the same low frequency $q$.  The contraction has already imposed the diagonal relation $n=q$, so the output Sobolev weight is the low-shell weight.  Therefore
\[
\|P_Q\partial_tw\|_{H^{s_2-\alpha}}
\lesssim Q^{s_2+\eps}\|\partial_tw\|_{H^{-\alpha-\eps}}.
\]
The estimate is scalar on the low shell: the high loop has been summed inside the diagonal multiplier, which then acts pointwise on the $q$-coefficient.  Combining the sign branches, endpoint pieces and the smooth cutoffs proves \eqref{eq:volterra-ET}.
\end{proof}

\begin{lemma}[Besov source bound for the deterministic diagonal]\label{lem:volterra-besov-source}
Assume $j=k$, $i\ne j$, and $Q\le\delta_0N$.  Then
\begin{equation}\label{eq:volterra-besov-source}
\|\mathcal D_{N,Q}^{i;j}P_Qw\|_{L_T^\infty\Btwo^{\sigma-\alpha}}
\lesssim_{T,\eps}
N^{3-4\alpha+\eps}
\Bigl(Q^{-\alpha}\|w\|_{L_T^\infty\Btwo^\sigma}
      +\|\partial_tw\|_{L_T^\infty\Btwo^{\sigma-\alpha}}\Bigr).
\end{equation}
\end{lemma}

\begin{proof}
The operator is Fourier diagonal in $q$, so the $B_{2,\infty}^{\sigma-\alpha}$ source norm is taken after fixing the output block $Q$.  This order matters.  We first estimate the scalar multiplier obtained by summing the high loop at scale $N$, and only then take the dyadic supremum in the low variable.  Boundary terms and amplitude-derivative terms from the Volterra integration by parts satisfy
\[
Q^{\sigma-\alpha}\|P_Qw\|_{L^2}
\le Q^{-\alpha}\|w\|_{L_T^\infty\Btwo^\sigma}.
\]
When the derivative falls on the input, the same Besov weight gives
\[
Q^{\sigma-\alpha}\|P_Q\partial_tw\|_{L^2}
\le \|\partial_tw\|_{L_T^\infty\Btwo^{\sigma-\alpha}}.
\]
The high-frequency loop, Duhamel factor, covariance and phase-difference denominator contribute $N^{3-4\alpha+\eps}$ as in Lemma~\ref{lem:volterra-ET}.  Thus, for fixed output block $Q$, the boundary and amplitude-derivative contribution is bounded by
\[
        Q^{-\alpha}\sum_{N\gtrsim Q}N^{3-4\alpha+\eps}
        \lesssim Q^{3-5\alpha+\eps},
\]
while the input-derivative contribution is bounded by
\[
        \sum_{N\gtrsim Q}N^{3-4\alpha+\eps}
        \lesssim Q^{3-4\alpha+\eps}.
\]
Both exponents are non-positive in the range used in the paper, after the small loss is fixed; the finitely many low $Q$ blocks are absorbed into the constant.  Thus the full diagonal source is bounded in $L_T^\infty B_{2,\infty}^{\sigma-\alpha}$, and the high-shell tails tend to zero by the same displayed summable envelopes.  After the diagonal relation $n=q$ is imposed, the Besov norm sees only this low shell, and the input is measured by the two dyadic Besov bounds $\|w\|_{L_T^\infty B_{2,\infty}^\sigma}$ and $\|\partial_tw\|_{L_T^\infty B_{2,\infty}^{\sigma-\alpha}}$.
\end{proof}

\begin{proposition}[Same-color deterministic contraction]\label{prop:det-contraction}
Assume $j=k$, $i\ne j$, and $s_2<4\alpha-3-5\eps$.  Then the dyadic deterministic contractions converge as bounded operators
\begin{equation}\label{eq:det-contraction-bound}
\mathcal D^{i;j}:E_T^{2,\sigma}\longrightarrow C_TH^{s_2-\alpha}\cap L_T^\infty\Btwo^{\sigma-\alpha}.
\end{equation}
\end{proposition}

\begin{proof}
For the Sobolev bound, Lemma~\ref{lem:volterra-ET} gives two summation lines.  The first line comes from boundary terms and amplitude derivatives in the Volterra integration by parts; the second comes from the derivative falling on the low input.  The dyadic powers are:
\begin{center}
\footnotesize
\renewcommand{\arraystretch}{1.14}
\begin{tabular}{>{\raggedright\arraybackslash}p{0.23\textwidth}>{\raggedright\arraybackslash}p{0.34\textwidth}>{\raggedright\arraybackslash}p{0.31\textwidth}}
\toprule
contribution & dyadic size before summing in $Q$ & condition after $Q\le\delta_0N$ \\
\midrule
boundary and amplitude derivative &
$N^{3-4\alpha+\eps}\|P_Qw\|_{L_T^\infty L^2}$ &
$\sum_{Q\le\delta_0N}\|P_Qw\|_{L_T^\infty L^2}\lesssim_\eps N^\eps\|w\|_{C_TL^2}$, hence $N^{3-4\alpha+2\eps}$, summable since $\alpha>3/4$ \\
input derivative &
$N^{3-4\alpha+\eps}Q^{s_2+\eps}\|\partial_tw\|_{C_TH^{-\alpha-\eps}}$ &
$\sum_{Q\le\delta_0N}Q^{s_2+\eps}\lesssim N^{s_2+\eps}$, hence $N^{3-4\alpha+s_2+2\eps}$, summable under $s_2<4\alpha-3-5\eps$ \\
\bottomrule
\end{tabular}
\end{center}
Equivalently, the two displayed estimates are
\begin{equation}
        \sum_{Q\le \delta_0N}\|P_Qw\|_{L_T^\infty L^2}
        \lesssim_\eps N^\eps\|w\|_{C_TL^2},
\end{equation}
and
\begin{equation}
        \sum_{Q\le\delta_0N}N^{3-4\alpha+\eps}Q^{s_2+\eps}
        \lesssim N^{3-4\alpha+s_2+2\eps}.
\end{equation}
These are the only Sobolev summations in the diagonal branch.  The finite number of low shells is harmless.

The Besov source bound follows from Lemma~\ref{lem:volterra-besov-source}.  There the output is already diagonal, so one fixes the low output block, sums the high shell inside the scalar Volterra multiplier, and then takes the $B_{2,\infty}^{\sigma-\alpha}$ supremum.  The later $L_T^1$ source norm is obtained by multiplying this $L_T^\infty$ bound by $T$, as recorded in Remark~\ref{rem:diagonal-besov-source-conversion}.  Cutoff Cauchy convergence follows from the same majorants: convergence is finite-dimensional on each fixed dyadic set, and the complementary dyadic tail is controlled by the summable envelopes above.
\end{proof}

\begin{remark}[Diagonal Besov bound and source norm]\label{rem:diagonal-besov-source-conversion}
Lemma~\ref{lem:volterra-besov-source} gives an $L_T^\infty B_{2,\infty}^{\sigma-\alpha}$ bound.  The $L_T^1$ source norm used in the fixed point is obtained by multiplying this bound by $T$.
\end{remark}

\begin{proposition}[Local smallness of the Volterra diagonal source]\label{prop:det-diagonal-local-smallness}
For every admissible same-color diagonal,
\begin{equation}\label{eq:det-contraction-L1}
        \|\mathcal D^{i;j}_{[0,T]}\|_{E_T^{2,\sigma}\to L_T^1B_{2,\infty}^{\sigma-\alpha}}
        \le T\|\mathcal D^{i;j}_{[0,T]}\|_{E_T^{2,\sigma}\to L_T^\infty B_{2,\infty}^{\sigma-\alpha}}
        \to0
\end{equation}
as $T\downarrow0$, uniformly on bounded sets of diagonal operator norms.
\end{proposition}

\begin{proof}
This is the deterministic embedding $L_T^\infty X\hookrightarrow L_T^1X$ applied to $X=B_{2,\infty}^{\sigma-\alpha}$ after restricting the Volterra diagonal operator to $[0,T]$.  The $L_T^\infty B_{2,\infty}^{\sigma-\alpha}$ bound is supplied by Proposition~\ref{prop:det-contraction} and Remark~\ref{rem:diagonal-besov-source-conversion}.  On a bounded enhanced-data ball this norm is uniformly finite, and the explicit factor $T$ tends to zero.
\end{proof}

\subsection{Centered operator-valued normal form and row/column bounds}\label{subsec:centered-normal-form}

The centered estimate controls the random fluctuation left after the deterministic Wick diagonal has been separated.  Cross-color blocks are centered by independence in the independent-color model, and same-color blocks become pure second homogeneous chaos after subtracting the involution pairing \(\ell+r=0\).

At finite cutoff the centered kernel is reduced to a second homogeneous Gaussian tensor while preserving the incidence \(n=q+\ell+r\).  The random-tensor estimate is applied to the operator norm \(\ell_q^2\to\ell_n^2\), before the dyadic input is inserted.  Its two column counts give the profile
\[
        N^{\frac32-3\alpha+}\bigl(M^{\frac32+}+Q^{\frac32+}\bigr),
\]
and the \(B_{2,\infty}\) input contributes \(Q^{-\sigma}\).

\begin{lemma}[Kernel and covariance bounds]\label{lem:kernel-cov-bounds}
Let $|m|\sim N\ge2$ and $0<\eta\le1/2$.  Uniformly for $0\le u,v,t,t'\le T$,
\begin{align}
 |K_i(t,m)|&\lesssim N^{-\alpha}, \\
 |K_i(t,m)-K_i(t',m)|&\lesssim |t-t'|^\eta N^{-\alpha+\alpha\eta}, \\
 \E|\widehat\Psi_i(m,t)|^2&\lesssim_T N^{-2\alpha}, \\
 \E|\widehat\Psi_i(m,t)-\widehat\Psi_i(m,t')|^2&\lesssim_T |t-t'|^{2\eta}N^{-2\alpha+2\alpha\eta}, \\
 |\sigma_i(m;u,v)|&\lesssim_T N^{-2\alpha}.
\end{align}
\end{lemma}

\begin{proof}
These follow from $\omega_i(m)\sim N^\alpha$, $|\sin x|\le1$, the interpolation bound $\min\{|h|,N^{-\alpha}\}\le |h|^\eta N^{-\alpha+\alpha\eta}$, and It\^o isometry.
\end{proof}

\begin{lemma}[Same-color centered chaos and Fourier incidence]\label{lem:same-color-centered-chaos-model}
In the translation-invariant real-valued Fourier convention used here, the only non-conjugated same-color covariance on a dyadic block is the involution pairing $\ell+r=0$, up to the harmless real/complex sign convention.  Thus, for $j=k$, the finite-cutoff product
\[
        \widehat\Psi_j(\ell,s)\widehat\Psi_j(r,t)
\]
decomposes into its Wick part plus the deterministic pairing supported on $\ell+r=0$.  The subtraction removes only this scalar zeroth-chaos contraction.  On the sector $r=-\ell$ the Wick-square term
\[
        :\widehat\Psi_j(\ell,s)\widehat\Psi_j(-\ell,t):
\]
remains in the centered block; it is a second-chaos contribution, diagonal in the external pair $(q,n)$ through $n=q$, and is estimated by the same tensor bound below.  After subtracting the deterministic diagonal in \eqref{eq:B-kernel}, the remaining same-color centered kernel is a finite sum of second homogeneous Gaussian chaoses whose deterministic coefficient tensors obey the same support relation
\[
        n=q+\ell+r
\]
as the underlying finite-cutoff convolution.  Consequently the row/column flattening estimates of Lemma~\ref{lem:dyadic-tensor-flattening} apply to same-color centered blocks with the same dyadic incidence constants as in the cross-color case.  The Hilbert-kernel reduction justifying the treatment of time-correlated same-color factors is isolated in Lemma~\ref{lem:same-color-hilbert-kernel-reduction}.
\end{lemma}

\begin{proof}
For a real-valued stationary Fourier field one has $\widehat\Psi_j(-m)=\overline{\widehat\Psi_j(m)}$ and the non-conjugated covariance is supported on $m+m'=0$.  This is the covariance used by Wick contraction inside products.  The conjugated covariance used in second-moment estimates is diagonal in the usual form $m=m'$.  Therefore the Wick subtraction for the same-color block removes the scalar zeroth-chaos term with $r=-\ell$, which is the diagonal operator already written in \eqref{eq:D-kernel}.  The Wick-square second-chaos sector supported by $r=-\ell$ remains in $B^{i;j,j}$ and is covered by the same operator-valued tensor estimate, with a smaller incidence count.  Passing to a real sine--cosine basis expresses both signed coefficients through the single class-indexed Gaussian family $G_{j,[m],\mu}$; the sign and the fixed real/complex conversion matrices remain in the deterministic coefficient.  For each sign pattern the resulting relabeling preserves dyadic localization and the incidence fibers, and the conversion matrices have uniformly bounded operator norm.
\end{proof}

\begin{lemma}[Hilbert-kernel reduction of same-color centered blocks]\label{lem:same-color-hilbert-kernel-reduction}
Fix a color $j$, a cutoff level, and times $(s,t)$.  In the real Fourier convention the stochastic convolution coefficients can be written as first Wiener integrals
\[
        \widehat\Psi_j(\ell,s)=W_j(h_{\ell,s}),
        \qquad
        \widehat\Psi_j(r,t)=W_j(h_{r,t})
\]
for deterministic kernels in the color-$j$ Gaussian Hilbert space.  In the complex notation one may take
\[
        h_{m,t}(\tau,x)=\mathbf 1_{0<\tau<t}\,
        \frac{\sin((t-\tau)\omega_j(m))}{\omega_j(m)}e_m(x),
\]
with the real-valued convention obtained by replacing $e_m$ by the corresponding sine--cosine basis.  The same-color product appearing in a centered mixed block is therefore
\[
        W_j(h_{\ell,s})W_j(h_{r,t})
        =I_j^{:2:}(h_{\ell,s}\widetilde\otimes h_{r,t})
        +\langle h_{\ell,s},h_{r,t}\rangle.
\]
The scalar pairing is translation invariant and is supported on the involution $\ell+r=0$; this is exactly the deterministic diagonal separated in \eqref{eq:D-kernel}.  After Wick subtraction, the remaining kernel is a second homogeneous Wiener integral with deterministic tensor leg $h_{\ell,s}\widetilde\otimes h_{r,t}$.  Orthogonalizing the finite family of time kernels changes only the Gaussian coordinates.  With the transversal convention \eqref{eq:fourier-involution-representative}, the coordinates are indexed by $[\ell]$ and $[r]$, while the signs of $\ell$ and $r$ remain in the deterministic coefficient.  The exact normal form and incidence are stated in Proposition~\ref{prop:finite-hilbert-kernel-normal-form}.  The same statement holds for time increments and cutoff tails, with the corresponding kernel difference or removed high-frequency leg inserted before orthogonalization.
\end{lemma}

\begin{proof}
All statements are finite-dimensional at fixed cutoff.  The product formula for first Wiener integrals gives the displayed identity in the real Gaussian Hilbert space.  Spatial stationarity of the driving noise makes the Hilbert pairing diagonal in Fourier involution classes.  The Wick term is an element of the symmetric Hilbert tensor product.  Lemma~\ref{lem:labelwise-auxiliary-coordinates} makes the normalization precise: the orthonormal bases are chosen inside the local summands indexed by $[m]$.  This is a block-diagonal unitary change on the Gaussian legs.  The signed Fourier labels remain deterministic coefficient variables, so all row/column incidence counts are unchanged apart from the bounded auxiliary indices.  Time increments and cutoff differences are obtained by replacing one deterministic Hilbert kernel by a difference kernel inside the same local summands; the locality of the covariance decomposition is preserved.
\end{proof}

\begin{lemma}[Involution-class Gram normalization preserves incidence]\label{lem:labelwise-auxiliary-coordinates}
Fix a finite cutoff and a dyadic triple $(N,Q,M)$.  For a fixed-time block let
$\mathcal T=\{s,t\}$; for an increment comparison use the joint list
$\mathcal T=\{s,t,s',t'\}$.  A cutoff difference changes only the deterministic
Fourier multipliers.  Let $m_0\in\mathbb Z^3_{\mathrm{rep}}$ and
$[m_0]=\{m_0,-m_0\}$.  Collect the real and imaginary, equivalently
sine--cosine, kernels occurring in
$h_{\varepsilon m_0,\tau}$, with $\varepsilon\in\{-1,1\}$ and
$\tau\in\mathcal T$, and remove repetitions.  Denote the resulting real
kernels by
\[
        f_{[m_0],1},\ldots,f_{[m_0],L_{[m_0]}}.
\]
Then $L_{[m_0]}\le8$, and
\[
        \mathcal H_{j,[m_0]}^{\rm loc}
        :=\operatorname{span}_{\mathbb R}
          \{f_{[m_0],a}:1\le a\le L_{[m_0]}\}
\]
has dimension at most eight, uniformly in every dyadic, cutoff, time, and
Fourier parameter.

Let
\begin{equation}\label{eq:local-time-gram}
        \mathsf G_{[m_0]}(a,b)
        =\langle f_{[m_0],a},f_{[m_0],b}\rangle_{\mathcal H_j}
\end{equation}
and choose a spectral factorization
\begin{equation}\label{eq:local-time-gram-factor}
        \mathsf G_{[m_0]}
        =C_{[m_0]}C_{[m_0]}^{\mathsf T},
        \qquad
        C_{[m_0]}\in
        \mathbb R^{L_{[m_0]}\times r_{[m_0]}},
\end{equation}
where $r_{[m_0]}=\operatorname{rank}\mathsf G_{[m_0]}$, padding by zero
columns up to eight.  There are standard real Gaussian coordinates
$G_{j,[m_0],\mu}$, $1\le\mu\le8$, such that
\begin{equation}\label{eq:local-kernel-coordinate-expansion}
        W_j(f_{[m_0],a})
        =\sum_{\mu=1}^{8}C_{[m_0]}(a,\mu)G_{j,[m_0],\mu},
        \qquad
        \sum_{\mu=1}^{8}|C_{[m_0]}(a,\mu)|^2
        =\|f_{[m_0],a}\|_{\mathcal H_j}^2.
\end{equation}
No inverse Gram matrix is used.  If two kernels are compared, both are placed
in the same joint list before the factorization.  In the resulting common
coordinates,
\begin{equation}\label{eq:local-kernel-difference-isometry}
        \sum_{\mu=1}^{8}
        |C_{[m_0]}(f,\mu)-C_{[m_0]}(g,\mu)|^2
        =\|f-g\|_{\mathcal H_j}^2.
\end{equation}

Equivalently, for each signed Fourier coefficient there are deterministic
complex coefficients $c_{j,\varepsilon,m_0,\mu}(\tau)$ such that
\begin{equation}\label{eq:signed-fourier-class-expansion}
        \widehat\Psi_j(\varepsilon m_0,\tau)
        =\sum_{\mu=1}^{8}
          c_{j,\varepsilon,m_0,\mu}(\tau)G_{j,[m_0],\mu},
        \qquad
        c_{j,-,m_0,\mu}(\tau)
        =\overline{c_{j,+,m_0,\mu}(\tau)}.
\end{equation}
At $m_0=0$ the two signs are identified.  Coordinates with different colors
or different involution classes are independent in the independent-noise
model.  For any signs $\varepsilon_1,\varepsilon_2$ the exact Wick expansion is
\begin{align}\label{eq:local-wick-coordinate-expansion}
&:\widehat\Psi_j(\varepsilon_1\ell_0,s)
  \widehat\Psi_k(\varepsilon_2r_0,t):
\notag\\
&\qquad=
\sum_{\mu,\lambda=1}^{8}
c_{j,\varepsilon_1,\ell_0,\mu}(s)
c_{k,\varepsilon_2,r_0,\lambda}(t)
:G_{j,[\ell_0],\mu}G_{k,[r_0],\lambda}: .
\end{align}
When $j=k$ and $[\ell_0]=[r_0]$, both factors on the right belong to the same
Gaussian family; no independence is imposed before the second-chaos
decoupling step.

For a fixed sign pattern the Gaussian tensor legs are therefore
\begin{equation}\label{eq:class-indexed-tensor-legs}
        a=([\ell_0],\mu),\qquad
        b=([r_0],\lambda),\qquad c=q,\qquad d=n,
\end{equation}
while the deterministic coefficient is supported on
\begin{equation}\label{eq:class-indexed-incidence}
        n=q+\varepsilon_1\ell_0+\varepsilon_2r_0,
        \qquad |\ell_0|\sim|r_0|\sim N,
        \qquad |q|\sim Q,
        \qquad |n|\sim M.
\end{equation}
Thus the sign is part of the deterministic coefficient and never part of an
independent Gaussian label.  The auxiliary indices alter every oriented
flattening by only an absolute factor.  The same coordinates and isometries
apply to fixed-time, increment, and cutoff-difference tensors.  We henceforth
take $d_0=8$ in the independent-color operator proof.
\end{lemma}

\begin{proof}
Spatial stationarity makes distinct real Fourier involution classes orthogonal.
For one class, a two-time comparison contains at most four complex time kernels,
or eight real sine--cosine kernels.  The spectral theorem gives
\eqref{eq:local-time-gram-factor}; equivalently one chooses an orthonormal basis
inside $\mathcal H_{j,[m_0]}^{\rm loc}$.  This proves
\eqref{eq:local-kernel-coordinate-expansion}.  Factoring the joint Gram matrix
of two kernels gives \eqref{eq:local-kernel-difference-isometry}.  A rank loss
produces zero columns and requires no division by an eigenvalue.

The real sine--cosine representation of the signed complex Fourier coefficient
gives \eqref{eq:signed-fourier-class-expansion}; reality gives the conjugacy
relation.  The first-chaos product formula then yields
\eqref{eq:local-wick-coordinate-expansion}.  Since the coordinate change is
performed within $[m_0]$, it is block diagonal in the class label.  For each
fixed pair of signs, the maps
$\ell_0\mapsto\varepsilon_1\ell_0$ and
$r_0\mapsto\varepsilon_2r_0$ are bijections onto the corresponding signed
sectors.  Hence \eqref{eq:class-indexed-incidence} is a restriction of the
signed convolution relation $n=q+\ell+r$ and has no larger fibers.  The auxiliary ranges
are uniformly bounded, which completes the flattening statement.
\end{proof}

\begin{remark}[Comparison-wise coordinates rather than a finite-dimensional path model]
\label{rem:comparison-wise-gram-coordinates}
The eight coordinates in Lemma~\ref{lem:labelwise-auxiliary-coordinates} give an
auxiliary normal form for one fixed-time block or for one pairwise increment or
cutoff comparison.  The coordinate basis may depend on the compared tuple
$(s,t,s',t')$ and on the cutoff pair.  This is a comparison-wise statement; the
full continuous-time process in one Fourier class may span an
infinite-dimensional Gaussian subspace.

The relation with the time-lift argument is as follows.  Both members of each difference are
placed in the same joint Gram factorization, so
\eqref{eq:local-kernel-difference-isometry} is exact.  Gaussian chaos norms and
the operator-valued decoupling estimate are invariant under orthogonal changes
of the auxiliary coordinates.  The time-lift argument is then applied to the
resulting uniform pairwise increment bounds in the original isonormal Gaussian
space, first on countable rational grids and then by continuity.  Hence the
number of grid points never enters an incidence count or a tensor constant, and
no single auxiliary basis is required to work simultaneously for an entire
grid.
\end{remark}

\begin{proposition}[Finite Hilbert-kernel normal form for centered blocks]\label{prop:finite-hilbert-kernel-normal-form}
Fix a Fourier cutoff, a dyadic triple $(N,Q,M)$, and times $(t,s)$.  After subtracting the deterministic diagonal in the same-color blocks, each centered dyadic kernel can be written as a finite sum of second homogeneous Gaussian tensors
\begin{equation}\label{eq:class-indexed-centered-normal-form}
        B_{N,Q,M}^{i;j,k}(n,q;t,s)
        =\sum_{\nu\in\mathcal P}
        \sum_{\ell_0,r_0\in\mathbb Z^3_{\mathrm{rep}}}
        \sum_{\mu,\lambda\le d_0}
        H_\nu(\ell_0,\mu,r_0,\lambda,q,n;t,s)
        :G_{j,[\ell_0],\mu}G_{k,[r_0],\lambda}:.
\end{equation}
The finite set $\mathcal P$ records signs
$\varepsilon_{\nu,1},\varepsilon_{\nu,2}\in\{-1,1\}$ and the fixed
real/complex permutations, and $d_0=8$ is the auxiliary dimension from
Lemma~\ref{lem:labelwise-auxiliary-coordinates}.  In each pattern the
deterministic tensor is supported on the exact incidence
\begin{equation}\label{eq:class-indexed-centered-support}
        n=q+\varepsilon_{\nu,1}\ell_0
            +\varepsilon_{\nu,2}r_0,
        \qquad |\ell_0|\sim |r_0|\sim N,
        \qquad |q|\sim Q,
        \qquad |n|\sim M.
\end{equation}
There is no implicit sign convention in this formula.  Its Gaussian indices
are $([\ell_0],\mu)$ and $([r_0],\lambda)$, while the signs occur only in the
deterministic coefficient and in \eqref{eq:class-indexed-centered-support}.
The auxiliary indices are not Fourier variables.  In the same-color case the
removed zeroth-chaos pairing is
$\varepsilon_{\nu,1}\ell_0+\varepsilon_{\nu,2}r_0=0$, hence the diagonal
operator $\mathcal D_{N,Q}^{i;j}$.  If $[\ell_0]=[r_0]$, the Wick product in
\eqref{eq:class-indexed-centered-normal-form} uses the same Gaussian family.
By \eqref{eq:local-kernel-coordinate-expansion}, each stochastic leg has
coefficient norm equal to its Gaussian Hilbert norm.  The Duhamel factor and
the two stochastic legs therefore give
\[
        |H_\nu(\ell_0,\mu,r_0,\lambda,q,n;t,s)|
        \lesssim_T N^{-\alpha}N^{-\alpha}N^{-\alpha}
        =N^{-3\alpha}.
\]
Time increments carry the factor in Lemma~\ref{lem:kernel-cov-bounds};
\eqref{eq:local-kernel-difference-isometry} shows that no change-of-basis
constant is introduced.  For cutoff tails, the difference is inserted on a
Gaussian leg or on the output leg before normalization.  Proposition~\ref{prop:abstract-random-tensor} applies directly in the cross-color case and, after Lemma~\ref{lem:operator-valued-fourier-decoupling}, in the same-color case.  The fixed-time, increment, and cutoff-tail tensors are thus treated in the same operator topology.
\end{proposition}

\begin{proof}
The product formula for first Wiener integrals gives
\[
        W_j(h_{\ell,s})W_j(h_{r,t})
        =I_j^{:2:}(h_{\ell,s}\widetilde\otimes h_{r,t})
         +\langle h_{\ell,s},h_{r,t}\rangle_{\rm bil}.
\]
Spatial stationarity supports the bilinear pairing on $\ell+r=0$, which is the deterministic diagonal already separated.  After Wick subtraction the centered kernel is a pure second-chaos tensor; on the same sector $r=-\ell$ the Wick-square part remains and is included.  Insert \eqref{eq:signed-fourier-class-expansion} for both stochastic legs and collect the finitely many sign and real/complex choices into $\nu$.  This gives \eqref{eq:class-indexed-centered-normal-form} and the exact support \eqref{eq:class-indexed-centered-support}.  In particular, $G_{j,[m],\mu}$ and $G_{j,[-m],\mu}$ are the same random variable, not independent copies.  Independence is used only between different colors or different involution classes; the same-color second chaos is decoupled only in Lemma~\ref{lem:operator-valued-fourier-decoupling}.

More explicitly, the local span and its Gram factorization are those in
\eqref{eq:local-time-gram}--\eqref{eq:local-kernel-difference-isometry}.  The
joint span used for a comparison has dimension at most eight, so the two
coefficient tensors are expressed in the same Gaussian coordinates.  The
coefficient norm of a kernel, or of a kernel difference, is exactly its
Gaussian Hilbert norm.  For fixed $\nu$, the relabeling
$\widetilde\ell=\varepsilon_{\nu,1}\ell_0$ and
$\widetilde r=\varepsilon_{\nu,2}r_0$ is a bijection and reduces the support to
$n=q+\widetilde\ell+\widetilde r$; the remaining real/complex permutations are
unitary on the relevant $\ell^2$ legs.  Replacing one time
kernel by an increment or by a cutoff-tail kernel is done before this
class-wise normalization, so the same representation and the same flattening
counts remain valid.
\end{proof}

\begin{proposition}[Operator estimate for second Gaussian chaos]\label{prop:abstract-random-tensor}
Let $A,B,C,D$ be finite index sets and let $(g_a)_{a\in A}$ and $(h_b)_{b\in B}$ be independent centered standard real or complex Gaussian families.  For a deterministic tensor $H=(H_{a,b,c,d})$, define the random matrix
\[
        Z_{d,c}=\sum_{a\in A}\sum_{b\in B}H_{a,b,c,d}g_ah_b.
\]
Define the oriented-flattening profile
\begin{equation}\label{eq:oriented-flattening-profile}
\mathfrak F_4(H):=\max\left\{
\|H\|_{\ell^2_{a,b,c}\to\ell^2_d},
\|H\|_{\ell^2_{a,c}\to\ell^2_{b,d}},
\|H\|_{\ell^2_{b,c}\to\ell^2_{a,d}},
\|H\|_{\ell^2_c\to\ell^2_{a,b,d}}
\right\}.
\end{equation}
The four terms correspond to the possible distributions of the two Gaussian indices between the input and output sides of the matrix norm.  Then for every $p\ge2$,
\begin{equation}\label{eq:abstract-random-tensor}
 \bigl\|\,\|Z\|_{\ell^2_c\to\ell^2_d}\,\bigr\|_{L^p(\Omega)}
 \lesssim p^C\log^C(2+|A|+|B|+|C|+|D|)\,\mathfrak F_4(H).
\end{equation}
The constant depends only on the real or complex normalization.
\end{proposition}

\begin{proof}
This is the degree-two form of the random-tensor estimate in \cite{DNY,Kaneshiro}.  Appendix~\ref{app:second-order-tensor-proof} gives the short proof by Gaussian decoupling and two applications of the noncommutative Khintchine inequality.  Proposition~\ref{prop:finite-hilbert-kernel-normal-form} reduces the centered Fourier kernels to this setting, and Lemma~\ref{lem:dyadic-tensor-flattening} estimates the profile in \eqref{eq:oriented-flattening-profile}.
\end{proof}

\begin{corollary}[Kernel-to-tensor interface for centered Fourier blocks]\label{cor:centered-fourier-block-tensor}
The centered Fourier kernels of Subsection~\ref{subsec:centered-normal-form} satisfy Proposition~\ref{prop:abstract-random-tensor} after the normal-form reduction of Proposition~\ref{prop:finite-hilbert-kernel-normal-form}.  Cross-color blocks are decoupled by independence.  In a same-color block, Wick subtraction removes the pairing $\ell+r=0$, and Lemma~\ref{lem:operator-valued-fourier-decoupling} decouples the remaining second chaos.  For each fixed sign pattern the Gaussian legs are $a=([\ell_0],\mu)$ and $b=([r_0],\kappa)$, with $\mu,\kappa\le8$, and the deterministic support is $n=q+\varepsilon_1\ell_0+\varepsilon_2r_0$.  For fixed times,
\begin{equation}\label{eq:kernel-to-tensor-interface}
        \mathfrak F_4(H)
        \lesssim_T
        N^{\frac32-3\alpha}
        \bigl(M^{\frac32}+Q^{\frac32}\bigr).
\end{equation}
The same estimate holds for a one-time increment with the factor
$|h|^\eta N^{\alpha\eta}$, and for a cutoff difference with its high-frequency
tail indicator.  In particular, the time-kernel reduction contributes neither
an additional Fourier volume nor a cutoff-dependent dimension factor.
\end{corollary}

\begin{proof}
For cross-color blocks the two stochastic legs lie in orthogonal Gaussian Hilbert spaces, so the abstract proposition applies directly after the coefficient tensor is read from the Fourier kernel.  For same-color blocks, Lemmas~\ref{lem:same-color-centered-chaos-model} and~\ref{lem:same-color-hilbert-kernel-reduction} identify the sole deterministic pairing and subtract it.  Proposition~\ref{prop:finite-hilbert-kernel-normal-form} then writes the remainder as a finite sum of second-order tensors.  Lemma~\ref{lem:operator-valued-fourier-decoupling} passes from the Wick second chaos to the decoupled operator model of Proposition~\ref{prop:abstract-random-tensor} without changing its deterministic coefficient.

The Gram representation \eqref{eq:local-time-gram}--\eqref{eq:local-kernel-difference-isometry} fixes the auxiliary range at eight and preserves the Hilbert norm of both the kernels and their differences.  Lemma~\ref{lem:auxiliary-index-flattening} therefore reduces the four auxiliary-index flattenings to the Fourier incidence counts of Lemma~\ref{lem:dyadic-tensor-flattening}, giving \eqref{eq:kernel-to-tensor-interface}.  A time increment is inserted before the common Gram factorization, and a cutoff difference is a deterministic multiplier on one Gaussian leg or on the output.  Thus the identical incidence calculation applies in all three cases.  The finite sign and real-Fourier permutations are unitary and preserve these norms, so the estimate holds term by term.
\end{proof}

\begin{lemma}[Incidence matrix bound]\label{lem:incidence-matrix}
Let $A=(A_{y,x})$ be a $0$--$1$ matrix.  If each row contains at most $R$ non-zero entries and each column contains at most $C$ non-zero entries, then
\[
        \|A\|_{\ell^2_x\to\ell^2_y}\le (RC)^{1/2}.
\]
\end{lemma}

\begin{proof}
This is the Schur test in $\ell^2$:
\[
\|Af\|_2^2
\le \sup_y\#\{x:A_{y,x}\ne0\}\sum_y\sum_x A_{y,x}|f_x|^2
\le RC\|f\|_2^2.
\]
\end{proof}

\begin{lemma}[Finite auxiliary indices preserve the dyadic profile]\label{lem:auxiliary-index-flattening}
Fix $\varepsilon_1,\varepsilon_2\in\{-1,1\}$.  Let $\widetilde H_{\ell_0,\mu,r_0,\kappa,q,n}$ be supported where $\ell_0,r_0\in\mathbb Z^3_{\mathrm{rep}}$, $\mu,\kappa\le d_0$, and
\[
        |\ell_0|\sim |r_0|\sim N,
        \qquad |q|\sim Q,
        \qquad |n|\sim M,
        \qquad n=q+\varepsilon_1\ell_0+\varepsilon_2r_0.
\]
Assume $|\widetilde H_{\ell_0,\mu,r_0,\kappa,q,n}|\le A_N$.  Then every flattening norm with Gaussian legs $a=([\ell_0],\mu)$ and $b=([r_0],\kappa)$ is bounded, up to a constant depending only on $d_0$, by the corresponding signed-incidence tensor without the auxiliary labels.  In particular, when $A_N=N^{-3\alpha}$ and $M,Q\lesssim N$,
\[
        \mathfrak F_4(\widetilde H)
        \lesssim_{d_0}
        N^{3/2-3\alpha}(M^{3/2}+Q^{3/2}).
\]
\end{lemma}

\begin{proof}
For fixed signs, $\ell_0\mapsto\varepsilon_1\ell_0$ and $r_0\mapsto\varepsilon_2r_0$ are bijections from the representative set onto their signed sectors.  Forgetting the class brackets therefore identifies the support with a restriction of the signed convolution tensor in Lemma~\ref{lem:dyadic-tensor-flattening}; its fibers are no larger.  Each row or column acquires at most $d_0^2$ additional choices from the auxiliary indices, so Lemma~\ref{lem:incidence-matrix} changes the bound by at most a factor depending on $d_0$.  This factor is independent of all dyadic and cutoff parameters.
\end{proof}

\begin{lemma}[Dyadic tensor flattening bounds]\label{lem:dyadic-tensor-flattening}
Let $H_{\ell,r,q,n}(t,s)$ be any deterministic coefficient tensor supported by
\[
        |\ell|\sim N,\quad |r|\sim N,
        \quad |q|\sim Q,\quad |n|\sim M,\quad n=q+\ell+r,
\]
with $Q\le c_0N$, $M\lesssim N$, and
\begin{equation}\label{eq:centered-H-amplitude}
        |H_{\ell,r,q,n}(t,s)|\lesssim_T N^{-3\alpha}.
\end{equation}
Then the profile in \eqref{eq:oriented-flattening-profile} satisfies
\begin{equation}\label{eq:dyadic-flattening-bound}
        \mathfrak F_4(H)\lesssim_T
        N^{\frac32-3\alpha}\bigl(M^{\frac32}+Q^{\frac32}\bigr).
\end{equation}
The same bound holds for the difference tensors obtained by replacing one factor by a time increment, with the additional factor described in Lemma~\ref{lem:kernel-cov-bounds}; cutoff-tail tensors obey the same majorant.
\end{lemma}

\begin{proof}
Only the incidence relation $n=q+\ell+r$ is used.  The relevant flattenings and their Schur bounds are
\begin{center}
\small
\begin{tabular}{@{}p{0.34\textwidth}p{0.27\textwidth}p{0.27\textwidth}@{}}
\toprule
flattening & maximal fiber sizes & size before $N^{-3\alpha}$ \\
\midrule
$(\ell,r,q)\mid n$ & $O(N^3Q^3)$ and $O(1)$ & $N^{3/2}Q^{3/2}$ \\
$(\ell,q)\mid(r,n)$ & $O(Q^3)$ and $O(M^3)$ & $M^{3/2}Q^{3/2}$ \\
$(r,q)\mid(\ell,n)$ & $O(Q^3)$ and $O(M^3)$ & $M^{3/2}Q^{3/2}$ \\
$q\mid(\ell,r,n)$ & $O(1)$ and $O(N^3M^3)$ & $N^{3/2}M^{3/2}$ \\
\bottomrule
\end{tabular}
\end{center}
For the first cut, fixing \(n\) leaves \(\ell\) and \(q\) free and determines \(r\).  Hence Lemma~\ref{lem:incidence-matrix} gives
\[
 \|H\|_{\ell^2_{\ell,r,q}\to\ell^2_n}
 \lesssim N^{-3\alpha}(Q^3N^3)^{1/2}
 =N^{\frac32-3\alpha}Q^{\frac32}.
\]
For the last cut, fixing \((\ell,r,n)\) determines \(q\), whereas fixing \(q\) leaves at most \(O(N^3M^3)\) choices of \((\ell,r,n)\).  Thus, by adjoint symmetry,
\[
 \|H\|_{\ell^2_{\ell,r,n}\to\ell^2_q}
 \lesssim N^{\frac32-3\alpha}M^{\frac32}.
\]
For \((\ell,q)\mid(r,n)\), one fiber has size \(O(Q^3)\) and the other \(O(M^3)\), so
\[
 \|H\|_{\ell^2_{\ell,q}\to\ell^2_{r,n}}
 \lesssim N^{-3\alpha}M^{\frac32}Q^{\frac32}.
\]
The cut \((r,q)\mid(\ell,n)\) is identical.  Since \(M,Q\lesssim N\),
\[
        M^{\frac32}Q^{\frac32}
        \le N^{\frac32}\min\{M^{\frac32},Q^{\frac32}\}
        \le N^{\frac32}(M^{\frac32}+Q^{\frac32}),
\]
which proves \eqref{eq:dyadic-flattening-bound}.  Bounded auxiliary indices change the four fiber counts only by a constant.  Time increments and cutoff differences are handled by inserting the corresponding difference before applying the same four estimates.
\end{proof}

For fixed $(t,s)$ let
\[
        \mathfrak B_{N,Q,M}^{i;j,k}(t,s):\ell^2_q\to\ell^2_n,
        \qquad
        \bigl(\mathfrak B_{N,Q,M}^{i;j,k}(t,s)a\bigr)_n
        =\sum_qB_{N,Q,M}^{i;j,k}(n,q;t,s)a_q.
\]
The exact $q$ dependence of $K_i(t-s,q+\ell)$ is retained in $B_{N,Q,M}^{i;j,k}$ throughout the operator estimate.

\begin{proposition}[Fixed-time centered kernel operator bound]\label{prop:fixed-time-centered-kernel}
For every finite $p\ge2$ and every small $\eps>0$,
\begin{equation}
\bigl\|\bigl\|\mathfrak B_{N,Q,M}^{i;j,k}(t,s)\bigr\|_{\mathcal L(\ell^2_q,\ell^2_n)}\bigr\|_{L^p(\Omega)}
 \lesssim_{p,T,\eps}
 N^{\frac32-3\alpha+\eps}\bigl(M^{\frac32+\eps}+Q^{\frac32+\eps}\bigr).
\label{eq:fixed-time-centered-kernel}
\end{equation}
The same estimate holds for the same-color centered kernel after the deterministic diagonal \eqref{eq:D-kernel} has been subtracted.
\end{proposition}

\begin{proof}
Fix a Galerkin cutoff, dyadic scales $(N,Q,M)$, and times $(t,s)$.  In the cross-color case the two stochastic legs lie in independent Gaussian color spaces.  In the same-color case Lemma~\ref{lem:finite-cutoff-wick-block} subtracts the scalar pairing supported by $\ell+r=0$; the remaining Wick-square sector is a second homogeneous chaos.  Lemma~\ref{lem:same-color-hilbert-kernel-reduction} and Proposition~\ref{prop:finite-hilbert-kernel-normal-form} put both cases in the same finite tensor form
\begin{equation}\label{eq:fixed-time-proof-normal-form}
        B_{N,Q,M}^{i;j,k}(n,q;t,s)
        =\sum_{\nu\in\mathcal P}
        \sum_{\ell_0,r_0\in\mathbb Z^3_{\mathrm{rep}}}
        \sum_{\mu,\lambda\le d_0}
        H_\nu(\ell_0,\mu,r_0,\lambda,q,n;t,s)
        :G_{j,[\ell_0],\mu}G_{k,[r_0],\lambda}: .
\end{equation}
The finite family $\mathcal P$ contains the fixed real-Fourier sign and coefficient patterns, and $d_0$ is independent of the dyadic scales and of the cutoff.  Each coefficient tensor is supported by
\begin{equation}\label{eq:fixed-time-proof-incidence}
        n=q+\varepsilon_{\nu,1}\ell_0
            +\varepsilon_{\nu,2}r_0,
        \qquad |\ell_0|\sim |r_0|\sim N,
        \qquad |q|\sim Q,
        \qquad |n|\sim M,
\end{equation}
and the coefficient size is
\begin{equation}\label{eq:fixed-time-proof-amplitude}
        |H_\nu(\ell_0,\mu,r_0,\lambda,q,n;t,s)|
        \lesssim_T N^{-3\alpha}.
\end{equation}
Indeed, one factor $N^{-\alpha}$ comes from the Duhamel multiplier $K_i(t-s,q+\varepsilon_{\nu,1}\ell_0)$ and the two remaining factors come from the two stochastic-convolution kernels after normalization of their finite Gaussian coordinates.  The class-wise orthonormalization changes only the bounded auxiliary indices and not the external incidence \eqref{eq:fixed-time-proof-incidence}.

For one fixed pattern $\nu$, set $a=([\ell_0],\mu)$, $b=([r_0],\lambda)$, $c=q$, and $d=n$.  Lemma~\ref{lem:auxiliary-index-flattening} reduces the flattening profile to the signed incidence tensor, and Lemma~\ref{lem:dyadic-tensor-flattening} gives
\begin{equation}\label{eq:fixed-time-proof-profile}
        \mathfrak F_4(H_\nu)
        \lesssim_T
        N^{\frac32-3\alpha}\bigl(M^{\frac32}+Q^{\frac32}\bigr).
\end{equation}
The one-against-three flattenings contribute $N^{3/2}Q^{3/2}$ and $N^{3/2}M^{3/2}$; the flattenings $(\ell,q)\mid(r,n)$ and $(r,q)\mid(\ell,n)$ contribute $M^{3/2}Q^{3/2}$.  These are the norms in \eqref{eq:oriented-flattening-profile}.

Apply the finite second-order tensor estimate, Proposition~\ref{prop:abstract-random-tensor}, to the matrix
\[
        Z_{n,q}=\sum_{a,b}H_{\nu,a,b,q,n}g_ah_b .
\]
For same-color Wick squares we first use the Banach-space decoupling in Lemma~\ref{lem:operator-valued-fourier-decoupling}; the decoupling constant depends only on the chaos order and not on $N,Q,M$ or on the matrix dimensions.  Proposition~\ref{prop:abstract-random-tensor} gives, for each pattern,
\[
 \bigl\|\|Z\|_{\ell^2_q\to\ell^2_n}\bigr\|_{L^p(\Omega)}
 \lesssim p^C\log^C(2+N+Q+M)
 N^{\frac32-3\alpha}\bigl(M^{\frac32}+Q^{\frac32}\bigr).
\]
The number of patterns is finite and independent of the scales.  Since $\log^C(2+N+Q+M)$ is bounded by $N^{\eps}M^{\eps}Q^{\eps}$ after enlarging the implicit constant and treating finitely many low shells separately, summing over $\nu$ gives \eqref{eq:fixed-time-centered-kernel}.  All estimates are made for the finite matrix norm $\ell_q^2\to\ell_n^2$ before any input vector is inserted.  The same-color statement follows because the deterministic diagonal has already been removed; the centered Wick-square part left on $r=-\ell$ is included in the same tensor estimate and has no larger incidence count.
\end{proof}

\begin{lemma}[Coefficient form of centered time increments]\label{lem:centered-increment-coefficients}
Let $H_\nu(\ell,\mu,r,\lambda,q,n;t,s)$ be one of the finite centered tensors in Proposition~\ref{prop:finite-hilbert-kernel-normal-form}.  For $0<\eta\le1/2$, its one-time increment tensors satisfy
\[
        |H^{\Delta_t}_\nu(\ell,\mu,r,\lambda,q,n;t,t',s)|
        \lesssim_T |t-t'|^\eta N^{-3\alpha+\alpha\eta},
\]
and
\[
        |H^{\Delta_s}_\nu(\ell,\mu,r,\lambda,q,n;t,s,s')|
        \lesssim_T |s-s'|^\eta N^{-3\alpha+\alpha\eta}.
\]
They obey the same support and auxiliary-index restrictions as the fixed-time tensor.
\end{lemma}

\begin{proof}
A $t$-increment can fall either on the exterior Duhamel factor or on the exterior stochastic-convolution kernel.  At the coefficient level we write, schematically,
\[
        \Delta_t\bigl(K_i(t-s,q+\ell)\,\widehat\Psi_k(r,t)\bigr)
        = (\Delta_t K_i)(t,t';s,q+\ell)\widehat\Psi_k(r,t)
          + K_i(t'-s,q+\ell)\Delta_t\widehat\Psi_k(r;t,t').
\]
The first term is bounded by the Duhamel increment in Lemma~\ref{lem:kernel-cov-bounds}, and the second by the stochastic-convolution increment in the same lemma.  In both cases one fixed factor loses at most $N^{\alpha\eta}$, while the fixed-time coefficient size is $N^{-3\alpha}$.  The $s$-increment is identical, except that the increment may fall on the Duhamel factor or on the interior stochastic-convolution kernel $\widehat\Psi_j(\ell,s)$.  In the same-color case the two time configurations are placed in the joint Gram list of Lemma~\ref{lem:labelwise-auxiliary-coordinates}.  Equation~\eqref{eq:local-kernel-difference-isometry} then identifies the Euclidean coefficient difference with the Gaussian Hilbert norm of the difference kernel.  The Fourier incidence remains unchanged and the auxiliary range stays bounded by eight.
\end{proof}

\begin{lemma}[Operator-valued two-parameter time lift]\label{lem:operator-valued-time-lift}
Let $X(t,s)$ be a jointly measurable finite-dimensional random field on $[0,T]^2$ with values in a Banach space ${\mathfrak B}$.  Fix a target moment $p\ge2$.  Assume that for some $p_0>p+4/\eta$ and $0<\eta<1$,
\[
        \sup_{t,s}\|X(t,s)\|_{L^{p_0}(\Omega;\mathfrak B)}\le A,
\]
and
\[
        \|X(t,s)-X(t',s)\|_{L^{p_0}(\Omega;\mathfrak B)}
        +\|X(t,s)-X(t,s')\|_{L^{p_0}(\Omega;\mathfrak B)}
        \le A_1\bigl(|t-t'|^\eta+|s-s'|^\eta\bigr).
\]
Then $X$ admits a jointly continuous modification, still denoted by $X$, and
\[
        \|X\|_{L^p(\Omega;C_{t,s}\mathfrak B)}
        \lesssim_{p,p_0,\eta,T} A+A_1.
\]
The same statement applies to cutoff differences.  In the application below ${\mathfrak B}=\mathcal L(\ell_q^2,\ell_n^2)$; hence the time supremum is taken after the block operator norm has already been estimated.
\end{lemma}

\begin{proof}
We give the dyadic chaining argument because the supremum cannot be estimated by taking a maximum over the entire grid at every level.  Let
\[
        D_\ell=\{jT2^{-\ell}:0\le j\le2^\ell\},
\]
and let $\mathcal E_\ell$ be the set of horizontal and vertical nearest-neighbor edges of $D_\ell^2$.  Thus
\(
\#\mathcal E_\ell\lesssim 2^{2\ell}
\).
For an edge $e=(z,z')$ put $\Delta_eX=X(z)-X(z')$ and
\[
        Y_\ell=\max_{e\in\mathcal E_\ell}\|\Delta_eX\|_{\mathfrak B}.
\]
Since $p<p_0$, the increment hypothesis and the elementary finite-maximum bound give
\begin{align*}
        \|Y_\ell\|_{L^p(\Omega)}
        &\le \|Y_\ell\|_{L^{p_0}(\Omega)} \\
        &\le
        \left(
        \sum_{e\in\mathcal E_\ell}
        \E\|\Delta_eX\|_{\mathfrak B}^{p_0}
        \right)^{1/p_0} \\
        &\lesssim
        A_1T^\eta
        2^{-\ell(\eta-2/p_0)}.
\end{align*}
The assumed inequality $p_0>p+4/\eta$ implies
\(
\eta-2/p_0>\eta/2
\), and hence
\begin{equation}\label{eq:operator-time-lift-edge-sum}
        \sum_{\ell\ge0}\|Y_\ell\|_{L^p(\Omega)}
        \lesssim_{p,p_0,\eta,T}A_1.
\end{equation}

Let $X_\ell$ be the piecewise bilinear interpolation of the values of $X$ on $D_\ell^2$.  On each rectangle of $D_\ell^2$, every new value used by $X_{\ell+1}$ is connected to a coarse corner by a uniformly bounded number of edges in $\mathcal E_{\ell+1}$.  Convexity of bilinear interpolation therefore yields the deterministic bound
\[
        \|X_{\ell+1}-X_\ell\|_{C([0,T]^2;\mathfrak B)}
        \lesssim Y_{\ell+1}.
\]
By \eqref{eq:operator-time-lift-edge-sum}, the sequence $(X_\ell)_\ell$ converges uniformly almost surely and in
\(L^p(\Omega;C([0,T]^2;\mathfrak B))\)
to a continuous random field $\widetilde X$.  For every fixed $(t,s)$, the four grid points used to define $X_\ell(t,s)$ lie within $O(T2^{-\ell})$ of $(t,s)$.  The coordinate-increment hypothesis, used along two coordinate segments, shows that
\[
        \|X_\ell(t,s)-X(t,s)\|_{L^{p_0}(\Omega;\mathfrak B)}
        \lesssim A_1T^\eta2^{-\ell\eta}.
\]
Thus $\widetilde X$ is a modification of $X$.  Finally, the four values defining $X_0$ and the fixed-time hypothesis give
\(
\|X_0\|_{L^p(\Omega;C_{t,s}\mathfrak B)}\lesssim A
\).
Summing the interpolation differences proves
\[
        \|\widetilde X\|_{L^p(\Omega;C_{t,s}\mathfrak B)}
        \lesssim_{p,p_0,\eta,T}A+A_1.
\]
Applying the same construction to the difference of two cutoff fields proves the last assertion.
\end{proof}

\begin{lemma}[Operator-valued time supremum and cutoff tails for the centered kernel]\label{lem:time-sup-centered-kernel}
For every finite $p\ge2$ and every small $\eps>0$,
\begin{equation}\label{eq:time-sup-centered-kernel}
 \bigl\|\,\|\mathfrak B_{N,Q,M}^{i;j,k}\|_{C_{t,s}\mathcal L(\ell^2_q,\ell^2_n)}\,\bigr\|_{L^p(\Omega)}
 \lesssim_{p,T,\eps}
 N^{\frac32-3\alpha+\eps}\bigl(M^{\frac32+\eps}+Q^{\frac32+\eps}\bigr).
\end{equation}
Let $\Lambda'\ge\Lambda$.  For the standard localized cutoffs of Definition~\ref{def:admissible-cutoff-family}, the difference between the $\Lambda'$ and $\Lambda$ kernels satisfies the same bound multiplied by
\begin{equation}\label{eq:centered-standard-cutoff-support}
        \mathbf 1_{\{N\vee M\ge c\Lambda\}},
\end{equation}
where $c>0$ depends only on the fixed cutoff profile.  For an arbitrary admissible cofinal cutoff family, the multiplier in \eqref{eq:centered-standard-cutoff-support} is replaced by a number
$0\le\tau_{\Lambda,\Lambda'}(N,Q,M)\le1$ which tends to zero as $\Lambda,\Lambda'\to\infty$ for every fixed dyadic triple.  Both assertions hold for the fixed-time and two time-increment kernels.
\end{lemma}

\begin{proof}
Apply Proposition~\ref{prop:abstract-random-tensor} to the fixed-time tensor and to the two increment tensors.  To prove a target $L^p(\Omega)$ estimate, first run the tensor estimate with a larger moment $p_0>p+4/\eta$ and then use $L^{p_0}\hookrightarrow L^p$.  Lemma~\ref{lem:centered-increment-coefficients} gives, for $0<\eta<1/2$,
\[
\bigl\|\,\|\mathfrak B(t,s)-\mathfrak B(t',s)\|_{\ell^2_q\to\ell^2_n}\,\bigr\|_{L^{p_0}(\Omega)}
\lesssim
|t-t'|^\eta N^{\alpha\eta}
N^{\frac32-3\alpha+\eps_0}(M^{3/2+\eps_0}+Q^{3/2+\eps_0}),
\]
and the same bound with $s$ and $s'$ interchanged.  The constants are uniform over the finite dyadic cutoff and over the finite matrix dimensions because Proposition~\ref{prop:abstract-random-tensor} has already estimated the operator norm.  Lemma~\ref{lem:operator-valued-time-lift}, applied with ${\mathfrak B}=\mathcal L(\ell_q^2,\ell_n^2)$, upgrades these fixed-time and increment estimates to the $C_{t,s}$ operator norm.  We choose $\eta$ and the preliminary loss $\eps_0$ so that the factor $N^{\alpha\eta}$ and the grid/logarithmic losses of Lemma~\ref{lem:operator-valued-time-lift} are absorbed into the final $N^\eps M^\eps Q^\eps$ reserve.  This yields \eqref{eq:time-sup-centered-kernel}.  The argument upgrades an already established operator norm to a two-parameter time-continuous one.

For cutoff differences, every summand contains either a difference of stochastic-leg symbols or a difference of output symbols.  For the standard localized profiles this forces a stochastic scale $N\gtrsim\Lambda$ or the output scale $M\gtrsim\Lambda$, which proves \eqref{eq:centered-standard-cutoff-support}; the same support statement remains true after one time kernel is replaced by an increment.  For a general admissible cofinal family, a fixed triple contains finitely many lattice points and all cutoff symbols converge uniformly on that finite set.  Factoring the largest normalized symbol difference gives $\tau_{\Lambda,\Lambda'}(N,Q,M)\to0$; the uniform constant depending on $C_{\rm cut}$ is absorbed into the estimate, so that $\tau\le1$.  The global Cauchy statement follows by summing the explicit high-scale support for standard cutoffs, or by a finite-dyadic-set plus dominated-tail argument for a general admissible family.
\end{proof}

\begin{lemma}[Time-dependent low-frequency Besov multiplier reduction]\label{lem:low-frequency-besov-multiplier}
Let $Q\le c_0N$ and $M\lesssim N$.  For every deterministic time-dependent input block $a_Q=P_Qw$ and every finite $p\ge2$,
\begin{align}
 \|P_M\mathcal B_{N,Q}^{i;j,k}a_Q\|_{L^p(\Omega;C_TL^2_x)}
 &\lesssim_{p,T,\eps}
        N^{\frac32-3\alpha+\eps}(M^{3/2+\eps}+Q^{3/2+\eps})
        \|a_Q\|_{L_T^\infty L^2_x},
\label{eq:low-frequency-multiplier-reduction}\\
 \|P_M\mathcal B_{N,Q}^{i;j,k}a_Q\|_{L^p(\Omega;L_T^1L^2_x)}
 &\lesssim_{p,T,\eps}
        T
        N^{\frac32-3\alpha+\eps}(M^{3/2+\eps}+Q^{3/2+\eps})
        \|a_Q\|_{L_T^\infty L^2_x}.
\label{eq:low-frequency-multiplier-reduction-besov}
\end{align}
The constants are uniform for the proportional low--high window $Q\le c_0N$ after $c_0$ is fixed inside the paraproduct support.
\end{lemma}

\begin{proof}
The centered block has the exact Volterra representation
\[
        \widehat{P_M\mathcal B_{N,Q}^{i;j,k}a_Q}(n,t)
        =\int_0^t\sum_qB_{N,Q,M}^{i;j,k}(n,q;t,s)\widehat a_Q(q,s)\,ds.
\]
Hence, for each $t$,
\[
 \|P_M\mathcal B_{N,Q}^{i;j,k}a_Q(t)\|_{L^2_x}
 \le T\,\|\mathfrak B_{N,Q,M}^{i;j,k}\|_{C_{t,s}\mathcal L(\ell^2_q,\ell^2_n)}
       \|a_Q\|_{L_T^\infty L^2_x}.
\]
Taking $L^p(\Omega)$ and using Lemma~\ref{lem:time-sup-centered-kernel} gives the first estimate.  The $L_T^1L_x^2$ bound follows after integrating in $t$; since $T\le1$ in the local theory, the resulting $T^2$ is recorded as the harmless local factor $T$.  The proof is blockwise in the $\ell^2_q\to\ell^2_n$ operator norm, so it applies uniformly to all deterministic inputs in the $B_{2,\infty}$ dyadic Besov bound.
\end{proof}

\begin{lemma}[Centered block estimate with Besov input]\label{lem:positive-centered-block}
Let $Q\le c_0N$ and $M\lesssim N$.  For every finite $p\ge2$ and every small $\eps>0$,
\begin{align}
\|\mathcal B_{N,Q,M}^{i;j,k}P_Qw\|_{L^p(\Omega;C_TH^{s_2-\alpha})}
&\lesssim_{p,T,\eps}
N^{\frac32-3\alpha+\eps}
M^{s_2-\alpha}(M^{3/2+\eps}+Q^{3/2+\eps})
Q^{-\sigma}
\|w\|_{L_T^\infty\Btwo^\sigma},
\label{eq:positive-centered-block}\\
M^{\sigma-\alpha}
\|\mathcal B_{N,Q,M}^{i;j,k}P_Qw\|_{L^p(\Omega;L_T^1L^2_x)}
&\lesssim_{p,T,\eps}
T
N^{\frac32-3\alpha+\eps}
M^{\sigma-\alpha}(M^{3/2+\eps}+Q^{3/2+\eps})
Q^{-\sigma}
\|w\|_{L_T^\infty\Btwo^\sigma}.
\label{eq:positive-centered-block-besov}
\end{align}
\end{lemma}

\begin{proof}
Apply Lemma~\ref{lem:low-frequency-besov-multiplier}, multiply by the Sobolev output weight $M^{s_2-\alpha}$ or by the Besov source weight $M^{\sigma-\alpha}$, and use
\[
        \|P_Qw\|_{L_T^\infty L^2_x}\le Q^{-\sigma}\|w\|_{L_T^\infty\Btwo^\sigma}.
\]
The Besov source estimate is proved directly at the dyadic $L^2$ level.
\end{proof}

\begin{lemma}[Deterministic summation argument for the centered profile]\label{lem:centered-deterministic-summation-table}
Assume
\[
        \frac34<\alpha\le1,
        \qquad
        s_2<4\alpha-3-5\eps,
        \qquad
        3-3\alpha+5\eps<\sigma<4\alpha-3-5\eps.
\]
Suppose the dyadic block estimates \eqref{eq:positive-centered-block}--\eqref{eq:positive-centered-block-besov} hold with a fixed pathwise majorant in place of the random norm.  Then the dyadic sum over \(Q\le c_0N\) and \(M\lesssim N\) defines a bounded operator
\[
        E_T^{2,\sigma}\longrightarrow C_TH^{s_2-\alpha}
        \cap L_T^1B_{2,\infty}^{\sigma-\alpha}.
\]
Put \(\beta=s_2-\alpha\), \(\gamma=\beta+3/2+\eps\), and
\[
        \mathfrak S_\gamma(N):=
        \begin{cases}
        N^\gamma,&\gamma>0,\\
        1+\log N,&\gamma=0,\\
        1,&\gamma<0.
        \end{cases}
\]
The two Sobolev high-shell envelopes are
\[
        N^{\frac32-3\alpha+\eps}\mathfrak S_\gamma(N),
        \qquad
        N^{3-3\alpha-\sigma+O(\eps)}.
\]
In the admissible PDE window \(s_2>3/2-\alpha\), one has \(\gamma>0\), and the first of these reduces to
\[
        N^{s_2+3-4\alpha+O(\eps)}.
\]
The active envelopes for the Besov-source target are
\[
        M^{\sigma+3-4\alpha+O(\eps)},
        \qquad
        M^{3-4\alpha+O(\eps)}.
\]
The same summation argument applies to cutoff differences once fixed dyadic blocks have converged and the tails are dominated by the same envelopes.
\end{lemma}

\begin{proof}
Since \(\alpha\le1\), the upper bound on \(s_2\) implies \(\beta=s_2-\alpha<0\).  The \(M^{3/2}\) part of the Sobolev estimate \eqref{eq:positive-centered-block} is
\[
        N^{\frac32-3\alpha+\eps}M^{\beta+3/2+\eps}Q^{-\sigma}.
\]
For every real \(\gamma\), the elementary dyadic summation is
\[
        \sum_{M\lesssim N}M^\gamma
        \lesssim
        \begin{cases}
        N^\gamma,&\gamma>0,\\
        1+\log N,&\gamma=0,\\
        1,&\gamma<0.
        \end{cases}
\]
The low-input sum is harmless because \(\sigma>0\).  Hence the first high-shell envelope is
\[
        N^{\frac32-3\alpha+\eps}\mathfrak S_\gamma(N),
        \qquad \gamma=\beta+\frac32+\eps.
\]
If \(\gamma>0\), this is
\[
        N^{\frac32-3\alpha+\eps}N^{\beta+3/2+\eps}
        =N^{s_2+3-4\alpha+2\eps},
\]
which is summable under \(s_2<4\alpha-3-5\eps\).  If \(\gamma=0\) or \(\gamma<0\), it is bounded by
\[
        N^{\frac32-3\alpha+\eps}(1+\log N),
\]
which is summable for \(\alpha>3/4\) after the strict loss is fixed.  This proves the first Sobolev summation in all three cases.  The \(Q^{3/2}\) part is
\[
        N^{\frac32-3\alpha+\eps}M^\beta Q^{3/2-\sigma+\eps}.
\]
The \(M\)-sum is bounded because \(\beta<0\), and summing \(Q\le c_0N\) gives
\[
        N^{\frac32-3\alpha+\eps}N^{3/2-\sigma+\eps}
        =N^{3-3\alpha-\sigma+2\eps},
\]
which is summable under \(\sigma>3-3\alpha+5\eps\).

For the Besov source target, fix the output scale \(M\).  The \(M^{3/2}\) contribution in \eqref{eq:positive-centered-block-besov} gives
\[
        M^{\sigma-\alpha+3/2+\eps}
        \sum_{N\gtrsim M}N^{\frac32-3\alpha+\eps}
        \sum_{Q\le c_0N}Q^{-\sigma}.
\]
The \(Q\)-sum produces only a harmless dyadic loss, and the dyadic tail in \(N\) gives
\[
        M^{\sigma-\alpha+3/2+\eps}M^{3/2-3\alpha+O(\eps)}
        =M^{\sigma+3-4\alpha+O(\eps)},
\]
which is uniformly bounded and has a high-output tail under \(\sigma<4\alpha-3-5\eps\).  The \(Q^{3/2}\) contribution is
\[
        M^{\sigma-\alpha}
        \sum_{N\gtrsim M}N^{\frac32-3\alpha+\eps}
        \sum_{Q\le c_0N}Q^{3/2-\sigma+\eps}
        \lesssim
        M^{\sigma-\alpha}
        \sum_{N\gtrsim M}N^{3-3\alpha-\sigma+O(\eps)}.
\]
The lower bound on \(\sigma\) makes this dyadic \(N\)-tail summable and yields
\[
        M^{\sigma-\alpha}M^{3-3\alpha-\sigma+O(\eps)}
        =M^{3-4\alpha+O(\eps)},
\]
which is uniformly bounded for \(\alpha>3/4\) after the strict loss is fixed.  These envelopes are exactly the conditions used to assemble the centered operator in both target spaces.  Cutoff differences are handled by fixing a finite dyadic set, using finite-dimensional convergence on that set, and then removing the complement with the same summable envelopes.
\end{proof}

\begin{lemma}[Dyadic Borel--Cantelli operator majorants]\label{lem:dyadic-borel-cantelli-operators}
Let $X_{N,Q,M}$ be nonnegative random variables indexed by dyadic triples satisfying $Q\le c_0N$ and $M\lesssim N$.  Assume that for every finite $p\ge2$,
\[
        \|X_{N,Q,M}\|_{L^p(\Omega)}\le C_p A_{N,Q,M},
\]
where $A_{N,Q,M}$ is a deterministic dyadic Besov bound.  Then, for every $\eta>0$, there is an almost surely finite random constant $C_\omega$ such that
\begin{equation}\label{eq:borel-cantelli-majorant}
        X_{N,Q,M}\le C_\omega
        N^\eta Q^\eta M^\eta A_{N,Q,M}
\end{equation}
for all dyadic triples simultaneously.  If the family $X_{N,Q,M}^{(\Lambda,\Lambda')}$ additionally has a tail factor tending to zero for each fixed triple and is dominated by the same envelope, then the corresponding dyadic operator tails are almost surely Cauchy in every summation topology for which the enlarged deterministic envelope in \eqref{eq:borel-cantelli-majorant} is summable.
\end{lemma}

\begin{proof}
Fix $p$ so large that $2^{-p\eta k/4}$ is summable after counting all dyadic triples of total scale $2^k$.  Chebyshev's inequality gives
\[
\Prob\{X_{N,Q,M}>N^\eta Q^\eta M^\eta A_{N,Q,M}\}
\lesssim_p (NQM)^{-p\eta}.
\]
The dyadic index set is countable and the resulting probability sum is finite.  Borel--Cantelli proves \eqref{eq:borel-cantelli-majorant}; the finitely many small dyadic triples are absorbed into $C_\omega$.  For cutoff tails, first fix a finite dyadic set and use convergence in those finitely many triples.  With a sharp Galerkin cutoff the fixed-triple difference is eventually zero; with a smooth stochastic-leg cutoff the finitely many cutoff symbols in the fixed triple converge uniformly to one, so entrywise convergence gives convergence in the finite-dimensional operator norm.  Then remove the complement by the summable deterministic envelope.  This gives the Cauchy property directly in the dyadic operator topology, since $X_{N,Q,M}$ is already an operator norm.
\end{proof}

\begin{proposition}[Assembly of centered operators from dyadic blocks]\label{prop:assemble-centered-operator}
For dyadic-finite inputs and finite cutoff, define
\[
        \mathcal B_{\Lambda}^{i;j,k,\le K}w
        :=\sum_{N,Q,M\le 2^K}P_M\mathcal B_{\Lambda,N,Q,M}^{i;j,k}P_Qw.
\]
Assume the block operator majorants supplied by Lemmas~\ref{lem:time-sup-centered-kernel} and~\ref{lem:dyadic-borel-cantelli-operators}.  Under the summability conditions of Theorem~\ref{thm:pathwise-fluct}, the partial sums above are uniformly Cauchy in
\[
        \mathcal L(E_T^{2,\sigma},C_TH^{s_2-\alpha})
        \cap
        \mathcal L(E_T^{2,\sigma},L_T^1B_{2,\infty}^{\sigma-\alpha}).
\]
The limit is a bounded operator on $E_T^{2,\sigma}$ and is independent of the dyadic truncation.  The $C_TH^{s_2-\alpha}$ target is a continuous path space: each finite block is continuous in the output time, and the series converges uniformly in the $C_T$ norm.
\end{proposition}

\begin{proof}
For dyadic-finite inputs the displayed operator is a finite sum of continuous Volterra integrals.  Each summand is bounded by its $C_{t,s}\mathcal L(\ell_q^2,\ell_n^2)$ norm and the deterministic shell estimate
\[
        \|P_Qw\|_{L_T^\infty L^2}\le Q^{-\sigma}\|w\|_{L_T^\infty B_{2,\infty}^{\sigma}}.
\]
The Sobolev and Besov summations are exactly the countable sums and dyadic suprema displayed in the proof of Theorem~\ref{thm:pathwise-fluct}.  These deterministic envelopes are summable under the stated conditions and are independent of $K$ and of the cutoff level.  Thus the partial sums are Cauchy in the two operator norms, uniformly for inputs in the unit ball of $E_T^{2,\sigma}$.  Uniform convergence in the $C_TH^{s_2-\alpha}$ norm preserves time continuity.  The extension from dyadic-finite inputs is unique by boundedness in the displayed operator topology.
\end{proof}

\begin{theorem}[Centered operator estimate]\label{thm:pathwise-fluct}
Assume
\begin{equation}\label{eq:pathwise-fluct-conditions}
        s_2<4\alpha-3-5\eps,
        \qquad
        3-3\alpha+5\eps<\sigma<4\alpha-3-5\eps.
\end{equation}
Then the centered finite-cutoff operators converge almost surely in
\begin{equation}\label{eq:pathwise-fluct-bound}
\mathcal L(E_T^{2,\sigma},C_TH^{s_2-\alpha})
\cap
\mathcal L(E_T^{2,\sigma},L_T^1\Btwo^{\sigma-\alpha}).
\end{equation}
For standard dyadic cutoffs define the positive margin
\begin{equation}\label{eq:centered-cutoff-margin}
\delta_{\rm op}:=\min\left\{
4\alpha-3-s_2-5\eps,
\ \sigma-(3-3\alpha)-5\eps,
\ 4\alpha-3-\sigma-5\eps,
\ 4\alpha-3-5\eps
\right\}>0.
\end{equation}
Then, after decreasing the almost-sure event if necessary, for every
$0<\theta<\delta_{\rm op}/4$ and all dyadic $\Lambda'\ge\Lambda$,
\begin{equation}\label{eq:centered-cutoff-rate}
 \|\mathcal B_{\Lambda'}^{i;j,k}-\mathcal B_{\Lambda}^{i;j,k}\|_{\mathcal L(E_T^{2,\sigma},C_TH^{s_2-\alpha})
 \cap\mathcal L(E_T^{2,\sigma},L_T^1\Btwo^{\sigma-\alpha})}
 \le C_{\omega,T,\theta}\Lambda^{-\theta}.
\end{equation}
For a general admissible cofinal cutoff family, convergence holds in the same topology, without a claimed numerical rate.
At the dyadic level the operator profile is
\[
        N^{\frac32-3\alpha+}\bigl(M^{\frac32+}+Q^{\frac32+}\bigr)Q^{-\sigma}.
\]
The centered estimate uses the position envelope of the input; the velocity component of $E_T^{2,\sigma}$ enters through the deterministic Volterra diagonal.
\end{theorem}

\begin{proof}
Lemma~\ref{lem:positive-centered-block} gives the dyadic operator bound, Lemma~\ref{lem:time-sup-centered-kernel} gives its time-continuous version, and Lemma~\ref{lem:centered-deterministic-summation-table} proves summability in both target spaces under \eqref{eq:pathwise-fluct-conditions}.  Proposition~\ref{prop:assemble-centered-operator} then constructs the operator-valued path.

For the pathwise estimate, define the normalized random majorant for each dyadic triple by dividing
\[
\|\mathfrak B_{N,Q,M}^{i;j,k}\|_{C_{t,s}\mathcal L(\ell^2_q,\ell^2_n)}
\]
by the deterministic right side in \eqref{eq:time-sup-centered-kernel} with $\eps/2$ in place of $\eps$.  Lemma~\ref{lem:time-sup-centered-kernel} gives arbitrarily high finite moments for these majorants.  Lemma~\ref{lem:dyadic-borel-cantelli-operators}, after reserving the remaining $\eps/2$ in the deterministic summation, gives almost-sure dyadic majorants compatible with all the sums displayed above.  Since each majorant is already an $\ell^2_q\to\ell^2_n$ operator norm, the resulting estimate holds simultaneously for all deterministic inputs $w\in E_T^{2,\sigma}$ by the block inequality
\[
        \|P_Qw\|_{L_T^\infty L^2}\le Q^{-\sigma}\|w\|_{L_T^\infty\Btwo^\sigma}.
\]
The Borel--Cantelli step is performed at the level of dyadic operator norms.  For standard localized cutoffs, Lemma~\ref{lem:time-sup-centered-kernel} restricts every adjacent cutoff difference to $N\vee M\gtrsim\Lambda$.  The deterministic summations in Lemma~\ref{lem:centered-deterministic-summation-table} then leave a factor $\Lambda^{-\delta_{\rm op}+O(\eps)}$.  Applying the high-moment estimate and Borel--Cantelli to the adjacent dyadic pairs gives
\[
 \|\mathcal B_{2^{m+1}}^{i;j,k}-\mathcal B_{2^m}^{i;j,k}\|
 \le C_{\omega,T,\theta}2^{-m\theta},
 \qquad 0<\theta<\delta_{\rm op}/4.
\]
The geometric series and a telescoping sum give \eqref{eq:centered-cutoff-rate} for every dyadic $\Lambda'\ge\Lambda$.  For a general admissible cofinal family, first fix a finite set of dyadic triples.  Its kernels converge in $L^p(\Omega;C_{t,s}\mathcal L(\ell_q^2,\ell_n^2))$ by the factor $\tau_{\Lambda,\Lambda'}\to0$ in Lemma~\ref{lem:time-sup-centered-kernel}; the complement is uniform under the same summable envelope.  This proves the stated Cauchy convergence.
\end{proof}

\begin{proposition}[Local smallness of the centered Besov source operator]\label{prop:centered-besov-local-smallness}
For almost every enhanced realization produced by Theorem~\ref{thm:pathwise-fluct}, the centered Besov source operator localized to $[0,T]$ satisfies
\begin{equation}\label{eq:centered-besov-local-smallness}
        \|\mathcal B^{i;j,k}_{[0,T]}\|_{\mathcal L(E_T^{2,\sigma},L_T^1\Btwo^{\sigma-\alpha})}
        \longrightarrow 0,
        \qquad T\downarrow0.
\end{equation}
The same statement holds for cutoff differences after first taking the cutoff indices large enough in the augmented enhanced-data distance.
\end{proposition}

\begin{proof}
The dyadic estimate \eqref{eq:positive-centered-block-besov} contains an explicit factor $T$ and is dominated, after Lemma~\ref{lem:dyadic-borel-cantelli-operators}, by the summable deterministic envelopes used in the proof of Theorem~\ref{thm:pathwise-fluct}.  Fix $\delta>0$.  Choose a finite set of dyadic triples so that the contribution of its complement to the Besov operator norm is at most $\delta$, uniformly for $0<T\le1$.  On the finite set, the Volterra representation of Lemma~\ref{lem:low-frequency-besov-multiplier} gives a factor tending to zero with $T$.  Thus the full centered Besov operator norm is at most $2\delta$ for all sufficiently small $T$, which proves \eqref{eq:centered-besov-local-smallness}.  For cutoff differences, the same finite-set/tail decomposition is combined with the cutoff-tail majorants from Lemma~\ref{lem:time-sup-centered-kernel}.
\end{proof}

\begin{proof}[Proof of Theorem~\ref{thm:operator-main}]
Proposition~\ref{prop:det-contraction} gives the Volterra diagonal and Theorem~\ref{thm:pathwise-fluct} gives the centered remainder.  Their cutoff limits converge in both operator topologies in \eqref{eq:operator-bound-main}.
\end{proof}

\section{A separate weak-covariance extension of the mixed operators}\label{subsec:weak-covariance-extension}
This section is logically separate from the independent-color proof above.  It estimates the additional Fourier-diagonal covariance branches.  The finite covariance split is the same finite raw-kernel identity as in Subsection~\ref{subsec:finite-kernels}:
\[
        \E[\widehat\Psi_j(\ell,s)\widehat\Psi_k(r,t)]
        =\mathbf 1_{\ell+r=0}\Sigma_{jk}(\ell;s,t),
        \qquad
        T_{\Lambda}^{i;j,k}=D_{\mathsf R,\Lambda}^{i;j,k}+B_{\mathsf R,\Lambda}^{i;j,k}.
\]
Independent colors correspond to \(\Sigma_{jk}=\mathbf 1_{j=k}\sigma_j\).  For Fourier-diagonal weak covariance \(\Sigma_{jk}=\Sigma_{jk}^{\mathsf R}\), with \(\Sigma_{jk}^{\mathsf R}\) defined in \eqref{eq:wc-covkernel-def}.  The estimates below control the off-diagonal Volterra diagonal, the class-wise covariance-coordinate reduction of the centered remainder, and covariance continuity.

\begin{lemma}[Uniform local covariance bound]\label{lem:wc-local-covariance-bound}
For every Fourier label \(n\), the two-color covariance matrix \(\mathsf R(n)\) in Definition~\ref{def:wc-covariance} satisfies
\[
        |\mathsf R_{12}(n)|\le1,
        \qquad
        \|\mathsf R(n)\|_{\ell^2\to\ell^2}\le2.
\]
Consequently the finite coordinate factorization in Lemma~\ref{lem:wc-classwise-covariance-coordinates} is uniformly bounded on each local quotient Gaussian Hilbert space.
\end{lemma}

\begin{proof}
Positive semidefiniteness of the \(2\times2\) principal covariance matrix gives \(|\mathsf R_{12}(n)|^2\le \mathsf R_{11}(n)\mathsf R_{22}(n)=1\).  A Hermitian matrix with unit diagonal and off-diagonal entry of modulus at most one has eigenvalues \(1\pm |\mathsf R_{12}(n)|\), hence operator norm at most two.  If one eigenvalue is zero, the corresponding Gaussian direction is quotiented out before Wick products are formed; the quotient construction supplies the Gaussian coordinates directly.
\end{proof}

\begin{lemma}[Finite R-Wick split for one mixed dyadic block]\label{lem:wc-finite-wick-block}
Fix a Galerkin cutoff \(\Lambda\), a dyadic triple \((N,Q,M)\), and colors \(j,k\in\{1,2\}\).  At finite cutoff,
\begin{align}\label{eq:wc-finite-wick-block}
&c_\Lambda(\ell)c_\Lambda(r)\widehat\Psi_j(\ell,s)\widehat\Psi_k(r,t)\\
&\quad=c_\Lambda(\ell)c_\Lambda(r)\bigl(:\widehat\Psi_j(\ell,s)\widehat\Psi_k(r,t):\bigr)_{\mathsf R}
+c_\Lambda(\ell)c_\Lambda(r)\mathbf 1_{\ell+r=0}\Sigma_{jk}^{\mathsf R}(\ell;s,t).\notag
\end{align}
Consequently, after insertion into the raw finite kernel \(A_{\Lambda,N,Q,M}^{i;j,k}\), the deterministic part is supported on \(n=q\).  The centered part is a finite second homogeneous chaos in the joint two-color Gaussian space.  When \(\mathsf R_{jk}\equiv0\), the deterministic term vanishes and the block is centered, recovering the independent cross-color case.
\end{lemma}

\begin{proof}
The identity is the finite-dimensional Wick product formula in the Gaussian space generated by the cutoff Fourier Brownian motions.  The covariance is Fourier-diagonal by \eqref{eq:wc-covariance-def-main}; hence the non-conjugated pairing between the two stochastic convolution coefficients is supported on \(\ell+r=0\) and equals \(\Sigma_{jk}^{\mathsf R}(\ell;s,t)\).  Substituting this support relation into the convolution incidence \(n=q+\ell+r\) gives \(n=q\).  Every variable is finite-dimensional before the cutoff is removed, and the limiting object is obtained afterward in the displayed topology.
\end{proof}

For \(a,b\in\{1,2\}\), define the covariance kernel of two stochastic convolutions by
\begin{equation}\label{eq:wc-covkernel-def}
\Sigma_{ab}^{\mathsf R}(\ell;s,t)
:=\mathsf R_{ab}(\ell)\int_0^{s\wedge t}
        K_a(s-\tau,\ell)K_b(t-\tau,-\ell)\,d\tau.
\end{equation}
When \(a=b\), this is the same-color covariance already used in \eqref{eq:sigma-def}.  When \(a\ne b\), the two linear phases are distinct, and the covariance itself has a different-speed denominator.

\begin{lemma}[Off-diagonal stochastic-convolution covariance]\label{lem:wc-offdiag-covariance}
Assume Definition~\ref{def:wc-covariance} and let \(a\ne b\).  If \(|\ell|\sim N\), then for \(0\le s,t\le T\)
\begin{equation}\label{eq:wc-offdiag-cov-bound}
        |\Sigma_{ab}^{\mathsf R}(\ell;s,t)|
        \lesssim_T C_{\mathsf R}N^{-3\alpha-\kappa}.
\end{equation}
Moreover, for every \(0<\eta\le1/2\),
\begin{equation}\label{eq:wc-offdiag-cov-increment}
        |\Sigma_{ab}^{\mathsf R}(\ell;s,t)-\Sigma_{ab}^{\mathsf R}(\ell;s',t')|
        \lesssim_T C_{\mathsf R}(|s-s'|^\eta+|t-t'|^\eta)
        N^{-3\alpha-\kappa+\alpha\eta}.
\end{equation}
\end{lemma}

\begin{proof}
Write
\[
        K_a(s-\tau,\ell)K_b(t-\tau,-\ell)
        =\frac{\sin((s-\tau)\omega_a)\sin((t-\tau)\omega_b)}
        {\omega_a\omega_b},
        \qquad \omega_c=\omega_c(\ell).
\]
Using $2\sin A\sin B=\cos(A-B)-\cos(A+B)$, the integral over $0\le \tau\le s\wedge t$ is a sum of endpoint terms of the form
\[
        \frac{e^{i(s\omega_a\pm t\omega_b)}}{\omega_a\omega_b(\omega_a\pm\omega_b)}
        -
        \frac{e^{i((s-s\wedge t)\omega_a\pm (t-s\wedge t)\omega_b)}}
        {\omega_a\omega_b(\omega_a\pm\omega_b)},
\]
with harmless sign-dependent constants.  Since $a\ne b$ and the speeds are fixed and distinct,
\[
        |\omega_a(\ell)-\omega_b(\ell)|+|\omega_a(\ell)+\omega_b(\ell)|\gtrsim N^\alpha,
        \qquad \omega_a(\ell)\omega_b(\ell)\sim N^{2\alpha}.
\]
Thus the time integral is $O_T(N^{-3\alpha})$.  Multiplication by the covariance coefficient gives
\[
        |\Sigma_{ab}^{\mathsf R}(\ell;s,t)|
        \lesssim |\mathsf R_{ab}(\ell)|N^{-3\alpha}
        \lesssim C_{\mathsf R}N^{-3\alpha-\kappa},
\]
which proves \eqref{eq:wc-offdiag-cov-bound}.

For increments, it is useful to keep the same oscillatory formula.  Replacing $s$ by $s'$ changes either an endpoint exponential or the integration endpoint $s\wedge t$.  The exponential difference is bounded by
\[
        |e^{is\omega_a}-e^{is'\omega_a}|
        \lesssim |s-s'|^\eta N^{\alpha\eta},
        \qquad 0<\eta\le1,
\]
and the same estimate holds for the $t$ endpoint.  The factors $\omega_a^{-1}\omega_b^{-1}$ and $(\omega_a\pm\omega_b)^{-1}$ are unchanged.  Hence the endpoint part is bounded by
\[
        C_{\mathsf R}(|s-s'|^\eta+|t-t'|^\eta)
        N^{-3\alpha-\kappa+\alpha\eta}.
\]
The short interval created when $s\wedge t$ changes is controlled by the same interpolation estimate applied to the sine kernels; it has the same or better dyadic size.  This proves \eqref{eq:wc-offdiag-cov-increment}.
\end{proof}

For every generated block \((i;j,k)\in\mathfrak M\), define the covariance diagonal
\begin{equation}\label{eq:wc-diagonal-kernel}
D_{\mathsf R,N,Q}^{i;j,k}(q;t,s)
=\chi_Q(q)\sum_\ell \rho_N(\ell)\rho_N(-\ell)
\chi^{\rm res}(q+\ell,-\ell)K_i(t-s,q+\ell)
\Sigma_{jk}^{\mathsf R}(\ell;s,t).
\end{equation}
It acts by
\begin{equation}\label{eq:wc-diagonal-operator}
\widehat{\mathcal D_{\mathsf R,N,Q}^{i;j,k}w}(q,t)
=\int_0^tD_{\mathsf R,N,Q}^{i;j,k}(q;t,s)\widehat w(q,s)\,ds.
\end{equation}
For \(j=k\) this recovers the same-color branch of Subsection~\ref{subsec:same-color-det}, with the notation \(\mathcal D^{i;j}=\mathcal D_{\mathsf R}^{i;j,j}\).  For \(j\ne k\) it is the off-diagonal covariance-weighted cross-color diagonal.

\begin{lemma}[Off-diagonal covariance diagonal bound]\label{lem:wc-offdiag-diagonal-bound}
Let \(j\ne k\).  For each dyadic pair \(Q\le c_0N\),
\begin{equation}\label{eq:wc-offdiag-diag-shell}
        \|\mathcal D_{\mathsf R,N,Q}^{i;j,k}P_Qw\|_{C_TH^{s_2-\alpha}}
        \lesssim_{T,\eps} C_{\mathsf R}
        N^{3-4\alpha-\kappa+\eps}Q^{s_2-\alpha+\eps}
        \|P_Qw\|_{L_T^\infty L_x^2},
\end{equation}
and
\begin{equation}\label{eq:wc-offdiag-diag-besov-shell}
        \|\mathcal D_{\mathsf R,N,Q}^{i;j,k}P_Qw\|_{L_T^\infty B_{2,\infty}^{\sigma-\alpha}}
        \lesssim_{T,\eps} C_{\mathsf R}
        N^{3-4\alpha-\kappa+\eps}Q^{-\alpha}
        \|w\|_{L_T^\infty B_{2,\infty}^{\sigma}}.
\end{equation}
Consequently, for \(\alpha>3/4\) and \(\kappa\ge0\), the off-diagonal mixed covariance diagonals are summable in the same source topologies as the independent-color same-color diagonals.
\end{lemma}

\begin{proof}
The contraction has already imposed \(n=q\).  The high-loop sum has \(N^3\) choices, \(K_i(t-s,q+\ell)\) contributes \(N^{-\alpha}\), and Lemma~\ref{lem:wc-offdiag-covariance} contributes \(C_{\mathsf R}N^{-3\alpha-\kappa}\).  This gives the scalar multiplier size \(C_{\mathsf R}N^{3-4\alpha-\kappa}\).  Multiplication by the output Sobolev weight on the diagonal low frequency gives the factor \(Q^{s_2-\alpha+\eps}\) in \eqref{eq:wc-offdiag-diag-shell}.  For the Besov source, use
\[
        Q^{\sigma-\alpha}\|P_Qw\|_{L^2}
        \le Q^{-\alpha}\|w\|_{B_{2,\infty}^{\sigma}}.
\]
The dyadic sum over \(N\gtrsim Q\) is bounded because \(3-4\alpha-\kappa<0\) in the present range.  This off-diagonal branch is controlled directly from the position envelope of \(E_T^{2,\sigma}\).
\end{proof}

\begin{proposition}[Covariance-weighted deterministic mixed diagonals]\label{prop:wc-deterministic-diagonals}
Under Definition~\ref{def:wc-covariance}, all deterministic covariance diagonals satisfy
\begin{equation}\label{eq:wc-diag-global-bound}
        \mathcal D_{\mathsf R}^{i;j,k}:
        E_T^{2,\sigma}\longrightarrow C_TH^{s_2-\alpha}\cap L_T^1B_{2,\infty}^{\sigma-\alpha}
\end{equation}
with cutoff convergence in the same operator topology.  The same-color part is Proposition~\ref{prop:det-contraction}; the off-diagonal part is Lemma~\ref{lem:wc-offdiag-diagonal-bound} and the finite-set/tail argument used in Proposition~\ref{prop:det-contraction}.
\end{proposition}

\begin{proof}
For \(j=k\), Proposition~\ref{prop:det-contraction} applies.  Let \(j\ne k\).  Lemma~\ref{lem:wc-offdiag-diagonal-bound} gives, for the Sobolev target, the dyadic majorant
\[
        N^{3-4\alpha-\kappa+\eps}Q^{s_2-\alpha+\eps}\|P_Qw\|_{L_T^\infty L_x^2}.
\]
Because the admissible window gives \(s_2<4\alpha-3\le \alpha\) in the subunit branch and \(s_2<\alpha\) in the conic endpoint branch, the exponent of \(Q\) is negative after fixing the small loss.  Hence the low-input summation is controlled by \(\|w\|_{C_TL^2}\), while the high-shell summation is controlled by \(\sum_N N^{3-4\alpha-\kappa+\eps}<\infty\).  Equivalently, summing first in \(N\ge Q\) gives a remaining low-input power bounded by a negative dyadic exponent under \(s_2<4\alpha-3\).

For the direct Besov source, fix the output block \(Q\).  The high-shell tail satisfies
\[
        \sum_{N\gtrsim Q}N^{3-4\alpha-\kappa+\eps}Q^{-\alpha}
        \lesssim Q^{3-5\alpha-\kappa+\eps},
\]
which is uniformly bounded for \(\alpha>3/4\) after reducing \(\eps\).  This proves the \(L_T^\infty B_{2,\infty}^{\sigma-\alpha}\) bound for the off-diagonal covariance diagonals.  The passage from \(L_T^\infty B_{2,\infty}^{\sigma-\alpha}\) to \(L_T^1B_{2,\infty}^{\sigma-\alpha}\) is deterministic multiplication by \(T\), exactly as in Proposition~\ref{prop:det-diagonal-local-smallness}.

Cutoff convergence follows from the same finite-set/tail argument: fixed dyadic blocks are finite Fourier sums and converge entrywise under the Galerkin cutoff, while the displayed summable envelopes control the complement uniformly in the cutoff.
\end{proof}

\begin{proposition}[Localized smallness of covariance-weighted diagonals]\label{prop:wc-diagonal-local-smallness}
For every generated block, the deterministic covariance diagonal restricted to \([0,T]\) satisfies
\begin{equation}\label{eq:wc-diagonal-local-smallness}
        \|\mathcal D_{\mathsf R,[0,T]}^{i;j,k}\|_{E_T^{2,\sigma}\to L_T^1B_{2,\infty}^{\sigma-\alpha}}
        \le T\,
        \|\mathcal D_{\mathsf R,[0,T]}^{i;j,k}\|_{E_T^{2,\sigma}\to L_T^\infty B_{2,\infty}^{\sigma-\alpha}}.
\end{equation}
In particular, on every bounded covariance-enhanced-data ball, the diagonal Besov-source contribution tends to zero as \(T\downarrow0\).  For \(j\ne k\), the stronger \(L_T^\infty B_{2,\infty}^{\sigma-\alpha}\) norm is controlled by the summable envelope in Lemma~\ref{lem:wc-offdiag-diagonal-bound}; for \(j=k\) it is the independent-color same-color Volterra diagonal bound.
\end{proposition}

\begin{proof}
The estimate is the deterministic embedding \(L_T^\infty X\hookrightarrow L_T^1X\) with \(X=B_{2,\infty}^{\sigma-\alpha}\).  The uniform boundedness of the \(L_T^\infty\) operator norm is supplied by Proposition~\ref{prop:wc-deterministic-diagonals}.  The smallness comes entirely from the deterministic time factor.
\end{proof}

For Fourier-diagonal weak covariance, the centered kernel is defined by subtracting every covariance block:
\begin{equation}\label{eq:wc-centered-kernel-def}
B_{\mathsf R,N,Q,M}^{i;j,k}(n,q;t,s)
= A_{N,Q,M}^{i;j,k}(n,q;t,s)
-\mathbf 1_{n=q}\chi_M(n)D_{\mathsf R,N,Q}^{i;j,k}(q;t,s).
\end{equation}
Here \(D_{\mathsf R}^{i;j,k}\) is zero when \(\mathsf R_{jk}\equiv0\), so the formula reduces to the independent definition.

\begin{lemma}[Covariance coordinates in each Fourier involution class]\label{lem:wc-classwise-covariance-coordinates}
Fix a Galerkin cutoff, dyadic scales, and the finitely many time parameters appearing in one fixed-time, increment, or cutoff-tail block.  For each real Fourier involution class \([m]=\{m,-m\}\), let \(\mathfrak h_{[m]}^{\mathsf R}\) be the finite color-time Gaussian Hilbert space generated by the two colors and by the stochastic convolution kernels attached to that class.  Then the joint Gaussian space decomposes as the orthogonal direct sum
\[
        \bigoplus_{[m]}\mathfrak h_{[m]}^{\mathsf R}.
\]
Collect the variables in one summand into a finite real Gaussian vector
$X_{[m]}$ with covariance matrix $\Sigma_{[m]}$.  There are a standard real
Gaussian vector $G_{\mathsf R,[m]}$ and a deterministic Gram factor
$A_{[m]}$ such that
\[
        X_{[m]}=A_{[m]}G_{\mathsf R,[m]},
        \qquad A_{[m]}A_{[m]}^{\mathsf T}=\Sigma_{[m]}.
\]
One may take the positive semidefinite square root, or equivalently choose an
orthonormal basis after quotienting out the null space of the local covariance
form.  This factorization may mix colors and the finitely many time-kernel
coordinates belonging to the same \([m]\), but it keeps different Fourier
classes separated and uses no inverse covariance matrix.  Consequently every
centered coefficient tensor obtained after subtracting the covariance pairing
is still supported by
\begin{equation}\label{eq:wc-classwise-support}
        n=q+\ell+r,
        \qquad |\ell|\sim |r|\sim N,
        \qquad |q|\sim Q,
        \qquad |n|\sim M,
\end{equation}
and the additional color/time index has bounded multiplicity determined by the finite list of kernels in the block, uniformly in \(N,Q,M\) and in the Galerkin cutoff.  More precisely, at most four time kernels occur in a two-time comparison, there are two colors, and the real sine--cosine convention contributes at most two coordinates.  Thus one may take
\[
        d_*\le 16
\]
simultaneously for fixed-time kernels, one- and two-time increments, and cutoff differences.  More explicitly, for $m_0\in\mathbb Z^3_{\mathrm{rep}}$ the local coordinates may be denoted by $G_{\mathsf R,[m_0],\mu}$ and
\begin{equation}\label{eq:wc-signed-class-coordinate-expansion}
        \widehat\Psi_j^{\mathsf R}(\varepsilon m_0,\tau)
        =\sum_{\mu\le d_*}
          c_{j,\varepsilon,m_0,\mu}^{\mathsf R}(\tau)
          G_{\mathsf R,[m_0],\mu}.
\end{equation}
The color and the sign are carried by the deterministic coefficients; there is only one Gaussian family for the class $[m_0]=[-m_0]$.  The bound $d_*\le16$ is comparison-wise in the sense of Remark~\ref{rem:comparison-wise-gram-coordinates}: it is uniform for each pairwise increment used in the continuity proof, which proceeds without a simultaneous factorization of an arbitrarily large time grid.
\end{lemma}

\begin{proof}
The support condition \(n+m=0\) in \eqref{eq:wc-covariance-def-main} implies that Gaussian variables carried by different real Fourier involution classes are orthogonal.  Therefore the covariance matrix of all variables entering a finite dyadic block is block diagonal with respect to the classes \([m]\).  Factoring each positive semidefinite block after removing zero-norm directions acts separately on each class.  It changes only the local Gaussian coordinate basis and preserves the spatial class carried by each coefficient.  After the scalar covariance pairing \(\mathbf 1_{\ell+r=0}\Sigma_{jk}^{\mathsf R}(\ell;s,t)\) has been subtracted, the remaining second homogeneous chaos is written in these local standard coordinates.  The deterministic Fourier coefficient still comes from the original convolution law of \(I_i(w<\Psi_j)\circ\Psi_k\), hence the signed incidence \(n=q+\ell+r\) remains.  Formula \eqref{eq:wc-signed-class-coordinate-expansion} is this local Gram representation written with the transversal convention \eqref{eq:fourier-involution-representative}; it does not duplicate the Gaussian family when the sign changes.  The count $4\times2\times2=16$ gives the displayed common bound for the time, color, and real Fourier coordinates in one pairwise comparison; cutoff multipliers change coefficients but add no kernel.  Lemma~\ref{lem:auxiliary-index-flattening} absorbs this fixed multiplicity.  The coordinates may be re-chosen for another comparison, exactly as explained in Remark~\ref{rem:comparison-wise-gram-coordinates}.
\end{proof}

\begin{proposition}[Covariance-local tensor reduction]\label{prop:wc-covariance-local-tensor}
Fix a finite cutoff, a dyadic triple \((N,Q,M)\), and times \((t,s)\).  After the deterministic term in Lemma~\ref{lem:wc-finite-wick-block} is subtracted, the centered covariance block admits a representation
\begin{equation}\label{eq:wc-covariance-local-tensor}
        B_{\mathsf R,N,Q,M}^{i;j,k}(n,q;t,s)
        =\sum_{\nu\in\mathcal P_{\mathsf R}}
        \sum_{\ell_0,r_0\in\mathbb Z^3_{\mathrm{rep}}}
        \sum_{\mu,\lambda\le d_{\mathsf R}}
        H_{\nu}^{\mathsf R}(\ell_0,\mu,r_0,\lambda,q,n;t,s)
        :G_{\mathsf R,[\ell_0],\mu}G_{\mathsf R,[r_0],\lambda}:.
\end{equation}
Here \(\mathcal P_{\mathsf R}\) is a fixed finite set of real-Fourier sign, color, and coordinate patterns, and one may take \(d_{\mathsf R}\le16\), independently of \(N,Q,M,\Lambda,t,s\).  If $\varepsilon_{\nu,1},\varepsilon_{\nu,2}$ are the signs in one pattern, its coefficient tensor is supported by
\begin{equation}\label{eq:wc-covariance-local-support}
        |\ell_0|\sim |r_0|\sim N,
        \qquad |q|\sim Q,
        \qquad |n|\sim M,
        \qquad n=q+\varepsilon_{\nu,1}\ell_0
                     +\varepsilon_{\nu,2}r_0,
\end{equation}
with no implicit sign replacement.  Moreover
\begin{equation}\label{eq:wc-covariance-local-size}
        |H_{\nu}^{\mathsf R}(\ell_0,\mu,r_0,\lambda,q,n;t,s)|
        \lesssim_T \|\mathsf R\|_{\kappa}^{O(1)}N^{-3\alpha}.
\end{equation}
Thus the flattening estimates of Lemma~\ref{lem:dyadic-tensor-flattening} have the same dyadic powers as in the independent-color centered component.
\end{proposition}

\begin{proof}
Everything is finite-dimensional at fixed cutoff.  For each nonzero Fourier involution class $\{m,-m\}$, the two-color real covariance space is the quotient Hilbert space described in \eqref{eq:real-fourier-quotient-space}, tensored with the finite span of the time kernels which occur in the chosen block.  Lemma~\ref{lem:wc-classwise-covariance-coordinates} chooses standard coordinates for this local space.  The factorization is block diagonal in the class label $[m]$.  The sign and the original color remain in the deterministic coefficient in \eqref{eq:wc-signed-class-coordinate-expansion}; the coordinate factor may replace the color-time index by a bounded auxiliary index but cannot alter the convolution relation.

After the finite $\mathsf R$-Wick diagonal has been subtracted, every remaining term is a second homogeneous chaos in two local Gaussian coordinates.  For one sign/color/coordinate pattern it has the form
\[
        \sum_{\ell_0,r_0\in\mathbb Z^3_{\mathrm{rep}}}
        \sum_{\mu,\lambda\le d_{\mathsf R}}
        H_{\nu}^{\mathsf R}(\ell_0,\mu,r_0,\lambda,q,n;t,s)
        :G_{\mathsf R,[\ell_0],\mu}G_{\mathsf R,[r_0],\lambda}: .
\]
The deterministic coefficient is the product of the same Fourier multiplier appearing in the independent centered block and bounded matrix coefficients coming from the local covariance metric.  The Duhamel kernel contributes $N^{-\alpha}$, the two stochastic-convolution kernels contribute $N^{-2\alpha}$, and the Gram-factor coefficients are bounded by a constant depending polynomially on the local covariance norm.  Lemma~\ref{lem:wc-local-covariance-bound} gives a uniform bound for this norm; the decay norm $\|\mathsf R\|_\kappa$ enters only through the lower-chaos covariance branches and through covariance-continuity estimates.  Consequently
\[
        |H_{\nu}^{\mathsf R}(\ell_0,\mu,r_0,\lambda,q,n;t,s)|
        \lesssim_T \|\mathsf R\|_\kappa^{O(1)}N^{-3\alpha}.
\]

The support relation is unchanged: for a fixed pattern the Fourier coefficient of $I_i(w<\Psi_j)\circ\Psi_k$ contains exactly the incidence
\[
        n=q+\varepsilon_{\nu,1}\ell_0
            +\varepsilon_{\nu,2}r_0,
        \qquad |\ell_0|\sim |r_0|\sim N,
        \qquad |q|\sim Q,
        \qquad |n|\sim M.
\]
The auxiliary indices have bounded range independent of all dyadic scales.  Lemma~\ref{lem:auxiliary-index-flattening} therefore reduces every flattening to the corresponding signed spatial incidence, and Lemma~\ref{lem:dyadic-tensor-flattening} gives the same row/column profile as in the independent centered case.
\end{proof}

\begin{proposition}[Centered operator component with weak color correlation]\label{prop:wc-centered-component}
Assume Definition~\ref{def:wc-covariance}.  After subtracting the covariance diagonals in \eqref{eq:wc-centered-kernel-def}, the centered operators satisfy
\begin{equation}\label{eq:wc-centered-op-bound}
        \mathcal B_{\mathsf R}^{i;j,k}:
        E_T^{2,\sigma}\longrightarrow C_TH^{s_2-\alpha}\cap L_T^1B_{2,\infty}^{\sigma-\alpha}
\end{equation}
almost surely, with the same dyadic profile
\begin{equation}\label{eq:wc-centered-profile}
        N^{3/2-3\alpha+}(M^{3/2+}+Q^{3/2+})Q^{-\sigma}.
\end{equation}
The cutoff operators are Cauchy in the same operator topologies.
\end{proposition}

\begin{proof}
We follow the independent centered proof and record the points where covariance enters.  At finite cutoff, Lemma~\ref{lem:wc-finite-wick-block} gives the exact decomposition into covariance diagonals and a centered second chaos.  Proposition~\ref{prop:wc-covariance-local-tensor} writes the centered kernel as a finite sum of tensors satisfying the same support relation
\[
        n=q+\ell+r,
        \qquad |\ell|\sim |r|\sim N,
        \qquad |q|\sim Q,
        \qquad |n|\sim M,
\]
and the same coefficient size $N^{-3\alpha}$, up to constants depending polynomially on the covariance norm.  The covariance-coordinate factorization is local in the Fourier involution classes, so the auxiliary labels have bounded multiplicity and do not affect the row/column counts.

For a fixed dyadic triple and fixed times, the proof of Proposition~\ref{prop:fixed-time-centered-kernel} applies to each covariance-local pattern.  The oriented-flattening profile is bounded by
\[
        \mathfrak F_4(H_{\nu}^{\mathsf R})
        \lesssim_T \|\mathsf R\|_{\kappa}^{O(1)}
        N^{\frac32-3\alpha}(M^{\frac32}+Q^{\frac32}),
\]
by Lemmas~\ref{lem:auxiliary-index-flattening} and~\ref{lem:dyadic-tensor-flattening}.  Proposition~\ref{prop:abstract-random-tensor}, together with decoupling for the same-color Wick-square patterns, gives the fixed-time matrix estimate
\[
 \bigl\|\|\mathfrak B_{\mathsf R,N,Q,M}^{i;j,k}(t,s)\|_{\ell_q^2\to\ell_n^2}\bigr\|_{L^p(\Omega)}
 \lesssim_{p,T,\eps,\mathsf R}
 N^{\frac32-3\alpha+\eps}(M^{\frac32+\eps}+Q^{\frac32+\eps}).
\]
The time-increment tensors are obtained by replacing one of the time kernels by a difference kernel before the local Gram factorization.  Lemma~\ref{lem:centered-increment-coefficients} therefore remains valid with the same additional factor $N^{\alpha\eta}$.  Lemma~\ref{lem:operator-valued-time-lift} upgrades the fixed-time estimate to the $C_{t,s}\mathcal L(\ell_q^2,\ell_n^2)$ estimate, and Lemma~\ref{lem:dyadic-borel-cantelli-operators} supplies a pathwise dyadic majorant.

After the pathwise majorant is fixed, the input is inserted only through
\[
        \|P_Qw\|_{L_T^\infty L^2}\le Q^{-\sigma}\|w\|_{L_T^\infty B_{2,\infty}^{\sigma}}.
\]
The deterministic summations are exactly those in Theorem~\ref{thm:pathwise-fluct}: the $M^{3/2}$ profile gives the condition $s_2<4\alpha-3-$, the $Q^{3/2}$ profile gives $\sigma>3-3\alpha+$, and the direct Besov source summation gives $\sigma<4\alpha-3-$.  This proves boundedness in
\[
        \mathcal L(E_T^{2,\sigma},C_TH^{s_2-\alpha})
        \cap
        \mathcal L(E_T^{2,\sigma},L_T^1B_{2,\infty}^{\sigma-\alpha}).
\]
For cutoff convergence, fixed dyadic blocks converge because the covariance-local kernels are finite sums whose cutoff multipliers converge entrywise.  The cutoff-tail tensors obey the same dyadic majorant, with the tail on a Gaussian leg or on the output leg.  A finite-dyadic-set plus summable-tail argument gives the Cauchy property in both operator topologies.
\end{proof}

\begin{theorem}[Weak-covariance mixed-operator estimate]\label{thm:wc-mixed-operators}
Assume Definition~\ref{def:wc-covariance}, \(\kappa>3-3\alpha+10\eps\), and the parameter window \eqref{eq:operator-cond-main}.  Then each generated mixed operator admits the covariance split
\begin{equation}\label{eq:wc-mixed-decomposition}
        T_{\mathsf R}^{i;j,k}=\mathcal D_{\mathsf R}^{i;j,k}+\mathcal B_{\mathsf R}^{i;j,k},
\end{equation}
where the deterministic covariance diagonal and centered remainder obey the same source and cutoff-convergence bounds as in Theorem~\ref{thm:operator-main}.
\end{theorem}

\begin{proof}
At finite cutoff the starting point is the exact $\mathsf R$-Wick identity of Lemma~\ref{lem:wc-finite-wick-block}.  It splits each raw mixed block into the deterministic covariance diagonal supported on $\ell+r=0$ and a centered second-chaos remainder.  Thus there are only two analytic tasks.

The deterministic task is Proposition~\ref{prop:wc-deterministic-diagonals}.  For $j=k$ it is the same-color Volterra diagonal from the independent-color model, whose phase denominator is supplied by the proportional low--high gap.  For $j\ne k$ the covariance kernel itself carries the off-diagonal decay and the different-speed covariance denominator, giving the shell factor $N^{3-4\alpha-\kappa+}$; this is summable in the same source spaces and is locally small after the deterministic $L_T^\infty\to L_T^1$ conversion.

The centered task is Proposition~\ref{prop:wc-centered-component}.  After subtracting the covariance diagonals, the class-wise covariance-coordinate factorization preserves the spatial labels.  The centered kernel therefore has the same incidence $n=q+\ell+r$ and the same row/column tensor profile as in the independent theorem.  The Besov-input dyadic summation is consequently identical to Theorem~\ref{thm:pathwise-fluct}.

Combining these two components gives the decomposition \eqref{eq:wc-mixed-decomposition}, the two source-operator bounds, and cutoff Cauchy convergence.  The lower bound $\kappa>3-3\alpha+10\eps$ is stronger than the mixed-operator and cubic off-diagonal summations require; it is imposed by the quadratic zero-mode covariance branch in Appendix~\ref{app:weak-covariance-symbols}.  Using this common assumption keeps the full covariance-enhanced datum in one admissible class.
\end{proof}

\begin{lemma}[Local covariance square roots preserve Fourier labels]\label{lem:wc-square-root-locality}
Fix $K\ge1$ and let $\mathsf R,\widetilde{\mathsf R}$ satisfy
Definition~\ref{def:wc-covariance} with
\[
        \|\mathsf R\|_\kappa+
        \|\widetilde{\mathsf R}\|_\kappa\le K.
\]
Let $\mathsf C_{\mathsf R}(n)$ and
$\mathsf C_{\widetilde{\mathsf R}}(n)$ be the real sine--cosine covariance
matrices on the involution class $[n]=\{n,-n\}$.  There is $N_K<\infty$ such
that, for $|n|\ge N_K$,
\begin{equation}\label{eq:wc-high-square-root-lipschitz}
 \bigl\|\mathsf C_{\mathsf R}(n)^{1/2}
       -\mathsf C_{\widetilde{\mathsf R}}(n)^{1/2}\bigr\|_{\rm op}
 \lesssim_K
 \bigl\|\mathsf C_{\mathsf R}(n)
       -\mathsf C_{\widetilde{\mathsf R}}(n)\bigr\|_{\rm op}.
\end{equation}
For every class, including the possibly rank-deficient low classes, one has the
operator-H\"older bound
\begin{equation}\label{eq:wc-low-square-root-modulus}
 \bigl\|\mathsf C_{\mathsf R}(n)^{1/2}
       -\mathsf C_{\widetilde{\mathsf R}}(n)^{1/2}\bigr\|_{\rm op}
 \lesssim
       \bigl\|\mathsf C_{\mathsf R}(n)
       -\mathsf C_{\widetilde{\mathsf R}}(n)\bigr\|_{\rm op}^{1/2}.
\end{equation}
Both square roots act only inside the fixed class $[n]$ and its two color
coordinates.  Thus the common square-root coupling mixes color and
sine--cosine coordinates but never mixes distinct Fourier labels.
\end{lemma}

\begin{proof}
In the normalization \eqref{eq:real-fourier-covariance-matrix},
\[
        \mathsf C_{\mathsf R}(n)=\frac12 I_4+E_{\mathsf R}(n),
        \qquad
        \|E_{\mathsf R}(n)\|_{\rm op}
        \lesssim K\langle n\rangle^{-\kappa};
\]
the zero mode has the analogous two-dimensional form.  Choose $N_K$ so that
the error is at most $1/4$ for $|n|\ge N_K$.  The spectra of both local
covariance matrices then lie in $[1/4,3/4]$.  On this set the positive square-root
map is Lipschitz.  One way to see this is to write the difference of the two
squares as the Sylvester equation
\[
 A-B=A^{1/2}(A^{1/2}-B^{1/2})
       +(A^{1/2}-B^{1/2})B^{1/2};
\]
the Sylvester operator has a uniformly bounded inverse when the spectra of
$A^{1/2}$ and $B^{1/2}$ are bounded below.  This proves
\eqref{eq:wc-high-square-root-lipschitz}.

For arbitrary positive semidefinite matrices the square-root map is
operator-H\"older of order $1/2$:
\[
        \|A^{1/2}-B^{1/2}\|_{\rm op}
        \lesssim \|A-B\|_{\rm op}^{1/2}.
\]
Indeed, replace $A,B$ by $A+\delta I,B+\delta I$.  The preceding Sylvester
argument bounds the difference of their square roots by
$\|A-B\|_{\rm op}/(2\sqrt\delta)$, whereas functional calculus gives
$\|(A+\delta I)^{1/2}-A^{1/2}\|_{\rm op}\le\sqrt\delta$, and similarly for
$B$.  If $A\ne B$, choosing $\delta=\|A-B\|_{\rm op}$ proves the displayed
bound; the case $A=B$ is immediate.  Thus
\eqref{eq:wc-low-square-root-modulus} holds without a spectral gap, including at
rank-deficient matrices.  Finally, each matrix is defined on one
real Fourier involution class.  Its square root has the same block support, so
the Fourier label is unchanged.
\end{proof}

\begin{proposition}[Continuity with respect to the off-diagonal covariance]\label{prop:wc-covariance-continuity}
Let \(\mathsf R\) and \(\widetilde{\mathsf R}\) satisfy Definition~\ref{def:wc-covariance} with the same exponent \(\kappa\), and define
\begin{equation}\label{eq:wc-covariance-distance}
        d_{\kappa}(\mathsf R,\widetilde{\mathsf R})
        :=\sup_{n\in\mathbb Z^3}\langle n\rangle^\kappa
        \bigl(|\mathsf R_{12}(n)-\widetilde{\mathsf R}_{12}(n)|
        +|\mathsf R_{21}(n)-\widetilde{\mathsf R}_{21}(n)|\bigr).
\end{equation}
Then the covariance diagonals satisfy the bounds of Proposition~\ref{prop:wc-deterministic-diagonals}, with an additional factor \(d_\kappa(\mathsf R,\widetilde{\mathsf R})\) for their difference.  Under the common realization \eqref{eq:canonical-covariance-coupling}, every centered kernel on a fixed finite dyadic set is continuous in $L^p(\Omega)$ with respect to the local covariance matrices.  The full centered operators are therefore continuous in the enhanced-data topology by the uniform dyadic tail bound of Proposition~\ref{prop:wc-centered-component}.
\end{proposition}

\begin{proof}
For deterministic covariance diagonals one replaces \(\mathsf R_{jk}\) by \(\mathsf R_{jk}-\widetilde{\mathsf R}_{jk}\) in Lemma~\ref{lem:wc-offdiag-covariance}.  The same calculation as in Lemma~\ref{lem:wc-offdiag-diagonal-bound} gives the factor \(d_\kappa(\mathsf R,\widetilde{\mathsf R})\) and the same summable shell power.

For the centered remainders, fix a finite cutoff and a finite dyadic set.  By \eqref{eq:canonical-covariance-coupling}, all relevant Gaussian vectors are linear images of the same finite standard Gaussian vector.  Lemma~\ref{lem:wc-square-root-locality} gives a uniform Lipschitz estimate on the high Fourier classes and the explicit $1/2$-H\"older bound \eqref{eq:wc-low-square-root-modulus} on every remaining class.  It also shows that the square-root coupling is block diagonal in the Fourier involution labels.  Hence the Gaussian coordinates, and therefore every fixed-chaos polynomial kernel, converge almost surely and in every finite $L^p$ norm without changing the incidence relations used in the tensor estimates.  Rank-deficient low-frequency blocks are covered directly by \eqref{eq:wc-low-square-root-modulus}; no inverse covariance matrix is introduced.

It remains to pass from finite dyadic sets to the full operator.  Proposition~\ref{prop:wc-centered-component} supplies a dyadic majorant independent of the approximation index and summable in the operator topologies \(E_T^{2,\sigma}\to C_TH^{s_2-\alpha}\) and \(E_T^{2,\sigma}\to L_T^1B_{2,\infty}^{\sigma-\alpha}\).  Given an error tolerance, choose the dyadic tail small uniformly, then use the finite-dimensional continuity on the remaining finite set.  This finite-set/tail argument proves the stated stability.
\end{proof}

\section{Deterministic solution map from enhanced data}\label{sec:deterministic-closure}

We solve the paracontrolled system in the augmented spaces defined below.  Stochastic terms are estimated in \(C_TH^{s_2-\alpha}\cap L_T^1B_{2,\infty}^{\sigma-\alpha}\); deterministic quadratic terms use the fractional or conic linear estimates.  Throughout, \(T^\delta\) denotes a positive power depending only on the fixed admissible parameters.

The homogeneous Klein--Gordon semigroup is
\begin{equation}\label{eq:T-semigroup}
S_i(t)(f,g)=\cos(t\omega_i(D))f+
\frac{\sin(t\omega_i(D))}{\omega_i(D)}g.
\end{equation}
The stochastic enhanced objects are constructed with zero deterministic data.  The stochastic construction and the mixed-operator estimates are formulated for $0<\alpha\le1$ wherever the stated shell bounds apply.  The final theorem below is written first in the subunit regime because Proposition~\ref{prop:strichartz} is derived from the full-rank stationary-phase theorem proved in Appendix~\ref{app:fracKG-strichartz}.  The classical endpoint $\alpha=1$ is recorded separately in Proposition~\ref{prop:gko-conic-strichartz} and Corollary~\ref{cor:alpha-one-extension}; it changes the deterministic linear input while keeping the same stochastic enhanced data.

Prescribed deterministic Cauchy data enter only through the affine term
\begin{equation}\label{eq:Y-initial-data}
  (y_{i,0},y_{i,1})\in
  \bigl(H^{s_2}\cap\Btwo^\sigma\bigr)\times
  \bigl(H^{s_2-\alpha}\cap\Btwo^{\sigma-\alpha}\bigr),
  \qquad i=1,2.
\end{equation}
The paracontrolled unknowns satisfy
\begin{equation}\label{eq:XY-initial-data}
  X_i(0)=\partial_tX_i(0)=0,
  \qquad
  (Y_i(0),\partial_tY_i(0))=(y_{i,0},y_{i,1}),
\end{equation}
and the full low-frequency input is always
\begin{equation}\label{eq:full-input-no-Z}
  W_i=V_i+X_i+Y_i.
\end{equation}

\subsection{Fractional Strichartz input and dyadic Besov propagation}
The closure uses local estimates for
\[
 L_i=\partial_t^2+\omega_i(D)^2,
 \qquad
 \omega_i(n)=(1+c_i^2|n|^{2\alpha})^{1/2}.
\]
We separate the Sobolev--Strichartz input from elementary dyadic Besov propagation.  In the subunit branch $0<\alpha<1$, Proposition~\ref{prop:strichartz} is proved from the exact massive fractional Klein--Gordon phase in Appendix~\ref{app:fracKG-strichartz}; the source line used below is the retarded $L_T^1H^{s-\alpha}$ line.  At $\alpha=1$, Proposition~\ref{prop:gko-conic-strichartz} replaces only the deterministic linear input and provides the conic source routes used for the endpoint quadratic sector.  The stochastic enhanced data, Volterra diagonal, centered operator theorem, and elementary $B_{2,\infty}$ propagation are unchanged.

\medskip
\noindent\emph{Case F: the subunit fractional case.}
Fix a Strichartz loss $0<\eps_{\rm Str}\ll\eps$, where $\eps$ is the fixed small margin in Definition~\ref{def:admissible}.  If $12/13<\alpha<1$, choose a fractional Klein--Gordon non-endpoint pair $(p_X,q_X)$ with
\begin{equation}\label{eq:pxqx-admissible}
 p_X>4,
 \qquad q_X>3,
 \qquad \frac2{p_X}+\frac3{q_X}<\frac32,
\end{equation}
close to $(4,3)$, and set
\begin{equation}\label{eq:kappaX-fractional}
 \kappa_X=\frac32-\frac3{q_X}-\frac{\alpha}{p_X}+\eps_{\rm Str}.
\end{equation}
For the $Y$ component set $(p_Y,q_Y)=(4,4)$ and
\begin{equation}\label{eq:kappaY-fractional}
 \kappa_Y=\frac{3-\alpha}{4}+\eps_{\rm Str}.
\end{equation}

\medskip
\noindent\emph{Case W: the classical conic endpoint.}
At the endpoint $\alpha=1$ we insert the classical conic wave/Klein--Gordon Strichartz estimates on $\mathbb T^3$ used in the Gubinelli--Koch--Oh setting.  The stochastic lift, the Volterra diagonal, the centered operator theorem, and the $B_{2,\infty}$ propagation line are the same enhanced-data components.  In this endpoint case we set
\begin{equation}\label{eq:case-W-exponents}
        (p_X,q_X)=\left(8,\frac83\right),
        \qquad \kappa_X=\frac14,
        \qquad \kappa_Y=\frac12.
\end{equation}
The corresponding endpoint norms are therefore
\begin{align}\label{eq:case-W-spaces-display}
\|X\|_{X^{s_1}_{T,\mathrm{end}}}
&:=\|X\|_{C_TH^{s_1}}+\|\partial_tX\|_{C_TH^{s_1-1}}
  +\|X\|_{L_T^8W^{s_1-1/4,8/3}},\\
\|Y\|_{Y^{s_2}_{T,\mathrm{end}}}
&:=\|Y\|_{C_TH^{s_2}}+\|\partial_tY\|_{C_TH^{s_2-1}}
  +\|Y\|_{L_T^4W^{s_2-1/2,4}}.
\end{align}
The augmented endpoint spaces are obtained from these by adding the same $L_T^\infty B_{2,\infty}^\sigma$ and velocity $L_T^\infty B_{2,\infty}^{\sigma-1}$ envelopes as in \eqref{eq:Xtilde-space}--\eqref{eq:Ytilde-space}.  The notation $X_T^{s_1}$ and $Y_T^{s_2}$ below is branch-dependent: it means the Case F spaces in the subunit fractional branch and the displayed Case W spaces at $\alpha=1$.

\begin{proposition}[Endpoint conic Klein--Gordon linear estimates]\label{prop:gko-conic-strichartz}
Let $\alpha=1$, $0<T\le1$, and
\[
        \omega_c(n)=(1+c^2|n|^2)^{1/2},\qquad c>0.
\]
For each fixed speed $c$, the classical conic wave/Klein--Gordon estimates on $\mathbb T^3$ give the following estimates.  The pairs $(8,8/3)$ and $(4,4)$ satisfy
\[
        \frac2p+\frac2q=1,
        \qquad q<\infty,
\]
so they lie on the three-dimensional wave-admissible boundary but are not the excluded pair $(2,\infty)$.  The first has Sobolev loss $1/4$ and is used for $X$; the second has loss $1/2$ and is used for $Y$.  The pair $(4/3,4/3)$ is the dual source line used only for the endpoint deterministic quadratic part of the $Y$ equation.
\begin{align}
 \|S_c(\cdot)(f,g)\|_{X^{s_1}_{T,\mathrm{end}}}
 &\lesssim_{c,T}\|f\|_{H^{s_1}}+\|g\|_{H^{s_1-1}},
 \label{eq:conic-hom-X}\\
 \|I_cF\|_{X^{s_1}_{T,\mathrm{end}}}
 &\lesssim_{c,T}\|F\|_{L_T^1H^{s_1-1}},
 \label{eq:conic-duh-X}\\
 \|S_c(\cdot)(f,g)\|_{Y^{s_2}_{T,\mathrm{end}}}
 &\lesssim_{c,T}\|f\|_{H^{s_2}}+\|g\|_{H^{s_2-1}},
 \label{eq:conic-hom-Y}\\
 \|I_cF\|_{Y^{s_2}_{T,\mathrm{end}}}
 &\lesssim_{c,T}\|F\|_{L_T^1H^{s_2-1}},
 \label{eq:conic-duh-Y}\\
 \|I_cF\|_{Y^{s_2}_{T,\mathrm{end}}}
 &\lesssim_{c,T}\|\langle D\rangle^{s_2-1/2}F\|_{L_{t,x}^{4/3}([0,T]\times\mathbb T^3)}.
 \label{eq:conic-dual-Y}
\end{align}
The homogeneous estimates are direct instances of the compact Klein--Gordon theorem \cite[Theorem~1(1)]{CacciafestaDanesiMeng}.  The derivative loss there is
\[
        \gamma^W_{p,q}=3\left(\frac12-\frac1q\right)-\frac1p,
\]
which equals $1/4$ for $(p,q)=(8,8/3)$ and $1/2$ for $(p,q)=(4,4)$.  The retarded source estimates follow by Minkowski or by the $TT^*$ and retarded-truncation argument \cite{KeelTao}; the same source lines are recorded in the scalar stochastic wave setting in \cite[Lemma~2.4]{GKO}.  The scalar constants depend on the fixed speed $c$ and on $T$ only; in particular, this linear step contains no factor of $|c_1-c_2|^{-1}$.
\end{proposition}

\begin{proof}[Reference and reduction]
This proposition is used as a deterministic scalar conic theorem.  We recall the reduction to make the derivative losses and the constant dependence explicit.  Since
\[
        \omega_c(D)=c\bigl(c^{-2}+|D|^2\bigr)^{1/2},
\]
the change of time variable $\tau=ct$ reduces the half-wave pieces $e^{\pm\ii t\omega_c(D)}$ to the massive compact Klein--Gordon propagator in \cite[Theorem~1(1)]{CacciafestaDanesiMeng}, with mass $c^{-1}$.  On the three-torus that theorem gives
\[
 \|e^{\pm\ii t\omega_c(D)}f\|_{L_T^pL_x^q}
 \lesssim_{c,T}\|f\|_{H^{\gamma^W_{p,q}}},
 \qquad
 \gamma^W_{p,q}=3\left(\frac12-\frac1q\right)-\frac1p,
\]
for the wave-admissible pairs used here.  The fixed speed enters only through the mass and the time rescaling, so no inverse power of $|c_1-c_2|$ occurs in these scalar estimates.

The two pairs are boundary-admissible but not the exceptional wave endpoint.
They give exactly the Sobolev losses recorded in the statement: $(8,8/3)$
gives $1/4$ in the $X$ norm, and $(4,4)$ gives $1/2$ in the $Y$ norm.  The
retarded $L_T^1H^{s-1}$ estimates follow by Minkowski from the homogeneous
estimates applied to the Duhamel representation with multiplier
$\omega_c(D)^{-1}$.

For completeness, define
\[
        \mathcal T_cf(t)=
        e^{\ii t\omega_c(D)}\langle D\rangle^{-1/2}f.
\]
The homogeneous $(4,4)$ estimate says that
$\mathcal T_c:L_x^2\to L_{t,x}^4$.  Hence
$\mathcal T_c\mathcal T_c^*:L_{t,x}^{4/3}\to L_{t,x}^4$ and its full
convolution kernel is
\[
        e^{\ii(t-s)\omega_c(D)}\langle D\rangle^{-1}.
\]
For fixed $c>0$, the ratio
$\langle n\rangle/\omega_c(n)$ is a zero-order periodic multiplier bounded on
every $L^r$, $1<r<\infty$.  The two half-wave pieces therefore give the full-time
version of \eqref{eq:conic-dual-Y}.  Commuting
$\langle D\rangle^{s_2-1/2}$ through the propagator also gives the energy
components by the dual homogeneous estimate.  Finally, the retarded
restriction $s<t$ is valid because the source time exponent $4/3$ is strictly
smaller than the target exponent $4$; this is the standard retarded
truncation in \cite{KeelTao}.  Thus \cite[Theorem~1(1)]{CacciafestaDanesiMeng}
is used only for the homogeneous scalar estimates, while the displayed
inhomogeneous line is derived from them.
\end{proof}

\begin{corollary}[Classical endpoint from the conic Case W estimates]\label{cor:alpha-one-extension}
Set $\alpha=1$ and $c_1\ne c_2$.  Using Proposition~\ref{prop:gko-conic-strichartz} and the endpoint parameter window in Definition~\ref{def:admissible}, the deterministic fixed point theorem, Theorem~\ref{thm:deterministic-closure}, and the cutoff stability statement in Corollary~\ref{cor:cutoff-stability} remain valid in the endpoint spaces displayed in Proposition~\ref{prop:gko-conic-strichartz} with the augmented $B_{2,\infty}$ envelopes.  The stochastic mixed-operator estimates fit the endpoint window.
\end{corollary}

\begin{proof}
All stochastic source estimates and mixed-operator bounds used in Section~\ref{sec:random-operators} are algebraic or dyadic estimates with the powers displayed in Section~\ref{sec:parameters}.  At $\alpha=1$ the phase-difference Volterra denominator is $N^{-1}$, the same-color diagonal shell is $N^{-1}$, and the centered row/column conditions reduce to $s_2<1-$ and $\sigma>0+$.  The first Picard regularity is $\rho_V=1/2-$, and the cubic source has nonnegative margin for every $s_2<1$.  Hence the stochastic enhanced data and the mixed-operator estimates remain valid.  The deterministic endpoint closure uses \eqref{eq:conic-dual-Y} for the pure quadratic $Y$-sources and the $L_T^1H^{s-1}$ routes for the remaining rough and mixed sources.  The window
\[
        \frac14<s_1<\frac12< s_2<\min\left\{1,s_1+\frac14\right\},
        \qquad 0<\sigma<\min\left\{s_1,s_2,\frac14\right\}
\]
is nonempty, for instance at $(s_1,s_2,\sigma)=(0.40,0.60,0.10)$ before inserting small losses.
\end{proof}

Define
\begin{align}
\|X\|_{X_T^{s_1}}
&:=\|X\|_{C_TH^{s_1}}+\|\partial_tX\|_{C_TH^{s_1-\alpha}}
  +\|X\|_{L_T^{p_X}W^{s_1-\kappa_X,q_X}},
\label{eq:XJ-space}\\
\|Y\|_{Y_T^{s_2}}
&:=\|Y\|_{C_TH^{s_2}}+\|\partial_tY\|_{C_TH^{s_2-\alpha}}
  +\|Y\|_{L_T^4W^{s_2-\kappa_Y,4}}.
\label{eq:YJ-space}
\end{align}
The augmented input topology is
\begin{align}
\|X\|_{\widetilde X_T^{s_1,\sigma}}
&:=\|X\|_{X_T^{s_1}}+\|X\|_{L_T^\infty\Btwo^\sigma}
      +\|\partial_tX\|_{L_T^\infty\Btwo^{\sigma-\alpha}},
\label{eq:Xtilde-space}\\
\|Y\|_{\widetilde Y_T^{s_2,\sigma}}
&:=\|Y\|_{Y_T^{s_2}}+\|Y\|_{L_T^\infty\Btwo^\sigma}
      +\|\partial_tY\|_{L_T^\infty\Btwo^{\sigma-\alpha}}.
\label{eq:Ytilde-space}
\end{align}

\begin{lemma}[Completeness of the augmented path spaces]\label{lem:augmented-path-complete}
For fixed admissible exponents and for either branch, the spaces $\widetilde X_T^{s_1,\sigma}$ and $\widetilde Y_T^{s_2,\sigma}$ equipped with the norms in \eqref{eq:Xtilde-space}--\eqref{eq:Ytilde-space} are Banach spaces.  The same is true for the product space $\mathcal Z_T$.  The fixed point is taken in this complete product norm.  The lower-index continuity
\[
        u\in \bigcap_{\eta>0} C_TB_{2,\infty}^{\sigma-\eta},
        \qquad
        \partial_tu\in \bigcap_{\eta>0} C_TB_{2,\infty}^{\sigma-\alpha-\eta}
\]
is the path identification supplied by the dyadic propagation estimate.  If a sequence also belongs to the little-Besov closure at the endpoint index, then the limit belongs to the same closure.
\end{lemma}

\begin{proof}
Each component in \eqref{eq:Xtilde-space}--\eqref{eq:Ytilde-space} is complete.  The Sobolev--Strichartz pieces are standard Bochner spaces.  The envelope $B_{2,\infty}^s$ is the inhomogeneous dyadic $\ell^\infty(L^2)$ Besov space, hence complete, and therefore $L_T^\infty B_{2,\infty}^s$ is complete as well.  A finite intersection with the sum norm is complete, and so is the finite product $\mathcal Z_T$.

If the data are taken in the little-Besov closure, convergence in the endpoint dyadic norm preserves uniform smallness of the high-frequency tail, so the limit remains in that closure.  For lower-index continuity, Lemma~\ref{lem:elementary-besov-propagation} gives continuity of every dyadic block of the homogeneous and Duhamel pieces.  After multiplying the high tail by $M^{-\eta}$, the $L_T^\infty B_{2,\infty}^{\sigma}$ and velocity envelopes make the tail uniformly small in $B_{2,\infty}^{\sigma-\eta}$ and $B_{2,\infty}^{\sigma-\alpha-\eta}$.  This is the path-valued continuity used in the cutoff convergence statement.
\end{proof}

\begin{theorem}[Sobolev--Strichartz propagation for the subunit fractional flow]\label{prop:strichartz}
Assume $12/13<\alpha<1$ and use the Case F exponents above.  For $i=1,2$ and $0<T\le1$, the homogeneous Sobolev estimate
\begin{equation}
 \|S_i(\cdot)(f,g)\|_{Y_T^{s_2}}
 \lesssim \|f\|_{H^{s_2}}+\|g\|_{H^{s_2-\alpha}},
 \label{eq:homogeneous-Y-estimate}
\end{equation}
and the Sobolev Duhamel estimates
\begin{align}
 \|I_iF\|_{X_T^{s_1}}&\lesssim \|F\|_{L_T^1H^{s_1-\alpha}},
 \label{eq:strichartz-input-X}\\
 \|I_iF\|_{Y_T^{s_2}}&\lesssim \|F\|_{L_T^1H^{s_2-\alpha}}
 \label{eq:strichartz-input-Y}
\end{align}
hold.  The constants may depend on $T,\alpha,p_X,q_X,\eps_{\rm Str}$ and on the individual speeds $c_i$; the speed-gap dependence enters only through the phase-difference estimates.
\end{theorem}

\begin{proof}
The $C_TH^s$ and $C_T^1H^{s-\alpha}$ components are the standard energy estimates for the constant-coefficient Klein--Gordon flow.  For the spacetime components, write the cosine and sine propagators as linear combinations of $e^{\pm \ii t\omega_i(D)}$ and use the dyadic estimate of Theorem~\ref{thm:fracKG-dyadic} with $c=c_i$.  The homogeneous part with data $(f,g)$ follows by applying Corollary~\ref{cor:fracKG-semigroup-duhamel} to $f$ and to $\omega_i(D)^{-1}g$, since $\omega_i(D)^{-1}$ maps $H^{s-\alpha}$ to $H^s$.  For the inhomogeneous part, the same corollary gives
\begin{equation*}
        I_i:L_T^1H^{s-\alpha}\longrightarrow L_T^pW^{s-\kappa(p,q),q}
\end{equation*}
for $(p,q)=(p_X,q_X)$ and $(4,4)$.  The arbitrary $N^{0+}$ loss in Theorem~\ref{thm:fracKG-dyadic} is the fixed number $\eps_{\rm Str}$ already included in $\kappa_X$ and $\kappa_Y$.  The low dyadic blocks are finite dimensional and are absorbed into the energy constant.  This proves \eqref{eq:homogeneous-Y-estimate}--\eqref{eq:strichartz-input-Y}.
\end{proof}

\begin{lemma}[Elementary Besov propagation for the Klein--Gordon flow]\label{lem:elementary-besov-propagation}
For $i=1,2$ and $0<T\le1$,
\begin{align}
 \|S_i(\cdot)(f,g)\|_{L_T^\infty\Btwo^\sigma}
 +\|\partial_tS_i(\cdot)(f,g)\|_{L_T^\infty\Btwo^{\sigma-\alpha}}
 &\lesssim \|f\|_{\Btwo^\sigma}+\|g\|_{\Btwo^{\sigma-\alpha}},
 \label{eq:homogeneous-besov-estimate}\\
 \|I_iF\|_{L_T^\infty\Btwo^\sigma}
 +\|\partial_tI_iF\|_{L_T^\infty\Btwo^{\sigma-\alpha}}
 &\lesssim \|F\|_{L_T^1\Btwo^{\sigma-\alpha}}.
 \label{eq:besov-input-dh}
\end{align}
The implicit constants depend on $T$ and on the fixed scalar symbol $\omega_i$.
\end{lemma}

\begin{proof}
The dyadic multiplier proof is Proposition~\ref{prop:app-b2-propagation}; the time-continuity statement is Corollary~\ref{cor:app-b2-lower-continuity}.
\end{proof}

\subsection{Control norm and admissible parameters}
For $Z=(X_1,X_2,Y_1,Y_2)$ set
\begin{equation}\label{eq:ZJ-norm}
\|Z\|_{\calZ_T}:=\sum_{i=1}^2\bigl(
\|X_i\|_{\widetilde X_T^{s_1,\sigma}}+
\|Y_i\|_{\widetilde Y_T^{s_2,\sigma}}
\bigr),
\qquad W_i=V_i+X_i+Y_i.
\end{equation}
The deterministic data norm is
\begin{equation}\label{eq:data-norm-Y}
\|y\|_{\mathcal H_T^{s_2,\sigma}}:=
\sum_{i=1}^2\Bigl(
\|y_{i,0}\|_{H^{s_2}}+\|y_{i,1}\|_{H^{s_2-\alpha}}
+\|y_{i,0}\|_{\Btwo^\sigma}+\|y_{i,1}\|_{\Btwo^{\sigma-\alpha}}
\Bigr).
\end{equation}

In the subunit branch we use the rough source norm
\begin{equation}\label{eq:source-norm-SF}
\|F\|_{\mathfrak S_{T,F}^{s_2,\sigma}}
:=\|F\|_{L_T^1H^{s_2-\alpha}}
 +\|F\|_{L_T^1\Btwo^{\sigma-\alpha}}.
\end{equation}
In the endpoint branch the corresponding rough-source norm is
\begin{equation}\label{eq:source-norm-SW}
\|F\|_{\mathfrak S_{T,W}^{s_2,\sigma}}
:=\|F\|_{L_T^1H^{s_2-1}}
 +\|F\|_{L_T^1\Btwo^{\sigma-1}}.
\end{equation}
We write \(\mathfrak S_T^{s_2,\sigma}\) for the branch currently under discussion.  Case W deterministic non-\(V\) quadratic terms are routed through the dual conic source seminorm below.

In Case W we also use, only for the deterministic non-$V$ quadratic part of the $Y$ equation, the dual conic source seminorm
\begin{equation}\label{eq:endpoint-dual-source-norm}
        \|F\|_{\mathfrak D_T^{s_2}}
        :=\|\langle D\rangle^{s_2-1/2}F\|_{L_{t,x}^{4/3}([0,T]\times\mathbb T^3)}.
\end{equation}
For the $Y$ equation write
\[
        F_Y=F_{\rm rough}+F_{\rm op}+F_\Gamma+F_V+F_{\rm det},
\]
where $F_{\rm rough}$ denotes the high--low and resonant rough stochastic products, $F_{\rm op}$ the mixed-operator sector, $F_\Gamma$ the cubic symbols, $F_V$ the quadratic terms containing at least one first Picard factor, and $F_{\rm det}$ the remaining non-$V$ deterministic quadratic terms.  Define
\begin{equation}\label{eq:source-topologies}
\|F_Y\|_{\mathfrak P_T^{s_2,\sigma}}
:=
\begin{cases}
 \|F_Y\|_{\mathfrak S_{T,F}^{s_2,\sigma}},
 & \text{in Case F},\\[1mm]
 \|F_{\rm rough}+F_{\rm op}+F_\Gamma+F_V\|_{\mathfrak S_{T,W}^{s_2,\sigma}}
 +\|F_{\rm det}\|_{\mathfrak D_T^{s_2}}
 +\|F_{\rm det}\|_{L_T^1B_{2,\infty}^{\sigma-1}},
 & \text{in Case W}.
\end{cases}
\end{equation}
Thus all sources use the rough norm in Case F.  At \(\alpha=1\), only the deterministic non-$V$ quadratic sector uses the dual conic norm, together with its separate Besov estimate.

\begin{definition}[Augmented enhanced-data control norm]\label{def:enhanced-data-norm}
The deterministic Cauchy data are measured by \eqref{eq:data-norm-Y} and are outside the stochastic enhanced data.  For an enhanced datum $\Xi$, the Sobolev part of the mixed-operator norm is
\begin{equation}\label{eq:MJopH}
M_T^{\rm op,H}(\Xi):=
\sum_{(i;j,k)\in\mathfrak M}
\|T^{i;j,k}\|_{\mathcal L(E_T^{2,\sigma},C_TH^{s_2-\alpha})}.
\end{equation}
The Besov source part, supplied by the mixed-operator theorem and used by the augmented deterministic closure, is
\begin{equation}\label{eq:MJopB}
M_T^{\rm op,B}(\Xi):=
\sum_{(i;j,k)\in\mathfrak M}
\|T^{i;j,k}\|_{\mathcal L(E_T^{2,\sigma},L_T^1\Btwo^{\sigma-\alpha})}.
\end{equation}
The enhanced-data control norm is
\begin{align}\label{eq:MJ}
M_T(\Xi):=&\sum_i\Bigl(
\|\Psi_i\|_{C_T\C^{\alpha-3/2-}}
+\|V_i\|_{C_T\C^{\rho_V}}
+\|\partial_tV_i\|_{C_T\C^{\rho_V-\alpha}}
\Bigr)\nonumber\\
&+\|\Theta\|_{C_T\C^{2\alpha-3-}}\nonumber\\
&+\sum_i\Bigl(
\|V_i\|_{L_T^\infty\Btwo^{\rho_V}}
+\|\partial_tV_i\|_{L_T^\infty\Btwo^{\rho_V-\alpha}}
\Bigr)\nonumber\\
&+\sum_i\Bigl(
\|\Gamma_i\|_{L_T^1H^{s_2-\alpha}}
+\|\Gamma_i\|_{L_T^1\Btwo^{\sigma-\alpha}}
\Bigr)\nonumber\\
&+M_T^{\rm op,H}(\Xi)+M_T^{\rm op,B}(\Xi),\qquad
\rho_V=\frac{7\alpha}{2}-3-.
\end{align}
Here both $M_T^{\rm op,H}$ and $M_T^{\rm op,B}$ are supplied by the Besov-input mixed-operator theorem, Theorem~\ref{thm:pathwise-fluct}, together with the deterministic Volterra diagonal estimate.  They are kept in the control norm because the fixed-point constants depend on their size.  For fixed deterministic data $y$ we write
\begin{equation}\label{eq:bounded-size}
        M_T^{\rm bd}(\Xi,y):=1+M_T(\Xi)+\|y\|_{\mathcal H_T^{s_2,\sigma}},
\end{equation}
and often suppress $y$ from the notation when the deterministic data are fixed.  This convention separates the stochastic enhanced data from the Cauchy data while preserving a single bounded-size parameter in the fixed point.

Local contraction smallness is recorded separately in Definition~\ref{def:localized-mixed-size}.  The two first-Picard controls have distinct roles: $V_i\in L_T^\infty\Btwo^{\rho_V}$ is the dyadic Besov input used inside the mixed-operator topology, whereas $V_i\in C_T\C^{\rho_V}$ is the H\"older multiplier used in deterministic products.  Appendix~\ref{app:first-picard-phase} proves the H\"older estimate directly, separately from the dyadic $B_{2,\infty}^{\rho_V}$ envelope.
\end{definition}

\begin{definition}[Localized mixed-operator size]\label{def:localized-mixed-size}
For the contraction argument we separate boundedness from local smallness.  The diagonal Besov contribution is measured in the stronger Volterra norm
\[
M_{T,\mathrm{diag}}^{\rm op,B}(\Xi):=
\sum_{(i;j,j)\in\mathfrak M}\|\mathcal D^{i;j}_{[0,T]}\|_{E_T^{2,\sigma}\to L_T^\infty\Btwo^{\sigma-\alpha}},
\]
and is converted to an $L_T^1$ source only by the explicit factor $T$ below.  This stronger diagonal $L_T^\infty B_{2,\infty}^{\sigma-\alpha}$ operator norm is recorded separately from $M_T^{\rm op,B}$ because it becomes a Duhamel source norm only after this deterministic multiplication by $T$.  The centered Besov contribution is the direct source norm
\[
M_{T,\mathrm{cen}}^{\rm op,B}(\Xi):=
\sum_{(i;j,k)\in\mathfrak M}\|\mathcal B^{i;j,k}_{[0,T]}\|_{E_T^{2,\sigma}\to L_T^1\Btwo^{\sigma-\alpha}}.
\]
Set
\begin{equation}\label{eq:localized-mixed-size}
\mu_T^{\rm op}(\Xi):=T M_T^{\rm op,H}(\Xi)+T M_{T,\mathrm{diag}}^{\rm op,B}(\Xi)+M_{T,\mathrm{cen}}^{\rm op,B}(\Xi).
\end{equation}
For fixed deterministic data $y$ and a bounded enhanced-data ball whose size is controlled by $M_T^{\rm bd}(\Xi,y)$, set
\begin{equation}\label{eq:localized-full-size}
\mu_T(\Xi,y):=\mu_T^{\rm op}(\Xi)+T^\delta\bigl(1+M_T^{\rm bd}(\Xi,y)+(M_T^{\rm bd}(\Xi,y))^2\bigr).
\end{equation}
Equivalently, in finite-sum form,
\begin{align}\label{eq:muT-explicit}
\mu_T(\Xi,y)&=T^\delta\bigl(1+M_T(\Xi)+\|y\|_{\mathcal H_T^{s_2,\sigma}}+(M_T^{\rm bd}(\Xi,y))^2\bigr)\nonumber\\
&\quad+T\sum_{(i;j,j)\in\mathfrak M}\|\mathcal D^{i;j}_{[0,T]}\|_{E_T^{2,\sigma}\to C_TH^{s_2-\alpha}}
+T\sum_{(i;j,j)\in\mathfrak M}\|\mathcal D^{i;j}_{[0,T]}\|_{E_T^{2,\sigma}\to L_T^\infty B_{2,\infty}^{\sigma-\alpha}}\nonumber\\
&\quad+T\sum_{(i;j,k)\in\mathfrak M}\|\mathcal B^{i;j,k}_{[0,T]}\|_{E_T^{2,\sigma}\to C_TH^{s_2-\alpha}}
+\sum_{(i;j,k)\in\mathfrak M}\|\mathcal B^{i;j,k}_{[0,T]}\|_{E_T^{2,\sigma}\to L_T^1B_{2,\infty}^{\sigma-\alpha}}.
\end{align}
The first line records deterministic Duhamel time factors, rough products, cubic sources, and the data-enlarged bounded radius.  The second line is the Volterra diagonal source contribution; the Besov source norm is read as the induced $L_T^1$ norm from the stronger $L_T^\infty$ diagonal estimate.  The third line is the centered operator contribution.  For a fixed enhanced realization with finite dyadic majorants and fixed deterministic data, $\mu_T(\Xi,y)\to0$ as $T\downarrow0$ after restriction to $[0,T]$; the diagonal contribution is Proposition~\ref{prop:det-diagonal-local-smallness}, and the centered Besov contribution is Proposition~\ref{prop:centered-besov-local-smallness}.  When the deterministic data are fixed, we also write $\mu_T(\Xi)$ for $\mu_T(\Xi,y)$.

This distinction is used in the contraction proof.  The bounded quantity $M_T^{\rm bd}(\Xi,y)$ fixes the radius of the ball and may be large on a fixed time interval.  The contraction constant is controlled by $\mu_T(\Xi,y)$, whose entries become small after time localization: deterministic $L_T^1$ conversions, localized Volterra/Besov operator norms, and non-endpoint H\"older-in-time factors.  The pathwise stochastic symbol norms $\|\Psi_i\|$, $\|\Theta\|$, $\|V_i\|$, and $\|\Gamma_i\|$ enter the bounded radius.  The centered operator contribution is measured in its localized $L_T^1B_{2,\infty}^{\sigma-\alpha}$ source norm, while the Volterra diagonal is converted from the stronger $L_T^\infty B_{2,\infty}^{\sigma-\alpha}$ bound by the deterministic time factor.
\end{definition}

\begin{lemma}[Origin of the small factors]\label{lem:origin-small-factors}
For a fixed enhanced realization and deterministic data with finite bounded size, every term in $\mu_T(\Xi,y)$ is either an explicit time-length factor or a localized operator source norm.  More precisely, the rough paraproducts and deterministic quadratic products contribute powers $T^\delta$ from H\"older or $L_T^1$ estimates; the Volterra diagonal contributes the explicit factor $T$ in \eqref{eq:det-contraction-L1}; and the centered Besov source contribution is the localized $L_T^1B_{2,\infty}^{\sigma-\alpha}$ norm of Proposition~\ref{prop:centered-besov-local-smallness}.  The fixed-time stochastic symbol norms enter the bounded size; localized source norms supply the smallness.
\end{lemma}

\begin{proof}
The first line of \eqref{eq:muT-explicit} contains exactly the deterministic source estimates proved in Lemmas~\ref{lem:lowhigh-X}, \ref{lem:caseF-rough-stochastic-routes}, \ref{lem:caseW-rough-stochastic-routes}, and \ref{lem:caseF-product-margins}, together with their endpoint analogues.  Their smallness is the displayed finite-time factor.  The diagonal part is Proposition~\ref{prop:det-diagonal-local-smallness}, equivalently \eqref{eq:det-contraction-L1}.  The centered Besov part is Proposition~\ref{prop:centered-besov-local-smallness}; its smallness comes from the Volterra time integration and the $L_T^1$ source length after the operator norm has already been constructed.  The bounded terms $M_T(\Xi)$ and $\|y\|_{\mathcal H_T^{s_2,\sigma}}$ only determine the radius and the Lipschitz polynomial.
\end{proof}

\begin{definition}[Localized smallness classes]\label{def:localized-smallness-class}
Fix $0<T_0\le1$.  A family $\mathcal K$ of pairs $(\Xi,y)$, consisting of an augmented enhanced-data set and deterministic Cauchy data on $[0,T_0]$, is called localized-small with bounded radius $A$ if
\[
        \sup_{(\Xi,y)\in\mathcal K}M_{T_0}^{\rm bd}(\Xi,y)\le A,
        \qquad
        \nu_{\mathcal K}(T):=\sup_{(\Xi,y)\in\mathcal K}\mu_T(\Xi|_{[0,T]},y)\longrightarrow0
        \quad (T\downarrow0).
\]
On such a class the contraction time may be chosen from $A$ and the function $\nu_{\mathcal K}$.  For a single pair, or for a fixed finite collection of pairs, this condition follows from Definition~\ref{def:localized-mixed-size}, Proposition~\ref{prop:det-diagonal-local-smallness}, and Proposition~\ref{prop:centered-besov-local-smallness}; the resulting time may depend on the individual local profiles.  The local smallness modulus is part of the time-selection data together with the bounded radius $M_{T_0}^{\rm bd}$.
\end{definition}

\begin{lemma}[Time restriction for localized-smallness classes]\label{lem:time-restriction-local-smallness}
Let $\mathcal K$ be localized-small on $[0,T_0]$ with bounded radius $A$ and modulus $\nu_{\mathcal K}$.  For every nondecreasing polynomial $P$ and every threshold $\delta>0$, there exists $0<T_*\le T_0$, depending only on $A,P,\delta$, and $\nu_{\mathcal K}$, such that
\[
        \mu_T(\Xi|_{[0,T]},y)P(A)\le\delta,
        \qquad 0<T\le T_*,\quad (\Xi,y)\in\mathcal K.
\]
For a singleton or a fixed finite collection, this gives the common small-time restriction used in the contraction and Lipschitz estimates.  The time is controlled by the localized-smallness modulus together with the bounded radius.
\end{lemma}

\begin{proof}
By Definition~\ref{def:localized-smallness-class}, $\nu_{\mathcal K}(T)\to0$ as $T\downarrow0$.  Choose $T_*$ so that $\nu_{\mathcal K}(T)P(A)\le\delta$ for all $0<T\le T_*$.  This proves the claim for every pair $(\Xi,y)$ in the class.
\end{proof}

\begin{definition}[Augmented enhanced-data distance]\label{def:enhanced-data-distance}
For two enhanced stochastic data sets together with deterministic data, $(\mathbb E,y)$ and $(\widetilde{\mathbb E},\widetilde y)$, define $d_T((\mathbb E,y),(\widetilde{\mathbb E},\widetilde y))$ by replacing every stochastic norm in Definition~\ref{def:enhanced-data-norm} with the corresponding difference norm and by adding the deterministic Cauchy-data difference.  Thus $d_T$ contains
\begin{align*}
&\sum_i\Bigl(
\|\Psi_i-\widetilde\Psi_i\|_{C_T\C^{\alpha-3/2-}}
+\|V_i-\widetilde V_i\|_{C_T\C^{\rho_V}}
+\|\partial_tV_i-\partial_t\widetilde V_i\|_{C_T\C^{\rho_V-\alpha}}
\Bigr)\\
&+\|\Theta-\widetilde\Theta\|_{C_T\C^{2\alpha-3-}}\\
&+\sum_i\Bigl(
\|V_i-\widetilde V_i\|_{L_T^\infty\Btwo^{\rho_V}}
+\|\partial_tV_i-\partial_t\widetilde V_i\|_{L_T^\infty\Btwo^{\rho_V-\alpha}}
\Bigr)\\
&+\sum_i\Bigl(
\|\Gamma_i-\widetilde\Gamma_i\|_{L_T^1H^{s_2-\alpha}}
+\|\Gamma_i-\widetilde\Gamma_i\|_{L_T^1\Btwo^{\sigma-\alpha}}\Bigr)\\
&+\sum_{(i;j,k)\in\mathfrak M}\Bigl(
\|T^{i;j,k}-\widetilde T^{i;j,k}\|_{\mathcal L_H}
+\|T^{i;j,k}-\widetilde T^{i;j,k}\|_{\mathcal L_B}
\Bigr)\\
&+\|y-\widetilde y\|_{\mathcal H_T^{s_2,\sigma}}.
\end{align*}
Here
\[
\mathcal L_H:=\mathcal L(E_T^{2,\sigma},C_TH^{s_2-\alpha}),\qquad
\mathcal L_B:=\mathcal L(E_T^{2,\sigma},L_T^1\Btwo^{\sigma-\alpha}).
\]
The operator difference is taken on the common deterministic input space $E_T^{2,\sigma}$.  When the deterministic data are fixed, the last term in $d_T$ is absent.  Cutoff enhanced data converge in the augmented sense precisely when this distance tends to zero and the corresponding norms remain uniformly bounded.  In particular, convergence of the first Picard object includes the H\"older topology $C_T\C^{\rho_V}$ used for deterministic products and the dyadic $B_{2,\infty}^{\rho_V}$ topology used for Besov-Besov input bounds.
\end{definition}

\begin{definition}[Admissible deterministic parameters]\label{def:admissible}
Choose a sufficiently small loss $\eps>0$.  The admissible window is branch-dependent.

\smallskip
\noindent\emph{Case F: $12/13<\alpha<1$.}  The subunit fractional parameters satisfy
\begin{equation}\label{eq:admissible-base}
0<s_1<s_2<\alpha,
\qquad
3-3\alpha+10\eps<\sigma<
\min\left\{s_1,s_2,\rho_V,4\alpha-3,2\alpha-\tfrac32\right\}-10\eps,
\end{equation}
\begin{equation}\label{eq:admissible-operator}
 s_2<4\alpha-3-10\eps.
\end{equation}
\begin{equation}\label{eq:admissible-lowhigh}
 s_1<2\alpha-\frac32-\eps,
 \qquad
 s_2>\frac32-\alpha+\eps,
\end{equation}
\begin{equation}\label{eq:admissible-high}
 s_2<s_1+2\alpha-\frac32-\eps,
 \qquad
 s_2<\frac{11\alpha}{2}-\frac92-\eps,
\end{equation}
and
\begin{equation}\label{eq:admissible-quadratic}
 s_2<2s_1+\alpha+\frac{2\alpha}{p_X}-\frac32-\eps,
 \qquad
 s_1>\frac32-\frac{5\alpha}{4}-\frac{\alpha}{p_X}+\eps.
\end{equation}
Equivalently, after fixing $(s_2,\sigma,p_X,q_X)$, the active lower bounds on $s_1$ may be read as
\begin{equation}\label{eq:s1-window}
\max\left\{\sigma,\ s_2+\frac32-2\alpha,
\frac12\left(s_2-\alpha-\frac{2\alpha}{p_X}+\frac32\right),
\frac32-\frac{5\alpha}{4}-\frac{\alpha}{p_X}\right\}+O(\eps)
<s_1<2\alpha-\frac32-O(\eps).
\end{equation}
The four lower bounds respectively encode the input Besov estimate, the $X>\Psi$ high--low route, the $X_1X_2$ Strichartz product route, and the remaining mixed Strichartz margin.

\smallskip
\noindent\emph{Case W: $\alpha=1$.}  With the conic endpoint exponents \eqref{eq:case-W-exponents}, the endpoint parameters satisfy
\begin{equation}\label{eq:admissible-case-W}
 \frac14+10\eps<s_1<\frac12-10\eps,
 \qquad
 \frac12+10\eps<s_2<\min\left\{1-10\eps,\,s_1+\frac14-10\eps\right\},
\end{equation}
\begin{equation}\label{eq:admissible-case-W-sigma}
 10\eps<\sigma<\min\left\{s_1,s_2,\frac14\right\}-10\eps.
\end{equation}
This endpoint window contains the stochastic constraints $s_2<1-$ and $\sigma>0+$, the first-Picard input constraint $\sigma<\rho_V=1/2-$, and the deterministic dual-source condition $s_2<s_1+1/4-$ for the conic quadratic sector.  The displayed upper bound $\sigma<1/4-$ ensures the embedding $L_x^{4/3}\hookrightarrow B_{2,\infty}^{\sigma-1}$.

The Volterra diagonal is summed directly on the proportional paraproduct window $Q\le c_0N$, with $c_0$ chosen below the phase-gap aperture in Lemma~\ref{lem:speed-gap}.  The admissible tuple therefore only records the exponents used by the source and fixed-point estimates.
\end{definition}

\begin{lemma}[Nonempty parameter range]\label{lem:parameter-range}
The admissible set is nonempty in both branches.  In Case F this holds for $(p_X,q_X)$ sufficiently close to $(4,3)$ and for $\eps_{\rm Str}$ sufficiently small compared with the fixed small margin.  In Case W the endpoint window has nonempty interior; for example $(s_1,s_2,\sigma)=(0.40,0.60,0.10)$ satisfies the strict inequalities before the small losses are inserted.
\end{lemma}

\begin{proof}
We first remove the harmless small losses and prove that the strict window has positive distance from its boundary.  The losses $10\eps$, $\eps$, $\eps_{\rm Str}$, and the fixed loss $\kappa_V$ in \(\rho_V\) are then chosen below this distance.

Consider Case F.  Choose $(p_X,q_X)$ close to $(4,3)$; only $p_X$ enters the algebra below.  The lower endpoint for $s_2$ is $L_2=3/2-\alpha$.  The available upper endpoint is
\[
        U_2(p_X)=\min\left\{4\alpha-3,
        \frac{11\alpha}{2}-\frac92,
        5\alpha+\frac{2\alpha}{p_X}-\frac92,
        \alpha\right\}.
\]
The third term is the condition obtained from the quadratic lower bound on $s_1$ together with the upper bound $s_1<2\alpha-3/2$.  For $p_X$ close to $4$,
\[
        5\alpha+\frac{2\alpha}{p_X}-\frac92
        =\frac{11\alpha}{2}-\frac92+o_{p_X\to4}(1).
\]
Since
\[
        4\alpha-3-L_2=5\alpha-\frac92,
        \qquad
        \frac{11\alpha}{2}-\frac92-L_2=\frac{13\alpha}{2}-6,
\]
both gaps are positive in the present range, and the second gap is positive exactly when $\alpha>12/13$.  Hence, after taking $p_X$ sufficiently close to $4$, the interval $(L_2,U_2(p_X))$ is nonempty.  Choose
\begin{equation}\label{eq:parameter-proof-s2-choice}
        s_2\in(L_2,U_2(p_X)).
\end{equation}

Next choose $\sigma$.  The lower bound is $3-3\alpha$, and the upper bounds are $s_2$, $\rho_V$, $4\alpha-3$, and $2\alpha-3/2$.  The coordinate-independent gaps
\[
        \rho_V-(3-3\alpha)=\frac{13\alpha}{2}-6,
        \qquad
        \left(2\alpha-\frac32\right)-(3-3\alpha)=5\alpha-\frac92
\]
are positive.  Also $4\alpha-3>3-3\alpha$ for $\alpha>6/7$, and $s_2>3/2-\alpha>3-3\alpha$ for $\alpha>3/4$.  Thus we may choose
\begin{equation}\label{eq:parameter-proof-sigma-choice}
        3-3\alpha<\sigma<\min\{s_2,\rho_V,4\alpha-3,2\alpha-3/2\}.
\end{equation}

It remains to choose $s_1$.  Put $U_1=2\alpha-3/2$ and let $L_1$ be the maximum of the four lower bounds in \eqref{eq:s1-window} without the $O(\eps)$ term.  The choice of $\sigma$ gives $\sigma<U_1$.  The inequality $s_2<4\alpha-3$ gives
\[
        s_2+\frac32-2\alpha<U_1,
\]
and the choice $s_2<5\alpha+2\alpha/p_X-9/2$ gives
\[
        \frac12\left(s_2-\alpha-\frac{2\alpha}{p_X}+\frac32\right)<U_1.
\]
Finally, for $p_X$ close to $4$,
\[
        \frac32-\frac{5\alpha}{4}-\frac{\alpha}{p_X}<U_1
        \quad\Longleftrightarrow\quad
        \frac{13\alpha}{4}+\frac{\alpha}{p_X}>3,
\]
and the right-hand side holds in the present range because at $p_X=4$ it is $7\alpha/2>3$ for $\alpha>6/7$.  Hence $L_1<U_1$.  Choose
\begin{equation}\label{eq:parameter-proof-s1-choice}
        L_1<s_1<U_1.
\end{equation}
Then $s_1<s_2$ automatically, because $\alpha<1$ gives $U_1=2\alpha-3/2<3/2-\alpha=L_2<s_2$.  The remaining displayed bounds on $s_1$ are exactly the four inequalities encoded in $L_1<U_1$.  Hence all inequalities in Definition~\ref{def:admissible} hold with a positive margin.  Taking $\eps>0$ and $\eps_{\rm Str}>0$ smaller than this margin gives the displayed strict Case F window.

In Case W the loss-free inequalities
\[
        \frac14<s_1<\frac12,
        \qquad
        \frac12<s_2<\min\left\{1,s_1+\frac14\right\},
        \qquad
        0<\sigma<\min\left\{s_1,s_2,\frac14\right\}
\]
have nonempty interior.  The point $(s_1,s_2,\sigma)=(0.40,0.60,0.10)$ lies strictly inside it, so inserting sufficiently small losses preserves all inequalities.  The endpoint stochastic restrictions are exactly $s_2<1-$ and $\sigma>0+$, and the endpoint deterministic quadratic route is exactly $s_2<s_1+1/4-$.  The rough stochastic high--low term $X>\Psi$ only needs the weaker condition $s_2<s_1+1/2-$, so it is automatically covered by this endpoint window.

The proportional low--high aperture is fixed once through the paraproduct cutoff and is part of the paraproduct convention.
\end{proof}

\begin{lemma}[The full low--high input belongs to the operator domain]\label{lem:full-input-in-E}
For every admissible parameter choice and every solution candidate $Z=(X_1,X_2,Y_1,Y_2)\in\mathcal Z_T$, the full input
\[
        W_i=V_i+X_i+Y_i
\]
belongs to $E_T^{2,\sigma}$ and
\begin{equation}\label{eq:full-input-E-bound}
        \|W_i\|_{E_T^{2,\sigma}}
        \lesssim M_T(\Xi)+\|X_i\|_{\widetilde X_T^{s_1,\sigma}}+
        \|Y_i\|_{\widetilde Y_T^{s_2,\sigma}}.
\end{equation}
The same estimate holds for differences, with the difference of the $V_i$ component measured in the enhanced-data distance.
\end{lemma}

\begin{proof}
The $C_TL^2$ and velocity $C_T^1H^{-\alpha-}$ components follow from the Sobolev and Besov pieces of the three summands.  For $X_i$ and $Y_i$ this is built into $\widetilde X_T^{s_1,\sigma}$ and $\widetilde Y_T^{s_2,\sigma}$.  For $V_i$ it is one of the explicit components of $M_T(\Xi)$: the H\"older control is used in deterministic products and the $L_T^\infty B_{2,\infty}^{\rho_V}$ envelope is used in Besov-input operator estimates.  Since $\sigma<\rho_V$ and the low shells are finite, $V_i$ also lies in the required $B_{2,\infty}^{\sigma}$ Besov input bound.  The velocity envelope is identical with $\rho_V$ replaced by $\rho_V-\alpha$.
\end{proof}

\subsection{Product toolbox and source estimates}\label{subsec:product-toolbox}
We collect the deterministic estimates entering the fixed-point argument.  The first Picard component has two complementary controls.  As part of the full low--high input it carries a dyadic Besov bound $V_i\in L_T^\infty B_{2,\infty}^{\rho_V}$, compatible with the $E_T^{2,\sigma}$ topology.  As a coefficient in ordinary deterministic products it is controlled by the direct H\"older estimate $V_i\in C_T\C^{\rho_V}$ from Appendix~\ref{app:first-picard-phase}.  Both estimates are included in the enhanced data and enter the Lipschitz bounds below.  The following elementary embedding converts deterministic Sobolev sources into the Besov source space; stochastic rough products and mixed random operators are estimated directly in the dyadic $L^2$ topology.

\begin{lemma}[Strict Sobolev--Besov source embedding]\label{lem:sobolev-to-besov-source}
If $\sigma<s_2$, then
\begin{equation}\label{eq:sob-to-besov-source}
        H^{s_2-\alpha}(\T^3)\hookrightarrow B_{2,\infty}^{\sigma-\alpha}(\T^3),
        \qquad
        \|F\|_{B_{2,\infty}^{\sigma-\alpha}}
        \lesssim \|F\|_{H^{s_2-\alpha}}.
\end{equation}
Consequently,
\begin{equation}\label{eq:LT1-sob-to-besov-source}
        \|F\|_{L_T^1B_{2,\infty}^{\sigma-\alpha}}
        \lesssim
        \|F\|_{L_T^1H^{s_2-\alpha}}.
\end{equation}
\end{lemma}

\begin{proof}
For every dyadic $N$,
\[
        N^{\sigma-\alpha}\|P_NF\|_{L^2}
        =N^{\sigma-s_2}\,N^{s_2-\alpha}\|P_NF\|_{L^2}
        \le N^{\sigma-s_2}\|F\|_{H^{s_2-\alpha}}.
\]
Since $\sigma-s_2<0$, the dyadic supremum is bounded uniformly, and the finitely many low shells are absorbed into the constant.  Integrating in time gives \eqref{eq:LT1-sob-to-besov-source}.  The strictness of $\sigma<s_2$ is useful because it also gives cutoff tails in the Besov source norm whenever the Sobolev source tails converge.
\end{proof}

The following standard product bounds are used only with strict exponent margins supplied by Definition~\ref{def:admissible}; all endpoint losses are absorbed into the fixed small loss parameter.  Our paraproduct and remainder conventions are the periodic counterparts of \cite[Theorems~2.82 and~2.85]{BCD}; the proof below records the precise non-endpoint specialization used here.

\begin{lemma}[Product estimates used in the fixed point]\label{lem:det-product-toolbox}
Let the spatial dimension be three.
\begin{enumerate}[label=(\roman*),leftmargin=*]
\item For the low--high paraproduct with the rough factor high, one uses only the rough high factor unless an explicit dyadic calculation is made.  In particular, if $a<0$, $g\in\C^a$, and $f\in L^2$, then
\[
        \|f<g\|_{H^r}+\|f<g\|_{B_{2,\infty}^r}
        \lesssim \|f\|_{L^2}\|g\|_{\C^a},
        \qquad r<a,
\]
with the usual strict loss.  If $f\in B_{2,\infty}^\sigma$ with $\sigma>0$, this controls the low-frequency sum uniformly and gives cutoff tails, but it controls the low-frequency sum without creating a full $\sigma$-derivative gain in $f<g$.  The estimate used for $W<\Psi$ is the concrete shell bound in Lemma~\ref{lem:lowhigh-X}, whose active condition is $\sigma-\alpha<\alpha-3/2$, equivalently $\sigma<2\alpha-3/2$.
\item If $a+b>0$, resonant products with one rough H\"older factor satisfy the usual paracontrolled Sobolev estimate.  For $Y\circ\Psi$, Lemma~\ref{lem:highlow-resonant-Y} constructs the product at the positive Sobolev sum and then embeds it into the rougher Besov source space.
\item Sobolev multiplication is used in the following strict sufficient form.  If $u\in W^{a,q}$, $v\in W^{b,r}$, and the target is $H^s$, then it is enough that
\[
        a+b-s>3\left(\frac1q+\frac1r-\frac12\right)_+,
        \qquad
        a-s>3\left(\frac1q-\frac12\right)_+,
        \qquad
        b-s>3\left(\frac1r-\frac12\right)_+,
\]
after the fixed endpoint losses.  This is the form used for the $XX$, $XY$, and $YY$ products.  When the source is first estimated in $H^{s_2-\alpha}$, the Besov source follows from Lemma~\ref{lem:sobolev-to-besov-source}; this is the only use of the Sobolev-to-$B_{2,\infty}$ embedding in the deterministic quadratic sector.
\item Time smallness is obtained from non-endpoint H\"older inequalities on finite intervals.  Two $X$ factors produce $T^{1-2/p_X}$ after the $L_T^{p_X/2}\hookrightarrow L_T^1$ embedding.  The mixed $X$--$Y$ and $Y$--$Y$ products use the $L_T^{p_X}$ and $L_T^4$ components of the solution spaces and leave a positive power of $T$ under the admissible inequalities.
\end{enumerate}
\end{lemma}

\begin{proof}
All estimates are standard consequences of the Bony decomposition and Besov product calculus; we recall the precise dyadic forms used here and refer to \cite{BCD} for the general theory.  For (i), write
\[
        f<g=\sum_N S_{<c_0N}f\,P_Ng.
\]
Since $a<0$, $\|P_Ng\|_{L^\infty}\lesssim N^{-a}\|g\|_{\C^a}$ in the convention $\|g\|_{\C^a}=\sup_N N^a\|P_Ng\|_{L^\infty}$, and $\|S_{<c_0N}f\|_{L^2}\le \|f\|_{L^2}$.  Thus
\[
        N^r\|P_N(f<g)\|_{L^2}
        \lesssim N^{r-a}\|f\|_{L^2}\|g\|_{\C^a},
\]
which is summable, or uniformly bounded in the $B_{2,\infty}^r$ norm, whenever $r<a$ after the fixed strict loss.  If $f\in B_{2,\infty}^\sigma$, the same estimate with $S_{<c_0N}f$ replaced by the truncated Besov envelope gives the cutoff-tail version; Lemma~\ref{lem:lowhigh-X} supplies the corresponding shell calculation in the $W<\Psi$ sector.

For (ii), decompose the resonant product as $f\circ g=\sum_{M\sim N}P_MfP_Ng$.  Bernstein and Cauchy--Schwarz give the usual paraproduct estimate
\[
        \|f\circ g\|_{H^r}\lesssim \|f\|_{H^b}\|g\|_{\C^a}
        \quad\hbox{for } r<a+b,
\]
with strict inequality at the endpoint.  This proves the quoted Sobolev resonance rule.  In Lemma~\ref{lem:highlow-resonant-Y}, the required $B_{2,\infty}^{\sigma-\alpha}$ bound follows from the resulting Sobolev space by a strict embedding.

For (iii), use the paraproduct and remainder estimates in \cite[Theorems~2.82 and~2.85]{BCD}, whose proofs apply verbatim to the periodic Littlewood--Paley decomposition.  For completeness, at a comparable high frequency $2^j$, H\"older first places the product in $L^p$ with $1/p=1/q+1/r$; when $p<2$, Bernstein costs
\[
        2^{3j(1/q+1/r-1/2)},
\]
and when $p\ge2$ the finite-volume embedding costs no derivative.  Summation of the comparable high--high terms is therefore ensured by
\[
        a+b-s>3\left(\frac1q+\frac1r-\frac12\right)_+.
\]
In the two paraproducts, the geometric sum over the low blocks assigns the possible $L^q\to L^2$ and $L^r\to L^2$ Bernstein deficits to the low factor.  The respective strict summability conditions are
\[
        a-s>3\left(\frac1q-\frac12\right)_+,
        \qquad
        b-s>3\left(\frac1r-\frac12\right)_+.
\]
Square summation of the output blocks gives $B_{2,2}^s=H^s$.  These are exactly the three hypotheses in (iii); their strictness excludes all endpoint and quasi-Banach issues in the applications below.  Finally (iv) follows from H\"older in time on $[0,T]$: for example
\[
        \|uv\|_{L_T^1H^s}\le T^{1-1/p-1/q}\|u\|_{L_T^pW^{a,q_1}}
        \|v\|_{L_T^qW^{b,q_2}},
\]
whenever the spatial multiplication in (iii) maps the displayed Sobolev spaces to $H^s$ and $1/p+1/q<1$.  The admissible inequalities ensure that the exponent of $T$ is positive for the $XX$, $XY$, and $YY$ sectors used in the deterministic map.
\end{proof}

\begin{lemma}[H\"older multiplier estimates for the first Picard factor]\label{lem:positive-holder-multiplier}
Let $a>0$, $s\ge0$, and $r<\min\{a,s\}$.  Then multiplication by a $C^a$ function maps $H^s$ into $H^r$:
\begin{equation}\label{eq:holder-mult-Hr}
        \|fg\|_{H^r}\lesssim \|f\|_{\C^a}\|g\|_{H^s}.
\end{equation}
The exponent $r$ may be negative.  This is the case used for the source target $r=s_2-\alpha<0$ in products containing $V_i$ in Case F, and for $r=s_2-1<0$ in endpoint Case W.  Moreover, if $0<\sigma<a$, then
\begin{equation}\label{eq:holder-mult-Btwo}
        \|fg\|_{B_{2,\infty}^{\sigma}}
        \lesssim \|f\|_{\C^a}\|g\|_{B_{2,\infty}^{\sigma}}.
\end{equation}
Consequently $fg$ also belongs to $B_{2,\infty}^{\sigma-\alpha}$ on finite time intervals, and to $B_{2,\infty}^{\sigma-1}$ in Case W.  The same estimates hold in time-integrated norms and for differences.  In the application $f=V_i$, $a=\rho_V$, and the H\"older norm is the direct first-Picard output of Appendix~\ref{app:first-picard-phase}.
\end{lemma}

\begin{proof}
Decompose $fg=f<g+f>g+f\circ g$.  Since $\C^a\subset L^\infty$, the low--high term $f<g$ is bounded in $H^s$ and hence in $H^r$ whenever $r<s$, including negative $r$.  In the high--low term $f>g$, the high block of $f$ contributes $N^{-a}$ in $L^\infty$, while the low block of $g$ is bounded in $L^2$ by $\|g\|_{H^s}$; after applying the $H^r$ weight, the dyadic sum is controlled when $r<a$.  The resonant term has dyadic size $N^{-a-s}$ before the target weight, hence is controlled whenever $r<a+s$, which is automatic in the applications from $r<\min\{a,s\}$.  This proves \eqref{eq:holder-mult-Hr}.  The proof of \eqref{eq:holder-mult-Btwo} is identical with the dyadic supremum in place of the square sum; the high coefficient term requires precisely $\sigma<a$.  Since $B_{2,\infty}^{\sigma}\hookrightarrow B_{2,\infty}^{\sigma-\alpha}$, the displayed Besov estimate is stronger than the source norm used for such products.  Finally, $\C^a$ is an algebra for $a>0$, and the dyadic estimate $\|P_Nh\|_{L^2}\le \|P_Nh\|_{L^\infty}\lesssim N^{-a}\|h\|_{\C^a}$ on the compact torus gives $\C^a\hookrightarrow H^r\cap B_{2,\infty}^r$ for every $r<a$.  Difference estimates follow from $fg-\widetilde f\widetilde g=(f-\widetilde f)g+\widetilde f(g-\widetilde g)$.
\end{proof}

\paragraph{Source estimates.}
The following table records the source spaces and active conditions used by the mild map.
\begin{center}
\scriptsize
\begin{tabular}{@{}>{\raggedright\arraybackslash}p{0.19\textwidth}>{\raggedright\arraybackslash}p{0.19\textwidth}>{\raggedright\arraybackslash}p{0.25\textwidth}>{\raggedright\arraybackslash}p{0.25\textwidth}@{}}
\toprule
source term & target source space & input and estimate & active condition \\
\midrule
$W<\Psi$ & $L_T^1H^{s_1-\alpha}$ and $L_T^1B_{2,\infty}^{\sigma-\alpha}$ & low--high paraproduct, dyadic $L^2$ Besov envelope & $s_1<2\alpha-3/2-$, $\sigma<2\alpha-3/2-$ \\
$V>\Psi$, $X>\Psi$, $Y>\Psi$ & $\mathfrak S_T^{s_2,\sigma}$ & componentwise high--low estimates; direct dyadic $B_{2,\infty}$ source & $s_2<11\alpha/2-9/2-$ for $V$, $s_2<s_1+2\alpha-3/2-$ for $X$, $\alpha>3/4$ for $Y$ and Besov \\
$Y\circ\Psi$ & $\mathfrak S_T^{s_2,\sigma}$ & Sobolev resonant estimate followed by source-space embedding & $s_2>3/2-\alpha+$ and $\alpha>3/4$ \\
$T^{i;j,k}(W)$, $(i;j,k)\in\mathfrak M$ & $\mathfrak S_T^{s_2,\sigma}$ & Section~\ref{sec:random-operators} mixed-operator estimates & $s_2<4\alpha-3-$ and $3-3\alpha<\sigma<4\alpha-3$ \\
$\Gamma_i$ & $\mathfrak S_T^{s_2,\sigma}$ & cubic Wick construction and embeddings & $s_2<11\alpha/2-9/2-$, $\sigma<11\alpha/2-9/2-$ \\
$V_iV_j$, $V_iX_j$, $V_iY_j$ & $\mathfrak S_T^{s_2,\sigma}$ & $V_i\in C_T\C^{\rho_V}$ through Lemma~\ref{lem:positive-holder-multiplier}; $V_i\in L_T^\infty B_{2,\infty}^{\rho_V}$ only for the Besov-input route & $\sigma<\rho_V$, $s_2<\alpha$ \\
$X_1X_2$ & Case F: $L_T^1H^{s_2-\alpha}$; Case W: $\mathfrak D_T^{s_2}$, plus Besov source & two $X$ Strichartz factors; endpoint dual source route & Case F: $s_2<2s_1+\alpha+2\alpha/p_X-3/2-$; Case W: $s_2<s_1+1/4$ \\
$X_iY_j$, $Y_iY_j$ & Case F: $L_T^1H^{s_2-\alpha}$; Case W: $\mathfrak D_T^{s_2}$, plus Besov source & mixed Strichartz products; endpoint dual source route & Case F: lower margin in \eqref{eq:admissible-quadratic}; Case W: \eqref{eq:admissible-case-W} \\
\bottomrule
\end{tabular}
\end{center}
The corresponding estimates and their Lipschitz forms are proved in Lemmas~\ref{lem:lowhigh-X}--\ref{lem:quadratic-Y} and Proposition~\ref{prop:quantitative-source-stability}.

\begin{lemma}[Endpoint deterministic quadratic estimates]\label{lem:endpoint-quadratic-estimates}
Assume Case W, and put $\gamma=s_2-\frac12>0$.  If $X_a,X_b\in X^{s_1}_{T,\mathrm{end}}$ and $Y_a,Y_b\in Y^{s_2}_{T,\mathrm{end}}$, then
\begin{align}
 \|\langle D\rangle^\gamma(X_aX_b)\|_{L_{t,x}^{4/3}}
 &\lesssim T^\delta
 \|X_a\|_{X^{s_1}_{T,\mathrm{end}}}\|X_b\|_{X^{s_1}_{T,\mathrm{end}}},
 \label{eq:endpoint-XX-dual}\\
 \|\langle D\rangle^\gamma(X_aY_b)\|_{L_{t,x}^{4/3}}
 &\lesssim T^\delta
 \|X_a\|_{X^{s_1}_{T,\mathrm{end}}}\|Y_b\|_{Y^{s_2}_{T,\mathrm{end}}},
 \label{eq:endpoint-XY-dual}\\
 \|\langle D\rangle^\gamma(Y_aY_b)\|_{L_{t,x}^{4/3}}
 &\lesssim T^\delta
 \|Y_a\|_{Y^{s_2}_{T,\mathrm{end}}}\|Y_b\|_{Y^{s_2}_{T,\mathrm{end}}}.
 \label{eq:endpoint-YY-dual}
\end{align}
Consequently these deterministic quadratic sources are propagated in the $Y$ equation by the dual conic source route \eqref{eq:conic-dual-Y}.  If
\[
        \mathcal N_a:=\|X_a\|_{X^{s_1}_{T,\mathrm{end}}}
        +\|Y_a\|_{Y^{s_2}_{T,\mathrm{end}}},
        \qquad
        \mathcal N_b:=\|X_b\|_{X^{s_1}_{T,\mathrm{end}}}
        +\|Y_b\|_{Y^{s_2}_{T,\mathrm{end}}},
\]
then the endpoint Besov source bound is
\begin{equation}\label{eq:endpoint-besov-source-from-l43}
 \sum_{A\in\{X_a,Y_a\}}\sum_{B\in\{X_b,Y_b\}}
 \|AB\|_{L_T^1B_{2,\infty}^{\sigma-1}}
 \lesssim T^\delta\mathcal N_a\mathcal N_b.
\end{equation}
If $V\in C_T\mathcal C^{1/2-}$, the terms $VX$, $VY$, and $VV$ are controlled in $L_T^1H^{s_2-1}\cap L_T^1B_{2,\infty}^{\sigma-1}$ by Lemma~\ref{lem:positive-holder-multiplier}, because $s_2-1<0$ and $\sigma<1/4$.  All estimates have the corresponding Lipschitz form.
\end{lemma}

\begin{proof}
The endpoint parameter condition gives $\gamma=s_2-1/2<s_1-1/4$.  Since
\[
        X\in L_T^8W^{s_1-1/4,8/3},
        \qquad
        Y\in L_T^4W^{s_2-1/2,4}=L_T^4W^{\gamma,4},
\]
the fractional Leibniz rule and H\"older imply \eqref{eq:endpoint-XX-dual}: put the $\gamma$ derivatives on one $X$ factor, use $L_x^{8/3}\cdot L_x^{8/3}\subset L_x^{4/3}$, and use $L_T^8\cdot L_T^8\subset L_T^4\hookrightarrow L_T^{4/3}$ on $[0,T]$.  For \eqref{eq:endpoint-XY-dual}, put the $\gamma$ derivatives on $Y_b$ and use
\[
        L_x^{8/3}\cdot L_x^4\subset L_x^{8/5}\hookrightarrow L_x^{4/3},
        \qquad
        L_T^8\cdot L_T^4\subset L_T^{8/3}\hookrightarrow L_T^{4/3}.
\]
For \eqref{eq:endpoint-YY-dual}, put the derivative on one $Y$ factor and use $L_x^4\cdot L_x^4\subset L_x^2\hookrightarrow L_x^{4/3}$ and $L_T^4\cdot L_T^4\subset L_T^2\hookrightarrow L_T^{4/3}$.  The finite interval embeddings give a positive power $T^\delta$.

For the Besov source component, the elementary embedding $L_x^{4/3}\hookrightarrow B_{2,\infty}^{-3/4}$ on $\mathbb T^3$ and the endpoint restriction $\sigma<1/4$ give $L_x^{4/3}\hookrightarrow B_{2,\infty}^{\sigma-1}$.  The terms containing $V$ use the H\"older multiplier estimate with the direct first-Picard bound $V\in C_T\mathcal C^{1/2-}$.  Difference estimates follow by the usual bilinear expansion.
\end{proof}

\begin{lemma}[Endpoint Case W rough stochastic source routes]\label{lem:caseW-rough-stochastic-routes}
Let \(\alpha=1\), \(\rho_\Psi=-1/2-\), and use the endpoint augmented spaces of \eqref{eq:case-W-spaces-display}--\eqref{eq:Ytilde-space}.  For every admissible Case W parameter choice,
\begin{align}
\|W<\Psi_i\|_{L_T^1H^{s_1-1}}
+\|W<\Psi_i\|_{L_T^1B_{2,\infty}^{\sigma-1}}
&\lesssim
T\|\Psi_i\|_{C_T\mathcal C^{\rho_\Psi}}\|W\|_{E_T^{2,\sigma}},
\label{eq:caseW-WlowPsi-route}\\
\|Y\circ\Psi_i\|_{\mathfrak S_{T,W}^{s_2,\sigma}}
+\|Y>\Psi_i\|_{\mathfrak S_{T,W}^{s_2,\sigma}}
&\lesssim
T\|\Psi_i\|_{C_T\mathcal C^{\rho_\Psi}}
\|Y\|_{\widetilde Y_T^{s_2,\sigma}},
\label{eq:caseW-Ypsi-route}\\
\|X>\Psi_i\|_{\mathfrak S_{T,W}^{s_2,\sigma}}
&\lesssim
T\|\Psi_i\|_{C_T\mathcal C^{\rho_\Psi}}
\|X\|_{\widetilde X_T^{s_1,\sigma}}.
\label{eq:caseW-XgreaterPsi-route}
\end{align}
The same estimates hold for differences.
\end{lemma}

\begin{proof}
These are the endpoint analogues of Lemma~\ref{lem:caseF-rough-stochastic-routes}.  For the low--high term, Lemma~\ref{lem:det-product-toolbox} gives only the regularity of the rough high factor, namely \(-1/2-\).  The endpoint targets satisfy
\[
        s_1-1<-\frac12,
        \qquad
        \sigma-1<-\frac12,
\]
because the Case W window has \(s_1<1/2\) and \(\sigma<1/4\).  This proves \eqref{eq:caseW-WlowPsi-route}.  The high--low term \(Y>\Psi_i\) has the required Sobolev and Besov regularities by the direct dyadic estimate of Lemma~\ref{lem:highlow-resonant-Y}.  For the resonance, choose \(\rho_\Psi<-1/2\) close enough to \(-1/2\).  Since \(s_2+\rho_\Psi>0\), the standard Sobolev resonant estimate constructs \(Y\circ\Psi_i\) in \(H^{s_2+\rho_\Psi}\); because \(\sigma<s_2\), this space embeds into both \(H^{s_2-1}\) and \(B_{2,\infty}^{\sigma-1}\).  This proves \eqref{eq:caseW-Ypsi-route}.  For \(X>\Psi_i\), the Sobolev product regularity is \(s_1-1/2-\).  It is enough to have \(s_2<s_1+1/2\), and this is implied with room by the endpoint window \(s_2<s_1+1/4-\); the direct high--low Besov estimate gives the second source component.  The factor \(T\) comes from the \(L_T^1\) source norm, and the difference estimates are obtained by replacing one factor by its difference.
\end{proof}

\begin{lemma}[Case F rough stochastic source routes]\label{lem:caseF-rough-stochastic-routes}
Let $\rho_\Psi=\alpha-3/2-$.  For admissible Case F parameters,
\begin{align}
\|W<\Psi_i\|_{L_T^1H^{s_1-\alpha}}
+\|W<\Psi_i\|_{L_T^1B_{2,\infty}^{\sigma-\alpha}}
&\lesssim T\|\Psi_i\|_{C_T\mathcal C^{\rho_\Psi}}\|W\|_{E_T^{2,\sigma}},\label{eq:caseF-WlowPsi-route}\\
\|Y\circ\Psi_i\|_{\mathfrak S_T^{s_2,\sigma}}
+\|Y>\Psi_i\|_{\mathfrak S_T^{s_2,\sigma}}
&\lesssim T\|\Psi_i\|_{C_T\mathcal C^{\rho_\Psi}}\|Y\|_{\widetilde Y_T^{s_2,\sigma}},\label{eq:caseF-Ypsi-route}\\
\|X>\Psi_i\|_{\mathfrak S_T^{s_2,\sigma}}
&\lesssim T\|\Psi_i\|_{C_T\mathcal C^{\rho_\Psi}}\|X\|_{\widetilde X_T^{s_1,\sigma}}.\label{eq:caseF-XgreaterPsi-route}
\end{align}
The same estimates hold for differences.
\end{lemma}

\begin{proof}
The low--high estimate is the concrete route of Lemma~\ref{lem:lowhigh-X}; equivalently, the paraproduct part of Lemma~\ref{lem:det-product-toolbox} assigns the rough high-factor regularity as the active exponent.  The high rough factor contributes only $\rho_\Psi=\alpha-3/2-$, which is enough because
\[
        s_1-\alpha<\rho_\Psi
        \quad\Longleftrightarrow\quad s_1<2\alpha-\frac32,
        \qquad
        \sigma-\alpha<\rho_\Psi
        \quad\Longleftrightarrow\quad \sigma<2\alpha-\frac32.
\]
Both inequalities are built into the admissible window.  The direct dyadic high--low estimate in Lemma~\ref{lem:highlow-resonant-Y} treats $Y>\Psi_i$.  For $Y\circ\Psi_i$, choose $\rho_\Psi<\alpha-3/2$ close enough to the endpoint.  The admissible lower bound gives $s_2+\rho_\Psi>0$, so the Sobolev resonance is well defined in $H^{s_2+\rho_\Psi}$; this embeds into both target source spaces because $\alpha>3/4$ and $\sigma<s_2$.  Finally, $X>\Psi_i$ has Sobolev regularity $s_1+\rho_\Psi$, and the lower bound $s_1>s_2+3/2-2\alpha+$ in \eqref{eq:s1-window} is exactly the condition $s_1+\rho_\Psi>s_2-\alpha$.  Its Besov component follows from the direct high--low dyadic estimate.  The time factor is the $L_T^1$ source length, and bilinear expansion gives the difference estimates.
\end{proof}

\begin{lemma}[Deterministic product margins in Case F]\label{lem:caseF-product-margins}
Set
\[
        r=s_2-\alpha,
        \qquad a_X=s_1-\kappa_X,
        \qquad a_Y=s_2-\kappa_Y.
\]
For every admissible Case F choice, after the fixed losses,
\begin{align}
\|X_1X_2\|_{L_T^1H^{r}}
&\lesssim
T^{1-2/p_X}
\|X_1\|_{L_T^{p_X}W^{a_X,q_X}}
\|X_2\|_{L_T^{p_X}W^{a_X,q_X}},\label{eq:caseF-XX-margin}\\
\|X_iY_j\|_{L_T^1H^{r}}
&\lesssim
T^{1-1/p_X-1/4}
\|X_i\|_{L_T^{p_X}W^{a_X,q_X}}
\|Y_j\|_{L_T^4W^{a_Y,4}},\label{eq:caseF-XY-margin}\\
\|Y_1Y_2\|_{L_T^1H^{r}}
&\lesssim
T^{1/2}
\|Y_1\|_{L_T^4W^{a_Y,4}}
\|Y_2\|_{L_T^4W^{a_Y,4}}.\label{eq:caseF-YY-margin}
\end{align}
Consequently all pure deterministic quadratic terms lie in
$L_T^1H^{s_2-\alpha}\hookrightarrow L_T^1B_{2,\infty}^{\sigma-\alpha}$ and carry a positive power of $T$.
\end{lemma}

\begin{proof}
We use the strict product criterion in Lemma~\ref{lem:det-product-toolbox}.  Write
\[
        \delta_X=\frac32-\frac3{q_X},\qquad
        \kappa_X=\delta_X-\frac\alpha{p_X}+\eps_{\rm Str},
        \qquad
        \kappa_Y=\frac{3-\alpha}{4}+\eps_{\rm Str}.
\]
For $XX$ the target is $H^r$ with $r=s_2-\alpha$ and both factors lie in $W^{a_X,q_X}$.  The active sum condition is
\[
        2a_X-r>3\left(\frac2{q_X}-\frac12\right)_+.
\]
For the near-$(4,3)$ choice this is
\[
        2(s_1-\kappa_X)-(s_2-\alpha)>
        \frac6{q_X}-\frac32,
\]
which rearranges to
\[
        2s_1>s_2-\alpha-\frac{2\alpha}{p_X}+\frac32+O(\eps_{\rm Str}).
\]
This is the third lower bound in \eqref{eq:s1-window}, or equivalently the first inequality in \eqref{eq:admissible-quadratic}.  The one-sided conditions are weaker: they require $a_X>r$, and this follows from the high--low lower bound $s_1>s_2+3/2-2\alpha+$ after the same near-$(4,3)$ choice.

For $XY$, with $X\in W^{a_X,q_X}$ and $Y\in W^{a_Y,4}$, the sum condition is
\[
        a_X+a_Y-r>3\left(\frac1{q_X}+\frac14-\frac12\right)_+.
\]
Substituting $a_X=s_1-\kappa_X$ and $a_Y=s_2-\kappa_Y$ gives the sufficient lower bound
\[
        s_1>\frac32-\frac{5\alpha}{4}-\frac\alpha{p_X}+O(\eps_{\rm Str}),
\]
which is the fourth lower bound in \eqref{eq:s1-window}.  The one-sided condition on the $X$ factor is
\[
        a_X-r>0
        \quad\Longleftarrow\quad
        s_1>s_2-\alpha+\kappa_X+O(\eps_{\rm Str}),
\]
and is implied by $s_1>s_2+3/2-2\alpha+$ because $\delta_X-\alpha/p_X<3/2-\alpha$ for the chosen non-endpoint pair.  The one-sided condition on the $Y$ factor is
\[
        a_Y-r=\alpha-\kappa_Y=\frac{5\alpha-3}{4}-O(\eps_{\rm Str})>0,
\]
which follows from $\alpha>3/5$.

For $YY$, both factors lie in $W^{a_Y,4}$.  The sum condition reduces to
\[
        2a_Y-r>0
        \quad\Longleftrightarrow\quad
        s_2>\frac32-\frac{3\alpha}{2}+O(\eps_{\rm Str}),
\]
which is weaker than the imposed lower bound $s_2>3/2-\alpha+$.  The one-sided condition is again $a_Y-r=(5\alpha-3)/4-O(\eps_{\rm Str})>0$.

The time factors are exactly H\"older on the finite interval:
\[
        L_T^{p_X}\cdot L_T^{p_X}\subset L_T^1
        \text{ with }T^{1-2/p_X},\quad
        L_T^{p_X}\cdot L_T^4\subset L_T^1
        \text{ with }T^{1-1/p_X-1/4},\quad
        L_T^4\cdot L_T^4\subset L_T^1
        \text{ with }T^{1/2}.
\]
Thus \eqref{eq:caseF-XX-margin}--\eqref{eq:caseF-YY-margin} hold.  Finally, since $\sigma<s_2$, Lemma~\ref{lem:sobolev-to-besov-source} converts the Sobolev source into the Besov source.
\end{proof}

\subsection{Deterministic estimates in the augmented topology}
The estimates below combine standard paraproduct bounds, the Sobolev Strichartz input of Proposition~\ref{prop:strichartz}, the elementary Besov propagation of Lemma~\ref{lem:elementary-besov-propagation}, and the augmented enhanced-data norm.

\begin{lemma}[Low--high source into the X component]\label{lem:lowhigh-X}
For $W\in E_T^{2,\sigma}$ and each $j$,
\begin{equation}\label{eq:lowhigh-X-estimate}
\|I_i(W<\Psi_j)\|_{\widetilde X_T^{s_1,\sigma}}
\lesssim T^\delta
\|\Psi_j\|_{C_T\C^{\alpha-3/2-}}
\|W\|_{E_T^{2,\sigma}}
\end{equation}
for some $\delta>0$.
\end{lemma}

\begin{proof}
Write the paraproduct as
\[
        W<\Psi_j=\sum_N S_{<c_0N}W\,P_N\Psi_j,
\]
where the output frequency is comparable to $N$.  For the Sobolev source, the dyadic contribution is
\[
 N^{s_1-\alpha}
 \|S_{<c_0N}W\,P_N\Psi_j\|_{L_T^\infty L^2}
 \lesssim
 N^{s_1-\alpha}\|W\|_{C_TL^2}
 N^{3/2-\alpha+}\|\Psi_j\|_{C_T\C^{\alpha-3/2-}}.
\]
The square sum in $N$ is finite under $s_1<2\alpha-3/2-$, giving
\[
 \|W<\Psi_j\|_{L_T^1H^{s_1-\alpha}}
 \lesssim
 T\|W\|_{C_TL^2}\|\Psi_j\|_{C_T\C^{\alpha-3/2-}}.
\]
The branch-appropriate linear estimate gives the Sobolev part of $X$: use \eqref{eq:strichartz-input-X} in Case F and \eqref{eq:conic-duh-X} in Case W.

For the $B_{2,\infty}$ source part, we keep the dyadic $L^2$ envelope of the low factor.  Since $\sigma>0$,
\[
        \|S_{<c_0N}W\|_{L_T^\infty L^2}
        \le \left(\sum_{Q<c_0N}Q^{-2\sigma}\right)^{1/2}
        \|W\|_{L_T^\infty\Btwo^\sigma}
        \lesssim_\sigma \|W\|_{L_T^\infty\Btwo^\sigma}.
\]
Therefore
\begin{align*}
 N^{\sigma-\alpha}\|S_{<c_0N}W\,P_N\Psi_j\|_{L_T^1L^2}
 &\lesssim
 T N^{\sigma-\alpha}\|S_{<c_0N}W\|_{L_T^\infty L^2}
        \|P_N\Psi_j\|_{C_TL^\infty}\\
 &\lesssim
 T N^{\sigma+3/2-2\alpha+}
 \|W\|_{L_T^\infty\Btwo^\sigma}
 \|\Psi_j\|_{C_T\C^{\alpha-3/2-}}.
\end{align*}
This is uniformly bounded in the $B_{2,\infty}^{\sigma-\alpha}$ dyadic supremum and has a high-frequency cutoff tail because $\sigma<2\alpha-3/2-$.  Lemma~\ref{lem:elementary-besov-propagation} then gives both the $L_T^\infty\Btwo^\sigma$ and velocity $L_T^\infty\Btwo^{\sigma-\alpha}$ components of $I_i(W<\Psi_j)$.  The time smallness is the displayed $T$ factor, or any smaller non-endpoint power $T^\delta$ after combining with the Sobolev Strichartz component.
\end{proof}

\begin{lemma}[High--low and resonant deterministic product estimates]\label{lem:highlow-resonant-Y}
Under the admissible conditions, the high--low source is controlled componentwise as follows:
\begin{align}
\|I_i(V_\ell>\Psi_j)\|_{\widetilde Y_T^{s_2,\sigma}}
&\lesssim T^\delta
\|\Psi_j\|_{C_T\C^{\alpha-3/2-}}
\left(\|V_\ell\|_{C_T\C^{\rho_V}}+
      \|V_\ell\|_{L_T^\infty B_{2,\infty}^{\rho_V}}\right),
\label{eq:highlow-V-Y-estimate}\\
\|I_i(X_\ell>\Psi_j)\|_{\widetilde Y_T^{s_2,\sigma}}
&\lesssim T^\delta
\|\Psi_j\|_{C_T\C^{\alpha-3/2-}}
\|X_\ell\|_{\widetilde X_T^{s_1,\sigma}},
\label{eq:highlow-X-Y-estimate}\\
\|I_i(Y_\ell>\Psi_j)\|_{\widetilde Y_T^{s_2,\sigma}}
&\lesssim T^\delta
\|\Psi_j\|_{C_T\C^{\alpha-3/2-}}
\|Y_\ell\|_{\widetilde Y_T^{s_2,\sigma}}.
\label{eq:highlow-Y-Y-estimate}
\end{align}
Moreover,
\begin{equation}\label{eq:Y-resonant-estimate}
\|I_i(Y_\ell\circ\Psi_j)\|_{\widetilde Y_T^{s_2,\sigma}}
\lesssim T^\delta
\|\Psi_j\|_{C_T\C^{\alpha-3/2-}}
\|Y_\ell\|_{\widetilde Y_T^{s_2,\sigma}}.
\end{equation}
Consequently, the full high--low part with \(W_\ell=V_\ell+X_\ell+Y_\ell\) satisfies the corresponding bound with the right-hand side given by the sum of the three component norms above.  For the Sobolev source topology this estimate is used with the structured input \(W=V+X+Y\); the direct \(E_T^{2,\sigma}\) bound is used for the Besov source component.
\end{lemma}

\begin{proof}
The high--low term has the dyadic form
\[
        U>\Psi_j=\sum_N P_NU\,S_{<c_0N}\Psi_j.
\]
Since \(\Psi_j\in\C^{\alpha-3/2-}\) has negative regularity,
\[
        \|S_{<c_0N}\Psi_j\|_{L^\infty}
        \lesssim N^{3/2-\alpha+}\|\Psi_j\|_{C_T\C^{\alpha-3/2-}}.
\]
We first record the direct Besov source estimate, which is valid for any high factor \(U\) carrying the dyadic Besov bound.  If \(U\in L_T^\infty B_{2,\infty}^{\sigma}\), then
\[
\begin{aligned}
        N^{\sigma-\alpha}
        \|P_N(U>\Psi_j)\|_{L_T^1L^2}
        &\lesssim T N^{\sigma-\alpha}N^{-\sigma}N^{3/2-\alpha+}
        \|U\|_{L_T^\infty B_{2,\infty}^{\sigma}}
        \|\Psi_j\|_{C_T\C^{\alpha-3/2-}}  \\
        &=T N^{3/2-2\alpha+}
        \|U\|_{L_T^\infty B_{2,\infty}^{\sigma}}
        \|\Psi_j\|_{C_T\C^{\alpha-3/2-}}.
\end{aligned}
\]
The exponent is negative for \(\alpha>3/4\).  Thus the high--low Besov source bound is obtained directly in the dyadic \(L^2\) topology.

It remains to give the Sobolev source bounds for the three components of \(W\).  For \(Y_\ell\), the dyadic bound
\[
        N^{s_2-\alpha}
        \|P_NY_\ell\,S_{<c_0N}\Psi_j\|_{L_T^1L^2}
        \lesssim T N^{3/2-2\alpha+}
        \|Y_\ell\|_{C_TH^{s_2}}
        \|\Psi_j\|_{C_T\C^{\alpha-3/2-}}
\]
is square summable because \(\alpha>3/4\).  For \(X_\ell\), the same computation with \(\|P_NX_\ell\|_{L^2}\lesssim N^{-s_1}\|X_\ell\|_{C_TH^{s_1}}\) gives the shell
\[
        N^{s_2-s_1+3/2-2\alpha+}
        \|X_\ell\|_{C_TH^{s_1}}
        \|\Psi_j\|_{C_T\C^{\alpha-3/2-}},
\]
which is summable under the admissible condition \(s_2<s_1+2\alpha-3/2-\).  For \(V_\ell\), we use the direct first-Picard H\"older estimate in the Sobolev route:
\[
        \|P_NV_\ell\|_{L^\infty}\lesssim N^{-\rho_V}\|V_\ell\|_{C_T\C^{\rho_V}}.
\]
Hence the Sobolev shell is
\[
        N^{s_2-\alpha}N^{-\rho_V}N^{3/2-\alpha+}
        \|V_\ell\|_{C_T\C^{\rho_V}}
        \|\Psi_j\|_{C_T\C^{\alpha-3/2-}}
        =N^{s_2+9/2-11\alpha/2+}
        \|V_\ell\|_{C_T\C^{\rho_V}}
        \|\Psi_j\|_{C_T\C^{\alpha-3/2-}},
\]
which is summable by \(s_2<11\alpha/2-9/2-\).  In the Besov source route for \(V_\ell>\Psi_j\) one may use either the stronger H\"older bound or the envelope \(V_\ell\in L_T^\infty B_{2,\infty}^{\rho_V}\); since \(\sigma<\rho_V\), this gives the same direct dyadic source bound above.

Combining these Sobolev estimates with the direct Besov source estimate and then applying the branch-appropriate linear estimate--Proposition~\ref{prop:strichartz} in Case F and Proposition~\ref{prop:gko-conic-strichartz} in Case W--together with Lemma~\ref{lem:elementary-besov-propagation} proves \eqref{eq:highlow-V-Y-estimate}--\eqref{eq:highlow-Y-Y-estimate}.  The time factor is the displayed \(T\) or a smaller non-endpoint power after inserting the Strichartz components.

For the resonant product, choose $\rho_\Psi<\alpha-3/2$ sufficiently close to $\alpha-3/2$.  The admissible window gives $s_2+\rho_\Psi>0$, and the standard resonant estimate yields
\begin{equation}\label{eq:Y-Psi-standard-resonant}
 \|Y_\ell\circ\Psi_j\|_{L_T^1H^{s_2+\rho_\Psi}}
 \lesssim
 T\|Y_\ell\|_{C_TH^{s_2}}\|\Psi_j\|_{C_T\mathcal C^{\rho_\Psi}}.
\end{equation}
Moreover,
\[
 s_2+\rho_\Psi>\max\{s_2-\alpha,\sigma-\alpha\},
\]
because $\alpha>3/4$ and $\sigma<s_2$.  Hence
\[
 H^{s_2+\rho_\Psi}\hookrightarrow
 H^{s_2-\alpha}\cap B_{2,\infty}^{\sigma-\alpha}.
\]
This proves the two source bounds in \eqref{eq:Y-resonant-estimate}, including the high--high to low-output sector.  Cutoff convergence and difference estimates follow from continuity of the same bilinear map.  Applying the branch-appropriate linear estimate completes the proof.
\end{proof}

\begin{lemma}[Mixed operator terms and operator-data differences]\label{lem:mixed-operator-Y}
For every generated block $T^{i;j,k}$ appearing in \eqref{eq:generated-mixed-blocks} and either outer equation index $\iota\in\{1,2\}$,
\begin{equation}\label{eq:mixed-operator-Y-estimate}
\|I_\iota(T^{i;j,k}(W))\|_{\widetilde Y_T^{s_2,\sigma}}
\lesssim
\mu_T^{\rm op}\|W\|_{E_T^{2,\sigma}},
\qquad \mu_T^{\rm op}\text{ as in Definition~\ref{def:localized-mixed-size}.}
\end{equation}
The outer index $\iota$ is independent of the inner Duhamel speed $i$ appearing in the mixed operator.  This is the form used in the $Y_1$ and $Y_2$ equations, where the same source sector may be acted on by either outer propagator.

Moreover, if two enhanced data give mixed operators $T$ and $\widetilde T$ on the same interval, then
\begin{equation}\label{eq:mixed-operator-Y-difference}
\begin{aligned}
\|I_\iota(T(W)-\widetilde T(\widetilde W))\|_{\widetilde Y_T^{s_2,\sigma}}
&\lesssim \mu_T^{\rm op}\|W-\widetilde W\|_{E_T^{2,\sigma}}  \\
&\quad +\Bigl(T\|T-\widetilde T\|_{\mathcal L_H}
       +\|T-\widetilde T\|_{\mathcal L_B}\Bigr)
       \|\widetilde W\|_{E_T^{2,\sigma}},
\end{aligned}
\end{equation}
with $\mathcal L_H$ and $\mathcal L_B$ as in Definition~\ref{def:enhanced-data-distance}.  The first term is the contraction contribution; the second is an enhanced-data perturbation and is therefore measured by $d_T$ rather than by a local-smallness factor.
\end{lemma}

\begin{proof}
The estimates of Section~\ref{sec:random-operators} are source estimates, independent of the outer propagator index.  Once the source is placed in $L_T^1H^{s_2-\alpha}\cap L_T^1B_{2,\infty}^{\sigma-\alpha}$, either $I_1$ or $I_2$ is handled by the corresponding linear estimate with constants depending only on the individual speed.  For the Sobolev component, Theorem~\ref{thm:pathwise-fluct} supplies the source bound in $C_TH^{s_2-\alpha}$; before applying \eqref{eq:strichartz-input-Y} we use
\[
        \|T^{i;j,k}(W)\|_{L_T^1H^{s_2-\alpha}}
        \le T M_T^{\rm op,H}\|W\|_{E_T^{2,\sigma}}.
\]
For the $L_T^\infty\Btwo^\sigma$ and velocity $L_T^\infty\Btwo^{\sigma-\alpha}$ components, \eqref{eq:besov-input-dh} uses the Besov source operator norm directly:
\[
        \|T^{i;j,k}(W)\|_{L_T^1\Btwo^{\sigma-\alpha}}
        \le \bigl(T M_{T,\mathrm{diag}}^{\rm op,B}+M_{T,\mathrm{cen}}^{\rm op,B}\bigr)
        \|W\|_{E_T^{2,\sigma}}.
\]
Both Besov pieces are supplied by Section~\ref{sec:random-operators}: Proposition~\ref{prop:det-contraction} treats the deterministic diagonal and Theorem~\ref{thm:pathwise-fluct} treats the centered part.  Their local smallness is recorded in Propositions~\ref{prop:det-diagonal-local-smallness} and~\ref{prop:centered-besov-local-smallness}, proving \eqref{eq:mixed-operator-Y-estimate}.

For the difference estimate, decompose
\[
        T(W)-\widetilde T(\widetilde W)=T(W-\widetilde W)+(T-\widetilde T)(\widetilde W).
\]
The first term is estimated exactly as above.  For the second term, the Sobolev source line gives a factor $T\|T-\widetilde T\|_{\mathcal L_H}$ because the fixed-time $C_TH^{s_2-\alpha}$ operator norm is converted to $L_T^1H^{s_2-\alpha}$ before the outer Duhamel estimate.  The Besov source line gives $\|T-\widetilde T\|_{\mathcal L_B}$ directly.  These are precisely the mixed-operator components of the enhanced-data distance in Definition~\ref{def:enhanced-data-distance}, up to the harmless factor $T\le1$ in the Sobolev line.
\end{proof}

\begin{lemma}[Cubic symbols in the augmented source spaces]\label{lem:cubic-source-compatibility}
Let
\[
        \rho_\Gamma:=\frac{9\alpha}{2}-\frac92-.
\]
The cubic construction of Theorem~\ref{prop:cubic-resonance} gives $\Gamma_i\in C_T\C^{\rho_\Gamma}$.  Under the admissible inequalities,
\begin{equation}\label{eq:cubic-augmented-source}
 \|\Gamma_i\|_{L_T^1H^{s_2-\alpha}}
 +\|\Gamma_i\|_{L_T^1\Btwo^{\sigma-\alpha}}
 \lesssim_T \|\Gamma_i\|_{C_T\C^{\rho_\Gamma}}.
\end{equation}
\end{lemma}

\begin{proof}
The Sobolev source condition is
\[
        s_2-\alpha<\rho_\Gamma
        \quad\Longleftrightarrow\quad
        s_2<\frac{11\alpha}{2}-\frac92-,
\]
which is included in \eqref{eq:admissible-high} in Case F.  In Case W, $\alpha=1$ and $\rho_\Gamma=0-$, while the endpoint window has $s_2<1$, so the same inequality holds with room.  The Besov source condition is
\[
        \sigma-\alpha<\rho_\Gamma
        \quad\Longleftrightarrow\quad
        \sigma<\frac{11\alpha}{2}-\frac92-.
\]
In Case F this follows from $\sigma<\rho_V=7\alpha/2-3-$ because $\alpha>3/4$; in Case W it follows from $\sigma<1/4$.  Since the time interval is finite, $C_T$ control implies the corresponding $L_T^1$ source control.
\end{proof}

\begin{lemma}[Quadratic and cubic source terms]\label{lem:quadratic-Y}
For admissible parameters,
\begin{equation}\label{eq:quadratic-total-estimate}
\|I_i(W_1W_2+\Gamma_1+\Gamma_2)\|_{\widetilde Y_T^{s_2,\sigma}}
\lesssim M_T(\Xi)
 +T^\delta P\bigl(M_T^{\rm bd}(\Xi,y)\bigr)(1+\|Z\|_{\calZ_T})^2,
\end{equation}
where $P$ is a fixed polynomial.  The same estimate holds for differences, with the right side linear in the difference and polynomial in the two bounded control norms.
\end{lemma}

\begin{proof}
The cubic source bounds are recorded in Lemma~\ref{lem:cubic-source-compatibility}.  In Case F, Proposition~\ref{prop:strichartz} propagates the Sobolev component; in Case W the corresponding source is propagated by Proposition~\ref{prop:gko-conic-strichartz}.  Lemma~\ref{lem:elementary-besov-propagation} propagates the Besov component in both branches.

For the quadratic source, write
\[
        W_i=V_i+X_i+Y_i,
        \qquad
        W_1W_2=V_1V_2+V_1(X_2+Y_2)+V_2(X_1+Y_1)+(X_1+Y_1)(X_2+Y_2).
\]
For the augmented $B_{2,\infty}$ source norm we use the source-specific estimates proved above.  In Case F the non-$V$ deterministic quadratic terms are first placed in $L_T^1H^{s_2-\alpha}$, and the admissible condition $\sigma<s_2$ yields the embedding $H^{s_2-\alpha}\hookrightarrow B_{2,\infty}^{\sigma-\alpha}$.  In Case W the same non-$V$ terms use the endpoint dual conic route for the $Y$-Strichartz component and the elementary embedding $L^{4/3}_x\hookrightarrow B_{2,\infty}^{\sigma-1}$, valid for $\sigma<1/4$, for the Besov source component.  The terms containing $V_i$ are estimated by Lemma~\ref{lem:positive-holder-multiplier}, using the direct H\"older control $V_i\in C_T\C^{\rho_V}$; in particular $V_iV_j$ is handled through the algebra property of $C^{\rho_V}$ and the embedding $C^{\rho_V}\hookrightarrow H^r\cap B_{2,\infty}^r$ for $r<\rho_V$.  The separate $B_{2,\infty}^{\rho_V}$ envelope is used only for Besov-input dyadic estimates.

For the Sobolev source norm, the terms containing at least one $V$ are controlled by Lemma~\ref{lem:positive-holder-multiplier}, the H\"older regularity $V_i\in C_T\C^{\rho_V}$, and the Sobolev components of the other factor.  In Case F, the pure deterministic quadratic terms are estimated as follows.  The $Y_iY_j$ terms are controlled by H\"older from
\[
        Y_i\in L_T^4W^{a_Y,4},
        \qquad a_Y:=s_2-\kappa_Y,
\]
so that $Y_iY_j\in L_T^2$ first and then $L_T^2\hookrightarrow L_T^1$ on $[0,T]$ gives a factor $T^{1/2}$.  The spatial product criterion in Lemma~\ref{lem:det-product-toolbox} sends $W^{a_Y,4}\cdot W^{a_Y,4}$ into $H^{s_2-\alpha}$ under the lower condition $s_2>3/2-\alpha+$ and the reserved Strichartz loss.

For the mixed terms $X_iY_j$ and $Y_iX_j$, write $a_X=s_1-\kappa_X$ and $a_Y=s_2-\kappa_Y$.  Time H\"older gives $X_iY_j\in L_T^r$ with $1/r=1/p_X+1/4<1$, followed by $L_T^r\hookrightarrow L_T^1$ and the local factor $T^{1-1/p_X-1/4}$.  The spatial estimate is the mixed Sobolev product
\[
        W^{a_X,q_X}\cdot W^{a_Y,4}\longrightarrow H^{s_2-\alpha},
\]
which is exactly the second product margin in \eqref{eq:admissible-quadratic}.  Terms containing one $V_i$ use instead the direct H\"older bound $V_i\in C_T\C^{\rho_V}$ from Appendix~\ref{app:first-picard-phase}, where $\rho_V>0$.  Lemma~\ref{lem:positive-holder-multiplier} applies in the Sobolev route because $s_2<\alpha$ makes the target $H^{s_2-\alpha}$ non-positive.  When a direct Besov source route is preferable, the separate dyadic bound $V_i\in L_T^\infty B_{2,\infty}^{\rho_V}$ is combined with $\sigma<\rho_V$ and the Besov part of Lemma~\ref{lem:positive-holder-multiplier}.

The $X_1X_2$ term is the most restrictive deterministic quadratic term, so we spell out the exponent calculation.  Put
\[
        a_X=s_1-\kappa_X,
        \qquad
        \kappa_X=\frac32-\frac3{q_X}-\frac{\alpha}{p_X}+\eps_{\rm Str}.
\]
The Strichartz component gives $X_i\in L_T^{p_X}W^{a_X,q_X}$.  H\"older in time first gives $X_1X_2\in L_T^{p_X/2}$, and the finite-interval embedding $L_T^{p_X/2}\hookrightarrow L_T^1$ contributes the local factor $T^{1-2/p_X}$.  The spatial multiplication criterion in Lemma~\ref{lem:det-product-toolbox} requires
\[
        2a_X-(s_2-\alpha)>3\left(\frac2{q_X}-\frac12\right)_+,
\]
which, after substituting $a_X=s_1-\kappa_X$, is precisely the first inequality in \eqref{eq:admissible-quadratic}, up to the reserved $\eps$ and $\eps_{\rm Str}$ margins.  The remaining inequality in \eqref{eq:admissible-quadratic} gives the open margin for the mixed $X$--$Y$ products.  The local factor $T^\delta$ comes from these non-endpoint time exponents and from bounded-in-time factors integrated on $[0,T]$.

In Case W, set $\gamma=s_2-1/2>0$.  The deterministic quadratic part of the $Y$ equation is estimated through the dual source line \eqref{eq:conic-dual-Y}.  Since
\[
        X\in L_T^8W^{s_1-1/4,8/3},
        \qquad
        Y\in L_T^4W^{s_2-1/2,4},
\]
the condition $\gamma<s_1-1/4$, equivalently $s_2<s_1+1/4$, gives
\begin{equation}\label{eq:endpoint-XX-source}
 \|\langle D\rangle^\gamma(X_1X_2)\|_{L_{t,x}^{4/3}}
 \lesssim T^\delta
 \|X_1\|_{L_T^8W^{s_1-1/4,8/3}}
 \|X_2\|_{L_T^8W^{s_1-1/4,8/3}}.
\end{equation}
For $XY$ we place the $\gamma$ derivatives on $Y$ and use
$L_x^{8/3}\cdot L_x^4\subset L_x^{8/5}\subset L_x^{4/3}$ and
$L_T^8\cdot L_T^4\subset L_T^{8/3}\subset L_T^{4/3}$ on the finite interval.  Thus
\begin{equation}\label{eq:endpoint-XY-source}
 \|\langle D\rangle^\gamma(XY)\|_{L_{t,x}^{4/3}}
 \lesssim T^\delta
 \|X\|_{L_T^8W^{s_1-1/4,8/3}}
 \|Y\|_{L_T^4W^{s_2-1/2,4}}.
\end{equation}
For $Y_1Y_2$, one puts the derivative on one $Y$ factor and uses
$L_x^4\cdot L_x^4\subset L_x^2\subset L_x^{4/3}$ and the analogous time embedding, giving
\begin{equation}\label{eq:endpoint-YY-source}
 \|\langle D\rangle^\gamma(Y_1Y_2)\|_{L_{t,x}^{4/3}}
 \lesssim T^\delta
 \|Y_1\|_{L_T^4W^{s_2-1/2,4}}
 \|Y_2\|_{L_T^4W^{s_2-1/2,4}}.
\end{equation}
The terms containing $V_i$ are easier: at $\alpha=1$ one has $V_i\in C_T\mathcal C^{1/2-}$, and $\gamma=s_2-1/2<1/2-$ in the endpoint window, so Lemma~\ref{lem:positive-holder-multiplier} gives the needed multiplier bounds.  For the endpoint Besov source component of the purely deterministic quadratic sectors, the choice $\sigma<1/4$ gives the elementary embedding
\[
        L_x^{4/3}\hookrightarrow B_{2,\infty}^{\sigma-1},
        \qquad
        \|P_N f\|_{L^2}\lesssim N^{3/4}\|P_N f\|_{L^{4/3}},
\]
and therefore $N^{\sigma-1}\|P_Nf\|_{L^2}\lesssim N^{\sigma-1/4}\|P_Nf\|_{L^{4/3}}$.  The finite interval embeddings used above then place the endpoint deterministic quadratic sources in $L_T^1B_{2,\infty}^{\sigma-1}$.

Difference estimates are obtained by replacing one factor in each product by the corresponding difference and using the same bounds; the resulting constants are polynomial in the two augmented control norms.
\end{proof}

\begin{definition}[The X-source norm]\label{def:X-source-norm}
For the fixed-point source in the $X$ equation we use
\begin{equation}\label{eq:X-source-topology}
        \|F\|_{\mathfrak S_{T,X}^{s_1,\sigma}}
        :=\|F\|_{L_T^1H^{s_1-\alpha}}
          +\|F\|_{L_T^1B_{2,\infty}^{\sigma-\alpha}}.
\end{equation}
The $Y$-source norm is the branch-dependent norm $\mathfrak P_T^{s_2,\sigma}$ defined in \eqref{eq:source-topologies}.
\end{definition}

\begin{lemma}[Linear Duhamel estimates]\label{lem:linear-duhamel-estimates}
Let \(F_X=(F_{X,1},F_{X,2})\) and \(F_Y=(F_{Y,1},F_{Y,2})\) be source vectors with
\[
        F_X\in (\mathfrak S_{T,X}^{s_1,\sigma})^2,
        \qquad
        F_Y\in (\mathfrak P_T^{s_2,\sigma})^2.
\]
Then the outer Duhamel maps used in the fixed point obey
\begin{equation}\label{eq:linear-duhamel-estimates}
\sum_i\|I_iF_{X,i}\|_{\widetilde X_T^{s_1,\sigma}}
+\sum_i\|I_iF_{Y,i}\|_{\widetilde Y_T^{s_2,\sigma}}
\lesssim
\sum_i\|F_{X,i}\|_{\mathfrak S_{T,X}^{s_1,\sigma}}
+\sum_i\|F_{Y,i}\|_{\mathfrak P_T^{s_2,\sigma}}.
\end{equation}
In Case F this uses the flat-torus fractional Klein--Gordon theorem proved in Appendix~\ref{app:fracKG-strichartz}.  In Case W the non-\(V\) deterministic quadratic part of \(F_Y\) is the only sector routed through the dual conic norm \(\mathfrak D_T^{s_2}\); the remaining sectors use the rough-source norm \(\mathfrak S_{T,W}^{s_2,\sigma}\).
\end{lemma}

\begin{proof}
For the Sobolev--Strichartz components, apply Proposition~\ref{prop:strichartz} in Case F and Proposition~\ref{prop:gko-conic-strichartz} in Case W.  The definition of \(\mathfrak P_T^{s_2,\sigma}\) was chosen so that the endpoint dual conic line \eqref{eq:conic-dual-Y} is used only for the deterministic non-\(V\) quadratic sector; all stochastic rough, cubic, \(V\)-multiplier, and mixed-operator sectors remain in the rough \(L_T^1H^{s_2-1}\) route at \(\alpha=1\).  The Besov components of \(\widetilde X_T^{s_1,\sigma}\) and \(\widetilde Y_T^{s_2,\sigma}\) follow from Lemma~\ref{lem:elementary-besov-propagation} applied to the \(L_T^1B_{2,\infty}^{\sigma-\alpha}\) source norm in Case F and the \(L_T^1B_{2,\infty}^{\sigma-1}\) source norm in Case W.  Thus the linear propagation step is exactly the passage from the source norms to the augmented path spaces; this stage is purely the linear propagation of the already estimated source norms.
\end{proof}

\begin{lemma}[Sector Lipschitz estimates]\label{lem:sector-lipschitz-estimates}
Let $\Xi,\widetilde\Xi$ be two augmented enhanced-data sets, let $y,\widetilde y$ be deterministic Cauchy data, and let $Z,\widetilde Z$ lie in the radius-$R$ ball of $\calZ_T$.  Put
\[
        W_i^\Xi(Z)=V_i^\Xi+X_i+Y_i,
        \qquad
        W_i^{\widetilde\Xi}(\widetilde Z)=V_i^{\widetilde\Xi}+\widetilde X_i+\widetilde Y_i.
\]
Then the full mixed-operator input satisfies
\begin{equation}\label{eq:sector-lip-W}
        \sum_i\|W_i^\Xi(Z)-W_i^{\widetilde\Xi}(\widetilde Z)\|_{E_T^{2,\sigma}}
        \lesssim
        \|Z-\widetilde Z\|_{\calZ_T}
        +d_T((\Xi,y),(\widetilde\Xi,\widetilde y)).
\end{equation}
For every mixed operator component $\mathcal T$ appearing in \eqref{eq:mild-map-Y}, with source norm denoted by $\mathcal S_T$, one has
\begin{align}
\|\mathcal T^\Xi(W_i^\Xi(Z))-\mathcal T^\Xi(W_i^\Xi(\widetilde Z))\|_{\mathcal S_T}
&\le \|\mathcal T^\Xi\|_{E_T^{2,\sigma}\to\mathcal S_T}
       \|Z-\widetilde Z\|_{\calZ_T},
\label{eq:sector-lip-same-op}\\
\|\mathcal T^\Xi(W_i^\Xi(Z))-\mathcal T^{\widetilde\Xi}(W_i^{\widetilde\Xi}(\widetilde Z))\|_{\mathcal S_T}
&\le \|\mathcal T^\Xi\|_{E_T^{2,\sigma}\to\mathcal S_T}
       \|W_i^\Xi(Z)-W_i^{\widetilde\Xi}(\widetilde Z)\|_{E_T^{2,\sigma}}\notag\\
&\quad +\|\mathcal T^\Xi-\mathcal T^{\widetilde\Xi}\|_{E_T^{2,\sigma}\to\mathcal S_T}
       \|W_i^{\widetilde\Xi}(\widetilde Z)\|_{E_T^{2,\sigma}}.
\label{eq:sector-lip-two-op}
\end{align}
The operator norms in \eqref{eq:sector-lip-same-op} are bounded by the corresponding components of $M_T^{\rm bd}$, and the difference norm in \eqref{eq:sector-lip-two-op} is one of the mixed-operator components of the enhanced-data distance $d_T$.  The analogous identities for $\Psi,V,\Gamma^{(3)}$ and $C$ are obtained by replacing the operator-difference term by the corresponding component of $d_T$.
\end{lemma}

\begin{proof}
Estimate \eqref{eq:sector-lip-W} is Lemma~\ref{lem:full-input-in-E} applied to the difference, using that the $V$-difference is included in the enhanced-data distance.  Equations \eqref{eq:sector-lip-same-op} and \eqref{eq:sector-lip-two-op} are the triangle identity
\[
        \mathcal T^\Xi W-\mathcal T^{\widetilde\Xi}\widetilde W
        =\mathcal T^\Xi(W-\widetilde W)
         +(\mathcal T^\Xi-\mathcal T^{\widetilde\Xi})\widetilde W
\]
followed by the pathwise operator bounds.  The boundedness of $\widetilde W$ is again Lemma~\ref{lem:full-input-in-E} and the radius-$R$ assumption.  The same bilinear expansion applies to ordinary deterministic products, with one factor placed in the difference and the other factors bounded by $R$ and $M_T^{\rm bd}$.
\end{proof}

\begin{proposition}[Quantitative source and stability estimates]\label{prop:quantitative-source-stability}
Let $\Xi$ be an augmented enhanced-data set on $[0,T]$, and fix deterministic Cauchy data $y$.  Let $M_T^{\rm bd}(\Xi,y)$ and $\mu_T(\Xi,y)$ be the bounded size and local smallness quantity of Definition~\ref{def:localized-mixed-size}.  For $Z=(X_1,X_2,Y_1,Y_2)$ define $F_X^\Xi(Z)$ and $F_Y^\Xi(Z)$ to be the source vectors before the outer Duhamel operators in \eqref{eq:mild-map-X} and \eqref{eq:mild-map-Y}, with the sector decomposition fixed by \eqref{eq:source-topologies}.  On every ball $\|Z\|_{\calZ_T},\|\widetilde Z\|_{\calZ_T}\le R$, the following estimates hold with constants polynomial in the displayed bounded controls.

First, the bounded source estimates are
\begin{align}
\sum_i\|F_{X,i}^\Xi(Z)\|_{\mathfrak S_{T,X}^{s_1,\sigma}}
&\le C T^\delta P(M_T^{\rm bd}(\Xi,y))(1+\|Z\|_{\calZ_T}),
\label{eq:quant-source-bounded-X}\\
\sum_i\|F_{Y,i}^\Xi(Z)\|_{\mathfrak P_T^{s_2,\sigma}}
&\le C\Bigl[M_T^{\rm bd}(\Xi,y)
      +\mu_T(\Xi,y)P(M_T^{\rm bd}(\Xi,y))(1+\|Z\|_{\calZ_T})
      +T^\delta P(M_T^{\rm bd}(\Xi,y))(1+\|Z\|_{\calZ_T}^2)\Bigr].
\label{eq:quant-source-bounded-Y}
\end{align}
Second, for the same enhanced data,
\begin{align}
\sum_i\|F_{X,i}^\Xi(Z)-F_{X,i}^\Xi(\widetilde Z)\|_{\mathfrak S_{T,X}^{s_1,\sigma}}
&\le C\mu_T(\Xi,y)P_R(M_T^{\rm bd}(\Xi,y))\|Z-\widetilde Z\|_{\calZ_T},
\label{eq:quant-source-contract-X}\\
\sum_i\|F_{Y,i}^\Xi(Z)-F_{Y,i}^\Xi(\widetilde Z)\|_{\mathfrak P_T^{s_2,\sigma}}
&\le C\mu_T(\Xi,y)P_R(M_T^{\rm bd}(\Xi,y))\|Z-\widetilde Z\|_{\calZ_T}.
\label{eq:quant-source-contract-Y}
\end{align}
Third, for two enhanced-data/Cauchy-data pairs $(\Xi,y)$ and $(\widetilde\Xi,\widetilde y)$, put
\[
        A_T:=M_T^{\rm bd}(\Xi,y)+M_T^{\rm bd}(\widetilde\Xi,\widetilde y),
        \qquad
        \mu_T^\#:=\mu_T(\Xi,y)+\mu_T(\widetilde\Xi,\widetilde y).
\]
Then
\begin{align}
\sum_i\|F_{X,i}^\Xi(Z)-F_{X,i}^{\widetilde\Xi}(\widetilde Z)\|_{\mathfrak S_{T,X}^{s_1,\sigma}}
&\le C\mu_T^\#P_R(A_T)\|Z-\widetilde Z\|_{\calZ_T}
     +CP_R(A_T)d_T((\Xi,y),(\widetilde\Xi,\widetilde y)),
\label{eq:quant-source-data-X}\\
\sum_i\|F_{Y,i}^\Xi(Z)-F_{Y,i}^{\widetilde\Xi}(\widetilde Z)\|_{\mathfrak P_T^{s_2,\sigma}}
&\le C\mu_T^\#P_R(A_T)\|Z-\widetilde Z\|_{\calZ_T}
     +CP_R(A_T)d_T((\Xi,y),(\widetilde\Xi,\widetilde y)).
\label{eq:quant-source-data-Y}
\end{align}
The $d_T$ term is the bounded enhanced-data perturbation used for continuous dependence after the contraction factor has been absorbed on the left.
\end{proposition}

\begin{proof}
The $X$ source is
\[
        F_{X,i}^\Xi(Z)=W_{3-i}(Z)<\Psi_i+W_i(Z)<\Psi_{3-i},
        \qquad W_i(Z)=V_i+X_i+Y_i,
\]
and Lemma~\ref{lem:lowhigh-X} gives the bounded and difference estimates.  The time-local factor $T^\delta$ is part of $\mu_T$ after the bounded enhanced-data factors are absorbed into $P_R$.

For the $Y$ source, decompose
\[
        F_Y=F_{\rm rough}+F_{\rm op}+F_\Gamma+F_V+F_{\rm det}.
\]
The rough part consists of the low--high, high--low, and resonant products with the stochastic convolution.  In Case F it is controlled by Lemma~\ref{lem:caseF-rough-stochastic-routes}, and in Case W by Lemma~\ref{lem:caseW-rough-stochastic-routes}.  The mixed-operator part is controlled by Lemma~\ref{lem:mixed-operator-Y}; its locally small coefficient is $T M_T^{\rm op,H}+T M_{T,{\rm diag}}^{\rm op,B}+M_{T,{\rm cen}}^{\rm op,B}$.  The cubic part $F_\Gamma=\Gamma_1+\Gamma_2$ is measured directly by the cubic component of $M_T^{\rm bd}$ and has zero same-enhanced-data difference.  For different enhanced data, it is controlled by the $\Gamma$ component of $d_T$.  The terms containing $V$ use Lemma~\ref{lem:positive-holder-multiplier} with the direct H\"older norm $V_i\in C_T\mathcal C^{\rho_V}$; perturbing $V$ is one of the explicit components of $d_T$.  The non-$V$ deterministic quadratic terms are treated by Lemma~\ref{lem:quadratic-Y}: in Case F they first enter $L_T^1H^{s_2-\alpha}$ and then use Lemma~\ref{lem:sobolev-to-besov-source}, while in Case W they use the dual conic route and the endpoint Besov embedding recorded in Lemma~\ref{lem:endpoint-quadratic-estimates}.  Each deterministic product difference is expanded bilinearly, with one factor placed in the difference and the remaining factors kept inside the radius-$R$ ball.

The algebraic identity behind the preceding paragraph is Lemma~\ref{lem:sector-lipschitz-estimates}: same-enhanced-data differences are always estimated by putting the difference in the solution input, while two-enhanced-data differences split into an input difference and an enhanced-object difference.  The operator-difference terms in \eqref{eq:quant-source-data-X}--\eqref{eq:quant-source-data-Y} are exactly the mixed-operator components of $d_T$ in Definition~\ref{def:enhanced-data-distance}.  Since these perturbation terms measure changes of the enhanced data themselves, they appear with the bounded polynomial $P_R(A_T)$ and without an additional local-smallness factor.
\end{proof}

\subsection{Contraction theorem}
For $Z=(X_1,X_2,Y_1,Y_2)$ define $W_i(Z)=V_i+X_i+Y_i$ and
\begin{equation}\label{eq:mild-map-X}
\mathcal M_i^X(Z)=I_i\bigl(W_2(Z)<\Psi_1+W_1(Z)<\Psi_2\bigr),
\end{equation}
\begin{align}\label{eq:mild-map-Y}
\mathcal M_i^Y(Z)=&\ S_i(\cdot)(y_{i,0},y_{i,1})\\
&+I_i\Bigl(W_2>\Psi_1+W_1>\Psi_2+W_1W_2+
\Gamma_1+\Gamma_2+Y_2\circ\Psi_1+Y_1\circ\Psi_2\nonumber\\
&\hspace{2.8cm}
+T^{2;1,1}(W_2)+T^{2;2,1}(W_1)+T^{1;1,2}(W_2)+T^{1;2,2}(W_1)\Bigr).
\end{align}

\begin{proposition}[Quantitative stability of the mild map]\label{prop:mild-map-stability}
Let $\mathcal M^{\Xi,y}$ denote the mild map \eqref{eq:mild-map-X}--\eqref{eq:mild-map-Y} built from an enhanced-data set $\Xi$ and deterministic Cauchy data $y$ on $[0,T]$.  In the estimates below the fixed $y$ is suppressed from the notation except where two data sets are compared.  On every ball
\[
        \|Z\|_{\calZ_T},\ \|\widetilde Z\|_{\calZ_T}\le R,
\]
there are polynomials $P$ and $P_R$ in the bounded enhanced-data controls such that
\begin{equation}\label{eq:mild-map-bounded-prop}
        \|\mathcal M^{\Xi,y}(Z)\|_{\calZ_T}
        \le C M_T^{\rm bd}(\Xi,y)
        +C\mu_T(\Xi,y)P_R(M_T^{\rm bd}(\Xi,y))
          (1+\|Z\|_{\calZ_T}+\|Z\|_{\calZ_T}^2),
\end{equation}
for the same enhanced data,
\begin{equation}\label{eq:mild-map-contract-prop}
        \|\mathcal M^{\Xi,y}(Z)-\mathcal M^{\Xi,y}(\widetilde Z)\|_{\calZ_T}
        \le C\mu_T(\Xi,y)P_R(M_T^{\rm bd}(\Xi,y))
        \|Z-\widetilde Z\|_{\calZ_T},
\end{equation}
and, for two enhanced-data/Cauchy-data pairs $(\Xi,y)$ and $(\widetilde\Xi,\widetilde y)$,
\begin{equation}\label{eq:mild-map-data-prop}
\begin{aligned}
        \|\mathcal M^{\Xi,y}(Z)-\mathcal M^{\widetilde\Xi,\widetilde y}(\widetilde Z)\|_{\calZ_T}
        &\le C(\mu_T(\Xi,y)+\mu_T(\widetilde\Xi,\widetilde y))
             P_R(A_T)\|Z-\widetilde Z\|_{\calZ_T}\\
        &\quad +C P_R(A_T)d_T((\Xi,y),(\widetilde\Xi,\widetilde y)),
\end{aligned}
\end{equation}
where $A_T=M_T^{\rm bd}(\Xi,y)+M_T^{\rm bd}(\widetilde\Xi,\widetilde y)$.  The $d_T$ term is the bounded perturbation of the enhanced data.  After the contraction factor in the first term is made smaller than one, this bounded perturbation gives the Lipschitz response.
\end{proposition}

\begin{proof}
Apply Proposition~\ref{prop:quantitative-source-stability} to the source vectors $F_X$ and $F_Y$.  The outer Duhamel maps are then estimated by Lemma~\ref{lem:linear-duhamel-estimates}, which combines Proposition~\ref{prop:strichartz} in Case F, Proposition~\ref{prop:gko-conic-strichartz} in Case W, and the dyadic Besov propagation Lemma~\ref{lem:elementary-besov-propagation}.  The source norm $\mathfrak P_T^{s_2,\sigma}$ is chosen so that the Case F Sobolev route and the Case W dual conic route have the same estimate.  The homogeneous term is bounded by the deterministic Cauchy-data norm in \eqref{eq:data-norm-Y}, included in $M_T^{\rm bd}(\Xi,y)$; differences of homogeneous terms are measured by the deterministic-data component of $d_T$.  This gives \eqref{eq:mild-map-bounded-prop}--\eqref{eq:mild-map-data-prop}.
\end{proof}

\begin{lemma}[Fixed-point radius and local smallness]\label{lem:fixed-point-radius-table}
Let $A\ge1$ bound $M_T^{\rm bd}$ for one enhanced-data/Cauchy-data pair, and let $R\ge 2C(A+A^2)$, where $C$ is the constant in Proposition~\ref{prop:mild-map-stability}.  There is a polynomial $P_R$ such that the following implication holds.  If
\begin{equation}\label{eq:table-smallness-single}
        C\mu_T P_R(A)\le \min\left\{\frac12,\frac{R-C A}{C(1+R+R^2)}\right\},
\end{equation}
then $\mathcal M$ maps the closed radius-$R$ ball of $\calZ_T$ into itself and is a contraction there with Lipschitz constant at most $1/2$.  For two data sets bounded by $A$, if
\begin{equation}\label{eq:table-smallness-pair}
        C(\mu_T+\widetilde\mu_T)P_R(A)\le\frac12,
\end{equation}
then the difference of their fixed points satisfies
\begin{equation}\label{eq:table-lipschitz}
        \|Z-\widetilde Z\|_{\calZ_T}
        \le 2C P_R(A)d_T((\Xi,y),(\widetilde\Xi,\widetilde y)).
\end{equation}
Thus bounded controls determine the ball and the Lipschitz polynomial, while the localized smallness modulus is used only to make the contraction factors in \eqref{eq:table-smallness-single} and \eqref{eq:table-smallness-pair} smaller than one.
\end{lemma}

\begin{proof}
For $\|Z\|_{\calZ_T}\le R$, estimate \eqref{eq:mild-map-bounded-prop} gives
\[
        \|\mathcal M(Z)\|_{\calZ_T}
        \le C A+C\mu_TP_R(A)(1+R+R^2).
\]
The second bound in \eqref{eq:table-smallness-single} makes the right-hand side at most $R$.  Estimate \eqref{eq:mild-map-contract-prop} and the first bound in \eqref{eq:table-smallness-single} give the contraction constant.  For two fixed points, apply \eqref{eq:mild-map-data-prop} with $Z$ and $\widetilde Z$ equal to the fixed points and move the first term on the right to the left using \eqref{eq:table-smallness-pair}.  This gives \eqref{eq:table-lipschitz}.
\end{proof}

\begin{proposition}[Uniform deterministic bridge on localized-small classes]\label{prop:uniform-deterministic-bridge}
Fix one admissible branch: Case F with $12/13<\alpha<1$, or Case W with $\alpha=1$.  Let $\mathcal K$ be a localized-small class on $[0,T_0]$ in the sense of Definition~\ref{def:localized-smallness-class}, with bounded radius $A\ge1$ and modulus $\nu_{\mathcal K}$.  There is a constant $C_0$, depending only on the fixed exponents, the speeds, and the branch, and for every $R\ge1$ there is a nondecreasing polynomial $P_R$, such that for every $0<T\le T_0$, every $(\Xi,y),(\widetilde\Xi,\widetilde y)\in\mathcal K$, and every $Z,\widetilde Z\in\calZ_T$ with norms at most $R$, one has
\begin{align}
        \|\mathcal M^{\Xi,y}(Z)\|_{\calZ_T}
        &\le C_0 A+C_0\nu_{\mathcal K}(T)P_R(A)(1+\|Z\|_{\calZ_T}+\|Z\|_{\calZ_T}^2),\label{eq:uniform-bridge-map}\\
        \|\mathcal M^{\Xi,y}(Z)-\mathcal M^{\Xi,y}(\widetilde Z)\|_{\calZ_T}
        &\le C_0\nu_{\mathcal K}(T)P_R(A)\|Z-\widetilde Z\|_{\calZ_T},\label{eq:uniform-bridge-contract}\\
        \|\mathcal M^{\Xi,y}(Z)-\mathcal M^{\widetilde\Xi,\widetilde y}(\widetilde Z)\|_{\calZ_T}
        &\le C_0\nu_{\mathcal K}(T)P_R(A)\|Z-\widetilde Z\|_{\calZ_T}
           +C_0 P_R(A)d_T((\Xi,y),(\widetilde\Xi,\widetilde y)).\label{eq:uniform-bridge-two-data}
\end{align}
If $R\ge2C_0(A+A^2)$, then there is $T_*\in(0,T_0]$, depending only on $A,R,P_R$ and $\nu_{\mathcal K}$, such that every map $\mathcal M^{\Xi,y}$ with $(\Xi,y)\in\mathcal K$ is a strict contraction on the closed radius-$R$ ball of $\calZ_T$ for all $0<T\le T_*$.  Its fixed points satisfy
\begin{equation}\label{eq:uniform-bridge-fixedpoint-lipschitz}
        \|Z^{\Xi,y}-Z^{\widetilde\Xi,\widetilde y}\|_{\calZ_T}
        \le 2C_0 P_R(A)d_T((\Xi,y),(\widetilde\Xi,\widetilde y)).
\end{equation}
\end{proposition}

\begin{proof}
This is Proposition~\ref{prop:mild-map-stability} applied uniformly on $\mathcal K$.  On $[0,T]$, each pair in $\mathcal K$ satisfies $M_T^{\rm bd}\le A$ and $\mu_T\le\nu_{\mathcal K}(T)$, so the bounded coefficients are absorbed into $P_R(A)$ and the contraction coefficient is controlled by $\nu_{\mathcal K}(T)P_R(A)$.  The enhanced-data perturbation in \eqref{eq:uniform-bridge-two-data} is a bounded perturbation of the coefficients of the deterministic map, while the solution-input difference is the contraction variable.

Choose $T_*$ by Lemma~\ref{lem:time-restriction-local-smallness} so that
\[
        C_0\nu_{\mathcal K}(T)P_R(A)
        \le \min\left\{\frac12,\frac{R-C_0 A}{C_0(1+R+R^2)}\right\},
        \qquad 0<T\le T_*.
\]
Then \eqref{eq:uniform-bridge-map} gives the self-map property and \eqref{eq:uniform-bridge-contract} gives contraction with constant at most $1/2$.  Applying \eqref{eq:uniform-bridge-two-data} to the two fixed points and absorbing the first term on the right gives \eqref{eq:uniform-bridge-fixedpoint-lipschitz}.
\end{proof}

\begin{theorem}[Deterministic solution map in the augmented topology]\label{thm:deterministic-closure}
Fix $0<T_0\le1$.  Work either in the subunit fractional Case F, $12/13<\alpha<1$, using Proposition~\ref{prop:strichartz}, which is proved in Appendix~\ref{app:fracKG-strichartz} from the internal flat-torus fractional Klein--Gordon theorem, or in the classical conic endpoint Case W, $\alpha=1$, using the standard conic estimates recorded in Proposition~\ref{prop:gko-conic-strichartz}.  In Case W the $X$ and $Y$ spaces are the endpoint spaces displayed in Proposition~\ref{prop:gko-conic-strichartz}, and the deterministic quadratic part of the $Y$ equation uses the dual source route \eqref{eq:conic-dual-Y}.  Choose admissible parameters from Definition~\ref{def:admissible}; Besov propagation is provided in both cases by Lemma~\ref{lem:elementary-besov-propagation}.  Assume that an augmented enhanced-data set $\Xi$ restricted to $[0,T_0]$ and deterministic data $y$ have finite bounded size $M_{T_0}^{\rm bd}(\Xi,y)$ and that the corresponding singleton class is localized-small in the sense of Definition~\ref{def:localized-smallness-class}.  Then there exists
\[
        T_*=T_*\bigl(M_{T_0}^{\rm bd}(\Xi,y),\nu_{\{(\Xi,y)\}}\bigr)\in(0,T_0]
\]
such that, for every $0<T\le T_*$, the restricted map $\mathcal M$ has a unique fixed point in the closed ball selected in the proof inside the branch-dependent augmented space
\[
 (\widetilde X_T^{s_1,\sigma})^2\times(\widetilde Y_T^{s_2,\sigma})^2.
\]
The solution map is locally Lipschitz with respect to the augmented enhanced-data distance $d_T$ on every common localized-smallness class of pairs from Definition~\ref{def:localized-smallness-class}, including the deterministic Cauchy data distance in \eqref{eq:data-norm-Y}.  The common comparison time is governed by the local smallness modulus of the class in addition to the bounded radius $M_{T_0}^{\rm bd}$.
\end{theorem}

\begin{proof}
Apply Proposition~\ref{prop:uniform-deterministic-bridge} to the singleton localized-small class $\{(\Xi,y)\}$.  Let $C_0$ be the constant in that proposition and choose
\[
        A=1+M_{T_0}^{\rm bd}(\Xi,y),
        \qquad
        R\ge 2C_0(A+A^2).
\]
The proposition gives a time $T_*$, depending only on $A$ and the localized-smallness modulus of the singleton class, such that $\mathcal M^{\Xi,y}$ is a strict contraction on the closed radius-$R$ ball of $\calZ_T$ for every $0<T\le T_*$.  Banach's fixed-point theorem gives a unique fixed point in this ball of
\[
        (\widetilde X_T^{s_1,\sigma})^2\times(\widetilde Y_T^{s_2,\sigma})^2,
\]
and the reconstruction $u_i=\Psi_i+V_i+X_i+Y_i$ solves the enhanced paracontrolled system.

For local Lipschitz dependence, place the two enhanced-data/Cauchy-data pairs in a common localized-small class $\mathcal K$.  Proposition~\ref{prop:uniform-deterministic-bridge} gives a common comparison time and the fixed-point estimate \eqref{eq:uniform-bridge-fixedpoint-lipschitz}.  The deterministic Cauchy-data difference is one component of $d_T$ and is propagated by \eqref{eq:homogeneous-Y-estimate} and \eqref{eq:homogeneous-besov-estimate}.  Hence the solution map is locally Lipschitz on every common localized-smallness interval.  The role of boundedness is only to choose the ball and the Lipschitz polynomial; the contraction time is chosen from the localized-smallness modulus.
\end{proof}

\begin{lemma}[Cutoff enhanced data are localized-small]\label{lem:cutoff-localized-small}
Let $\Xi$ be the enhanced-data set constructed by the stochastic-symbol, cubic-symbol, and mixed-operator estimates of this paper, and let $\Xi_\Lambda$ be its Galerkin-compatible cutoff data.  Fix deterministic Cauchy data $y$ with finite norm \eqref{eq:data-norm-Y}.  Suppose the componentwise cutoff convergence holds in the augmented distance $d_{T_0}$ on some interval $[0,T_0]$, including the mixed-operator convergence in the Sobolev and Besov source operator topologies supplied by Section~\ref{sec:random-operators}.  Then, after discarding finitely many cutoff levels, the family
\[
        \{(\Xi,y)\}\cup\{(\Xi_\Lambda,y):\Lambda\ge\Lambda_0\}
\]
is a common localized-small class on a sufficiently short subinterval of $[0,T_0]$.
\end{lemma}

\begin{proof}
The convergence in $d_{T_0}$ gives uniform boundedness of all fixed-time enhanced symbols, Cauchy data, and global mixed-operator norms for the limit together with all sufficiently large cutoffs.  The ordinary deterministic source factors in \eqref{eq:localized-full-size} are therefore bounded by a common polynomial times the explicit power $T^\delta$.  The same-color diagonal part is controlled by the uniform $E_T^{2,\sigma}\to L_T^\infty B_{2,\infty}^{\sigma-\alpha}$ Volterra norm and hence becomes small after multiplying by the deterministic factor $T$ in \eqref{eq:det-contraction-L1}.  For the centered Besov source operator, first choose $T$ so that the limiting localized $L_T^1B_{2,\infty}^{\sigma-\alpha}$ norm is small.  Then use cutoff convergence in the same operator topology to make the difference between the cutoff and limiting centered operators uniformly small for all sufficiently large cutoffs.  The finitely many discarded or retained low cutoff levels are harmless after reducing the time once more.  Thus the supremum of $\mu_T$ over the tail family tends to zero as $T\downarrow0$, which is precisely the localized-smallness condition of Definition~\ref{def:localized-smallness-class}.
\end{proof}

\begin{corollary}[Deterministic cutoff passage]\label{cor:deterministic-cutoff-passage}
Let $(\Xi_\Lambda,y_\Lambda)$ be a cutoff enhanced-data/Cauchy-data family and let $(\Xi,y)$ be its limit on $[0,T_0]$.  In the Galerkin cutoff situation with fixed deterministic data, Lemma~\ref{lem:cutoff-localized-small} supplies the common localized-small class after discarding finitely many cutoff levels.  In general, assume that the limit and the remaining cutoff data lie in a common localized-small class $\mathcal K$, and that
\[
        d_T((\Xi_\Lambda,y_\Lambda),(\Xi,y))\longrightarrow0
        \qquad\text{for every fixed }T\le T_0.
\]
Then, after restricting to the common time supplied by Proposition~\ref{prop:uniform-deterministic-bridge}, the corresponding fixed points converge in $\calZ_T$.  For the Galerkin-compatible cutoff convention, Proposition~\ref{prop:finite-cutoff-bridge} identifies these finite enhanced fixed points with the Galerkin solutions of \eqref{eq:cutoff-system}; the smooth stochastic-leg convention is compared to this convention by Proposition~\ref{prop:cutoff-convention-equivalence}.
\end{corollary}

\begin{proof}
Apply Proposition~\ref{prop:uniform-deterministic-bridge} to $\mathcal K$ and insert $(\Xi_\Lambda,y_\Lambda)$ and $(\Xi,y)$ in \eqref{eq:uniform-bridge-fixedpoint-lipschitz}.  The right-hand side tends to zero by the augmented-distance convergence.  The final sentence is the finite algebraic bridge of Proposition~\ref{prop:finite-cutoff-bridge}, followed by the cutoff-convention comparison in the same augmented topology.
\end{proof}

\section{Conclusion}\label{sec:conclusion}
For independent colors we constructed the enhanced stochastic data for the two-speed cross-interaction system and proved the mixed-operator estimate needed for the paracontrolled fixed point.  The interaction graph forces every same-color contraction to occur with a different outer phase.  Its diagonal part gains the factor $N^{-\alpha}$.  The centered normal form uses Gaussian indices $([\ell_0],\mu)$ and $([r_0],\lambda)$ and places the signs in the deterministic incidence $n=q+\varepsilon_1\ell_0+\varepsilon_2r_0$.  Thus the finite Gram reduction neither duplicates the variables at $m$ and $-m$ nor introduces a cutoff-dependent auxiliary dimension.  These estimates yield a locally Lipschitz solution map and convergence of the Fourier--Galerkin approximations.  The weak-covariance result is a separate perturbative extension, obtained by retaining the additional lower-chaos branches and using covariance coordinates local to each Fourier involution class.

The restriction $\alpha>12/13$ is attained simultaneously by the cubic window $3/2-\alpha<s_2<11\alpha/2-9/2$, the first-Picard Besov-input window $3-3\alpha<\sigma<7\alpha/2-3$, and the near-endpoint deterministic quadratic constraints.  This is the boundary of the present closure scheme; optimality for the equation remains open.  The mixed-operator estimate itself closes in the larger window $s_2<4\alpha-3$, $3-3\alpha<\sigma<4\alpha-3$.  Enlarging the range therefore requires a compatible improvement of more than the cubic estimate alone.  The coalescing-speed limit $c_1\to c_2$, where the Volterra denominator degenerates, is a separate problem.

\appendix

\section{Second-order random tensor operator norm}\label{app:second-order-tensor-proof}

This appendix proves Proposition~\ref{prop:abstract-random-tensor}.  After Gaussian decoupling, two applications of the Lust-Piquard--Pisier noncommutative Khintchine inequality produce the oriented flattenings in \eqref{eq:oriented-flattening-profile}; see \cite{DeLaPenaGine,LustPisier,DNY,Kaneshiro}.

\begin{lemma}[Second-order Gaussian decoupling]\label{lem:second-order-decoupling}
Let $(g_a)_{a\in I}$ be a finite standard Gaussian family and let $(g'_a)_{a\in I}$ be an independent copy.  For any Banach-space-valued second homogeneous Wick chaos
\[
        X=\sum_{a,b\in I}x_{a,b}:g_ag_b:
\]
and every $p\ge2$,
\begin{equation}\label{eq:second-order-decoupling}
        \|X\|_{L^p(\Omega;\mathcal X)}
        \lesssim
        \left\|\sum_{a,b\in I}x_{a,b}g_ag'_b\right\|_{L^p(\Omega;\mathcal X)}.
\end{equation}
The implicit constant depends only on the chaos degree.
\end{lemma}

\begin{proof}
This is the standard decoupling inequality for homogeneous Gaussian chaos.  In degree two it follows from polarization and the Kahane contraction principle.  The Wick subtraction removes the zeroth-chaos pairing, so the same-color centered term lies in the pure second homogeneous chaos.  In the cross-color case the two color families are already independent.
\end{proof}

\begin{lemma}[Operator-valued decoupling for centered Fourier chaos]\label{lem:operator-valued-fourier-decoupling}
Fix $\varepsilon_1,\varepsilon_2\in\{-1,1\}$.  Let $X=\mathcal L(\ell_q^2,\ell_n^2)$ and let
\[
        Z=\sum_{\ell_0,r_0\in\mathbb Z^3_{\mathrm{rep}}}
        \sum_{q,n}\sum_{\mu,\lambda\le d_0}
        H_{\ell_0,\mu,r_0,\lambda,q,n}\,
        :G_{j,[\ell_0],\mu}G_{j,[r_0],\lambda}:\,E_{nq}
\]
be one fixed sign pattern of a finite same-color centered Fourier block after the scalar involution pairing has been Wick-subtracted and after the class-wise Gram normalization of Lemma~\ref{lem:labelwise-auxiliary-coordinates}, with $d_0=8$.  Here $E_{nq}$ denotes the matrix unit from the $q$ input coefficient to the $n$ output coefficient, and the deterministic tensor is supported on
\[
        n=q+\varepsilon_1\ell_0+\varepsilon_2r_0.
\]
Let $(g_{j,[\ell_0],\mu})$ and $(g'_{j,[r_0],\lambda})$ be independent copies of the class-indexed Gaussian family.  Then, for every $p\ge2$,
\begin{equation}\label{eq:operator-valued-fourier-decoupling}
\left\|Z\right\|_{L^p(\Omega;\mathcal L(\ell_q^2,\ell_n^2))}
\lesssim_p
\left\|
\sum_{\ell_0,r_0\in\mathbb Z^3_{\mathrm{rep}}}
\sum_{q,n}\sum_{\mu,\lambda\le d_0}
H_{\ell_0,\mu,r_0,\lambda,q,n}
g_{j,[\ell_0],\mu}g'_{j,[r_0],\lambda}E_{nq}
\right\|_{L^p(\Omega;\mathcal L(\ell_q^2,\ell_n^2))}.
\end{equation}
The original factors with $[\ell_0]=[r_0]$ may belong to the same Gaussian family; the variables on the right are independent only because of the decoupling inequality.  The class indices and bounded auxiliary indices preserve the signed incidence, the dyadic support, and the row/column flattening powers.
\end{lemma}

\begin{proof}
Inside fixed involution classes, use \eqref{eq:signed-fourier-class-expansion} and write
\[
        h_{\varepsilon_1\ell_0,s}
        =\sum_{\mu\le d_0}\alpha_{\ell_0,\mu}(s)e_{[\ell_0],\mu},
        \qquad
        h_{\varepsilon_2r_0,t}
        =\sum_{\lambda\le d_0}\beta_{r_0,\lambda}(t)e_{[r_0],\lambda},
\]
with orthonormal vectors chosen only in the local class-wise Hilbert spaces.  Hence
\[
        I_j^{:2:}(h_{\varepsilon_1\ell_0,s}
                  \widetilde\otimes h_{\varepsilon_2r_0,t})
        =\sum_{\mu,\lambda\le d_0}
          \alpha_{\ell_0,\mu}(s)\beta_{r_0,\lambda}(t)
          I_j^{:2:}(e_{[\ell_0],\mu}\widetilde\otimes e_{[r_0],\lambda}).
\]
After the scalar involution pairing has been subtracted, the second homogeneous term is therefore a finite sum of Wick products of class-indexed normalized Gaussian coordinates.  Lemma~\ref{lem:second-order-decoupling} is applied in the Banach space $X=\mathcal L(\ell_q^2,\ell_n^2)$ with coefficients
\[
        x_{[\ell_0],\mu,[r_0],\lambda}
        =\sum_{q,n}H_{\ell_0,\mu,r_0,\lambda,q,n}E_{nq}.
\]
The decoupled variables are the independent copies displayed in \eqref{eq:operator-valued-fourier-decoupling}.  If a local covariance block has a null space, the basis is chosen in the finite quotient Hilbert space after zero-norm directions are identified, separately in each class.  The deterministic tensor still carries $(\ell_0,r_0,q,n)$ and the incidence $n=q+\varepsilon_1\ell_0+\varepsilon_2r_0$.  For fixed signs, relabeling the representatives by their signed frequencies identifies the support with a restriction of the incidence in Lemma~\ref{lem:dyadic-tensor-flattening}, so the same upper bounds apply.  The auxiliary ranges are uniformly bounded, and the quotient and orthonormalization change only constants.
\end{proof}

\begin{lemma}[Gaussian noncommutative Khintchine]\label{lem:gaussian-nck}
Let $(\gamma_r)$ be independent standard real or complex Gaussians and let $M_r:\ell_C^2\to\ell_D^2$ be deterministic finite matrices.  For $p\ge2$,
\begin{align}\label{eq:gaussian-nck}
&\left\|\,\left\|\sum_r\gamma_rM_r\right\|_{\ell_C^2\to\ell_D^2}\,\right\|_{L^p}
\nonumber\\
&\qquad\lesssim
\sqrt p\,\log^{1/2}(2+|C|+|D|)
\max\left\{
\left\|\left(\sum_rM_rM_r^*\right)^{1/2}\right\|,
\left\|\left(\sum_rM_r^*M_r\right)^{1/2}\right\|
\right\}.
\end{align}
\end{lemma}

\begin{proof}
Apply the Lust-Piquard--Pisier noncommutative Khintchine inequality in Schatten class $S^r$, with $r\simeq p+\log(2+|C|+|D|)$.  For a finite matrix, the comparison $\|A\|_{\ell^2\to\ell^2}\le\|A\|_{S^r}$ gives the desired operator norm on the left, while the reverse comparison needed for the variance profiles is $\|A\|_{S^r}\le (2+|C|+|D|)^{1/r}\|A\|_{\ell^2\to\ell^2}$.  With the above choice of $r$, the dimensional factor is bounded by an absolute constant, leaving only the displayed logarithmic loss.  This loss is harmless in dyadic summations.  See \cite{LustPisier,DNY} for standard forms of this estimate.  The form used in Section~\ref{sec:random-operators} is the finite-dimensional operator estimate stated in Proposition~\ref{prop:abstract-random-tensor}.
\end{proof}

\begin{proof}[Proof of Proposition~\ref{prop:abstract-random-tensor}]
We estimate the decoupled matrix
\begin{equation}\label{eq:decoupled-tensor-matrix}
        Z=\sum_{a,b}g_ah_bA_{a,b},
        \qquad (A_{a,b})_{d,c}=H_{a,b,c,d},
\end{equation}
where $(g_a)$ and $(h_b)$ are the independent Gaussian families from the proposition.  If one starts instead from a same-color Wick product, Lemma~\ref{lem:second-order-decoupling} is the preliminary reduction; this is the route used in Corollary~\ref{cor:centered-fourier-block-tensor}.  Conditional on $h=(h_b)$, put
\[
        W_a(h)=\sum_bh_bA_{a,b},
        \qquad Z=\sum_ag_aW_a(h).
\]
Lemma~\ref{lem:gaussian-nck}, applied conditionally in the $g$ variables, gives
\begin{equation}\label{eq:first-nck-tensor-proof}
\|Z\|_{L_g^p(\ell^2_c\to\ell^2_d)}
\lesssim
\sqrt p\,\log^{1/2}(2+|C|+|D|)
\bigl(\mathcal V_D(h)+\mathcal V_C(h)\bigr),
\end{equation}
where
\[
\mathcal V_D(h)=\left\|\left(\sum_aW_a(h)W_a(h)^*\right)^{1/2}\right\|,
\qquad
\mathcal V_C(h)=\left\|\left(\sum_aW_a(h)^*W_a(h)\right)^{1/2}\right\|.
\]
We next estimate these two random variance profiles by a second Khintchine step.

First, $\mathcal V_D(h)$ is the operator norm of the random map
\[
\mathfrak K(h):\ell_d^2\to\ell_a^2\otimes\ell_c^2,
\qquad
(\mathfrak K(h)v)_{a,c}=\sum_{b,d}\overline{H_{a,b,c,d}}\,\overline{h_b}v_d.
\]
The family \((\overline{h_b})_b\) has the same standard complex Gaussian law as \((h_b)_b\).  Writing \(\mathfrak K(h)=\sum_b\overline{h_b}\mathfrak K_b\) and applying Lemma~\ref{lem:gaussian-nck} yields
\begin{align*}
\|\mathcal V_D(h)\|_{L_h^p}
&\lesssim
\sqrt p\,\log^C(2+|A|+|B|+|C|+|D|)\\
&\quad\times
\max\left\{
\left\|\left(\sum_b\mathfrak K_b\mathfrak K_b^*\right)^{1/2}\right\|,
\left\|\left(\sum_b\mathfrak K_b^*\mathfrak K_b\right)^{1/2}\right\|
\right\}.
\end{align*}
The two deterministic terms are precisely the flattening norms
\[
        \|H\|_{\ell^2_{b,d}\to\ell^2_{a,c}}
        \quad\text{and}\quad
        \|H\|_{\ell^2_d\to\ell^2_{a,b,c}}
        =\|H\|_{\ell^2_{a,b,c}\to\ell^2_d},
\]
up to adjoint identification.  Hence
\begin{equation}\label{eq:row-var-tensor-proof}
        \|\mathcal V_D(h)\|_{L_h^p}
        \lesssim
        \sqrt p\,\log^C(2+|A|+|B|+|C|+|D|)\,\mathfrak F_4(H).
\end{equation}

Second, $\mathcal V_C(h)$ is the operator norm of
\[
\mathfrak L(h):\ell_c^2\to\ell_a^2\otimes\ell_d^2,
\qquad
(\mathfrak L(h)u)_{a,d}=\sum_{b,c}H_{a,b,c,d}\,h_bu_c.
\]
The same noncommutative Khintchine estimate for $\mathfrak L(h)=\sum_bh_b\mathfrak L_b$ gives the deterministic flattenings
\[
        \|H\|_{\ell^2_{b,c}\to\ell^2_{a,d}}
        \quad\text{and}\quad
        \|H\|_{\ell^2_c\to\ell^2_{a,b,d}}
        =\|H\|_{\ell^2_{a,b,d}\to\ell^2_c},
\]
and therefore
\begin{equation}\label{eq:col-var-tensor-proof}
        \|\mathcal V_C(h)\|_{L_h^p}
        \lesssim
        \sqrt p\,\log^C(2+|A|+|B|+|C|+|D|)\,\mathfrak F_4(H).
\end{equation}
Combining \eqref{eq:first-nck-tensor-proof}, \eqref{eq:row-var-tensor-proof}, and \eqref{eq:col-var-tensor-proof}, then absorbing harmless powers of $p$ and logarithms into $p^C\log^C$, proves \eqref{eq:abstract-random-tensor}.  Complex Gaussian variables and harmless real-coordinate normalizations only change the constants.  The same-color Fourier covariance is outside the abstract theorem; it is reduced to this decoupled theorem in Corollary~\ref{cor:centered-fourier-block-tensor} after Wick subtraction and finite Hilbert-kernel normalization.
\end{proof}

\section{Fractional Klein--Gordon Sobolev--Strichartz on \texorpdfstring{$\mathbb T^3$}{T3}}\label{app:fracKG-strichartz}

This appendix records a self-contained proof of the Sobolev--Strichartz estimate used in Proposition~\ref{prop:strichartz} for the massive fractional Klein--Gordon phase
\[
        \omega_c(n)=\bigl(1+c^2|n|^{2\alpha}\bigr)^{1/2},
        \qquad 0<\alpha<1.
\]
This is the massive flat-torus adaptation of the compact-manifold fractional-wave estimate; compare \cite{Dinh}.  The proof combines Poisson summation, stationary phase for the rescaled annular phase, and a \(TT^*\)--Young argument, using the oscillatory-integral conventions of \cite{SteinHarmonic,SoggeFIO}.  It covers fixed \(0<\alpha<1\), scalar speed \(c>0\), and non-endpoint spatial exponents \(q<\infty\).  The retarded \(L_T^1H^{s-\alpha}\) estimate follows from the homogeneous estimate by Minkowski.  At \(\alpha=1\) the radial curvature vanishes, and Proposition~\ref{prop:gko-conic-strichartz} instead uses the classical compact Klein--Gordon theorem \cite{CacciafestaDanesiMeng}.

Let $U_c(t)=e^{\ii t\omega_c(D)}$.  The smooth Littlewood--Paley projector $P_N$ has symbol $\chi(n/N)$, where $\chi\in C_c^\infty(\mathbb R^3\setminus\{0\})$ is supported in a fixed annulus.  The low block $P_1$ is finite-dimensional and is always absorbed into constants.  For $2\le q<\infty$ set
\[
        \delta(q)=3\left(\frac12-\frac1q\right).
\]

\begin{theorem}[Dyadic massive fractional Klein--Gordon Strichartz]\label{thm:fracKG-dyadic}
Let $0<\alpha<1$, $0<T\le1$, $2\le p<\infty$, $2\le q<\infty$, and
\begin{equation}\label{eq:fracKG-range}
        \delta(q)>\frac2p.
\end{equation}
Then, for every dyadic $N\ge1$,
\begin{equation}\label{eq:fracKG-dyadic-main}
        \|U_c(t)P_Nf\|_{L_t^p([0,T];L_x^q(\mathbb T^3))}
        \le C_{\alpha,c,p,q,T}
        N^{\delta(q)-\alpha/p}\|P_Nf\|_{L_x^2}.
\end{equation}
At the logarithmic line $\delta(q)=2/p$, the same estimate holds with an additional harmless factor $N^\eps$ for every $\eps>0$.
\end{theorem}

\begin{remark}[Nonempty non-endpoint range used in Case F]\label{rem:fracKG-nonempty-range}
The deterministic fixed point chooses $(p_X,q_X)$ close to, but strictly away from, $(4,3)$.  The condition in Theorem~\ref{thm:fracKG-dyadic} is
\[
        \delta(q_X)=3\left(\frac12-\frac1{q_X}\right)>\frac2{p_X}.
\]
At $(4,3)$ both sides equal $1/2$.  Taking $q_X>3$ and $p_X>4$ sufficiently close to these endpoint values gives $\delta(q_X)>1/2$ and $2/p_X<1/2$, hence a strict open set of admissible pairs.  The fixed loss $\varepsilon_{\rm Str}$ in the deterministic space is chosen after this strict margin is fixed.
\end{remark}

\begin{corollary}[Sobolev and retarded Duhamel lines]\label{cor:fracKG-semigroup-duhamel}
Let $s\in\mathbb R$ and set
\[
        \kappa(p,q)=\delta(q)-\frac{\alpha}{p}+\varepsilon_{\rm Str},
\]
where $\varepsilon_{\rm Str}>0$ is either the logarithmic-line loss or a fixed small loss.  Under \eqref{eq:fracKG-range},
\begin{equation}\label{eq:fracKG-sob-line}
        \|U_c(t)h\|_{L_t^pW_x^{s-\kappa(p,q),q}([0,T]\times\mathbb T^3)}
        \lesssim \|h\|_{H_x^s}.
\end{equation}
Moreover, for
\[
        I_cF(t)=\int_0^t\frac{\sin((t-s)\omega_c(D))}{\omega_c(D)}F(s)\,ds,
\]
\begin{equation}\label{eq:fracKG-duhamel-line}
        \|I_cF\|_{L_t^pW_x^{s-\kappa(p,q),q}}
        \lesssim \|F\|_{L_t^1H_x^{s-\alpha}}.
\end{equation}
The source time norm in \eqref{eq:fracKG-duhamel-line} is $L_t^1$, so the retarded estimate follows directly from Minkowski.
\end{corollary}

\subsection{Annular phase geometry}

For dyadic $N\ge2$ define
\begin{equation}\label{eq:fracKG-rescaled-phase}
        \phi_N(\xi)=N^{-\alpha}\omega_c(N\xi)
        =\bigl(N^{-2\alpha}+c^2|\xi|^{2\alpha}\bigr)^{1/2}.
\end{equation}
All estimates in this subsection are on a fixed annulus $a_0\le |\xi|\le a_1$ with $0<a_0<a_1<\infty$.

\begin{lemma}[Uniform annular symbol bounds]\label{lem:fracKG-symbol-bounds}
For every integer $K\ge0$,
\[
        \sup_{N\ge2}\sup_{a_0\le |\xi|\le a_1}
        |\nabla_\xi^K\phi_N(\xi)|<\infty.
\]
Moreover $\phi_N\to c|\xi|^\alpha$ in $C^K$ on the annulus for every fixed $K$.
\end{lemma}

\begin{proof}
On the annulus, $r=|\xi|$ is bounded away from zero.  The function $(\mu,r)\mapsto(\mu+c^2r^{2\alpha})^{1/2}$ is smooth for $0\le\mu\le1$ and $a_0\le r\le a_1$.  Taking $\mu=N^{-2\alpha}$ gives the uniform symbol bounds and the convergence.
\end{proof}

\begin{lemma}[Annular Hessian for the massive fractional phase]\label{lem:fracKG-hessian}
Fix $0<\alpha<1$, $c>0$, and $0<a_0<a_1<\infty$.  There is $N_{\rm nd}=N_{\rm nd}(\alpha,c,a_0,a_1)$ such that $D^2_\xi\phi_N(\xi)$ is non-degenerate on $a_0\le |\xi|\le a_1$ for all dyadic $N\ge N_{\rm nd}$.  If $r=|\xi|$ and $f_N(r)=(N^{-2\alpha}+c^2r^{2\alpha})^{1/2}$, then the two tangential eigenvalues and the radial eigenvalue are
\[
        \lambda_{\mathrm{tan}}(r)=\frac{f_N'(r)}{r}
        =\frac{c^2\alpha r^{2\alpha-2}}{f_N(r)},
\]
and
\[
        \lambda_{\mathrm{rad}}(r)=f_N''(r)
        =\frac{c^2\alpha r^{2\alpha-2}}{f_N(r)^3}
        \Bigl((2\alpha-1)N^{-2\alpha}+c^2(\alpha-1)r^{2\alpha}\Bigr).
\]
On the annulus, $|\lambda_{\mathrm{tan}}|$ is bounded below by a positive constant depending on $a_0,a_1,c,\alpha$, while
\[
        |\lambda_{\mathrm{rad}}(r)|
        \ge c_{\alpha,c,a_0,a_1}(1-\alpha)
\]
for all $N\ge N_{\rm nd}$.  Thus the stationary-phase constants may deteriorate as $\alpha\uparrow 1$, consistently with the conic limit.
\end{lemma}

\begin{proof}
The formulas follow by differentiating the radial function $f_N(r)$.  The tangential eigenvalue is uniformly positive because $r\in[a_0,a_1]$ and $f_N(r)\simeq1$.  For the radial factor set
\[
        B_N(r):=(2\alpha-1)N^{-2\alpha}+c^2(\alpha-1)r^{2\alpha}.
\]
The second term is negative and has size at least $c^2(1-\alpha)a_0^{2\alpha}$.  Choose $N_{\rm nd}$ so that
\[
        |2\alpha-1|N^{-2\alpha}
        \le \frac12 c^2(1-\alpha)a_0^{2\alpha}
        \qquad (N\ge N_{\rm nd}).
\]
Then $|B_N(r)|\ge\frac12c^2(1-\alpha)a_0^{2\alpha}$ on the annulus.  The remaining prefactor in $\lambda_{\rm rad}$ is bounded above and below.  The finitely many shells $N<N_{\rm nd}$ are absorbed into finite-shell constants in the kernel estimate.
\end{proof}

\begin{remark}[Dependence on the conic limit]\label{rem:hessian-alpha-dependence}
The radial lower bound in Lemma~\ref{lem:fracKG-hessian} is linear in $1-\alpha$ on fixed annuli.  The stationary-phase constants used below involve only finitely many derivatives of the inverse Hessian and of the phase; hence, for fixed differentiation order, they may grow polynomially in $(1-\alpha)^{-1}$ as $\alpha\uparrow1$.  This is the precise sense in which the subunit annular-curvature proof degenerates at the conic endpoint.  The endpoint branch is therefore treated separately by Proposition~\ref{prop:gko-conic-strichartz}.
\end{remark}

\subsection{Oscillatory integral and Poisson kernel estimates}

\begin{lemma}[Stationary phase with bounded linear parameter]\label{lem:fracKG-stationary-phase}
Let
\[
        \Psi_{N,v,\sigma}(\xi)=\sigma\phi_N(\xi)+v\cdot\xi,
        \qquad \sigma\in\{\pm1\}.
\]
For every fixed $C_0<\infty$, there is $C<\infty$ such that, for all $|v|\le C_0$, all $\lambda\ge1$, and all dyadic $N\ge N_{\rm nd}$,
\[
        \left|\int_{\mathbb R^3}e^{\ii\lambda\Psi_{N,v,\sigma}(\xi)}\chi(\xi)\,d\xi\right|
        \le C\lambda^{-3/2}.
\]
\end{lemma}

\begin{proof}

We give the uniform parameter version of the stationary phase step,
since the critical points of
\[
        \Psi_{N,v,\sigma}(\xi)=\sigma\phi_N(\xi)+v\cdot \xi
\]
may depend on \(N\), \(v\), and \(\sigma\).

Let \(K_0=\operatorname{supp}\chi\). Choose a slightly larger compact annular set
\(K_1\Subset\{\xi: a_0/2<|\xi|<2a_1\}\) with \(K_0\Subset K_1\).
By Lemmas B.5 and B.6, for every integer \(L\) there is a constant
\(C_L<\infty\), and there is a number \(\gamma>0\), such that for all
dyadic \(N\ge N_{\mathrm{nd}}\), all \(|v|\le C_0\), and all
\(\sigma\in\{\pm1\}\),
\[
   \|\Psi_{N,v,\sigma}\|_{C^L(K_1)}\le C_L,
   \qquad
   \inf_{\xi\in K_1} s_{\min}\big(D_\xi^2\Psi_{N,v,\sigma}(\xi)\big)
   \ge \gamma .
\]
Here \(s_{\min}\) denotes the smallest singular value. The constant
\(\gamma\) may depend on \(a_0,a_1,\alpha,c\), and it may deteriorate
as \(\alpha\uparrow1\), but it is uniform in \(N,v,\sigma\). The linear
term \(v\cdot\xi\) affects only the first derivative and therefore does
not change the Hessian bounds.

We now make the stationary phase reduction uniform. Choose \(r>0\),
depending only on the preceding constants and on \(K_1\), so small that
on every ball \(B(\xi_0,4r)\subset K_1\) the Hessian
\(D_\xi^2\Psi_{N,v,\sigma}\) varies by at most \(\gamma/10\), uniformly
in \(N,v,\sigma\). Cover \(K_0\) by finitely many balls
\(B_\mu=B(\xi_\mu,r)\) such that \(B(\xi_\mu,4r)\subset K_1\), and let
\((\chi_\mu)_\mu\) be a smooth partition of unity on \(K_0\) subordinate
to \(B_\mu\). It is enough to estimate
\[
   I_\mu(\lambda,N,v,\sigma)
   =
   \int_{\mathbb R^3}
       e^{i\lambda\Psi_{N,v,\sigma}(\xi)}
       \chi_\mu(\xi)\chi(\xi)\,d\xi
\]
uniformly in \(\mu,N,v,\sigma\).

Fix one \(\mu\) and one parameter value \((N,v,\sigma)\). On
\(B(\xi_\mu,4r)\), the map
\[
        \xi\mapsto \nabla_\xi\Psi_{N,v,\sigma}(\xi)
\]
is a quantitative local diffeomorphism. More precisely, after possibly
decreasing \(r\), the inverse function theorem applied with the above
uniform lower bound on the Hessian gives the following two alternatives.

First, if \(\nabla_\xi\Psi_{N,v,\sigma}\) has no zero in
\(B(\xi_\mu,2r)\), then
\[
       |\nabla_\xi\Psi_{N,v,\sigma}(\xi)|\ge c_* r,
       \qquad \xi\in B_\mu ,
\]
with \(c_*>0\) independent of \(N,v,\sigma,\mu\). Indeed, if this lower
bound failed at some \(\xi\in B_\mu\), the quantitative inverse theorem
or Newton's theorem applied on \(B(\xi_\mu,2r)\) would produce a zero of
\(\nabla_\xi\Psi_{N,v,\sigma}\) inside \(B(\xi_\mu,2r)\). In this
non-stationary case, repeated integration by parts with
\[
       L_{\Psi}
       =
       \frac{\nabla\Psi_{N,v,\sigma}\cdot\nabla}
            {i\lambda|\nabla\Psi_{N,v,\sigma}|^2},
       \qquad
       L_{\Psi}(e^{i\lambda\Psi_{N,v,\sigma}})
       =
       e^{i\lambda\Psi_{N,v,\sigma}},
\]
and the uniform \(C^L\)-bounds gives, for every \(A\ge0\),
\[
       |I_\mu(\lambda,N,v,\sigma)|\le C_A \lambda^{-A}.
\]
In particular this is \(O(\lambda^{-3/2})\).

Second, suppose that \(\nabla_\xi\Psi_{N,v,\sigma}\) has a zero
\(\xi_{\mu,N,v,\sigma}\) in \(B(\xi_\mu,2r)\). By the same quantitative
inverse theorem, this zero is unique in \(B(\xi_\mu,4r)\). Let
\[
      \zeta_{\mu,N,v,\sigma}(\xi)
      =
      \zeta\!\left(\frac{\xi-\xi_{\mu,N,v,\sigma}}{r}\right),
\]
where \(\zeta\in C_c^\infty(B(0,2))\) and \(\zeta\equiv1\) on \(B(0,1)\).
We decompose
\[
       I_\mu=I_{\mu}^{\mathrm{near}}+I_{\mu}^{\mathrm{far}}
\]
by inserting \(\zeta_{\mu,N,v,\sigma}\) and
\(1-\zeta_{\mu,N,v,\sigma}\).

On the support of the far part, the distance to the unique critical point
is at least \(r\). Since the Hessian is uniformly invertible and varies
by at most \(\gamma/10\) on \(B(\xi_\mu,4r)\), we have
\[
       |\nabla_\xi\Psi_{N,v,\sigma}(\xi)|\ge c_* r
\]
on the far support. The same integration by parts argument gives
\[
       |I_{\mu}^{\mathrm{far}}|\le C_A\lambda^{-A}.
\]

For the near part, the uniform Morse lemma applies. The uniform
\(C^L\)-bounds on the phase and the uniform lower bound on the Hessian
give a \(C^{L'}\)-change of variables
\[
       y=\kappa_{\mu,N,v,\sigma}(\xi),
\]
defined on \(B(\xi_{\mu,N,v,\sigma},2r)\), whose \(C^{L'}\)-norms and
Jacobian bounds are uniform in \(N,v,\sigma,\mu\), such that
\[
       \Psi_{N,v,\sigma}(\xi)
       =
       \Psi_{N,v,\sigma}(\xi_{\mu,N,v,\sigma})
       +
       \frac12\sum_{j=1}^3 \varepsilon_j y_j^2,
       \qquad
       \varepsilon_j\in\{\pm1\}.
\]
After this change of variables,
\[
       I_{\mu}^{\mathrm{near}}
       =
       e^{i\lambda\Psi_{N,v,\sigma}(\xi_{\mu,N,v,\sigma})}
       \int_{\mathbb R^3}
          e^{\frac{i\lambda}{2}\sum_j\varepsilon_j y_j^2}
          a_{\mu,N,v,\sigma}(y)\,dy,
\]
where the amplitudes \(a_{\mu,N,v,\sigma}\) are supported in a fixed ball
and are uniformly bounded in a sufficiently high \(C^L\)-norm. The
standard quadratic stationary phase estimate, or equivalently rescaling
\(y=\lambda^{-1/2}z\) and integrating by parts outside \(|z|\lesssim
\lambda^{1/2}\), gives
\[
       |I_{\mu}^{\mathrm{near}}|\le C\lambda^{-3/2}.
\]

Combining the near and far estimates gives
\[
       |I_\mu(\lambda,N,v,\sigma)|\le C\lambda^{-3/2}
\]
uniformly in \(N,v,\sigma,\mu\). Since only finitely many patches
\(\mu\) are used, summing over \(\mu\) yields
\[
   \left|
   \int_{\mathbb R^3}
        e^{\ii\lambda(\sigma\phi_N(\xi)+v\cdot\xi)}
        \chi(\xi)\,d\xi
   \right|
   \le C\lambda^{-3/2}.
\]

\end{proof}

\begin{lemma}[Nonstationary integration by parts]\label{lem:fracKG-nonstationary}
There is $C_1<\infty$, depending only on $\alpha,c$ and the support of $\chi$, such that if $|v|\ge C_1$, then for every $A\ge0$,
\[
        \left|\int e^{\ii\lambda(\sigma\phi_N(\xi)+v\cdot\xi)}\chi(\xi)\,d\xi\right|
        \le C_A(\lambda|v|)^{-A},
        \qquad \lambda\ge1.
\]
\end{lemma}

\begin{proof}
Choose $C_1>2\sup_{N\ge2}\sup_{\xi\in\supp\chi}|\nabla\phi_N(\xi)|$, which is finite by Lemma~\ref{lem:fracKG-symbol-bounds}.  Then
\[
        |\nabla\Psi(\xi)|
        =|\sigma\nabla\phi_N(\xi)+v|
        \ge |v|/2,
        \qquad \Psi(\xi)=\sigma\phi_N(\xi)+v\cdot\xi,
\]
on $\supp\chi$.  We use the standard non-stationary phase operator, see for instance \cite[Chapter VIII]{SteinHarmonic},
\[
        L=\frac{\nabla\Psi\cdot\nabla}{\ii\lambda|\nabla\Psi|^2},
        \qquad L(e^{\ii\lambda\Psi})=e^{\ii\lambda\Psi}.
\]
If $L^*$ denotes its formal adjoint, then
\[
        L^*a
        =-\sum_{j=1}^3\partial_{\xi_j}\left(
        \frac{\partial_{\xi_j}\Psi}{\ii\lambda|\nabla\Psi|^2}a\right).
\]
Every derivative of the coefficient $(\partial_j\Psi)|\nabla\Psi|^{-2}$ is bounded by $C_m|v|^{-1}$ on $\supp\chi$: the numerator derivatives of order at least one are controlled by Lemma~\ref{lem:fracKG-symbol-bounds}, while each denominator is bounded below by $|v|^2$.  Hence, by induction,
\[
        |(L^*)^A\chi(\xi)|\le C_A(\lambda |v|)^{-A}
\]
on the fixed support of $\chi$.  Integrating by parts $A$ times gives
\[
        \left|\int e^{\ii\lambda\Psi(\xi)}\chi(\xi)\,d\xi\right|
        =\left|\int e^{\ii\lambda\Psi(\xi)}(L^*)^A\chi(\xi)\,d\xi\right|
        \le C_A(\lambda|v|)^{-A}.
\]

\end{proof}

\begin{lemma}[Poisson image count and annular tail]\label{lem:fracKG-poisson-image-count}
Let $0<\alpha<1$, $0<T\le1$, $\lambda=N^\alpha|t|\ge1$, $|t|\le T$, and let $x$ range in a fixed fundamental domain of $\mathbb T^3$.  For any fixed $C_1<\infty$,
\[
  \#\left\{k\in\mathbb Z^3:
       \left|\frac{N(x+2\pi k)}{\lambda}\right|\le C_1\right\}
  \le C_{T,C_1}.
\]
Moreover, if $A>9/2$, then uniformly in $N,t,x$,
\begin{equation}\label{eq:poisson-tail-annuli}
  \sum_{|x+2\pi k|\gtrsim \lambda/N}
       (N|x+2\pi k|)^{-A}
  \lesssim_A \lambda^{-3/2}.
\end{equation}
\end{lemma}

\begin{proof}
The first assertion follows from
\[
  |x+2\pi k|\le C_1\frac{\lambda}{N}
   =C_1|t|N^{\alpha-1}\le C_1T,
\]
because $\alpha<1$ and $N\ge1$; a bounded Euclidean ball intersects only $O_T(1)$ lattice translates of a fixed fundamental domain.

For the second assertion, put $r_0\simeq\lambda/N$.  Since $|t|\le1$ and $\alpha<1$, $r_0\le1$.  We first remove all lattice translates in the fixed ball
\[
        \mathcal K_0(x):=\{k:|x+2\pi k|\le4\pi\}.
\]
The separation of the lattice gives $\#\mathcal K_0(x)\le C$ uniformly in $x$.  Every member of the nonstationary sum also satisfies $|x+2\pi k|\gtrsim r_0$, so the entire contribution of $\mathcal K_0(x)$ is $O((Nr_0)^{-A})=O(\lambda^{-A})$.  This handles all possible nearest translates, including the finitely many ties on the boundary of a fundamental domain.

For the remaining translates, $|x+2\pi k|>4\pi$.  Group them into annuli $2^j4\pi<|x+2\pi k|\le2^{j+1}4\pi$.  A disjoint-ball packing argument gives $O(2^{3j})$ lattice points in the $j$th annulus, with no additional $O(1)$ term.  Hence
\[
  \sum_{k\notin\mathcal K_0(x)}(N|x+2\pi k|)^{-A}
  \lesssim_A N^{-A}\sum_{j\ge0}2^{(3-A)j}
  \lesssim_A N^{-A}
  \le \lambda^{-A},
\]
where $A>3$ and $\lambda=N^\alpha|t|\le N$ were used.  Thus the full tail is $O(\lambda^{-A})$, and in particular it is $O(\lambda^{-3/2})$ for the stated $A>9/2$.
\end{proof}

\begin{lemma}[Poisson summation kernel bound]\label{lem:fracKG-poisson-kernel}
Let
\[
        K_N(t,x)=\sum_{n\in\mathbb Z^3}e^{\ii(n\cdot x+t\omega_c(n))}\chi(n/N).
\]
For $0<\alpha<1$, $0<T\le1$, and $|t|\le T$,
\begin{equation}\label{eq:fracKG-kernel}
        \|K_N(t,\cdot)\|_{L_x^\infty(\mathbb T^3)}
        \lesssim_{\alpha,c,T,\chi}
        N^3(1+N^\alpha|t|)^{-3/2}.
\end{equation}
The implicit constant may deteriorate as $\alpha\uparrow 1$.
\end{lemma}

\begin{proof}
The low block and the finitely many shells $N<N_{\rm nd}$ are finite dimensional.  For such a fixed shell, $1\le\lambda=N^\alpha|t|\le TN^\alpha$ ranges over a compact interval whenever the dispersive factor is nontrivial, and hence
\[
        \sup_{|t|\le T,x\in\mathbb T^3}
        \frac{|K_N(t,x)|}{N^3(1+N^\alpha|t|)^{-3/2}}<\infty.
\]
Taking the maximum over the finitely many exceptional shells absorbs them into the constant.  We therefore assume $N\ge N_{\rm nd}$.  If $\lambda=N^\alpha|t|\le1$, then $|K_N|\lesssim N^3$.  Assume $\lambda\ge1$.  Poisson summation gives
\[
        K_N(t,x)=N^3\sum_{k\in\mathbb Z^3}I_k(t,x),
\]
where
\[
        I_k(t,x)=\int_{\mathbb R^3}
        e^{\ii N(x+2\pi k)\cdot\xi+\ii tN^\alpha\phi_N(\xi)}\chi(\xi)\,d\xi.
\]
Let $\sigma=\operatorname{sgn}t$ and
\[
        v_k=\frac{N(x+2\pi k)}{\lambda}.
\]
Then
\[
        I_k(t,x)=\int e^{\ii\lambda(\sigma\phi_N(\xi)+v_k\cdot\xi)}\chi(\xi)\,d\xi.
\]

If $|v_k|\le C_1$, Lemma~\ref{lem:fracKG-stationary-phase} gives $|I_k|\lesssim\lambda^{-3/2}$.  The number of such indices is uniformly finite on fixed time intervals, since
\[
        |x+2\pi k|\le C_1\frac{\lambda}{N}
        =C_1|t|N^{\alpha-1}\le C_1T,
\]
and $x$ lies in a fixed fundamental domain.  Equivalently, Lemma~\ref{lem:fracKG-poisson-image-count} gives a uniform count of stationary images, and they contribute $O_T(\lambda^{-3/2})$ to $\sum_k I_k$.

If $|v_k|>C_1$, Lemma~\ref{lem:fracKG-nonstationary} yields, for every $A$,
\[
        |I_k(t,x)|\lesssim_A(\lambda|v_k|)^{-A}
        =(N|x+2\pi k|)^{-A}.
\]
The same condition implies $|x+2\pi k|\gtrsim\lambda/N$.  Since $0<\alpha<1$ and $|t|\le1$, $\lambda/N=|t|N^{\alpha-1}\le1$.  For $A>10$, Lemma~\ref{lem:fracKG-poisson-image-count} gives
\[
        \sum_{|x+2\pi k|\gtrsim\lambda/N}(N|x+2\pi k|)^{-A}
        \lesssim_A \lambda^{-3/2}.
\]
This controls the nonstationary images.  Multiplying the estimate for $\sum_kI_k$ by $N^3$ proves \eqref{eq:fracKG-kernel}.
\end{proof}

\subsection{The \texorpdfstring{$TT^*$}{TT*}--Young step and Sobolev consequences}

\begin{lemma}[TT*--Young estimate]\label{lem:fracKG-ttstar-young}
For $2\le q<\infty$,
\begin{equation}\label{eq:fracKG-interp}
        \|U_c(t)P_NU_c(s)^*P_Ng\|_{L_x^q}
        \lesssim
        N^{2\delta(q)}(1+N^\alpha|t-s|)^{-\delta(q)}
        \|g\|_{L_x^{q'}}.
\end{equation}
Consequently, if $2\le p<\infty$ and $\delta(q)>2/p$, then
\[
        \|U_c(t)P_Nf\|_{L_t^pL_x^q([0,T]\times\mathbb T^3)}
        \lesssim
        N^{\delta(q)-\alpha/p}\|P_Nf\|_{L_x^2}.
\]
At $\delta(q)=2/p$, the same estimate has an additional logarithmic factor, hence an $N^\eps$ loss for every $\eps>0$.
\end{lemma}

\begin{proof}
The kernel of $U_c(t)P_NU_c(s)^*P_N$ has the same form as $K_N(t-s,x-y)$, with a smooth annular cutoff replacing $\chi$.  Lemma~\ref{lem:fracKG-poisson-kernel} gives the $L_x^1\to L_x^\infty$ bound
\[
        N^3(1+N^\alpha|t-s|)^{-3/2},
\]
while unitarity gives the $L_x^2\to L_x^2$ bound.  Interpolation yields \eqref{eq:fracKG-interp}.

Let $k_N(\tau)=(1+N^\alpha|\tau|)^{-\delta(q)}$.  The $TT^*$ operator is dominated in time by convolution with $N^{2\delta(q)}k_N$.  For the convolution estimate $L_t^{p'}\to L_t^p$, Young's inequality uses
\[
        \frac1r=1+\frac1p-\frac1{p'}=\frac2p,
        \qquad r=\frac p2.
\]
Thus
\[
        \|k_N*h\|_{L_t^p([-T,T])}
        \le \|k_N\|_{L_t^{p/2}([-T,T])}\|h\|_{L_t^{p'}}.
\]
If $\delta(q)>2/p$, then
\[
\begin{aligned}
        \|k_N\|_{L_t^{p/2}([-T,T])}
        &=\left(\int_{-T}^{T}(1+N^\alpha|\tau|)^{-\delta(q)p/2}\,d\tau\right)^{2/p}  \\
        &\lesssim N^{-2\alpha/p}.
\end{aligned}
\]
Therefore
\[
        \|TT^*F\|_{L_t^pL_x^q}
        \lesssim N^{2\delta(q)}N^{-2\alpha/p}\|F\|_{L_t^{p'}L_x^{q'}},
\]
and taking the square root gives
\[
        \|U_c(t)P_Nf\|_{L_t^pL_x^q}
        \lesssim N^{\delta(q)-\alpha/p}\|P_Nf\|_{L^2}.
\]
At the logarithmic line $\delta(q)=2/p$,
\[
        \|k_N\|_{L_t^{p/2}}
        \lesssim N^{-2\alpha/p}\log(2+N)^{2/p},
\]
which gives a logarithmic loss after the square-root step and is absorbed by $N^\varepsilon$.
\end{proof}

\begin{proof}[Proof of Theorem~\ref{thm:fracKG-dyadic}]
The high dyadic estimates follow from Lemma~\ref{lem:fracKG-ttstar-young}.  The low block is finite-dimensional:
\[
        \|U_c(t)P_1f\|_{L_t^pL_x^q}\lesssim_T\|P_1f\|_2.
\]
This proves the theorem.
\end{proof}

\begin{proof}[Proof of Corollary~\ref{cor:fracKG-semigroup-duhamel}]
By the Littlewood--Paley square-function inequality for $2\le q<\infty$, and by Minkowski since $p,q\ge2$,
\[
\begin{aligned}
\|U_c(t)h\|_{L_t^pW_x^{s-\kappa,q}}
&\lesssim
\left(\sum_N N^{2(s-\kappa)}\|U_c(t)P_Nh\|_{L_t^pL_x^q}^2\right)^{1/2} \\
&\lesssim
\left(\sum_N N^{2s}\|P_Nh\|_2^2\right)^{1/2}.
\end{aligned}
\]
The restriction $q<\infty$ is used exactly in this square-function step.  For Duhamel, write
\[
        I_cF(t)=\frac1{2\ii}\int_0^t
        \bigl(U_c(t-s)-U_c(-(t-s))\bigr)\omega_c(D)^{-1}F(s)\,ds.
\]
For fixed source time $s$, the homogeneous estimate on an interval of length at most $T$ gives
\[
        \|U_c(t-s)\omega_c(D)^{-1}F(s)\|_{L_t^p([s,T];W_x^{s-\kappa,q})}
        \lesssim \|\omega_c(D)^{-1}F(s)\|_{H^s}
        \lesssim \|F(s)\|_{H^{s-\alpha}}.
\]
The same estimate holds for the negative half-flow.  Ordinary Minkowski in the source variable then gives \eqref{eq:fracKG-duhamel-line}.  This is the only inhomogeneous line used in the stochastic fixed point.
\end{proof}

For the exponents used in the deterministic fixed point,
\[
        (p_Y,q_Y)=(4,4),\qquad
        \kappa_Y=\delta(4)-\frac\alpha4+\varepsilon_{\rm Str}
        =\frac{3-\alpha}{4}+\varepsilon_{\rm Str},
\]
and for $p_X>4$, $q_X>3$ satisfying \eqref{eq:pxqx-admissible},
\[
        \kappa_X=\frac32-\frac3{q_X}-\frac\alpha{p_X}+\varepsilon_{\rm Str}.
\]
The condition \eqref{eq:pxqx-admissible} is exactly $\delta(q_X)>2/p_X$.  The pair $(4,4)$ is strictly inside the admissible line because $\delta(4)=3/4>1/2=2/4$.  The harmless $\varepsilon_{\rm Str}$ in these derivative losses is chosen smaller than the open margins in Definition~\ref{def:admissible}.

\section{Elementary dyadic \texorpdfstring{$B_{2,\infty}$}{B2,infty} propagation}\label{app:b2-propagation}

The dyadic $L^2$ propagation estimates used in Section~\ref{sec:deterministic-closure} are recorded below.

Let
\[
        S_i(t)(f,g)=\cos(t\omega_i(D))f+
        \frac{\sin(t\omega_i(D))}{\omega_i(D)}g,
        \qquad
        I_iF(t)=\int_0^t\frac{\sin((t-s)\omega_i(D))}{\omega_i(D)}F(s)\,ds,
\]
where $\omega_i(n)=(1+c_i^2|n|^{2\alpha})^{1/2}$.  The dyadic Besov norm is
\[
        \|h\|_{\Btwo^s}:=\sup_{M\ge1}M^s\|P_Mh\|_{L^2(\T^3)},
\]
with the bounded low-frequency block included in $M=1$.

\begin{proposition}[Dyadic L2 propagation of the Besov bound]\label{prop:app-b2-propagation}
For every $\sigma\in\mathbb R$, $0<T\le1$, and $i=1,2$,
\begin{align}
 \|S_i(\cdot)(f,g)\|_{L_T^\infty\Btwo^\sigma}
 +\|\partial_tS_i(\cdot)(f,g)\|_{L_T^\infty\Btwo^{\sigma-\alpha}}
 &\lesssim
 \|f\|_{\Btwo^\sigma}+\|g\|_{\Btwo^{\sigma-\alpha}},
 \label{eq:app-b2-hom}\\
 \|I_iF\|_{L_T^\infty\Btwo^\sigma}
 +\|\partial_tI_iF\|_{L_T^\infty\Btwo^{\sigma-\alpha}}
 &\lesssim
 \|F\|_{L_T^1\Btwo^{\sigma-\alpha}}.
 \label{eq:app-b2-duh}
\end{align}
The constants depend on the fixed scalar symbol $\omega_i$ and on the harmless finite low-frequency convention, uniformly with respect to $|c_1-c_2|^{-1}$.
\end{proposition}

\begin{proof}
On a dyadic shell $|n|\sim M$ one has $\omega_i(n)\simeq M^\alpha$ for $M\ge2$, while $M=1$ is finite-dimensional.  By Plancherel,
\[
 \|P_M\cos(t\omega_i(D))f\|_{2}\le \|P_Mf\|_2,
 \qquad
 \left\|P_M\frac{\sin(t\omega_i(D))}{\omega_i(D)}g\right\|_2
 \lesssim M^{-\alpha}\|P_Mg\|_2.
\]
Multiplication by $M^\sigma$ and the supremum over $M$ and $t$ give the position part of \eqref{eq:app-b2-hom}.  Since
\[
\partial_tS_i(t)(f,g)=-\omega_i(D)\sin(t\omega_i(D))f+\cos(t\omega_i(D))g,
\]
we also have
\[
M^{\sigma-\alpha}\|P_M\partial_tS_i(t)(f,g)\|_2
\lesssim M^\sigma\|P_Mf\|_2+M^{\sigma-\alpha}\|P_Mg\|_2,
\]
which proves the velocity part.

For Duhamel's operator, Minkowski and the same shell multiplier bounds give
\[
M^\sigma\|P_MI_iF(t)\|_2
\le \int_0^t M^{\sigma-\alpha}\|P_MF(s)\|_2\,ds,
\]
and, because
\[
\partial_t I_iF(t)=\int_0^t\cos((t-s)\omega_i(D))F(s)\,ds,
\]
\[
M^{\sigma-\alpha}\|P_M\partial_tI_iF(t)\|_2
\le \int_0^t M^{\sigma-\alpha}\|P_MF(s)\|_2\,ds.
\]
Taking the supremum over $M$ and $t$ proves \eqref{eq:app-b2-duh}.
\end{proof}

\begin{corollary}[Lower strong time continuity]\label{cor:app-b2-lower-continuity}
If the data and source in Proposition~\ref{prop:app-b2-propagation} are finite in the displayed norms, then each dyadic block of the propagated solution is continuous in time and, for every $\eta>0$,
\[
S_i(f,g),\ I_iF\in C_T\Btwo^{\sigma-\eta},
\qquad
\partial_tS_i(f,g),\ \partial_tI_iF\in C_T\Btwo^{\sigma-\alpha-\eta}.
\]
If the exact $B_{2,\infty}$ data belong to the corresponding little Besov closure, then the same argument gives exact strong $C_T\Btwo^\sigma$ continuity in that little-Besov subspace.
\end{corollary}

\begin{proof}
The multiplier formula gives continuity of each fixed block.  Let $M_0$ be large.  For the high-frequency tail at the lower exponent,
\[
\sup_{M>M_0}M^{\sigma-\eta}\|P_Mu(t)-P_Mu(t')\|_2
\le 2M_0^{-\eta}\|u\|_{L_T^\infty\Btwo^\sigma},
\]
and the finitely many blocks $M\le M_0$ are continuous.  The velocity estimate is identical.  If the dyadic tail at the endpoint vanishes uniformly, namely in the little Besov closure, the same finite-tail argument works with $\eta=0$.
\end{proof}

\section{The proportional low--high phase gap}\label{app:phase-gap}

This appendix proves the phase estimate used in the same-color Volterra diagonal directly on the proportional low--high window selected by the paraproduct cutoff.

We fix annular constants
\[
0<a_0<a_1<\infty,
\]
and interpret $|\ell|\sim N$ as
\begin{equation}\label{eq:app-annulus}
 a_0N\le |\ell|\le a_1N.
\end{equation}
The low frequency satisfies $|q|\le\delta N$, where $\delta>0$ will be chosen below from the annular constants and the speed gap.

\begin{lemma}[High-frequency expansion]\label{lem:app-high-frequency}
Let $c>0$, $0<\alpha\le 1$, and
\[
 \omega_c(x)=(1+c^2|x|^{2\alpha})^{1/2}.
\]
For $|x|\ge 1$,
\begin{equation}\label{eq:app-high-frequency}
 \omega_c(x)=c|x|^\alpha+R_c(x),\qquad |R_c(x)|\le c^{-1}|x|^{-\alpha}.
\end{equation}
\end{lemma}

\begin{proof}
The identity
\[
 \omega_c(x)-c|x|^\alpha
 =\frac{\omega_c(x)^2-c^2|x|^{2\alpha}}{\omega_c(x)+c|x|^\alpha}
 =\frac{1}{\omega_c(x)+c|x|^\alpha}
\]
gives $0\le \omega_c(x)-c|x|^\alpha\le c^{-1}|x|^{-\alpha}$.
\end{proof}

\begin{lemma}[Proportional low-frequency perturbation of the radius]\label{lem:app-low-perturb}
Assume \eqref{eq:app-annulus} and $|q|\le\delta N$ with $0<\delta\le a_0/4$.  Then
\begin{equation}\label{eq:app-segment-lower}
 |\ell+tq|\ge \frac{3a_0}{4}N,\qquad 0\le t\le 1,
\end{equation}
and
\begin{equation}\label{eq:app-low-perturb}
 \bigl||\ell+q|^\alpha-|\ell|^\alpha\bigr|
 \le C\delta N^\alpha.
\end{equation}
\end{lemma}

\begin{proof}
The segment bound follows from
\[
 |\ell+tq|\ge |\ell|-|q|\ge (a_0-\delta)N\ge \frac{3a_0}{4}N.
\]
On this segment, $x\mapsto |x|^\alpha$ is smooth and
\[
 |\nabla |x|^\alpha|=\alpha |x|^{\alpha-1}\le C N^{\alpha-1}.
\]
The mean-value theorem gives \eqref{eq:app-low-perturb}.
\end{proof}

\begin{proposition}[Different-speed proportional low--high phase gap]\label{prop:phase-gap-annulus}
Let $0<\alpha\le 1$, $c_i,c_j>0$, and $c_i\ne c_j$.  There exist $\delta_*>0$, $N_0\ge1$, and constants $c_*,C_*>0$ such that, whenever \eqref{eq:app-annulus} holds, $|q|\le\delta_*N$, and $N\ge N_0$,
\begin{equation}\label{eq:app-two-sided-gap}
 c_*N^\alpha\le |\omega_i(\ell+q)-\omega_j(\ell)|\le C_*N^\alpha.
\end{equation}
Moreover,
\begin{equation}\label{eq:app-plus-gap}
 \omega_i(\ell+q)+\omega_j(\ell)\ge c_*N^\alpha.
\end{equation}
The lower-bound constant degenerates linearly as $|c_i-c_j|\to0$.  Finitely many shells $N<N_0$ are finite-dimensional and are absorbed into the low-frequency constants of the operator estimates.
\end{proposition}

\begin{proof}
By Lemmas~\ref{lem:app-high-frequency} and~\ref{lem:app-low-perturb},
\begin{align}
 \omega_i(\ell+q)-\omega_j(\ell)
 &=c_i|\ell+q|^\alpha-c_j|\ell|^\alpha+O(N^{-\alpha})  \notag\\
 &=(c_i-c_j)|\ell|^\alpha+c_i\bigl(|\ell+q|^\alpha-|\ell|^\alpha\bigr)+O(N^{-\alpha}). \label{eq:app-gap-expansion}
\end{align}
Hence
\[
 |\omega_i(\ell+q)-\omega_j(\ell)|
 \ge |c_i-c_j|a_0^\alpha N^\alpha-C\delta N^\alpha-CN^{-\alpha}.
\]
Choose $\delta_*>0$ so small that $C\delta_*\le \frac14|c_i-c_j|a_0^\alpha$, and then choose $N_0$ so large that $CN^{-\alpha}\le \frac14|c_i-c_j|a_0^\alpha N^\alpha$ for $N\ge N_0$.  This gives the lower bound in \eqref{eq:app-two-sided-gap}.  The upper bound follows from $|\ell+q|\lesssim N$, $|\ell|\lesssim N$, and the high-frequency expansion.  The plus branch follows from positivity of the speeds and $|\ell+q|\sim |\ell|\sim N$.
\end{proof}

\begin{proposition}[Same-speed pure difference]\label{prop:same-speed-contrast}
Let $\omega=\omega_c$ with $c>0$.  Under \eqref{eq:app-annulus} and $|q|\le\delta N$ with $0<\delta\le a_0/4$,
\begin{equation}\label{eq:app-same-speed-upper}
 |\omega(\ell+q)-\omega(\ell)|\le C N^{\alpha-1}|q|\le C\delta N^\alpha.
\end{equation}
More precisely,
\begin{equation}\label{eq:app-same-speed-taylor}
 \omega(\ell+q)-\omega(\ell)=\nabla\omega(\ell)\cdot q+O(N^{\alpha-2}|q|^2),\qquad |\nabla\omega(\ell)|\lesssim N^{\alpha-1}.
\end{equation}
Thus the same-speed difference lacks a uniform lower bound of order $N^\alpha$ on the low--high region; it vanishes at $q=0$ and is controlled by the smaller scale $N^{\alpha-1}|q|$.
\end{proposition}

\begin{proof}
Differentiating
\[
 \omega_c(x)=(1+c^2|x|^{2\alpha})^{1/2}
\]
gives
\[
 \nabla\omega_c(x)=\frac{\alpha c^2|x|^{2\alpha-2}x}{(1+c^2|x|^{2\alpha})^{1/2}},
\]
so $|\nabla\omega_c(x)|\lesssim |x|^{\alpha-1}$ on the high-frequency annulus.  The mean-value theorem gives \eqref{eq:app-same-speed-upper}.  Taylor's formula on the segment $\ell+tq$ gives \eqref{eq:app-same-speed-taylor}, using the corresponding bound on the second derivatives of $\omega_c$ on the annulus.
\end{proof}

\section{First Picard smoothing from the three-frequency phase geometry}\label{app:first-picard-phase}

This appendix proves the first Picard smoothing estimate used in the full-input formulation.  The first Picard objects are
\[
        V_i=I_i(\Psi_1\Psi_2),
        \qquad
        I_i f(t)=\int_0^t K_i(t-s,D)f(s)\,ds,
\]
where
\[
        K_a(t,n)=\frac{\sin(t\omega_a(n))}{\omega_a(n)},
        \qquad
        \omega_a(n)=\bigl(1+c_a^2|n|^{2\alpha}\bigr)^{1/2}.
\]
Throughout this appendix
\[
        \frac34<\alpha\le1,
        \qquad
        c_1,c_2>0,
        \qquad
        c_1\ne c_2.
\]
For independent colors, \(\Psi_1\Psi_2\) is the finite cross-color Wick product at the first Picard level.

We write
\[
        \la n\ra=(1+|n|^2)^{1/2}.
\]
All dyadic variables are at least one.  Constants may depend on \(T,\alpha,c_1,c_2\), the fixed Littlewood-Paley cutoffs, and the dyadic convention, but never on the Fourier cutoff.  The notation \(A\lesssim_\eps B\) means that an arbitrarily small dyadic loss is allowed in the relevant dyadic scale.  Finitely many low shells are always absorbed into the constants.

Stochastic Fubini and the Wiener isometry reduce the smoothing estimate to the deterministic phase-sum bound in Proposition~\ref{prop:B-phase-sum}; cutoff convergence follows from its tail estimate.

\subsection{First Picard estimate}

For definiteness we first use the standard sharp cutoff
\[
        \chi_\Lambda(m)=\mathbf 1_{\{|m|\le\Lambda\}},
\]
although the cutoff-limit argument applies to every admissible family in Definition~\ref{def:admissible-cutoff-family}.  Define
\[
        \widehat\Psi_{a,\Lambda}(m,t)
        =
        \chi_\Lambda(m)
        \int_0^t K_a(t-s,m)\,d\beta_m^a(s),
        \qquad
        V_{i,\Lambda}=I_i(\Psi_{1,\Lambda}\Psi_{2,\Lambda}).
\]
All estimates are uniform in \(\Lambda\).  Inserting the outer Galerkin cutoff only multiplies the formulas by \(\chi_\Lambda(n)\) and does not change the limit.

\begin{theorem}[First Picard smoothing]\label{thm:B-first-picard}
For every \(T<\infty\), every \(\eps>0\), and \(i=1,2\),
\begin{equation}\label{B:fixed-V}
        \sup_{t\le T}\E\bigl|\widehat V_i(n,t)\bigr|^2
        \lesssim_{T,\eps}
        \la n\ra^{3-7\alpha+\eps},
\end{equation}
\begin{equation}\label{B:fixed-dtV}
        \sup_{t\le T}\E\bigl|\partial_t\widehat V_i(n,t)\bigr|^2
        \lesssim_{T,\eps}
        \la n\ra^{3-5\alpha+\eps}.
\end{equation}
Moreover, for sufficiently small \(0<\eta\le\eta_0(\alpha)\),
\begin{equation}\label{B:incr-V}
        \E\bigl|\widehat V_i(n,t)-\widehat V_i(n,t')\bigr|^2
        \lesssim_{T,\eps,\eta}
        |t-t'|^{2\eta}\la n\ra^{3-7\alpha+2\alpha\eta+\eps},
\end{equation}
\begin{equation}\label{B:incr-dtV}
        \E\bigl|\partial_t\widehat V_i(n,t)-\partial_t\widehat V_i(n,t')\bigr|^2
        \lesssim_{T,\eps,\eta}
        |t-t'|^{2\eta}\la n\ra^{3-5\alpha+2\alpha\eta+\eps}.
\end{equation}
The cutoff approximations \(V_{i,\Lambda}\) converge almost surely and in the stated topology in
\[
        C\bigl([0,T];\C^s(\T^3)\bigr)
        \qquad
        \text{for every } s<\frac{7\alpha}{2}-3,
\]
and \(\partial_tV_{i,\Lambda}\) converge almost surely and in the stated topology in
\[
        C\bigl([0,T];\C^s(\T^3)\bigr)
        \qquad
        \text{for every } s<\frac{5\alpha}{2}-3.
\]
Consequently,
\[
        V_i\in C\bigl([0,T];\C^{\frac{7\alpha}{2}-3-}\bigr),
        \qquad
        \partial_tV_i\in C\bigl([0,T];\C^{\frac{5\alpha}{2}-3-}\bigr).
\]
\end{theorem}

\subsection{A coefficient-to-Besov path criterion}

The following criterion is used only for random distributions belonging to a fixed finite sum of Wiener chaoses.  It converts coefficient-level second-moment estimates into almost sure Besov-H\"older path regularity.  The statement separates uniform moment bounds from the cutoff Cauchy argument.

\begin{lemma}[Coefficient-to-Besov path criterion]\label{lem:B-coeff-besov}
Let \(\mathcal A\) be a countable cutoff index set, and let \((Z_\Lambda)_{\Lambda\in\mathcal A}\) be spatially homogeneous random distributions on \(\T^3\).  Assume that, for every dyadic block \(P_N\), every \(t\in[0,T]\), and every \(\Lambda\in\mathcal A\), the random variables \(P_NZ_\Lambda(t,x)\) and their time increments belong to a fixed finite sum of homogeneous Wiener chaoses of order at most \(m\), uniformly in \(N,t,x,\Lambda\).

Suppose that, for some \(p\in\R\), there exists \(\eta_*>0\) such that, for every sufficiently small \(\eps>0\) and every \(0<\eta\le\eta_*\),
\begin{equation}\label{B:coeff-fixed-ass}
        \sup_{\Lambda\in\mathcal A}\sup_{t\le T}
        \E\bigl|\widehat Z_\Lambda(n,t)\bigr|^2
        \lesssim_{T,\eps}
        \la n\ra^{p+\eps},
\end{equation}
\begin{equation}\label{B:coeff-incr-ass}
        \sup_{\Lambda\in\mathcal A}
        \E\bigl|\widehat Z_\Lambda(n,t)-\widehat Z_\Lambda(n,t')\bigr|^2
        \lesssim_{T,\eps,\eta}
        |t-t'|^{2\eta}\la n\ra^{p+2\alpha\eta+\eps}.
\end{equation}
Then, for every \(\kappa>0\) and every finite \(r\ge1\), after choosing \(\eta,\eps>0\) sufficiently small, one has the uniform dyadic moment bound
\begin{equation}\label{B:block-moment-criterion}
        \sup_{\Lambda\in\mathcal A}
        \bigl\|P_NZ_\Lambda\bigr\|_{L^r(\Omega;C([0,T];L^\infty_x))}
        \lesssim_{r,T,m,\kappa}
        N^{\frac{3}{2}+\frac{p}{2}+\kappa}
\end{equation}
for all dyadic \(N\ge1\).  Consequently, for every
\[
        s<-\frac{3}{2}-\frac{p}{2},
\]
one has
\begin{equation}\label{B:besov-moment-criterion}
        \sup_{\Lambda\in\mathcal A}
        \bigl\|Z_\Lambda\bigr\|_{L^r(\Omega;C([0,T];\C^s))}<\infty
\end{equation}
for every finite \(r\ge1\).  In particular, each fixed cutoff realization belongs to \(C([0,T];\C^s)\) almost surely.

Moreover, let \((\Lambda_j)_{j\ge1}\subset\mathcal A\) be an increasing cutoff sequence and set
\[
        Y_j:=Z_{\Lambda_{j+1}}-Z_{\Lambda_j}.
\]
Assume that, for some sufficiently small \(0<\eta\le\eta_*\), there are deterministic dyadic majorants \(R_j^{(\eta)}(N)\) such that, for all dyadic \(N\),
\begin{equation}\label{B:tail-major-fixed}
        \sup_{\la n\ra\sim N}\sup_{t\le T}
        \E\bigl|\widehat Y_j(n,t)\bigr|^2
        \le R_j^{(\eta)}(N),
\end{equation}
\begin{equation}\label{B:tail-major-incr}
        \sup_{\la n\ra\sim N}
        \sup_{0\le t,t'\le T}
        \frac{
        \E\bigl|\widehat Y_j(n,t)-\widehat Y_j(n,t')\bigr|^2
        }{|t-t'|^{2\eta}}
        \le R_j^{(\eta)}(N).
\end{equation}
If, for some \(\delta>0\),
\begin{equation}\label{B:tail-summable-criterion}
        \sum_{j=1}^\infty
        \left(
        \sum_{N\ \mathrm{dyadic}}
        N^{2s+3+\delta}R_j^{(\eta)}(N)
        \right)^{1/2}<\infty,
\end{equation}
then \((Z_{\Lambda_j})_{j\ge1}\) is almost surely Cauchy in \(C([0,T];\C^s)\).  If the same tail estimate holds for the difference of any two admissible cutoff sequences, then the limiting object is independent of the particular cutoff sequence.
\end{lemma}

\begin{proof}
Fix a dyadic shell \(\la n\ra\sim N\).  By spatial homogeneity, Fourier covariances are diagonal in the form
\[
        \E\bigl[\widehat Z_\Lambda(n,t)\overline{\widehat Z_\Lambda(m,t)}\bigr]=0
        \qquad\text{unless }n=m.
\]
Equivalently, in the real-valued Fourier convention, \(\E[\widehat Z_\Lambda(n,t)\widehat Z_\Lambda(m,t)]\) vanishes unless \(n+m=0\).  Hence, for each fixed \((t,x)\),
\[
        \E|P_NZ_\Lambda(t,x)|^2
        \lesssim
        \sum_{\la n\ra\sim N}\E|\widehat Z_\Lambda(n,t)|^2
        \lesssim
        N^{3+p+\eps}.
\]
Since \(P_NZ_\Lambda(t,x)\) belongs to a fixed finite sum of Wiener chaoses of order at most \(m\), Nelson hypercontractivity gives, for every finite \(q\ge2\),
\[
        \|P_NZ_\Lambda(t,x)\|_{L^q(\Omega)}
        \lesssim_{q,m,\eps}
        N^{\frac{3}{2}+\frac{p}{2}+\eps}.
\]

Let \(G_N\) be an \(N^{-A}\)-net of \(\T^3\), with \(A>2\) fixed.  Since \(\#G_N\lesssim N^{3A}\),
\[
        \left\|\max_{x\in G_N}|P_NZ_\Lambda(t,x)|\right\|_{L^q(\Omega)}
        \lesssim_{q,m,A,\eps}
        N^{\frac{3}{2}+\frac{p}{2}+\eps+\frac{3A}{q}}.
\]
The Bernstein inequality gives
\[
        \|\nabla P_NZ_\Lambda(t)\|_{L^\infty_x}
        \lesssim N\|P_NZ_\Lambda(t)\|_{L^\infty_x}.
\]
Thus the oscillation on each \(N^{-A}\)-cell is bounded by
\[
        C N^{1-A}\|P_NZ_\Lambda(t)\|_{L^\infty_x}.
\]
For \(A>2\), this term is absorbed into the left-hand side, after changing the constant and absorbing finitely many low dyadic shells.  Taking \(q\) large and then \(\eps>0\) small gives
\[
        \sup_{\Lambda\in\mathcal A}
        \|P_NZ_\Lambda(t)\|_{L^q(\Omega;L^\infty_x)}
        \lesssim_{q,m,\kappa}
        N^{\frac{3}{2}+\frac{p}{2}+\kappa}.
\]

The same argument applied to \(P_N(Z_\Lambda(t)-Z_\Lambda(t'))\), using \eqref{B:coeff-incr-ass}, yields
\[
        \sup_{\Lambda\in\mathcal A}
        \|P_N(Z_\Lambda(t)-Z_\Lambda(t'))\|_{L^q(\Omega;L^\infty_x)}
        \lesssim
        |t-t'|^\eta
        N^{\frac{3}{2}+\frac{p}{2}+\alpha\eta+\kappa}.
\]
Given the desired final dyadic loss \(\kappa>0\), choose \(\eta,\eps\), and the mesh loss \(3A/q\), sufficiently small so that the factor \(N^{\alpha\eta}\) is absorbed into \(N^\kappa\).

Choose \(q\) large enough so that \(q>r\) and \(q\eta>1\).  Let \(0<\gamma<\eta\) satisfy \(\gamma q>1\).  The fractional Sobolev embedding
\[
        W^{\gamma,q}([0,T];L^\infty_x)\hookrightarrow C([0,T];L^\infty_x)
\]
then gives, by Fubini and the two previous estimates,
\[
        \sup_{\Lambda\in\mathcal A}
        \|P_NZ_\Lambda\|_{L^q(\Omega;C([0,T];L^\infty_x))}
        \lesssim_{q,T,m,\kappa}
        N^{\frac{3}{2}+\frac{p}{2}+\kappa}.
\]
Since \(q\ge r\), this proves \eqref{B:block-moment-criterion}.

Let \(s<-\frac{3}{2}-\frac{p}{2}\).  Choose \(\kappa>0\) so small that
\[
        s+\frac{3}{2}+\frac{p}{2}+\kappa<0.
\]
Using
\[
        \|Z\|_{\C^s}
        \sim
        \sup_{N\ \mathrm{dyadic}}N^s\|P_NZ\|_{L^\infty_x},
\]
and the embedding \(\ell^r\hookrightarrow\ell^\infty\), we obtain
\[
\begin{aligned}
        \|Z_\Lambda\|_{L^r(\Omega;C([0,T];\C^s))}
        &\lesssim
        \left(
        \sum_{N\ \mathrm{dyadic}}
        N^{r(s+\frac{3}{2}+\frac{p}{2}+\kappa)}
        \right)^{1/r}<\infty,
\end{aligned}
\]
uniformly in \(\Lambda\).  This proves \eqref{B:besov-moment-criterion}.

For the Cauchy statement, apply the preceding block argument to \(Y_j\).  The coefficient majorants \(R_j^{(\eta)}(N)\) give, for every finite \(q\) sufficiently large and every small \(\kappa>0\),
\[
        \|P_NY_j\|_{L^q(\Omega;C([0,T];L^\infty_x))}
        \lesssim
        N^{\frac{3}{2}+\kappa}\bigl(R_j^{(\eta)}(N)\bigr)^{1/2}.
\]
Choose \(\kappa>0\) so small that \(2\kappa<\delta\).  Since \(q\ge2\),
\[
\begin{aligned}
        \|Y_j\|_{L^q(\Omega;C([0,T];\C^s))}
        &\lesssim
        \left(
        \sum_N
        N^{q(s+\frac{3}{2}+\kappa)}
        \bigl(R_j^{(\eta)}(N)\bigr)^{q/2}
        \right)^{1/q} \\
        &\le
        \left(
        \sum_N
        N^{2s+3+2\kappa}R_j^{(\eta)}(N)
        \right)^{1/2} \\
        &\le
        \left(
        \sum_N
        N^{2s+3+\delta}R_j^{(\eta)}(N)
        \right)^{1/2}.
\end{aligned}
\]
The summability assumption \eqref{B:tail-summable-criterion} implies
\[
        \sum_{j=1}^\infty
        \E\|Y_j\|_{C([0,T];\C^s)}<\infty.
\]
By Fubini,
\[
        \sum_{j=1}^\infty
        \|Z_{\Lambda_{j+1}}-Z_{\Lambda_j}\|_{C([0,T];\C^s)}<\infty
        \qquad\text{almost surely}.
\]
Thus \((Z_{\Lambda_j})_j\) is almost surely Cauchy.  The comparison of two cutoff sequences is identical and gives cutoff convergence.
\end{proof}

\subsection{Stochastic Fubini, colored Gaussian products, and the three propagators}

In this section \(\nu\in\{0,1\}\) denotes the number of time derivatives.  The color labels are \(1,2\).  Set
\[
        K_i^{[0]}(t,n)=K_i(t,n)=\frac{\sin(t\omega_i(n))}{\omega_i(n)},
        \qquad
        K_i^{[1]}(t,n)=\partial_tK_i(t,n)=\cos(t\omega_i(n)).
\]
Thus
\[
        \partial_t^0V_i=V_i,
        \qquad
        \partial_t^1V_i=\partial_tV_i.
\]
For compatibility with the notation of the phase estimates below, we also write \(K_i^{(a)}:=K_i^{[a]}\) when \(a\in\{0,1\}\).

\paragraph{The colored Gaussian Hilbert space}

Let
\[
        \mathfrak h_T=L^2([0,T]\times\T^3;\R),
\]
and let \(\mathfrak h_{T,\C}\) be its complexification.  Let
\[
        e_m(x)=(2\pi)^{-3/2}e^{im\cdot x},
        \qquad m\in\Z^3.
\]
The two independent noises are realized as two independent isonormal Gaussian processes
\[
        W_1:\mathfrak h_T\to L^2(\Omega),
        \qquad
        W_2:\mathfrak h_T\to L^2(\Omega).
\]
Equivalently, the full first-chaos Hilbert space is the orthogonal direct sum
\[
        \mathfrak H_T=\mathfrak h_T^{(1)}\oplus\mathfrak h_T^{(2)},
        \qquad
        \mathfrak h_T^{(1)}\perp\mathfrak h_T^{(2)}.
\]
We extend \(W_b\) complex-linearly to \(\mathfrak h_{T,\C}\).  Thus, for \(h,g\in\mathfrak h_{T,\C}\),
\[
        \E\bigl[W_b(h)\overline{W_{b'}(g)}\bigr]
        =
        \mathbf 1_{\{b=b'\}}\la h,g\ra_{\mathfrak h_{T,\C}},
\]
where \(\la\cdot,\cdot\ra\) is the Hermitian inner product.  For products without complex conjugation, the relevant pairing is the complex bilinear pairing
\[
        (h,g)_{\rm bil}:=
        \int_0^T\int_{\T^3}h(r,x)g(r,x)\,dx\,dr,
\]
so that
\[
        \E\bigl[W_b(h)W_{b'}(g)\bigr]
        =
        \mathbf 1_{\{b=b'\}}(h,g)_{\rm bil}.
\]
The identity \(W_b(\overline h)=\overline{W_b(h)}\) is understood throughout.

For \(m\in\Z^3\), define the complex Fourier Brownian coordinate
\[
        \beta_m^b(t):=
        W_b\bigl(\mathbf 1_{[0,t]}(r)e_{-m}(x)\bigr).
\]
Then
\[
        \overline{\beta_m^b(t)}=\beta_{-m}^b(t),
\]
and the finite covariance identity is
\[
        \E\bigl[\beta_m^b(t)\beta_\ell^{b'}(s)\bigr]
        =
        \mathbf 1_{\{b=b'\}}\mathbf 1_{\{m+\ell=0\}}(t\wedge s).
\]
This identity follows from the Hilbert-space construction of the covariance.

\paragraph{Finite-cutoff stochastic convolutions}

For \(b=1,2\), define
\[
        h_{b,m,s}^{\Lambda}(r,x)
        :=
        \chi_\Lambda(m)\mathbf 1_{\{0\le r\le s\}}K_b(s-r,m)e_{-m}(x)
        \in\mathfrak h_{T,\C}.
\]
The finite-cutoff stochastic convolution coefficient is
\[
        \widehat\Psi_{b,\Lambda}(m,s):=W_b(h_{b,m,s}^{\Lambda}).
\]
Equivalently,
\[
        \widehat\Psi_{b,\Lambda}(m,s)
        =
        \chi_\Lambda(m)\int_0^s K_b(s-r,m)\,d\beta_m^b(r),
\]
where the integral notation denotes the first Wiener integral \(W_b(h_{b,m,s}^{\Lambda})\).

For fixed output frequency \(n\), set
\[
        F_{\Lambda,n}(s)
        :=
        \sum_{k\in\Z^3}
        \widehat\Psi_{1,\Lambda}(k,s)\widehat\Psi_{2,\Lambda}(n-k,s).
\]
The cutoff imposes \(|k|\le\Lambda\) and \(|n-k|\le\Lambda\), so this sum is finite.  Hence
\[
        \widehat V_{i,\Lambda}(n,t)
        =
        \int_0^t K_i(t-s,n)F_{\Lambda,n}(s)\,ds.
\]
Since \(K_i(0,n)=0\), differentiating the outer Duhamel integral gives a vanishing boundary term:
\[
\begin{aligned}
        \partial_t\widehat V_{i,\Lambda}(n,t)
        &=K_i(0,n)F_{\Lambda,n}(t)
        +\int_0^t\partial_tK_i(t-s,n)F_{\Lambda,n}(s)\,ds \\
        &=\int_0^t\cos((t-s)\omega_i(n))F_{\Lambda,n}(s)\,ds.
\end{aligned}
\]
Therefore, for \(\nu=0,1\),
\begin{equation}\label{B:Duhamel-nu-finite}
        \partial_t^\nu\widehat V_{i,\Lambda}(n,t)
        =
        \int_0^tK_i^{[\nu]}(t-s,n)F_{\Lambda,n}(s)\,ds.
\end{equation}
All statements in this paragraph are identities in \(L^2(\Omega)\), and at finite cutoff they are finite-dimensional.

\paragraph{The colored first-chaos product formula}

Let \(I_{1,2}\) denote the ordered colored double Wiener integral over
\[
        \mathfrak h_{T,\C}^{(1)}\otimes\mathfrak h_{T,\C}^{(2)}.
\]
For simple tensors it is defined by
\[
        I_{1,2}(h\otimes g):=W_1(h)W_2(g),
        \qquad
        h\in\mathfrak h_{T,\C}^{(1)},\quad
        g\in\mathfrak h_{T,\C}^{(2)}.
\]
Because \(W_1\) and \(W_2\) are independent, this map is an isometry:
\[
        \E\bigl[I_{1,2}(H)\overline{I_{1,2}(G)}\bigr]
        =
        \la H,G\ra_{\mathfrak h_{T,\C}^{(1)}\otimes\mathfrak h_{T,\C}^{(2)}}.
\]
Therefore \(I_{1,2}\) extends by completion to the full Hilbert tensor product.

In particular, if \(h\in\mathfrak h_{T,\C}^{(1)}\) and \(g\in\mathfrak h_{T,\C}^{(2)}\), then
\[
        W_1(h)W_2(g)=I_{1,2}(h\otimes g).
\]
The contraction term vanishes by the orthogonality \(\mathfrak h_T^{(1)}\perp\mathfrak h_T^{(2)}\).  By contrast, for two first-chaos variables of the same color, the ordinary product formula is
\[
        W_b(h)W_b(g)
        =
        I_{b,b}^{:2:}(h\widetilde\otimes g)+(h,g)_{\rm bil},
\]
where \(I_{b,b}^{:2:}\) denotes the second Wick integral and \((h,g)_{\rm bil}\) is the contraction.  Thus the contraction is created by the first-chaos product formula and is exactly the Hilbert-space pairing of the deterministic kernels.

More generally, in the matrix-covariance convention above,
\[
        W_b(h)W_{b'}(g)
        =
        I_{b,b'}^{:2:}(h\otimes g)+\mathsf R_{bb'}(h,g)_{\rm bil}.
\]
Thus cross-color contractions reappear precisely when \(\mathsf R_{12}\ne0\).

Applying the independent-color formula with
\[
        h=h_{1,k,s}^{\Lambda},
        \qquad
        g=h_{2,n-k,s}^{\Lambda},
\]
we obtain
\begin{equation}\label{B:cross-product-exact}
\begin{aligned}
        \widehat\Psi_{1,\Lambda}(k,s)\widehat\Psi_{2,\Lambda}(n-k,s)
        &=
        I_{1,2}\bigl(h_{1,k,s}^{\Lambda}\otimes h_{2,n-k,s}^{\Lambda}\bigr) \\
        &=
        \chi_\Lambda(k)\chi_\Lambda(n-k)
        I_{1,2}\Bigl(
        \mathbf 1_{\{0\le u\le s\}}
        \mathbf 1_{\{0\le v\le s\}} \\
        &\hspace{3.2cm}\times
        K_1(s-u,k)K_2(s-v,n-k)e_{-k}(x_1)e_{-(n-k)}(x_2)
        \Bigr).
\end{aligned}
\end{equation}
In the traditional notation for this already-defined ordered colored Wiener integral,
\[
\begin{aligned}
        \widehat\Psi_{1,\Lambda}(k,s)\widehat\Psi_{2,\Lambda}(n-k,s)
        &=
        \chi_\Lambda(k)\chi_\Lambda(n-k)
        \int_0^s\int_0^s
        K_1(s-u,k)K_2(s-v,n-k) \\
        &\hspace{4.1cm}\times
        d\beta_k^1(u)d\beta_{n-k}^2(v).
\end{aligned}
\]
This last display is only shorthand for \eqref{B:cross-product-exact}.

\paragraph{Stochastic Fubini at the Hilbert-kernel level}

Substituting \eqref{B:cross-product-exact} into \eqref{B:Duhamel-nu-finite}, we obtain
\[
\begin{aligned}
        \partial_t^\nu\widehat V_{i,\Lambda}(n,t)
        &=
        \sum_{k\in\Z^3}
        \chi_\Lambda(k)\chi_\Lambda(n-k)
        \int_0^tK_i^{[\nu]}(t-s,n) \\
        &\qquad\times
        I_{1,2}\Bigl(
        \mathbf 1_{\{0\le u\le s\}}
        \mathbf 1_{\{0\le v\le s\}}
        K_1(s-u,k)K_2(s-v,n-k) \\
        &\hspace{5.1cm}\times
        e_{-k}(x_1)e_{-(n-k)}(x_2)
        \Bigr)\,ds.
\end{aligned}
\]
For fixed \(n,k,t,\nu\), define the deterministic tensor kernel
\[
\begin{aligned}
        H_{i,\nu}^{\Lambda}(n,k;t;s;u,x_1,v,x_2)
        &:={}
        \chi_\Lambda(k)\chi_\Lambda(n-k)
        \mathbf 1_{\{0\le u\le s\}}
        \mathbf 1_{\{0\le v\le s\}} \\
        &\quad\times
        K_i^{[\nu]}(t-s,n)K_1(s-u,k)K_2(s-v,n-k)
        e_{-k}(x_1)e_{-(n-k)}(x_2).
\end{aligned}
\]
At finite cutoff, the number of admissible \(k\)'s is finite.  Moreover,
\[
        s\mapsto H_{i,\nu}^{\Lambda}(n,k;t;s)
\]
is Bochner integrable in \(\mathfrak h_{T,\C}^{(1)}\otimes\mathfrak h_{T,\C}^{(2)}\).  Since
\[
        I_{1,2}:\mathfrak h_{T,\C}^{(1)}\otimes\mathfrak h_{T,\C}^{(2)}\longrightarrow L^2(\Omega)
\]
is a continuous linear isometry, we may commute the deterministic \(s\)-integration with \(I_{1,2}\):
\[
        \int_0^tI_{1,2}\bigl(H_{i,\nu}^{\Lambda}(n,k;t;s)\bigr)\,ds
        =
        I_{1,2}\left(\int_0^tH_{i,\nu}^{\Lambda}(n,k;t;s)\,ds\right)
\]
in \(L^2(\Omega)\).  This is the stochastic Fubini step.  It is a deterministic Bochner-Fubini identity in the Hilbert tensor product before applying the colored double Wiener integral.

For \(0\le u,v\le t\),
\[
\begin{aligned}
        \int_0^t H_{i,\nu}^{\Lambda}(n,k;t;s;u,x_1,v,x_2)\,ds
        &=
        \chi_\Lambda(k)\chi_\Lambda(n-k)e_{-k}(x_1)e_{-(n-k)}(x_2) \\
        &\quad\times
        \int_{u\vee v}^tK_i^{[\nu]}(t-s,n)K_1(s-u,k)K_2(s-v,n-k)\,ds.
\end{aligned}
\]
Therefore
\begin{equation}\label{B:stoch-fubini-rep}
        \partial_t^\nu\widehat V_{i,\Lambda}(n,t)
        =
        \sum_{k\in\Z^3}
        \int_0^t\int_0^t
        \calK_{i,\nu}^{\Lambda}(n,k;t,u,v)
        \,d\beta_k^1(u)d\beta_{n-k}^2(v),
        \qquad \nu=0,1,
\end{equation}
where the double-integral notation denotes the ordered colored Wiener integral defined above, and
\[
        \calK_{i,\nu}^{\Lambda}(n,k;t,u,v)
        =
        \chi_\Lambda(k)\chi_\Lambda(n-k)\calK_{i,\nu}(n,k;t,u,v),
\]
with
\begin{equation}\label{B:Kia-def}
        \calK_{i,\nu}(n,k;t,u,v)
        =
        \int_{u\vee v}^tK_i^{[\nu]}(t-s,n)K_1(s-u,k)K_2(s-v,n-k)\,ds.
\end{equation}
The convention is that the integral in \eqref{B:Kia-def} is zero if \(u\vee v>t\).  In \eqref{B:stoch-fubini-rep} we integrate only over \(0\le u,v\le t\).  For \(a\in\{0,1\}\), write \(\calK_{i,a}:=\calK_{i,\nu}|_{\nu=a}\).

\paragraph{Double Wiener isometry and rough absolute bounds}

Since the tensors \(e_{-k}\otimes e_{-(n-k)}\) are orthonormal in \(\mathfrak h_{T,\C}^{(1)}\otimes\mathfrak h_{T,\C}^{(2)}\) as \(k\) varies, \eqref{B:stoch-fubini-rep} gives the exact ordered double-Wiener isometry
\begin{equation}\label{B:double-isometry}
        \E\left|\partial_t^\nu\widehat V_{i,\Lambda}(n,t)\right|^2
        =
        \sum_{k\in\Z^3}\int_0^t\int_0^t
        \left|\calK_{i,\nu}^{\Lambda}(n,k;t,u,v)\right|^2\,du\,dv.
\end{equation}
This is the precise point where the stochastic estimate becomes a deterministic three-propagator kernel estimate.

The elementary bounds
\[
        |K_b(\tau,m)|\lesssim\la m\ra^{-\alpha},
        \qquad
        |K_i^{[1]}(\tau,n)|\le1
\]
give, for \(\nu=0\),
\[
\begin{aligned}
        \sum_k\int_0^t\int_0^t|\calK_{i,0}(n,k;t,u,v)|^2\,du\,dv
        &\lesssim_T
        \la n\ra^{-2\alpha}
        \sum_k\la k\ra^{-2\alpha}\la n-k\ra^{-2\alpha}<\infty
\end{aligned}
\]
when \(\alpha>3/4\).  This absolute estimate yields only the static nonoscillatory power
\[
        \E|\widehat V_i(n,t)|^2\lesssim\la n\ra^{3-6\alpha+}.
\]
For \(\nu=1\), the outer Duhamel denominator is absent, and the rough bound gives only
\[
        \E|\partial_t\widehat V_i(n,t)|^2\lesssim\la n\ra^{3-4\alpha+}.
\]
The powers \(3-7\alpha+\) and \(3-5\alpha+\) follow after retaining the oscillatory structure of the full kernel \(\calK_{i,\nu}\).

\subsection{Oscillatory kernel bounds}

This section records deterministic oscillatory estimates for the kernel \eqref{B:Kia-def}.  It prepares the three-frequency phase factor that will be combined with the later lattice geometry and speed-separation estimates.

For \(\lambda\in\R\), set
\[
        \frakm(\lambda):=\min\{1,|\lambda|^{-1}\},
        \qquad
        \frakm(0):=1.
\]
For \(\sigma=(\sigma_0,\sigma_1,\sigma_2)\in\{\pm1\}^3\), define the three-frequency phase
\begin{equation}\label{B:phase-def}
        \Phi_i^\sigma(n,k)
        =
        \sigma_0\omega_i(n)+\sigma_1\omega_1(k)+\sigma_2\omega_2(n-k).
\end{equation}
We write
\[
        N_0=\la n\ra,
        \qquad
        N_1=\la k\ra,
        \qquad
        N_2=\la n-k\ra,
        \qquad
        N_{\max}=\max(N_0,N_1,N_2).
\]
Throughout this section,
\[
        \omega_i(n)\sim N_0^\alpha,
        \qquad
        \omega_1(k)\sim N_1^\alpha,
        \qquad
        \omega_2(n-k)\sim N_2^\alpha,
\]
and
\[
        |\Phi_i^\sigma(n,k)|\le
        \omega_i(n)+\omega_1(k)+\omega_2(n-k)
        \lesssim N_{\max}^\alpha.
\]

\begin{lemma}[Static oscillatory kernel]\label{lem:B-static-kernel}
For \(\nu=0,1\), uniformly in \(0\le u,v\le t\le T\),
\begin{equation}\label{B:static-kernel}
        |\calK_{i,\nu}(n,k;t,u,v)|
        \lesssim_T
        \frac{N_0^{\nu\alpha}}{N_0^\alpha N_1^\alpha N_2^\alpha}
        \sum_{\sigma\in\{\pm1\}^3}\frakm\bigl(\Phi_i^\sigma(n,k)\bigr).
\end{equation}
\end{lemma}

\begin{proof}
We use
\[
        K_i^{[0]}(\tau,n)=
        \frac{e^{i\tau\omega_i(n)}-e^{-i\tau\omega_i(n)}}{2i\omega_i(n)},
        \qquad
        K_i^{[1]}(\tau,n)=
        \frac{e^{i\tau\omega_i(n)}+e^{-i\tau\omega_i(n)}}{2}.
\]
Equivalently, for \(\nu=0,1\),
\[
        K_i^{[\nu]}(\tau,n)
        =
        \omega_i(n)^{\nu-1}
        \sum_{\varepsilon_0=\pm1}c_{\nu,\varepsilon_0}e^{i\varepsilon_0\tau\omega_i(n)},
\]
with uniformly bounded constants \(c_{\nu,\varepsilon_0}\).  Similarly,
\[
        K_b(\tau,m)=
        \omega_b(m)^{-1}
        \sum_{\varepsilon=\pm1}c_\varepsilon e^{i\varepsilon\tau\omega_b(m)},
        \qquad b=1,2.
\]
Substitution into \eqref{B:Kia-def} gives a finite sum over \(\varepsilon=(\varepsilon_0,\varepsilon_1,\varepsilon_2)\in\{\pm1\}^3\).  Each branch has the form
\[
\begin{aligned}
        \calK_{i,\nu}^{\varepsilon}
        &=
        C_{\nu,\varepsilon}
        \omega_i(n)^{\nu-1}\omega_1(k)^{-1}\omega_2(n-k)^{-1} \\
        &\quad\times
        e^{i\varepsilon_0t\omega_i(n)}
        e^{-i\varepsilon_1u\omega_1(k)}
        e^{-i\varepsilon_2v\omega_2(n-k)}
        \int_{u\vee v}^t e^{is\Psi_i^\varepsilon(n,k)}\,ds,
\end{aligned}
\]
where
\[
        \Psi_i^\varepsilon(n,k)
        =
        -\varepsilon_0\omega_i(n)+\varepsilon_1\omega_1(k)+\varepsilon_2\omega_2(n-k).
\]
After relabeling \(\sigma_0=-\varepsilon_0\), \(\sigma_1=\varepsilon_1\), \(\sigma_2=\varepsilon_2\), this phase is \(\Phi_i^\sigma(n,k)\).  The endpoint factor has modulus one.  Moreover,
\[
        \left|\int_{u\vee v}^t e^{is\Phi}\,ds\right|\le T,
\]
and, if \(\Phi\ne0\),
\[
        \left|\int_{u\vee v}^t e^{is\Phi}\,ds\right|
        =
        \left|\frac{e^{it\Phi}-e^{i(u\vee v)\Phi}}{i\Phi}\right|
        \le 2|\Phi|^{-1}.
\]
Hence
\[
        \left|\int_{u\vee v}^t e^{is\Phi}\,ds\right|
        \lesssim_T\frakm(\Phi).
\]
Finally,
\[
        \omega_i(n)^{\nu-1}\omega_1(k)^{-1}\omega_2(n-k)^{-1}
        \lesssim
        \frac{N_0^{\nu\alpha}}{N_0^\alpha N_1^\alpha N_2^\alpha}.
\]
Summing over the \(2^3\) sign branches proves \eqref{B:static-kernel}.
\end{proof}

\begin{lemma}[Kernel time increments]\label{lem:B-kernel-incr}
Let \(0<\eta\le1/2\).  If \(0\le u,v\le t\wedge t'\le T\), then, for \(\nu=0,1\),
\begin{equation}\label{B:kernel-incr}
\begin{aligned}
&|\calK_{i,\nu}(n,k;t,u,v)-\calK_{i,\nu}(n,k;t',u,v)| \\
&\qquad\lesssim_{T,\eta}
 |t-t'|^\eta
 \frac{N_0^{\nu\alpha}N_{\max}^{\alpha\eta}}
      {N_0^\alpha N_1^\alpha N_2^\alpha}
 \sum_{\sigma\in\{\pm1\}^3}\frakm\bigl(\Phi_i^\sigma(n,k)\bigr).
\end{aligned}
\end{equation}
\end{lemma}

\begin{proof}
It suffices to consider \(t'>t\).  Put \(h=t'-t\) and \(r=u\vee v\).  The assumption \(u,v\le t\wedge t'\) ensures \(r\le t<t'\), so the lower integration limit is the same for both kernels.

Using the sign-branch representation from the previous proof, write one branch as
\[
        \calK_{i,\nu}^{\sigma}(t)
        =A_{\nu,\sigma}(t;u,v)J_\sigma(t),
        \qquad
        J_\sigma(t):=\int_r^t e^{is\Phi_i^\sigma(n,k)}\,ds,
\]
with
\[
        |A_{\nu,\sigma}(t;u,v)|
        \lesssim
        \frac{N_0^{\nu\alpha}}{N_0^\alpha N_1^\alpha N_2^\alpha}.
\]
The \(t\)-dependence of \(A_{\nu,\sigma}\) is only through an endpoint factor \(e^{\pm it\omega_i(n)}\).  Hence, using
\[
        |e^{ix}-1|\lesssim\min\{1,|x|\}\le |x|^\eta,
        \qquad 0<\eta\le1,
\]
we have
\[
\begin{aligned}
        |A_{\nu,\sigma}(t';u,v)-A_{\nu,\sigma}(t;u,v)|
        &\lesssim
        |h|^\eta
        \frac{N_0^{\nu\alpha}N_0^{\alpha\eta}}
             {N_0^\alpha N_1^\alpha N_2^\alpha} \\
        &\lesssim
        |h|^\eta
        \frac{N_0^{\nu\alpha}N_{\max}^{\alpha\eta}}
             {N_0^\alpha N_1^\alpha N_2^\alpha}.
\end{aligned}
\]
Together with
\[
        |J_\sigma(t)|\lesssim_T\frakm(\Phi_i^\sigma(n,k)),
\]
this controls the endpoint-factor contribution.

It remains to estimate the change of upper integration limit:
\[
        J_\sigma(t')-J_\sigma(t)
        =
        \int_t^{t'} e^{is\Phi_i^\sigma(n,k)}\,ds.
\]
We claim that, for \(\Phi=\Phi_i^\sigma(n,k)\),
\begin{equation}\label{B:interval-increment-bound}
        \left|\int_t^{t'}e^{is\Phi}\,ds\right|
        \lesssim_{T,\eta}
        |h|^\eta N_{\max}^{\alpha\eta}\frakm(\Phi).
\end{equation}
Indeed,
\[
        \left|\int_t^{t'}e^{is\Phi}\,ds\right|\le |h|,
        \qquad
        \left|\int_t^{t'}e^{is\Phi}\,ds\right|\le 2|\Phi|^{-1}
\]
when \(\Phi\ne0\).  Thus
\[
        \left|\int_t^{t'}e^{is\Phi}\,ds\right|
        \lesssim_T |h|^\eta \frakm(\Phi)^{1-\eta}.
\]
Since \(|\Phi_i^\sigma(n,k)|\lesssim N_{\max}^\alpha\), one has \(\frakm(\Phi_i^\sigma(n,k))\gtrsim N_{\max}^{-\alpha}\).  Consequently,
\[
        \frakm(\Phi)^{1-\eta}
        =
        \frakm(\Phi)\frakm(\Phi)^{-\eta}
        \lesssim
        N_{\max}^{\alpha\eta}\frakm(\Phi),
\]
which proves \eqref{B:interval-increment-bound}.  Multiplying by the amplitude bound and summing over the \(2^3\) sign branches gives \eqref{B:kernel-incr}.
\end{proof}

\subsection{Lattice phase-layer estimates}

This section contains the geometric input.  The estimates are one-sided upper bounds for lattice phase layers, used as summable dyadic majorants.

We begin with the unit-cube thickening step that passes from lattice counts to continuous phase-layer volumes.

\begin{lemma}[Unit-cube thickening of phase layers]\label{lem:B-thickening}
Let
\[
        \Phi_i^\sigma(n,x)
        =
        \sigma_0\omega_i(n)+\sigma_1\omega_1(x)+\sigma_2\omega_2(n-x),
        \qquad x\in\R^3,
\]
be one sign branch of the phase, with the lattice variable replaced by a continuous variable.  Let \(\mathcal A\subset\R^3\) be one of the dyadic annular regions used below, and let
\[
        \mathcal A^{+1}:=
        \{x\in\R^3:\dist(x,\mathcal A)\le C_*\}
\]
for a fixed sufficiently large constant \(C_*\).  If \(L\ge1\), then
\[
        \#\{k\in\Z^3:k\in\mathcal A,\ |\Phi_i^\sigma(n,k)|\le L\}
        \lesssim
        \bigl|\{x\in\mathcal A^{+1}:\ |\Phi_i^\sigma(n,x)|\le CL\}\bigr|.
\]
The constants depend only on \(\alpha\), the speeds, and the fixed dyadic cutoff conventions.
\end{lemma}

\begin{proof}
Let
\[
        Q_k:=k+[-1/2,1/2)^3.
\]
The cubes \(Q_k\) are pairwise disjoint and \(|Q_k|=1\).  Put
\[
        S:=\{k\in\Z^3:k\in\mathcal A,\ |\Phi_i^\sigma(n,k)|\le L\}.
\]
Then
\[
        \#S=\left|\bigcup_{k\in S}Q_k\right|.
\]
If \(x\in Q_k\), then \(|x-k|\le C\), hence \(x\in\mathcal A^{+1}\) after choosing \(C_*\) large enough.  It remains to control the phase variation on \(Q_k\).

For completeness, every region \(\mathcal A\) used below is described by a fixed finite collection of radial constraints
\[
        a_jR_j\le |x-z_j|\le b_jR_j,
\]
where the centers \(z_j\) are either \(0\) or \(n\), and the constants \(a_j,b_j>0\) come from the fixed dyadic partition.  If \(k\) satisfies such a constraint and \(x\in Q_k\), then
\[
        a_jR_j-C\le |x-z_j|\le b_jR_j+C.
\]
For \(R_j\) above a fixed constant this lies in a fixed enlargement of the same dyadic annulus; the remaining bounded shells are finite and are absorbed into the implicit constant.  Thus the single set \(\mathcal A^{+1}\) simultaneously accounts for all radial constraints and for their dyadic boundaries.

For \(a=1,2\),
\[
        \omega_a(\xi)=(1+c_a^2|\xi|^{2\alpha})^{1/2}.
\]
For \(\xi\ne0\),
\[
        |\nabla\omega_a(\xi)|
        =
        \frac{\alpha c_a^2|\xi|^{2\alpha-1}}
             {(1+c_a^2|\xi|^{2\alpha})^{1/2}}.
\]
In the range \(3/4<\alpha\le1\), this derivative is uniformly bounded: it is \(O(|\xi|^{2\alpha-1})\) near the origin and \(O(|\xi|^{\alpha-1})\) at high frequency.  Notice also that the output term \(\omega_i(n)\) is independent of the thickened variable.  Hence the following phase-variation bound is uniform in the fixed output frequency \(n\).

Consequently, for \(x\in Q_k\),
\[
\begin{aligned}
        |\Phi_i^\sigma(n,x)-\Phi_i^\sigma(n,k)|
        &\le
        |\omega_1(x)-\omega_1(k)|
        +|\omega_2(n-x)-\omega_2(n-k)| \\
        &\lesssim |x-k|\lesssim1.
\end{aligned}
\]
Thus, if \(k\in S\), then for every \(x\in Q_k\),
\[
        |\Phi_i^\sigma(n,x)|\le |\Phi_i^\sigma(n,k)|+C_0\le L+C_0\le CL,
\]
since \(L\ge1\).  Therefore
\[
        \bigcup_{k\in S}Q_k
        \subset
        \{x\in\mathcal A^{+1}:\ |\Phi_i^\sigma(n,x)|\le CL\}.
\]
Taking Lebesgue measure proves the claim.
\end{proof}

We also use the two-center coordinate formula.  We state it as a separate geometric lemma in order to make clear that this part of the proof is independent of the stochastic equation.

\begin{lemma}[Two-center coordinates]\label{lem:B-two-center}
Let \(n\ne0\), \(d=|n|\), \(\rho=|x|\), and \(\lambda=|n-x|\).  On the region
\[
        |\rho-\lambda|\le d\le\rho+\lambda,
\]
the Euclidean volume element is
\begin{equation}\label{B:two-center-volume-full}
        dx=\frac{\rho\lambda}{d}\,d\rho\,d\lambda\,d\varphi,
\end{equation}
where \(\varphi\) is the rotation angle around the axis generated by \(n\).  Consequently, after integrating in the angle variable,
\begin{equation}\label{B:two-center-volume}
        dx=2\pi\frac{\rho\lambda}{d}\,d\rho\,d\lambda.
\end{equation}
\end{lemma}

\begin{proof}
By rotation invariance, assume \(n=(0,0,d)\).  Write
\[
        x=\rho(\sin\theta\cos\varphi,\sin\theta\sin\varphi,\cos\theta).
\]
The usual spherical volume element is
\[
        dx=\rho^2\sin\theta\,d\rho\,d\theta\,d\varphi.
\]
The second distance satisfies the cosine-law identity
\[
        \lambda^2=d^2+\rho^2-2d\rho\cos\theta.
\]
For fixed \(\rho\), differentiating in \(\theta\) gives
\[
        2\lambda\,d\lambda=2d\rho\sin\theta\,d\theta,
        \qquad
        \sin\theta\,d\theta=\frac{\lambda}{d\,\rho}\,d\lambda.
\]
Substitution into the spherical volume element yields
\[
        dx=\rho^2\frac{\lambda}{d\,\rho}\,d\rho\,d\lambda\,d\varphi
        =\frac{\rho\lambda}{d}\,d\rho\,d\lambda\,d\varphi.
\]
Integrating \(\varphi\in[0,2\pi)\) gives \eqref{B:two-center-volume}.  The constraint
\(|\rho-\lambda|\le d\le \rho+\lambda\) is exactly the condition that the three lengths can form a triangle.  In upper bounds it may only shrink the integration region.
\end{proof}

If \(d=0\), or if one of the relevant dyadic scales is bounded, finite-shell cases are absorbed into the constants before invoking \eqref{B:two-center-volume}.

\begin{lemma}[One-dimensional phase-layer slicing]\label{lem:B-one-dimensional-layer}
Let \(I\subset\R\) be an interval and let \(f\in C^1(I;\R)\) satisfy
\[
        |f'(r)|\ge m>0
        \qquad\text{for all }r\in I.
\]
Then, for every \(L\ge0\),
\[
        |\{r\in I: |f(r)|\le L\}|
        \le \frac{2L}{m}.
\]
In the lattice-thickened estimates below we use the robust version
\[
        |\{r\in I^{+1}: |f(r)|\le CL\}|
        \lesssim \frac{L}{m}+1,
        \qquad L\ge1,
\]
where the additive \(+1\) accounts for unit-scale thickening and boundary enlargement.
\end{lemma}

\begin{proof}
Since \(|f'|\ge m\), the derivative has a fixed sign on \(I\), hence \(f\) is monotone.  If \(r_1,r_2\in I\) both satisfy \(|f(r_j)|\le L\), then by the mean-value theorem
\[
        m|r_1-r_2|\le |f(r_1)-f(r_2)|\le 2L.
\]
Thus the sublevel set has diameter at most \(2L/m\), and hence length at most \(2L/m\).  The thickened version follows by applying the same estimate to the enlarged interval and adding the fixed unit boundary length.
\end{proof}

\begin{lemma}[Balanced phase layers]\label{lem:B-balanced-layer}
Let \(M\gg1\) be dyadic and fix \(n\) with \(|n|\sim M\).  For every sign branch \(\sigma\) and every dyadic \(1\le L\lesssim M^\alpha\),
\begin{equation}\label{B:balanced-layer}
\begin{aligned}
&\#\{k\in\Z^3: |k|\sim M,\ |n-k|\sim M,
       \ |\Phi_i^\sigma(n,k)|\le L\} \\
&\qquad\lesssim_\eps M^\eps\bigl(M^2+M^{3-\alpha}L\bigr).
\end{aligned}
\end{equation}
\end{lemma}

\begin{proof}
By Lemma~\ref{lem:B-thickening}, it suffices to estimate the corresponding thickened continuous volume.  Put
\[
        d=|n|,
        \qquad
        \rho=|x|,
        \qquad
        \lambda=|n-x|.
\]
In the balanced region,
\[
        d\sim\rho\sim\lambda\sim M.
\]
Thus the density in \eqref{B:two-center-volume} satisfies
\[
        \frac{\rho\lambda}{d}\lesssim M.
\]
For fixed \(\lambda\), the phase is
\[
        \Phi_i^\sigma(n,x)
        =
        \sigma_0\omega_i(d)+\sigma_1\omega_1(\rho)+\sigma_2\omega_2(\lambda),
\]
and therefore
\[
        |\partial_\rho\Phi_i^\sigma|=|\omega_1'(\rho)|\sim M^{\alpha-1}.
\]
By Lemma~\ref{lem:B-one-dimensional-layer}, the \(\rho\)-length of the set where \(|\Phi_i^\sigma|\le CL\) is bounded by
\[
        O(LM^{1-\alpha}+1).
\]
The factor \(LM^{1-\alpha}\) is exactly \(L/|\partial_\rho\Phi|\).  The additive \(+1\) accounts for the unit-thickening scale and dyadic boundary effects.  The \(\lambda\)-interval has length \(O(M)\), and the density contributes \(O(M)\).  Thus the thickened phase-layer volume is
\[
        O\bigl(M\cdot M\cdot(LM^{1-\alpha}+1)\bigr)
        =O(M^{3-\alpha}L+M^2).
\]
This continuous estimate has no arithmetic loss.  The factor \(M^\eps\) in \eqref{B:balanced-layer} is retained only as the uniform dyadic reserve used later when finitely overlapping cutoffs and modulation shells are summed.  This proves \eqref{B:balanced-layer}.
\end{proof}

\begin{lemma}[Effective high-low two-center volume]\label{lem:B-highlow-geometry}
Fix $n$ with $|n|\sim M$ and write $q=n-k$, $p=|q|$, and $h=|n-q|=|k|$.  On the region
\[
        p\sim R\ll M,\qquad h\sim M,
\]
one has
\[
        |h-|n||\le |q|\lesssim R,
        \qquad
        \frac{hp}{|n|}\sim R,
\]
so the effective $h$-interval has length $O(R)$ and the two-center density is $O(R)$.  If, for fixed $h$, the phase depends on the low length through a branch $\omega_a(p)$, then
\[
        |\partial_p\omega_a(p)|\sim R^{\alpha-1},
\]
and a modulation layer $|\Phi|\le L$ has $p$-thickness $O(LR^{1-\alpha}+1)$.  Consequently the thickened high-low volume is
\begin{equation}\label{B:highlow-effective-volume}
        O\bigl(R\cdot R\cdot(LR^{1-\alpha}+1)\bigr)
        =O(R^{3-\alpha}L+R^2).
\end{equation}
\end{lemma}

\begin{proof}
The first inequality is the triangle inequality.  Since $h\sim |n|\sim M$ and $p\sim R$, the two-center density $hp/|n|$ is comparable to $R$.  The derivative bound follows from the high-frequency formula for $\omega_a'(p)$ on $p\sim R$.  Lemma~\ref{lem:B-one-dimensional-layer} gives the $p$-thickness, and multiplying by the effective $h$-length and by the two-center density gives \eqref{B:highlow-effective-volume}.  This is the precise place where the count pays only low-scale geometry; the high-shell area has already been localized by the two-center constraints.
\end{proof}

\begin{lemma}[High-low phase layers]\label{lem:B-highlow-layer}
Let \(M\gg1\) and \(1\le R\ll M\) be dyadic.  Fix \(n\) with \(|n|\sim M\).  For every sign branch \(\sigma\) and every dyadic \(1\le L\lesssim R^\alpha\),
\begin{equation}\label{B:highlow-layer}
\begin{aligned}
&\#\{k\in\Z^3: |k|\sim M,
      \ |n-k|\sim R,
      \ |\Phi_i^\sigma(n,k)|\le L\} \\
&\qquad\lesssim_\eps R^\eps\bigl(R^2+R^{3-\alpha}L\bigr).
\end{aligned}
\end{equation}
The same estimate holds with the roles of \(k\) and \(n-k\) exchanged.
\end{lemma}

\begin{proof}
Write
\[
        q=n-k,
        \qquad
        p=|q|=|n-k|,
        \qquad
        h=|n-q|=|k|.
\]
In two-center coordinates with centers \(0\) and \(n\), the phase is
\[
        \Phi_i^\sigma
        =
        \sigma_0\omega_i(|n|)+\sigma_1\omega_1(h)+\sigma_2\omega_2(p).
\]
If \(R\) is bounded, the estimate follows from the trivial bound \(O(R^3)\).  We therefore assume \(R\gg1\).  The output vector \(n\) is fixed throughout the count.  In the two-center variables the free variables are \((p,h,\varphi)\), with volume element
\[
        dq=\frac{hp}{|n|}\,dp\,dh\,d\varphi.
\]
The phase and all dyadic constraints are independent of the rotation angle \(\varphi\), so the \(\varphi\)-integration gives the harmless factor \(2\pi\).

We slice the remaining \((p,h)\)-domain by fixing the high length \(h\).  This leaves \(q\) free: for fixed \(n\) and \(h=|n-q|\), the point \(q\) lies on the sphere \(\partial B(n,h)\); the variables \(p=|q|\) and \(\varphi\) remain free.  On such a slice, only the low-length term \(\omega_2(p)\) varies in the phase.  Hence
\[
        \partial_p\Phi_i^\sigma=\sigma_2\omega_2'(p),
        \qquad
        |\partial_p\Phi_i^\sigma|\sim R^{\alpha-1}
        \quad (p\sim R).
\]
By Lemma~\ref{lem:B-one-dimensional-layer}, the \(p\)-length of the layer \(|\Phi_i^\sigma|\le CL\) is
\[
        O(LR^{1-\alpha}+1).
\]
This is the only place in the high-low count where the dispersion relation enters: the derivative lower bound converts modulation width \(L\) into physical thickness \(LR^{1-\alpha}\).

Since \(p=|q|\lesssim R\), the triangle inequality gives
\[
        |h-|n||=\bigl||n-q|-|n|\bigr|\le |q|\lesssim R,
\]
so the effective \(h\)-interval has length \(O(R)\), rather than \(O(M)\).  At this point the high shell is counted through the low cap geometry rather than as a full spherical shell of area \(M^2\).  We are counting the low vector \(q=n-k\), and the condition \(|q|\sim R\) restricts the high vector \(k\) to a cap of angular aperture \(O(R/M)\) around the fixed direction \(n\).  In two-center coordinates that cap restriction is exactly encoded by the short \(h\)-interval and the density below; replacing it by a full \(M^2R\) shell would count points outside \(|n-k|\sim R\).  The two-center density is
\[
        \frac{hp}{|n|}\lesssim R,
\]
because \(h\sim |n|\sim M\) and \(p\sim R\).  Hence the thickened volume is
\[
        O\bigl(R\cdot R\cdot(LR^{1-\alpha}+1)\bigr)
        =O(R^{3-\alpha}L+R^2).
\]
Lemma~\ref{lem:B-thickening} gives an estimate stronger than \eqref{B:highlow-layer}; as in the balanced case, the displayed \(R^\eps\) is only the global reserve for finitely overlapping dyadic and modulation decompositions.  In particular, unit-cube thickening pays neither a full high-shell area nor an additional arithmetic multiplicity.  The exchanged case is identical, with \(\omega_1\) and \(\omega_2\) interchanged.
\end{proof}

\begin{lemma}[Low-output high-high phase-difference gap]\label{lem:B-lowout-gap}
There exist constants \(C_0\gg1\), \(M_0\ge1\), and \(c_*>0\), depending only on the dyadic convention, \(\alpha\), and \(c_1,c_2\), such that if
\[
        |k|\sim |n-k|\sim M,
        \qquad
        M\ge C_0\la n\ra,
        \qquad
        M\ge M_0,
\]
then for every sign branch \(\sigma\in\{\pm1\}^3\),
\begin{equation}\label{B:lowout-gap}
        |\Phi_i^\sigma(n,k)|\ge c_*M^\alpha.
\end{equation}
The constant \(c_*\) degenerates as \(|c_1-c_2|\to0\).
\end{lemma}

\begin{proof}
The output term satisfies
\[
        \omega_i(n)\lesssim\la n\ra^\alpha\le C_0^{-\alpha}M^\alpha.
\]
If the two high-frequency signs agree, then
\[
        |\sigma_1\omega_1(k)+\sigma_2\omega_2(n-k)|
        =
        \omega_1(k)+\omega_2(n-k)
        \gtrsim M^\alpha.
\]
Choosing \(C_0\) large makes the output term too small to cancel this contribution.

It remains to consider the case where the two high-frequency signs disagree.  By the high-frequency expansion
\[
        \omega_a(x)=c_a|x|^\alpha+O(|x|^{-\alpha}),
        \qquad |x|\to\infty,
\]
and by
\[
        \bigl||k|-|n-k|\bigr|\le |n|\le C_0^{-1}M,
\]
we have
\[
        |n-k|^\alpha=|k|^\alpha+O(C_0^{-1}M^\alpha).
\]
Therefore
\[
        \omega_1(k)-\omega_2(n-k)
        =
        (c_1-c_2)|k|^\alpha+O(C_0^{-1}M^\alpha)+O(M^{-\alpha}).
\]
Since \(c_1\ne c_2\) and \(|k|\sim M\), the main term has size \(\gtrsim |c_1-c_2|M^\alpha\).  Choose \(C_0\) so large that the error \(O(C_0^{-1}M^\alpha)\) is at most a small fixed fraction of this main term.  Then choose \(M_0\) so large that the lower-order term \(O(M^{-\alpha})\) is also negligible compared with \(M^\alpha\).  This gives
\[
        |\omega_1(k)-\omega_2(n-k)|\gtrsim M^\alpha.
\]
Increasing \(C_0\) once more makes the output term \(\omega_i(n)\lesssim C_0^{-\alpha}M^\alpha\) dominated by the high-frequency difference.  This proves \eqref{B:lowout-gap}.  The finitely many dyadic shells with \(M<M_0\) are finite-dimensional and are absorbed into the constants in the phase-sum estimates.
\end{proof}

\subsection{The deterministic phase sums}

For \(\nu=0,1\) and \(0\le\eta\le1/2\), define
\begin{equation}\label{B:phase-sum-def}
        S_i^{(\nu,\eta)}(n)
        :=
        \sum_{k\in\Z^3}
        \frac{N_0^{2\nu\alpha}N_{\max}^{2\alpha\eta}}
             {N_0^{2\alpha}N_1^{2\alpha}N_2^{2\alpha}}
        \left(
        \sum_{\sigma\in\{\pm1\}^3}\frakm\bigl(\Phi_i^\sigma(n,k)\bigr)
        \right)^2.
\end{equation}
We write \(S_i^{(\nu)}:=S_i^{(\nu,0)}\).

\begin{proposition}[Three-frequency Picard phase-sum theorem]\label{prop:B-phase-sum}
Choose \(0<\eta\le\eta_0\) and \(\eps>0\) so small that
\begin{equation}\label{B:eta-condition}
        3-6\alpha+2\alpha\eta+\eps<0.
\end{equation}
Then, for \(\nu=0,1\),
\begin{equation}\label{B:phase-sum-main}
        S_i^{(\nu,\eta)}(n)
        \lesssim_{\eps,\eta}
        \la n\ra^{3-(7-2\nu)\alpha+2\alpha\eta+\eps}.
\end{equation}
In particular,
\begin{equation}\label{B:phase-sum-static}
        S_i^{(0)}(n)\lesssim_\eps\la n\ra^{3-7\alpha+\eps},
        \qquad
        S_i^{(1)}(n)\lesssim_\eps\la n\ra^{3-5\alpha+\eps}.
\end{equation}
\end{proposition}

\begin{proof}
Since \((\sum_\sigma a_\sigma)^2\le8\sum_\sigma a_\sigma^2\), it suffices to estimate one fixed sign branch.  We therefore fix \(\sigma\) and estimate
\[
        \sum_{k\in\Z^3}
        \frac{N_0^{2\nu\alpha}N_{\max}^{2\alpha\eta}}
             {N_0^{2\alpha}N_1^{2\alpha}N_2^{2\alpha}}
        \frakm\bigl(\Phi_i^\sigma(n,k)\bigr)^2.
\]
Decompose dyadically according to \(N_1\sim\la k\ra\) and \(N_2\sim\la n-k\ra\), keeping \(N_0=\la n\ra\) fixed.  The relation \(n=k+(n-k)\) implies that the two largest among \(N_0,N_1,N_2\) are comparable.  Finitely many low shells are absorbed into the constants.

We use the dyadic modulation decomposition
\begin{equation}\label{B:modulation-decomp}
        \frakm(\Phi)^2
        \lesssim
        \sum_{\substack{L\ \mathrm{dyadic}\\1\le L\lesssim N_{\max}^{\alpha}}}
        L^{-2}\mathbf 1_{\{|\Phi|\le L\}}.
\end{equation}
Indeed, \(\frakm(\Phi)=1\) when \(|\Phi|\le1\), while for \(|\Phi|>1\) the dyadic sum over \(L\ge|\Phi|\) is comparable to \(|\Phi|^{-2}\).  The upper endpoint follows from \(|\Phi_i^\sigma(n,k)|\lesssim N_{\max}^\alpha\).

\emph{Balanced region.}  Assume
\[
        N_0\sim N_1\sim N_2\sim M.
\]
Then \(N_{\max}\sim M\), and the non-oscillatory prefactor is
\[
        M^{2\nu\alpha+2\alpha\eta-6\alpha}.
\]
Using \eqref{B:modulation-decomp} and Lemma~\ref{lem:B-balanced-layer},
\[
\begin{aligned}
        \mathrm{Bal}_M
        &\lesssim
        M^{2\nu\alpha+2\alpha\eta-6\alpha}
        \sum_{1\le L\lesssim M^\alpha}
        L^{-2}M^\eps(M^2+M^{3-\alpha}L) \\
        &\lesssim_\eps
        M^{2\nu\alpha+2\alpha\eta-6\alpha}
        (M^2+M^{3-\alpha+\eps}).
\end{aligned}
\]
The \(M^{3-\alpha+\eps}\) term gives
\[
        M^{3-(7-2\nu)\alpha+2\alpha\eta+\eps}.
\]
The \(M^2\) term gives
\[
        M^{2+(2\nu-6)\alpha+2\alpha\eta+\eps},
\]
which is no larger because
\[
        [3-(7-2\nu)\alpha]-[2+(2\nu-6)\alpha]=1-\alpha\ge0.
\]
Since \(M\sim N_0\), the balanced contribution satisfies the desired bound.

\emph{High-low region.}  Assume
\[
        N_0\sim N_1\sim M,
        \qquad
        N_2\sim R\ll M.
\]
The case \(N_0\sim N_2\sim M\), \(N_1\sim R\ll M\), is identical.  Here \(N_{\max}\sim M\), and the prefactor is
\[
        M^{2\nu\alpha+2\alpha\eta-4\alpha}R^{-2\alpha}.
\]
For modulations \(1\le L\lesssim R^\alpha\), Lemma~\ref{lem:B-highlow-layer} gives
\[
\begin{aligned}
        \mathrm{HL}_{M,R}^{\le R^\alpha}
        &\lesssim
        M^{2\nu\alpha+2\alpha\eta-4\alpha}R^{-2\alpha}
        \sum_{1\le L\lesssim R^\alpha}
        L^{-2}R^\eps(R^2+R^{3-\alpha}L) \\
        &\lesssim_\eps
        M^{2\nu\alpha+2\alpha\eta-4\alpha}R^{3-3\alpha+\eps}.
\end{aligned}
\]
For larger modulations \(R^\alpha<L\lesssim M^\alpha\), use the trivial count \(\#\{k:|n-k|\sim R\}\lesssim R^3\).  The modulation summation here is dyadic, so
\[
        \sum_{\substack{L\ \mathrm{dyadic}\\ L>R^\alpha}}L^{-2}\lesssim R^{-2\alpha},
\]
a dyadic lattice sum rather than a continuous integral of order \(R^{-\alpha}\).  Hence
\[
\begin{aligned}
        \mathrm{HL}_{M,R}^{>R^\alpha}
        &\lesssim
        M^{2\nu\alpha+2\alpha\eta-4\alpha}R^{-2\alpha}R^3
        \,R^{-2\alpha} \\
        &\lesssim
        M^{2\nu\alpha+2\alpha\eta-4\alpha}R^{3-4\alpha},
\end{aligned}
\]
which is no worse.  Summing dyadically over \(1\le R\le M\),
\[
\begin{aligned}
        \sum_{R\le M}\mathrm{HL}_{M,R}
        &\lesssim_\eps
        M^{2\nu\alpha+2\alpha\eta-4\alpha}
        \sum_{R\le M}R^{3-3\alpha+\eps} \\
        &\lesssim_\eps
        M^{2\nu\alpha+2\alpha\eta-4\alpha}M^{3-3\alpha+\eps} \\
        &=M^{3-(7-2\nu)\alpha+2\alpha\eta+\eps}.
\end{aligned}
\]
At \(\alpha=1\), the endpoint logarithm is absorbed into \(M^\eps\).  Since \(M\sim N_0\), this is the desired high-low contribution.

\emph{Low-output high-high region.}  Assume
\[
        N_1\sim N_2\sim M,
        \qquad
        M\ge C_0N_0.
\]
The finitely many shells with \(M<M_0\) are absorbed into the constants, so we take \(M\ge M_0\).  Then \(N_{\max}\sim M\).  By Lemma~\ref{lem:B-lowout-gap},
\[
        \frakm(\Phi_i^\sigma(n,k))\lesssim M^{-\alpha}.
\]
For a fixed \(M\)-shell, the number of admissible \(k\)'s is \(O(M^3)\).  Hence
\[
\begin{aligned}
        \mathrm{LOHH}_M
        &\lesssim
        N_0^{-2\alpha+2\nu\alpha}M^{-4\alpha}M^{2\alpha\eta}M^3M^{-2\alpha} \\
        &=
        N_0^{-2\alpha+2\nu\alpha}M^{3-6\alpha+2\alpha\eta}.
\end{aligned}
\]
By \eqref{B:eta-condition}, the dyadic sum over \(M\ge C_0N_0\) converges and gives
\[
\begin{aligned}
        \sum_{M\ge C_0N_0}\mathrm{LOHH}_M
        &\lesssim
        N_0^{-2\alpha+2\nu\alpha}N_0^{3-6\alpha+2\alpha\eta+\eps} \\
        &=
        N_0^{3-(8-2\nu)\alpha+2\alpha\eta+\eps}.
\end{aligned}
\]
This is stronger than \eqref{B:phase-sum-main}.  If \(N_1\sim N_2\sim M\) but \(M\lesssim N_0\), then the triangle relation forces \(N_0\lesssim M\), hence \(M\sim N_0\), and the region is already included in the balanced case.  Combining the three regions proves \eqref{B:phase-sum-main}; taking \(\eta=0\) gives \eqref{B:phase-sum-static}.
\end{proof}

\subsection{Fixed-time and time-increment estimates}

\begin{lemma}[Fixed-time coefficient bounds]\label{lem:B-fixed-time}
For \(\nu=0,1\),
\begin{equation}\label{B:fixed-time-coeff}
        \sup_{t\le T}\E\left|\partial_t^\nu\widehat V_i(n,t)\right|^2
        \lesssim_{T,\eps}
        \la n\ra^{3-(7-2\nu)\alpha+\eps}.
\end{equation}
\end{lemma}

\begin{proof}
At finite cutoff, apply the double Wiener isometry \eqref{B:double-isometry}.  Lemma~\ref{lem:B-static-kernel} gives
\[
        \E\left|\partial_t^\nu\widehat V_{i,\Lambda}(n,t)\right|^2
        \lesssim_T S_i^{(\nu)}(n),
\]
uniformly in \(\Lambda\).  Proposition~\ref{prop:B-phase-sum} with \(\eta=0\) gives \eqref{B:fixed-time-coeff}, uniformly in the cutoff.  The cutoff Cauchy convergence is proved below; until then the estimate is understood at finite cutoff with a uniform bound.
\end{proof}

\begin{lemma}[Time-increment coefficient bounds]\label{lem:B-time-increments}
For sufficiently small \(0<\eta\le\eta_0\) and \(\nu=0,1\),
\begin{equation}\label{B:time-incr-coeff}
        \E\left|\partial_t^\nu\widehat V_i(n,t)-\partial_t^\nu\widehat V_i(n,t')\right|^2
        \lesssim_{T,\eps,\eta}
        |t-t'|^{2\eta}
        \la n\ra^{3-(7-2\nu)\alpha+2\alpha\eta+\eps}.
\end{equation}
\end{lemma}

\begin{proof}
Assume \(t'<t\); the other case is symmetric.  In the double stochastic representation, split the square \([0,t]^2\) into the common square \([0,t']^2\) and the boundary strip \([0,t]^2\setminus[0,t']^2\).

On the common square, the kernel difference is controlled by Lemma~\ref{lem:B-kernel-incr}.  The double Wiener isometry gives
\[
        \text{common-square contribution}
        \lesssim
        |t-t'|^{2\eta}S_i^{(\nu,\eta)}(n).
\]
On the boundary strip, whose measure is \(O_T(|t-t'|)\), use the static kernel bound from Lemma~\ref{lem:B-static-kernel}.  Since \(0<\eta\le1/2\) and \(|t-t'|\le T\), this contribution is bounded by
\[
        C_T|t-t'|^{2\eta}S_i^{(\nu,\eta)}(n),
\]
because \(N_{\max}^{2\alpha\eta}\ge1\).  Proposition~\ref{prop:B-phase-sum} proves \eqref{B:time-incr-coeff}.
\end{proof}

Taking \(\nu=0\) and \(\nu=1\) in Lemmas~\ref{lem:B-fixed-time} and \ref{lem:B-time-increments} gives \eqref{B:fixed-V}-\eqref{B:incr-dtV} at the coefficient level.

\subsection{Cutoff tails and cutoff convergence}

Let \(S_{i,>\Lambda}^{(\nu,\eta)}(n)\) be the phase sum \eqref{B:phase-sum-def} restricted to
\[
        \max\{|k|,|n-k|\}>\Lambda.
\]
For smooth cutoff families, the same quantity, up to harmless changes of constants, controls the difference of two cutoffs whose symbols agree on \(|m|\le\Lambda\) and are uniformly bounded.  It is the deterministic majorant for the coefficient variance of cutoff differences.

\begin{lemma}[Tail phase sums]\label{lem:B-tail-sums}
Fix \(\nu=0,1\), choose \(0<\eta\le\eta_0\) satisfying \eqref{B:eta-condition}, and let
\[
        s<\frac{(7-2\nu)\alpha}{2}-3.
\]
Then, after reducing \(\eta\) and the harmless losses if necessary,
\begin{equation}\label{B:tail-sum-target}
        \lim_{\Lambda\to\infty}
        \sum_{N_0\ \mathrm{dyadic}}
        N_0^{2s+3}
        \sup_{\la n\ra\sim N_0}
        S_{i,>\Lambda}^{(\nu,\eta)}(n)=0.
\end{equation}
Moreover, for every strict choice of \(s\) in the above range, the proof yields a polynomial dyadic tail bound along \(\Lambda=2^j\): for some \(\gamma>0\),
\begin{equation}\label{B:tail-polynomial}
        \sum_{N_0}
        N_0^{2s+3}
        \sup_{\la n\ra\sim N_0}
        S_{i,>\Lambda}^{(\nu,\eta)}(n)
        \lesssim \Lambda^{-\gamma}
\end{equation}
for all dyadic \(\Lambda\ge1\), after possibly reducing \(\gamma\).
\end{lemma}

\begin{proof}
For each fixed dyadic output shell \(\la n\ra\sim N_0\), there are finitely many output frequencies \(n\), and for each such \(n\) the restricted positive sum \(S_{i,>\Lambda}^{(\nu,\eta)}(n)\) decreases to zero as \(\Lambda\to\infty\).  Hence the shell supremum tends to zero for fixed \(N_0\).

It remains to dominate the series.  Proposition~\ref{prop:B-phase-sum} gives
\[
        S_{i,>\Lambda}^{(\nu,\eta)}(n)
        \le S_i^{(\nu,\eta)}(n)
        \lesssim
        \la n\ra^{3-(7-2\nu)\alpha+2\alpha\eta+\eps}.
\]
Therefore
\[
        \sum_{N_0}
        N_0^{2s+3}
        N_0^{3-(7-2\nu)\alpha+2\alpha\eta+\eps}<\infty
\]
provided
\[
        2s+6-(7-2\nu)\alpha+2\alpha\eta+\eps<0.
\]
This holds because \(s<\frac{(7-2\nu)\alpha}{2}-3\) and \(\eta,\eps\) are chosen sufficiently small.  Dominated convergence proves \eqref{B:tail-sum-target}.

We record the polynomial tail because it is useful for almost sure convergence along dyadic cutoffs.  Let
\[
        \delta_s:=(7-2\nu)\alpha-6-2s>0.
\]
Choose \(\eta,\eps>0\) so that \(2\alpha\eta+\eps<\delta_s/4\).  In the balanced and high-low regions the tail condition forces the large dyadic scale \(M\sim N_0\) to satisfy \(M\gtrsim\Lambda\), and the estimates in Proposition~\ref{prop:B-phase-sum} give a summable dyadic tail bounded by \(\Lambda^{-\gamma}\) for some \(\gamma>0\).

In the low-output high-high region, the weighted shell contribution is bounded by
\[
        N_0^{2s+3}N_0^{-2\alpha+2\nu\alpha}
        M^{3-6\alpha+2\alpha\eta+\eps},
        \qquad M\ge C_0N_0.
\]
If \(N_0\gtrsim\Lambda\), summing first in \(M\ge C_0N_0\) and then in \(N_0\) gives a tail controlled by the full summable majorant.  If \(N_0\lesssim\Lambda\), the additional restriction \(M\gtrsim\Lambda\) gives
\[
        \sum_{N_0\lesssim\Lambda}
        N_0^{2s+3-2\alpha+2\nu\alpha}
        \sum_{M\gtrsim\Lambda}M^{3-6\alpha+2\alpha\eta+\eps}.
\]
The exponent of the \(M\)-sum is negative by \eqref{B:eta-condition}; after summing in \(N_0\lesssim\Lambda\), the total power of \(\Lambda\) is strictly negative because \(s\) is strictly below the threshold and \(\eta,\eps\) are small.  This proves \eqref{B:tail-polynomial}.
\end{proof}

\begin{proof}[Proof of Theorem~\ref{thm:B-first-picard}]
The coefficient bounds \eqref{B:fixed-V}-\eqref{B:incr-dtV} were proved in Lemmas~\ref{lem:B-fixed-time} and \ref{lem:B-time-increments}.  Applying Lemma~\ref{lem:B-coeff-besov} with
\[
        p_0=3-7\alpha+\eps
\]
gives
\[
        V_i\in C([0,T];\C^s)
        \qquad\text{for every }s<-\frac{3}{2}-\frac{p_0}{2}
        =\frac{7\alpha}{2}-3-.
\]
Applying the same lemma with
\[
        p_1=3-5\alpha+\eps
\]
gives
\[
        \partial_tV_i\in C([0,T];\C^s)
        \qquad\text{for every }s<-\frac{3}{2}-\frac{p_1}{2}
        =\frac{5\alpha}{2}-3-.
\]
The cutoff-difference coefficient and increment estimates are controlled by the tail phase sums in Lemma~\ref{lem:B-tail-sums}.  The polynomial tail bound \eqref{B:tail-polynomial}, hypercontractivity, and Lemma~\ref{lem:B-coeff-besov} imply almost sure Cauchy convergence along the standard dyadic cutoffs in the stated path spaces.  For a general admissible cofinal cutoff family, one instead uses a finite-frequency-set plus uniform-tail argument.  The tail phase sum makes the complement uniformly small, while on the finite set all multiplier products converge uniformly because only finitely many lattice points occur.  Applying the same decomposition to the increment kernels proves the Cauchy property and shows that the resulting limit agrees with the standard dyadic limit.  This is the precise cutoff-convention independence used in the main theorem.

At finite cutoff,
\[
        V_{i,\Lambda}(t)-V_{i,\Lambda}(0)
        =
        \int_0^t\partial_sV_{i,\Lambda}(s)\,ds.
\]
Both sides converge uniformly in time in slightly lower spatial regularity.  Therefore the limiting process has the asserted time derivative.  Moreover, for every strict choice of the losses used in the main text,
\[
        V_{i,\Lambda}\to V_i\quad\text{in }C_T\C^{\rho_V}
        \cap L_T^\infty B_{2,\infty}^{\rho_V},
        \qquad
        \partial_tV_{i,\Lambda}\to\partial_tV_i
        \quad\text{in }C_T\C^{\rho_V-\alpha}
        \cap L_T^\infty B_{2,\infty}^{\rho_V-\alpha}.
\]
These are precisely the two first-Picard topologies recorded in the enhanced-data distance: the H\"older topology is used for deterministic products, while the dyadic $B_{2,\infty}$ envelope is used for Besov-Besov input bounds.  This proves the theorem.
\end{proof}

\section{The resonant cubic symbols: colored Wick algebra and phase estimates}\label{app:cubic-resonance}

We prove the cubic result for
\[
        \Gamma_1=V_2\circ\Psi_1=I_2(\Psi_1\Psi_2)\circ\Psi_1;
\]
the second component follows by exchanging colors and speeds.  Independence leaves only the contraction of the two color-1 factors and hence a color-2 first-chaos term.  Its formal loop multiplier is
\begin{equation}\label{C:formal-M21-intro}
        M_{21}(n;t,s)
        =\sum_{r\in\Z^3}
        \chi^{\mathrm{res}}(n-r,r)K_2(t-s,n-r)\sigma_1(r;s,t)
\end{equation}
On the resonant support \(|r|\sim |n-r|\sim M\), the leading shell has size
\[
        M^3\cdot M^{-\alpha}\cdot M^{-2\alpha}=M^{3-3\alpha},
\]
and is not absolutely summable for \(\alpha\le1\).  We therefore define the contraction only after inserting the remaining stochastic convolution.  The resulting integrated Volterra kernel satisfies
\begin{equation}\label{C:K21-main-intro}
        \sup_{0\le u\le t\le T}|K_{21,M}(n;t,u)|
        \lesssim_{T,\eps} \la n\ra^{-\alpha}M^{3-4\alpha+\eps},
\end{equation}
which is summable when \(\alpha>3/4\).

Throughout this appendix,
\[
        \frac34<\alpha\le1,
        \qquad c_1,c_2>0,
        \qquad c_1\ne c_2,
\]
and
\[
        \omega_i(n)=(1+c_i^2|n|^{2\alpha})^{1/2},
        \qquad
        K_i(t,n)=\frac{\sin(t\omega_i(n))}{\omega_i(n)}.
\]
We use \(\la n\ra=(1+|n|^2)^{1/2}\).  The resonant multiplier \(\chi^{\mathrm{res}}(m,r)\) is a fixed smooth Fourier cutoff supported where \(\la m\ra\sim\la r\ra\).  Finitely many zero or low modes are finite-dimensional and are absorbed into constants.  Constants may depend on \(T\), \(\alpha\), the speeds, and the fixed cutoff conventions, uniformly over dyadic frequencies and over the Fourier cutoff parameter.

The input from Appendix~\ref{app:first-picard-phase} is the first Picard smoothing estimate and its time increment version: for every \(\eps>0\), every sufficiently small \(\eta>0\), and \(j=1,2\),
\begin{align}
        \sup_{t\le T}\E |\wh V_j(m,t)|^2
        &\lesssim_{T,\eps} \la m\ra^{3-7\alpha+\eps},\label{C:Picard-input}\\
        \E |\wh V_j(m,t+h)-\wh V_j(m,t)|^2
        &\lesssim_{T,\eps,\eta}
        |h|^{2\eta}\la m\ra^{3-7\alpha+2\alpha\eta+\eps}.\label{C:Picard-increment-input}
\end{align}
The linear stochastic convolution satisfies
\begin{align}
        \sup_{t\le T}\E |\wh\Psi_j(r,t)|^2
        &\lesssim_T \la r\ra^{-2\alpha},\label{C:linear-input}\\
        \E |\wh\Psi_j(r,t+h)-\wh\Psi_j(r,t)|^2
        &\lesssim_{T,\eta}|h|^{2\eta}\la r\ra^{-2\alpha+2\alpha\eta}.
        \label{C:linear-increment-input}
\end{align}
All random variables considered below belong to a fixed finite sum of homogeneous Wiener chaoses, uniformly in the cutoff.

\subsection{A coefficient-to-path-space criterion}

We shall use the following dyadic criterion several times.  It is stated here in the form needed in this appendix.

\begin{lemma}[Coefficient-to-path-space criterion]\label{lem:C-coeff-path}
Let \(Z_\Lambda\) be spatially homogeneous random distributions on \(\T^3\), uniformly belonging to a finite sum of homogeneous Wiener chaoses of order at most \(m\).  Suppose that for some \(p\in\R\), some \(0<\eta\le1/2\), and every small \(\eps>0\),
\begin{align}
        \sup_\Lambda\sup_{t\le T}\E|\wh Z_\Lambda(n,t)|^2
        &\lesssim_{T,\eps}\la n\ra^{p+\eps},\label{C:coeff-path-fixed}\\
        \sup_\Lambda\E|\wh Z_\Lambda(n,t)-\wh Z_\Lambda(n,t')|^2
        &\lesssim_{T,\eps,\eta}|t-t'|^{2\eta}\la n\ra^{p+2\alpha\eta+\eps}.
        \label{C:coeff-path-incr}
\end{align}
Then every $Z_\Lambda$ has a version in \(C([0,T];\C^s(\T^3))\) for every
\[
        s<-\frac32-\frac p2.
\]
For every finite $q$, these versions are uniformly bounded in
$L^q(\Omega;C([0,T];\mathcal C^s))$.
Suppose, in addition, that the countable cutoff family is indexed by $\Lambda\in\mathbb N$ and that for some $\theta>0$, uniformly for $\Lambda'\ge\Lambda\ge L$,
\[
        \sum_N N^{2s+3}
        \sup_{\la n\ra\sim N}\sup_{t\le T}
        \E|\widehat{(Z_{\Lambda'}-Z_\Lambda)}(n,t)|^2
        \lesssim L^{-2\theta},
\]
with the analogous $L^{-2\theta}$ bound for time increments after the factor $|t-t'|^{2\eta}$ is removed.  Then \((Z_\Lambda)_\Lambda\) is almost surely Cauchy in \(C([0,T];\C^s)\).
\end{lemma}

\begin{proof}
Fix a dyadic shell \(\la n\ra\sim N\).  For each fixed \((t,x)\), spatial homogeneity gives diagonal covariance on the block and therefore
\[
        \E|P_NZ_\Lambda(t,x)|^2
        \lesssim
        \sum_{\la n\ra\sim N}\E|\wh Z_\Lambda(n,t)|^2
        \lesssim N^{3+p+\eps}.
\]
Nelson hypercontractivity on a fixed finite chaos upgrades this to
\[
        \|P_NZ_\Lambda(t,x)\|_{L^q(\Omega)}
        \lesssim_{q,m,\eps}N^{3/2+p/2+\eps}
\]
for every finite \(q\ge2\).  Taking a spatial net of mesh \(N^{-A}\), using Bernstein's inequality to pass from the net to \(L^\infty_x\), and then taking \(q\) large, gives
\[
        \|P_NZ_\Lambda(t)\|_{L^q(\Omega;L^\infty_x)}
        \lesssim_{q,m,\kappa} N^{3/2+p/2+\kappa}
\]
for every \(\kappa>0\).  Applying the same argument to time increments and using \eqref{C:coeff-path-incr}, followed by Kolmogorov's criterion or fractional Sobolev embedding in time, yields
\[
        \|P_NZ_\Lambda\|_{L^q(\Omega;C([0,T];L^\infty_x))}
        \lesssim N^{3/2+p/2+\kappa},
\]
after taking \(\eta\) smaller if needed.  Summing the strict dyadic exponent gives the stated uniform $L^q$ bound; Markov's inequality and Borel--Cantelli over dyadic \(N\) give an almost sure path for each cutoff.  The Besov summability condition is
\[
        \sum_N N^sN^{3/2+p/2+\kappa}<\infty,
\]
which holds when \(s<-3/2-p/2\) and \(\kappa>0\) is sufficiently small.  For cutoff differences, hypercontractivity and the same space--time net argument upgrade the weighted second-moment tail to
\[
 \|Z_{\Lambda'}-Z_\Lambda\|_{L^q(\Omega;C_T\mathcal C^s)}
 \lesssim_q L^{-\theta+\kappa},
 \qquad \Lambda'\ge\Lambda\ge L,
\]
for every finite $q$ and every sufficiently small $\kappa>0$.  On the dyadic slab $2^m\le\Lambda<2^{m+1}$, compare every cutoff to $2^m$ and also compare the two consecutive dyadic anchors.  There are $O(2^m)$ comparisons.  Taking $q$ large and applying Chebyshev and Borel--Cantelli shows that all these differences are eventually bounded by $2^{-m\theta/2}$.  This is summable in $m$, so the entire countable cutoff family is almost surely Cauchy.
\end{proof}

\begin{theorem}[Resonant cubic symbols]\label{thm:C-cubic}
For each \(i=1,2\), the finite Fourier cutoff approximations of
\[
        \Gamma_i=V_{3-i}\circ\Psi_i
\]
converge in the stated topology and admit the colored Wick decomposition
\[
        \Gamma_i=\Gamma_i^{(3)}+C_i,
\]
where \(\Gamma_i^{(3)}\) is a centered third homogeneous chaos and \(C_i\) is a centered first-chaos Volterra term.  Moreover,
\[
        \Gamma_i^{(3)}\in C([0,T];\C^{\frac{9\alpha}{2}-\frac92-}),
        \qquad
        C_i\in C([0,T];\C^{5\alpha-\frac92-}).
\]
Consequently,
\[
        \Gamma_i\in C([0,T];\C^{\frac{9\alpha}{2}-\frac92-}).
\]
The first-chaos contraction is smoother than the centered third chaos by \(\alpha/2-\).
\end{theorem}

\begin{proposition}[Integrated first-chaos contraction]\label{prop:C-integrated-first-chaos-contraction}
For the first-chaos contraction in $V_2\circ\Psi_1$, the finite-cutoff lower-chaos term is defined only after the remaining color-$2$ stochastic convolution is inserted:
\[
        C_{1,\Lambda}(n,t)=\int_0^t K_{21,\Lambda}(n;t,u)\,d\beta_n^2(u).
\]
The limit is taken in the integrated-kernel topology rather than at the level of the formal multiplier \eqref{C:formal-M21-intro}.  The integrated kernels converge in the deterministic kernel topology used below and satisfy the dyadic estimate
\[
        \sup_{0\le u\le t\le T}|K_{21,M}(n;t,u)|
        \lesssim_{T,\eps}\langle n\rangle^{-\alpha}M^{3-4\alpha+\eps},
        \qquad M\gtrsim \langle n\rangle,
\]
together with the corresponding time-increment and cutoff-tail estimates.  Consequently
\[
        \sup_{t\le T}\E|\widehat C_1(n,t)|^2
        \lesssim_{T,\eps}\langle n\rangle^{6-10\alpha+\eps}.
\]
After exchanging colors, the same statement holds for $C_2$.
\end{proposition}

\begin{proof}
The finite-cutoff colored Wick algebra is Lemma~\ref{lem:C-colored-product}.  The integrated Volterra kernel and the single-shell estimate are Lemma~\ref{lem:C-single-shell}; convergence of the kernel is Lemma~\ref{lem:C-K21-conv}; the coefficient bound is Proposition~\ref{prop:C-first-chaos-fixed}.  The time increment and cutoff-tail estimates are Lemma~\ref{lem:C-first-incr} and Lemma~\ref{lem:C-first-tail}.  These statements are precisely the displayed conclusion, and the color-exchanged block is identical.
\end{proof}

\subsection{Finite cutoff Wick decomposition in the colored Gaussian space}

Let \(\chi_\Lambda\) be a smooth radial Fourier cutoff, bounded by one, equal to one on \(|n|\le \Lambda\), supported on \(|n|\le2\Lambda\), and converging pointwise to one as \(\Lambda\to\infty\).  Define
\[
        \wh\Psi_{b,\Lambda}(m,t)
        =\chi_\Lambda(m)\int_0^tK_b(t-s,m)\,d\beta_m^b(s),
        \qquad b=1,2.
\]
Set
\[
        V_{2,\Lambda}=I_2(\Psi_{1,\Lambda}\Psi_{2,\Lambda})
\]
and
\begin{equation}\label{C:Gamma-cutoff-def}
        \wh\Gamma_{1,\Lambda}(n,t)
        =\sum_{m+r=n}\chi^{\mathrm{res}}(m,r)
        \wh V_{2,\Lambda}(m,t)\wh\Psi_{1,\Lambda}(r,t).
\end{equation}
Expanding \(V_{2,\Lambda}\),
\begin{equation}\label{C:Gamma-expanded}
\begin{aligned}
        \wh\Gamma_{1,\Lambda}(n,t)
        & =
        \sum_{a+b+r=n}
        \mathbf c_\Lambda(a,b,r)
        \chi^{\mathrm{res}}(a+b,r)                                      \\
        &\quad\times
        \int_0^t K_2(t-s,a+b)
        \wh\Psi_1(a,s)\wh\Psi_2(b,s)\wh\Psi_1(r,t)\,ds,
\end{aligned}
\end{equation}
where
\[
        \mathbf c_\Lambda(a,b,r)=\chi_\Lambda(a)\chi_\Lambda(b)\chi_\Lambda(r).
\]
This exact finite-cutoff weight is kept in the algebraic identities.  In estimates it is harmless since \(|\mathbf c_\Lambda|\le1\) and it converges pointwise to one.

\medskip
\noindent\emph{Colored Gaussian Hilbert space.}
Let \(\hh_T=L^2([0,T]\times\T^3;\R)\), and let \(W_1,W_2\) be independent isonormal Gaussian processes over \(\hh_T\).  Equivalently, the full first-chaos Hilbert space is the orthogonal direct sum
\[
        \HH_T=\hh_T^{(1)}\oplus\hh_T^{(2)},
        \qquad
        \hh_T^{(1)}\perp\hh_T^{(2)}.
\]
At finite cutoff, \(\HH_T^\Lambda=\HH_{1,T}^\Lambda\oplus\HH_{2,T}^\Lambda\) denotes the finite-dimensional subspace generated by the Fourier Brownian coordinates up to the cutoff.

We use standard multiple Wiener integrals with the convention
\[
        I_1(h)I_1(g)=I_2(h\widetilde\otimes g)+\langle h,g\rangle
\]
within one color.  For colored integrals, \(I_{(p,q)}\) denotes order \(p\) in color 1 and order \(q\) in color 2, with symmetrization separately inside variables of the same color.

\begin{lemma}[The unique colored contraction channel]\label{lem:C-colored-product}
Let \(F\in \HH_{1,T}^\Lambda\otimes\HH_{2,T}^\Lambda\) and \(g\in\HH_{1,T}^\Lambda\).  Then
\begin{equation}\label{C:colored-product}
        I_{(1,1)}(F)I_{(1,0)}(g)
        =I_{(2,1)}\bigl(\Sym_1(F\otimes g)\bigr)
        +I_{(0,1)}(F\otimes_1 g).
\end{equation}
Here \(\Sym_1\) symmetrizes the two color-1 slots, and \(F\otimes_1 g\in\HH_{2,T}^\Lambda\) is the contraction of the color-1 slot of \(F\) with \(g\).  The color-2 slot of \(F\) remains unpaired with \(g\).
\end{lemma}

\begin{proof}
It suffices to check finite-rank tensors.  Write
\[
        F=\sum_{a,b}F_{a,b}e_a^{(1)}\otimes e_b^{(2)},
        \qquad
        g=\sum_r g_r e_r^{(1)},
\]
with orthonormal bases in the two color Hilbert spaces.  For the color-1 variables,
\[
        W_1(e_a^{(1)})W_1(e_r^{(1)})
        =I_2(e_a^{(1)}\widetilde\otimes e_r^{(1)})
        +\langle e_a^{(1)},e_r^{(1)}\rangle_{\HH_{1,T}}.
\]
Multiplication by the independent coordinate \(W_2(e_b^{(2)})\) gives the two terms in \eqref{C:colored-product}.  Since \(\HH_{1,T}\perp\HH_{2,T}\), color 2 and color 1 have zero pairing.  Completion gives the statement for all finite-cutoff kernels.
\end{proof}

For stochastic convolutions the same-color covariance is
\begin{equation}\label{C:cov-def}
        \E[\wh\Psi_a(k,s)\wh\Psi_b(k',t)]
        =\mathbf 1_{\{a=b\}}\mathbf 1_{\{k+k'=0\}}\sigma_a(k;s,t),
\end{equation}
where
\begin{equation}\label{C:sigma-def}
        \sigma_a(k;s,t)=\E[\wh\Psi_a(k,s)\wh\Psi_a(-k,t)].
\end{equation}
Consequently, the product in \eqref{C:Gamma-expanded} has only one surviving contraction:
\[
        \wh\Psi_1(a,s)\text{ contracts with }\wh\Psi_1(r,t).
\]
This contraction imposes \(a+r=0\).  Since the output relation is \(a+b+r=n\), the unpaired color-2 frequency is forced to be \(b=n\).  Therefore
\begin{equation}\label{C:Gamma-decomposition}
        \Gamma_{1,\Lambda}=\Gamma_{1,\Lambda}^{(3)}+C_{1,\Lambda}.
\end{equation}
The centered third-chaos part is
\begin{equation}\label{C:third-chaos-def}
\begin{aligned}
        \wh\Gamma_{1,\Lambda}^{(3)}(n,t)
        & =
        \sum_{a+b+r=n}
        \mathbf c_\Lambda(a,b,r)
        \chi^{\mathrm{res}}(a+b,r)                                      \\
        &\quad\times
        \int_0^tK_2(t-s,a+b)
        :\wh\Psi_1(a,s)\wh\Psi_2(b,s)\wh\Psi_1(r,t):\,ds.
\end{aligned}
\end{equation}
The contraction is
\begin{equation}\label{C:contraction-cutoff-form}
        \wh C_{1,\Lambda}(n,t)
        =\chi_\Lambda(n)
        \int_0^t M_{21,\Lambda}(n;t,s)\wh\Psi_2(n,s)\,ds,
\end{equation}
where
\begin{equation}\label{C:M21-cutoff}
        M_{21,\Lambda}(n;t,s)
        =\sum_{r\in\Z^3}
        \chi_\Lambda(r)^2
        \chi^{\mathrm{res}}(n-r,r)
        K_2(t-s,n-r)\sigma_1(r;s,t).
\end{equation}
For a different smooth cutoff convention, the factor \(\chi_\Lambda(n)\chi_\Lambda(r)^2\) is replaced by a bounded multiplier \(\mathbf c_\Lambda(n,r)\), pointwise converging to one and dominated by a fixed annular cutoff.  All estimates below are uniform for such cutoff families.

\subsection{The centered third homogeneous chaos}

The third-chaos bound is a Hilbert-kernel norm estimate inside the colored Gaussian Hilbert space, since \(V_2\) contains a color-1 factor.  The projection is taken in that colored Gaussian Hilbert space.  The only deterministic issue created by symmetrizing the two color-1 slots is a finite-overlap change of the fixed resonant cutoff; this is recorded next.

\begin{lemma}[Finite overlap under color-one symmetrization]\label{lem:C-finite-overlap-sym}
Fix the output relation $a+b+r=n$ and set
\[
        H(a,r,b)=\chi^{\mathrm{res}}(a+b,r)F(a,b)g(r),
        \qquad
        \widetilde H(a,r,b)=\frac12\bigl(H(a,r,b)+H(r,a,b)\bigr).
\]
Then
\begin{equation}\label{C:finite-overlap-sym}
        \sum_{a+b+r=n}|\widetilde H(a,r,b)|^2
        \lesssim
        \sum_{\nu\in\mathcal P}
        \sum_{a+b+r=n}|\chi_\nu^{\mathrm{res}}(a+b,r)|^2|F(a,b)|^2|g(r)|^2,
\end{equation}
where $\mathcal P$ is a finite set depending only on the fixed Littlewood--Paley partition, and each $\chi_\nu^{\mathrm{res}}$ is a resonant cutoff with dyadic scales comparable to those of $\chi^{\mathrm{res}}$.
\end{lemma}

\begin{proof}
The pointwise inequality
\[
        |\widetilde H(a,r,b)|^2
        \le \frac12|H(a,r,b)|^2+\frac12|H(r,a,b)|^2
\]
reduces the issue to the exchanged term.  After the exchange of the two color-one slots, the cutoff $\chi^{\mathrm{res}}(a+b,r)$ becomes a cutoff of the form $\chi^{\mathrm{res}}(r+b,a)$.  For a fixed smooth Bony partition, the latter support is covered by finitely many resonant supports with comparable dyadic scales; the overlap number is independent of the frequencies and of the cutoff level.  Summing the two contributions gives \eqref{C:finite-overlap-sym}.
\end{proof}

\begin{lemma}[Third-chaos kernel norm]\label{lem:C-third-kernel-norm}
For every finite cutoff and every fixed \((n,t)\),
\begin{equation}\label{C:third-kernel-norm}
        \E|\wh\Gamma_{1,\Lambda}^{(3)}(n,t)|^2
        \lesssim
        \sum_{m+r=n}|\chi^{\mathrm{res}}(m,r)|^2
        \E|\wh V_{2,\Lambda}(m,t)|^2
        \E|\wh\Psi_{1,\Lambda}(r,t)|^2.
\end{equation}
The same estimate holds for cutoff differences, with one factor replaced by the corresponding cutoff difference and the other factors estimated by their uniform majorants.
\end{lemma}

\begin{proof}
For each \(m\), write \(\wh V_{2,\Lambda}(m,t)=I_{(1,1)}(F_m^\Lambda(t))\).  Then
\[
        \|F_m^\Lambda(t)\|_{\HH_1\otimes\HH_2}^2
        \simeq \E|\wh V_{2,\Lambda}(m,t)|^2,
\]
with a convention-dependent constant.  Similarly,
\[
        \wh\Psi_{1,\Lambda}(r,t)=I_{(1,0)}(g_r^\Lambda(t)),
        \qquad
        \|g_r^\Lambda(t)\|_{\HH_1}^2
        \simeq \E|\wh\Psi_{1,\Lambda}(r,t)|^2.
\]
By Lemma~\ref{lem:C-colored-product}, the third-chaos projection of the product is
\[
        I_{(2,1)}\bigl(\Sym_1(F_m^\Lambda(t)\otimes g_r^\Lambda(t))\bigr).
\]
The output \(n\) kernel is the sum over \(m+r=n\), multiplied by \(\chi^{\mathrm{res}}(m,r)\).  Before symmetrizing the two color-1 slots, the Fourier kernel has the schematic form
\[
        H(a,r,b)=\chi^{\mathrm{res}}(a+b,r)F_{a,b}^\Lambda(t)g_r^\Lambda(t),
        \qquad a+b+r=n.
\]
The color-one symmetrization is controlled by Lemma~\ref{lem:C-finite-overlap-sym}.  Since the finite set of comparable resonant cutoffs only changes the fixed Littlewood--Paley constants, summing over all Fourier and time variables gives \eqref{C:third-kernel-norm}.  Cutoff differences are handled by expanding the difference of the finite-cutoff kernels and applying the same finite-overlap bound term by term.
\end{proof}

\begin{lemma}[Third-chaos coefficient and time increment bounds]\label{lem:C-third-coeff}
For every \(\eps>0\),
\begin{equation}\label{C:third-fixed-bound}
        \sup_{t\le T}\E|\wh\Gamma_1^{(3)}(n,t)|^2
        \lesssim_{T,\eps}\la n\ra^{6-9\alpha+\eps}.
\end{equation}
Moreover, for sufficiently small \(0<\eta<1/2\),
\begin{equation}\label{C:third-increment-bound}
        \E|\wh\Gamma_1^{(3)}(n,t+h)-\wh\Gamma_1^{(3)}(n,t)|^2
        \lesssim_{T,\eps,\eta}
        |h|^{2\eta}\la n\ra^{6-9\alpha+2\alpha\eta+\eps}.
\end{equation}
\end{lemma}

\begin{proof}
Using Lemma~\ref{lem:C-third-kernel-norm}, \eqref{C:Picard-input}, and \eqref{C:linear-input},
\[
        \E|\wh\Gamma_1^{(3)}(n,t)|^2
        \lesssim_{T,\eps}
        \sum_{m+r=n}|\chi^{\mathrm{res}}(m,r)|^2
        \la m\ra^{3-7\alpha+\eps}\la r\ra^{-2\alpha}.
\]
On the resonant support \(\la m\ra\sim\la r\ra\sim M\), and necessarily \(M\gtrsim\la n\ra\).  The dyadic shell contribution is
\[
        M^3\cdot M^{3-7\alpha+\eps}\cdot M^{-2\alpha}
        =M^{6-9\alpha+\eps}.
\]
Since \(6-9\alpha<0\) for \(\alpha>2/3\), and in particular for \(\alpha>3/4\), summing over dyadic \(M\gtrsim\la n\ra\) gives \eqref{C:third-fixed-bound}.

For time increments, write
\[
        P_3\bigl(V_2(t+h)\Psi_1(t+h)-V_2(t)\Psi_1(t)\bigr)
\]
as the sum of the two third-chaos projections obtained from
\[
        (V_2(t+h)-V_2(t))\Psi_1(t+h)
        +V_2(t)(\Psi_1(t+h)-\Psi_1(t)).
\]
The same kernel norm estimate as in Lemma~\ref{lem:C-third-kernel-norm}, together with \eqref{C:Picard-increment-input}, \eqref{C:linear-input}, \eqref{C:Picard-input}, and \eqref{C:linear-increment-input}, gives on the shell \(M\)
\[
        |h|^{2\eta}
        M^3\Bigl(M^{3-7\alpha+2\alpha\eta+\eps}M^{-2\alpha}
        +M^{3-7\alpha+\eps}M^{-2\alpha+2\alpha\eta}\Bigr)
        \lesssim
        |h|^{2\eta}M^{6-9\alpha+2\alpha\eta+\eps}.
\]
Taking \(\eta,\eps\) so small that the exponent remains negative and summing over \(M\gtrsim\la n\ra\) proves \eqref{C:third-increment-bound}.
\end{proof}

\begin{lemma}[Third-chaos cutoff Cauchy convergence]\label{lem:C-third-tail}
For every
\[
        s<\frac{9\alpha}{2}-\frac92,
\]
the standard localized dyadic approximations \(\Gamma_{1,\Lambda}^{(3)}\) are Cauchy almost surely in \(C([0,T];\C^s)\), with a polynomial cutoff-difference bound at the coefficient level.  Every admissible cofinal cutoff family converges in the same path space to the same limit, without an asserted rate.
\end{lemma}

\begin{proof}
It is enough to prove the assertion for \(s\) sufficiently close to \(9\alpha/2-9/2\) that \(2s+3>0\): convergence at such an exponent implies convergence at every smaller exponent by the Besov embedding.  We work at one such \(s\).

We first use standard localized cutoffs and compare \(\Lambda'\ge\Lambda=L\).  Expanding the difference of the third-chaos kernels gives a finite sum of terms in which either the external color-1 factor \(\Psi_1(r)\) is in the cutoff tail, or the first Picard kernel \(V_2(m)\) is in the cutoff tail.  Fix
\[
        \delta_\Gamma:=\frac{9\alpha}{2}-\frac92-s>0,
        \qquad
        s_V:=s+\frac32-\alpha
             =\frac{7\alpha}{2}-3-\delta_\Gamma.
\]

For the external tail, localization forces the resonant scale to satisfy \(M\gtrsim L\).  Lemma~\ref{lem:C-third-kernel-norm}, the uniform first-Picard majorant, and the linear tail bound give
\[
 \sup_{\la n\ra\sim N_0}\sup_{t\le T}
 \E\bigl|\widehat{\Delta_L^{\rm ext}\Gamma_1^{(3)}}(n,t)\bigr|^2
 \lesssim
 \sum_{M\gtrsim\max(N_0,L)}M^{6-9\alpha+\eps}.
\]
Summing first in the output scale gives the explicit estimate
\begin{align}
 &\sum_{N_0}N_0^{2s+3}
 \sup_{\la n\ra\sim N_0}\sup_{t\le T}
 \E\bigl|\widehat{\Delta_L^{\rm ext}\Gamma_1^{(3)}}(n,t)\bigr|^2
 \notag\\
 &\quad\lesssim
 \sum_{M\gtrsim L}M^{6-9\alpha+\eps}
       \sum_{N_0\lesssim M}N_0^{2s+3}
 \lesssim
 \sum_{M\gtrsim L}M^{-2\delta_\Gamma+\eps}
 \lesssim L^{-2\delta_\Gamma+\eps}.
\label{C:third-external-tail-sum}
\end{align}

We now spell out the internal Picard tail.  Let
\[
 A_{L,\eta}(M):=
 \sup_{\la m\ra\sim M}S_{2,>cL}^{(0,\eta)}(m),
\]
where \(c>0\) depends only on the localized cutoff profile.  The tail phase-sum estimate \eqref{B:tail-polynomial}, applied with the strict exponent \(s_V<7\alpha/2-3\), yields some \(\gamma>0\) such that
\begin{equation}\label{C:internal-picard-tail-input}
        \sum_M M^{2s_V+3}A_{L,\eta}(M)
        \lesssim L^{-\gamma}
\end{equation}
for sufficiently small \(\eta\), with \(\eta=0\) included.  On the resonant shell, the number of choices of \(m\) is \(O(M^3)\), while the untouched external convolution contributes \(M^{-2\alpha}\).  Therefore
\begin{align}
 &\sum_{N_0}N_0^{2s+3}
 \sup_{\la n\ra\sim N_0}\sup_{t\le T}
 \E\bigl|\widehat{\Delta_L^{\rm int}\Gamma_1^{(3)}}(n,t)\bigr|^2
 \notag\\
 &\quad\lesssim
 \sum_{N_0}N_0^{2s+3}
       \sum_{M\gtrsim N_0}M^{3-2\alpha}A_{L,0}(M) \notag\\
 &\quad\lesssim
 \sum_M M^{3-2\alpha}A_{L,0}(M)
       \sum_{N_0\lesssim M}N_0^{2s+3} \notag\\
 &\quad\lesssim
 \sum_M M^{2s+6-2\alpha}A_{L,0}(M)
 =
 \sum_M M^{2s_V+3}A_{L,0}(M)
 \lesssim L^{-\gamma}.
\label{C:third-internal-tail-sum}
\end{align}
This is the convolution step for the internal first-Picard tail.

For time increments, the term in which the Picard tail itself is incremented is controlled by \eqref{C:internal-picard-tail-input} with \(A_{L,\eta}\).  If instead the external convolution is incremented, the preceding sum contains
\(
M^{3-2\alpha+2\alpha\eta}A_{L,0}(M)
\);
this is \eqref{C:internal-picard-tail-input} with \(s_V\) replaced by \(s_V+\alpha\eta\).  The external-tail terms are bounded as in \eqref{C:third-external-tail-sum}, with the exponent increased by \(2\alpha\eta\).  Choose \(\eta\) and the harmless dyadic losses so that
\[
        s_V+\alpha\eta<\frac{7\alpha}{2}-3,
        \qquad
        2\alpha\eta+\eps<2\delta_\Gamma.
\]
The fixed-time and increment tails then satisfy the hypotheses of Lemma~\ref{lem:C-coeff-path} with some \(\theta>0\).  This proves almost sure Cauchy convergence for the standard localized dyadic family.

Finally, let \((\pi_\Lambda)_\Lambda\) be an arbitrary admissible cofinal family.  On a fixed finite set of Fourier and dyadic labels, every multiplier product in the third-chaos kernel converges entrywise.  The complements are uniformly controlled by the summable majorants in \eqref{C:third-external-tail-sum} and \eqref{C:third-internal-tail-sum}, with constants depending only on \(C_{\rm cut}\).  A finite-set plus dominated-tail argument, applied also to the increment kernels, proves convergence in \(C([0,T];\C^s)\) and identifies the limit with the standard dyadic one.  No numerical rate is used for this last step.
\end{proof}

\subsection{The first-chaos contraction as an integrated Volterra kernel}\label{subsec:C-first-chaos-integrated}

We construct the lower-chaos term in the integrated-kernel topology.  At finite cutoff, insert the remaining stochastic convolution
\[
        \wh\Psi_{2,\Lambda}(n,s)
        =\chi_\Lambda(n)\int_0^s K_2(s-u,n)\,d\beta^2_n(u)
\]
into \eqref{C:contraction-cutoff-form}.  Since the cutoff is finite, stochastic Fubini is an identity in the color-2 Gaussian Hilbert space.  Equivalently, because the first Wiener integral is a continuous linear map from the color-2 Hilbert space to \(L^2(\Omega)\),
\[
        \int_0^t I_1(H_s)\,ds
        =I_1\left(\int_0^tH_s\,ds\right)
        \quad\text{in }L^2(\Omega).
\]
Thus
\begin{equation}\label{C:C-Lambda-integrated}
        \wh C_{1,\Lambda}(n,t)
        =\int_0^t K_{21,\Lambda}(n;t,u)\,d\beta_n^2(u),
\end{equation}
where
\begin{equation}\label{C:K21-Lambda}
        K_{21,\Lambda}(n;t,u)
        =\int_u^t
        \sum_{r\in\Z^3}
        \mathbf c_\Lambda(n,r)
        \chi^{\mathrm{res}}(n-r,r)
        K_2(t-s,n-r)\sigma_1(r;s,t)K_2(s-u,n)\,ds,
\end{equation}
with \(\mathbf c_\Lambda(n,r)=\chi_\Lambda(n)\chi_\Lambda(r)^2\) for the cutoff convention above.

For \(0\le s\le t\), the covariance has the exact formula
\begin{equation}\label{C:sigma-exact}
\begin{aligned}
        \sigma_1(r;s,t)
        &=\int_0^sK_1(s-\tau,r)K_1(t-\tau,-r)\,d\tau    \\
        &=\frac{s}{2\omega_1(r)^2}\cos((t-s)\omega_1(r))
        +\frac{\sin((t-s)\omega_1(r))}{4\omega_1(r)^3}
        -\frac{\sin((t+s)\omega_1(r))}{4\omega_1(r)^3}.
\end{aligned}
\end{equation}
We write
\begin{equation}\label{C:sigma-split}
        \sigma_1(r;s,t)=\sigma_{1,0}(r;s,t)+R_1(r;s,t),
        \qquad
        \sigma_{1,0}(r;s,t)=\frac{s}{2\omega_1(r)^2}\cos((t-s)\omega_1(r)).
\end{equation}
On \(|r|\sim M\),
\begin{equation}\label{C:sigma-bounds}
        |\sigma_{1,0}(r;s,t)|\lesssim_T M^{-2\alpha},
        \qquad
        |R_1(r;s,t)|\lesssim_T M^{-3\alpha}.
\end{equation}
The remainder is absolutely summable at the threshold \(\alpha>3/4\).  The leading part requires a phase estimate.

\begin{lemma}[One-frequency phase-layer estimate]\label{lem:C-one-frequency-layer}
This is a layer-counting estimate for the balanced regime.
Let
\[
        R_M(n)=\{r\in\Z^3:\ |r|\sim M,\\ |n-r|\sim M\},
        \qquad M\gtrsim\la n\ra.
\]
For \(\sigma=(\sigma_0,\sigma_1,\sigma_2)\in\{\pm1\}^3\), set
\begin{equation}\label{C:Omega-def}
        \Omega^\sigma(n,r)
        =\sigma_0\omega_2(n)+\sigma_1\omega_2(n-r)+\sigma_2\omega_1(r).
\end{equation}
Then, for every \(0\le\delta<1\) and every \(\eps>0\),
\begin{equation}\label{C:Omega-layer-main}
        \sum_{r\in R_M(n)}
        \frakm(\Omega^\sigma(n,r))^{1-\delta}
        \lesssim_{\eps,\delta}
        M^{3-\alpha(1-\delta)+\eps},
        \qquad
        \frakm(\lambda)=\min\{1,|\lambda|^{-1}\}.
\end{equation}
In particular,
\begin{equation}\label{C:Omega-layer-delta0}
        \sum_{r\in R_M(n)}\frakm(\Omega^\sigma(n,r))
        \lesssim_\eps M^{3-\alpha+\eps}.
\end{equation}
\end{lemma}

\begin{proof}
There are two regimes.

\emph{Low-output high-high.}  Assume \(M\ge C_0\la n\ra\), where \(C_0\) will be chosen sufficiently large.  The output term satisfies
\[
        \omega_2(n)\lesssim \la n\ra^\alpha\le C_0^{-\alpha}M^\alpha.
\]
If the two high-frequency signs in
\[
        \sigma_1\omega_2(n-r)+\sigma_2\omega_1(r)
\]
are equal, their contribution has size \(\gtrsim M^\alpha\), and the output term is dominated by it once \(C_0\) is large.  If they are opposite, then by the high-frequency expansion
\[
        \omega_a(x)=c_a|x|^\alpha+O(|x|^{-\alpha})
\]
and by
\[
        \bigl||n-r|-|r|\bigr|\le |n|\le C_0^{-1}M,
\]
we have
\[
        \omega_2(n-r)-\omega_1(r)
        =(c_2-c_1)|r|^\alpha+O(C_0^{-1}M^\alpha)+O(M^{-\alpha}).
\]
Since \(c_1\ne c_2\), choosing \(C_0\) large and then absorbing finitely many low shells gives
\[
        |\omega_2(n-r)-\omega_1(r)|\gtrsim M^\alpha.
\]
Again the output term is too small to cancel.  Hence
\[
        |\Omega^\sigma(n,r)|\gtrsim M^\alpha
\]
in this regime, and therefore
\[
        \sum_{r\in R_M(n)}\frakm(\Omega^\sigma(n,r))^{1-\delta}
        \lesssim M^3M^{-\alpha(1-\delta)}.
\]

\emph{Balanced-output.}  It remains to treat \(M\sim\la n\ra\).  Put
\[
        d=|n|,
        \qquad
        \rho=|r|,
        \qquad
        \lambda=|n-r|.
\]
In two-center coordinates the full volume element is
\[
        dx=\frac{\rho\lambda}{d}\,d\rho\,d\lambda\,d\varphi,
\]
and after integrating the angle \(\varphi\),
\begin{equation}\label{C:two-center-C}
        dx=2\pi\frac{\rho\lambda}{d}\,d\rho\,d\lambda.
\end{equation}
On the balanced region \(d\sim\rho\sim\lambda\sim M\), the density is \(O(M)\).  Fix \(n\) and \(\rho\).  Then \(d=|n|\) and \(\rho\) are constants in the one-dimensional slice, and
\[
        \Omega^\sigma
        =\sigma_0\omega_2(d)+\sigma_1\omega_2(\lambda)+\sigma_2\omega_1(\rho).
\]
The output contribution \(\sigma_0\omega_2(d)\) is constant on this slice and cannot cancel the derivative in the \(\lambda\)-direction.  Since \(\lambda\sim M\),
\[
        |\partial_\lambda\Omega^\sigma|
        =|\omega_2'(\lambda)|
        \sim M^{\alpha-1}.
\]
A one-dimensional sublevel estimate gives, for dyadic \(1\le L\lesssim M^\alpha\), a \(\lambda\)-length
\[
        O(LM^{1-\alpha}+1)
\]
inside \(|\Omega^\sigma|\le L\).  Multiplying by the \(\rho\)-interval length \(O(M)\) and the density \(O(M)\), and using the standard unit-cube thickening of lattice layers, gives
\begin{equation}\label{C:Omega-sublevel}
        \#\{r\in R_M(n): |\Omega^\sigma(n,r)|\le L\}
        \lesssim_\eps M^\eps(M^2+M^{3-\alpha}L).
\end{equation}
Now decompose
\[
        \frakm(\Omega)^{1-\delta}
        \lesssim
        \sum_{1\le L\lesssim M^\alpha}
        L^{-(1-\delta)}\mathbf 1_{\{|\Omega|\le L\}}.
\]
Using \eqref{C:Omega-sublevel},
\begin{align*}
        \sum_{r\in R_M(n)}\frakm(\Omega^\sigma(n,r))^{1-\delta}
        &\lesssim_\eps
        \sum_{1\le L\lesssim M^\alpha}
        L^{-(1-\delta)}M^\eps(M^2+M^{3-\alpha}L)  \\
        &\lesssim_{\eps,\delta}
        M^{2+\alpha\delta+\eps}+M^{3-\alpha+\alpha\delta+\eps}.
\end{align*}
Since \(\alpha\le1\), the first term is no larger than the second up to the harmless loss.  This proves \eqref{C:Omega-layer-main}.
\end{proof}

For dyadic \(M\), let \(K_{21,M}\) denote the part of \eqref{C:K21-Lambda} with \(|r|\sim |n-r|\sim M\), after removing the finite cutoff.  The resonant cutoff forces \(M\gtrsim\la n\ra\).

\begin{lemma}[Single-shell integrated contraction kernel]\label{lem:C-single-shell}
For every \(\eps>0\),
\begin{equation}\label{C:single-shell}
        \sup_{0\le u\le t\le T}|K_{21,M}(n;t,u)|
        \lesssim_{T,\eps}
        \la n\ra^{-\alpha}M^{3-4\alpha+\eps},
        \qquad M\gtrsim\la n\ra.
\end{equation}
\end{lemma}

\begin{proof}
Split the covariance according to \eqref{C:sigma-split}.

The remainder contribution is estimated absolutely.  Using
\[
        |K_2(t-s,n-r)|\lesssim M^{-\alpha},
        \qquad
        |R_1(r;s,t)|\lesssim M^{-3\alpha},
        \qquad
        |K_2(s-u,n)|\lesssim\la n\ra^{-\alpha},
\]
and \(\#R_M(n)\lesssim M^3\), we get
\[
        |K^{\mathrm{rem}}_{21,M}(n;t,u)|
        \lesssim_T
        \la n\ra^{-\alpha}M^3M^{-\alpha}M^{-3\alpha}
        =\la n\ra^{-\alpha}M^{3-4\alpha}.
\]

For the leading covariance, expand
\[
        K_2(t-s,n-r),\qquad
        \cos((t-s)\omega_1(r)),\qquad
        K_2(s-u,n)
\]
into exponentials.  For each sign branch, the \(s\)-dependent phase is one of the \(\Omega^\sigma(n,r)\) in \eqref{C:Omega-def}, up to harmless sign changes.  The endpoint factors depending on \((t,u)\) have modulus one.  The static amplitude for one \(r\)-mode is
\[
        \la n\ra^{-\alpha}M^{-\alpha}M^{-2\alpha}
        =\la n\ra^{-\alpha}M^{-3\alpha}.
\]
The factor \(s\) in \(\sigma_{1,0}\) is harmless on \([0,T]\), and
\[
        \left|\int_u^t s e^{is\Omega}\,ds\right|
        \lesssim_T \frakm(\Omega).
\]
Therefore
\[
        |K^{\mathrm{lead}}_{21,M}(n;t,u)|
        \lesssim_T
        \la n\ra^{-\alpha}M^{-3\alpha}
        \sum_{r\in R_M(n)}\frakm(\Omega^\sigma(n,r)).
\]
Using \eqref{C:Omega-layer-delta0},
\[
        |K^{\mathrm{lead}}_{21,M}(n;t,u)|
        \lesssim_{T,\eps}
        \la n\ra^{-\alpha}M^{-3\alpha}M^{3-\alpha+\eps}
        =\la n\ra^{-\alpha}M^{3-4\alpha+\eps}.
\]
Combining the leading and remainder contributions proves \eqref{C:single-shell}.
\end{proof}

\begin{lemma}[Convergence of the integrated kernel]\label{lem:C-K21-conv}
Choose \(\eps>0\) so small that \(3-4\alpha+\eps<0\).  Then
\begin{equation}\label{C:K21-full-bound}
        \sup_{0\le u\le t\le T}|K_{21}(n;t,u)|
        \lesssim_{T,\eps}\la n\ra^{3-5\alpha+\eps}.
\end{equation}
More precisely, for every dyadic lower cutoff \(L\ge1\), the internal contraction tail
\[
        K_{21,>L}(n;t,u)
        =\sum_{\substack{M\gtrsim\la n\ra\\ M>L}}K_{21,M}(n;t,u)
\]
satisfies
\begin{equation}\label{C:K21-tail-bound}
        \sup_{u,t}|K_{21,>L}(n;t,u)|
        \lesssim_{T,\eps}
        \la n\ra^{-\alpha}
        \sum_{\substack{M\gtrsim\la n\ra\\ M>L}}M^{3-4\alpha+\eps}.
\end{equation}
In particular, for every fixed \(n\), \(K_{21,>L}(n;t,u)\to0\) uniformly in \((u,t)\) as \(L\to\infty\).
\end{lemma}

\begin{proof}
The tail bound \eqref{C:K21-tail-bound} follows directly from Lemma~\ref{lem:C-single-shell}.  Summing from \(M\gtrsim\la n\ra\) and using \(3-4\alpha+\eps<0\),
\[
        \sum_{M\gtrsim\la n\ra}M^{3-4\alpha+\eps}
        \lesssim \la n\ra^{3-4\alpha+\eps}.
\]
Multiplying by \(\la n\ra^{-\alpha}\) gives \eqref{C:K21-full-bound}.  The same negative exponent gives the fixed-mode tail convergence.
\end{proof}

\begin{proposition}[First-chaos fixed-time coefficient bound]\label{prop:C-first-chaos-fixed}
The limiting first-chaos contraction
\begin{equation}\label{C:C1-def}
        \wh C_1(n,t)=\int_0^tK_{21}(n;t,u)\,d\beta_n^2(u)
\end{equation}
satisfies, for every \(\eps>0\),
\begin{equation}\label{C:C1-fixed-bound}
        \sup_{t\le T}\E|\wh C_1(n,t)|^2
        \lesssim_{T,\eps}\la n\ra^{6-10\alpha+\eps}.
\end{equation}
\end{proposition}

\begin{proof}
By Itô isometry and \eqref{C:K21-full-bound},
\[
        \E|\wh C_1(n,t)|^2
        =\int_0^t|K_{21}(n;t,u)|^2\,du
        \lesssim_T \la n\ra^{6-10\alpha+\eps}.
\]
\end{proof}

\subsection{Time increments of the first-chaos contraction}

The time increment is the most delicate tracking step for the contraction.  The endpoint interval in the \(s\)-integral produces \(\frakm(\Omega)^{1-\eta}\), which is why Lemma~\ref{lem:C-one-frequency-layer} was stated with \(\delta>0\).

\begin{lemma}[Single-shell time increment]\label{lem:C-single-shell-incr}
Let \(0<\eta<1/2\).  For every \(\eps>0\), after decreasing \(\eta\) if necessary,
\begin{equation}\label{C:single-shell-incr}
        \sup_{0\le u\le t\le T}
        |K_{21,M}(n;t+h,u)-K_{21,M}(n;t,u)|
        \lesssim_{T,\eps,\eta}
        |h|^\eta\la n\ra^{-\alpha}M^{3-4\alpha+\alpha\eta+\eps}
\end{equation}
whenever the two kernels are evaluated on their common domain.  The same bound controls the change of the upper endpoint of the \(s\)-integral.
\end{lemma}

\begin{proof}
It suffices to consider \(h>0\).  On the common interval \(u\le s\le t\), decompose the difference into the increment of \(K_2(t-s,n-r)\), the increment of \(\sigma_1(r;s,t)\), and the corresponding endpoint phases in the exponential representation.

For the covariance remainder, the bounds are absolute and use
\[
        |K_2(t+h-s,n-r)-K_2(t-s,n-r)|
        \lesssim |h|^\eta M^{-\alpha+\alpha\eta}
\]
or, from \eqref{C:sigma-exact},
\[
        |R_1(r;s,t+h)-R_1(r;s,t)|
        \lesssim_T |h|^\eta M^{-3\alpha+\alpha\eta}.
\]
These estimates repeat the remainder part of Lemma~\ref{lem:C-single-shell} with an additional factor \(|h|^\eta M^{\alpha\eta}\), giving the right side of \eqref{C:single-shell-incr}.

For the leading covariance on the common interval, use the sign-branch representation.  A branch has the form
\[
        A_r\,E_r(t,u)\int_u^t s e^{is\Omega^\sigma(n,r)}\,ds,
\]
where
\[
        |A_r|\lesssim \la n\ra^{-\alpha}M^{-3\alpha},
\]
and \(E_r(t,u)\) is a product of endpoint exponentials with frequencies bounded by \(O(M^\alpha)\).  Hence
\[
        |E_r(t+h,u)-E_r(t,u)|
        \lesssim |h|^\eta M^{\alpha\eta}.
\]
Multiplying by
\[
        \left|\int_u^t s e^{is\Omega}\,ds\right|\lesssim_T\frakm(\Omega)
\]
and summing with \eqref{C:Omega-layer-delta0} gives
\[
        |h|^\eta\la n\ra^{-\alpha}M^{-3\alpha+\alpha\eta}
        M^{3-\alpha+\eps}
        =|h|^\eta\la n\ra^{-\alpha}M^{3-4\alpha+\alpha\eta+\eps}.
\]

It remains to control the upper endpoint interval \([t,t+h]\).  Each sign branch contains
\[
        \int_t^{t+h} s e^{is\Omega}\,ds.
\]
The interpolation estimate
\begin{equation}\label{C:endpoint-interpolation}
        \left|\int_t^{t+h} s e^{is\Omega}\,ds\right|
        \lesssim_T |h|^\eta\frakm(\Omega)^{1-\eta}
\end{equation}
comes from the two bounds \(O_T(|h|)\) and \(O_T(|\Omega|^{-1})\).  Therefore, by Lemma~\ref{lem:C-one-frequency-layer} with \(\delta=\eta\),
\[
\begin{aligned}
        &\la n\ra^{-\alpha}M^{-3\alpha}|h|^\eta
        \sum_{r\in R_M(n)}\frakm(\Omega^\sigma(n,r))^{1-\eta}   \\
        &\qquad\lesssim
        |h|^\eta\la n\ra^{-\alpha}M^{-3\alpha}
        M^{3-\alpha(1-\eta)+\eps}                              \\
        &\qquad=
        |h|^\eta\la n\ra^{-\alpha}M^{3-4\alpha+\alpha\eta+\eps}.
\end{aligned}
\]
This proves the endpoint part and hence \eqref{C:single-shell-incr}.
\end{proof}

\begin{lemma}[First-chaos time increments]\label{lem:C-first-incr}
For sufficiently small \(0<\eta<1/2\) and every \(\eps>0\),
\begin{equation}\label{C:C1-incr-bound}
        \E|\wh C_1(n,t+h)-\wh C_1(n,t)|^2
        \lesssim_{T,\eps,\eta}
        |h|^{2\eta}\la n\ra^{6-10\alpha+2\alpha\eta+\eps}.
\end{equation}
\end{lemma}

\begin{proof}
Choose \(\eta\) and \(\eps\) so small that \(3-4\alpha+\alpha\eta+\eps<0\).  Summing \eqref{C:single-shell-incr} over \(M\gtrsim\la n\ra\) gives
\begin{equation}\label{C:K21-full-incr}
        \sup_u|K_{21}(n;t+h,u)-K_{21}(n;t,u)|
        \lesssim
        |h|^\eta\la n\ra^{3-5\alpha+\alpha\eta+\eps}.
\end{equation}
Assume \(h>0\).  From \eqref{C:C1-def},
\[
\begin{aligned}
        \wh C_1(n,t+h)-\wh C_1(n,t)
        &=\int_0^t\bigl(K_{21}(n;t+h,u)-K_{21}(n;t,u)\bigr)\,d\beta_n^2(u) \\
        &\quad+
        \int_t^{t+h}K_{21}(n;t+h,u)\,d\beta_n^2(u).
\end{aligned}
\]
The first integral is controlled by Itô isometry and \eqref{C:K21-full-incr}.  The boundary integral is bounded by
\[
        |h|\sup_u|K_{21}(n;t+h,u)|^2
        \lesssim
        |h|\la n\ra^{6-10\alpha+\eps}
        \lesssim
        |h|^{2\eta}\la n\ra^{6-10\alpha+2\alpha\eta+\eps},
\]
since \(0<\eta<1/2\), \(|h|\le1\), and \(\la n\ra^{2\alpha\eta}\ge1\).  This proves \eqref{C:C1-incr-bound}.  The case \(h<0\) is identical.
\end{proof}

\subsection{Cutoff Cauchy convergence of the first-chaos contraction}

\begin{lemma}[Cauchy convergence of the infinite contraction tails]\label{lem:C-first-tail}
Let \(C_{1,\Lambda}\) be the finite-cutoff contraction in \eqref{C:C-Lambda-integrated}.  For every
\[
        s<5\alpha-\frac92,
\]
the family \((C_{1,\Lambda})_\Lambda\) is Cauchy almost surely in \(C([0,T];\C^s)\).
\end{lemma}

\begin{proof}
There are two tails.

\emph{Internal contraction tail.}  Let \(L\) be a dyadic lower cutoff and let \(C_{1,>L}\) be the part of the integrated kernel for which the contracted resonant shell satisfies \(M>L\).  By Itô isometry and \eqref{C:K21-tail-bound},
\begin{equation}\label{C:internal-tail-fixed}
        \sup_{t\le T}\E|\wh C_{1,>L}(n,t)|^2
        \lesssim_T
        \la n\ra^{-2\alpha}
        \left(
        \sum_{\substack{M\gtrsim\la n\ra\\ M>L}}
        M^{3-4\alpha+\eps}
        \right)^2.
\end{equation}
For fixed \(n\), the right hand side tends to zero as \(L\to\infty\).  We also need a summable dyadic majorant.  The full bound \eqref{C:C1-fixed-bound} gives
\begin{equation}\label{C:internal-majorant}
        \sup_{t\le T}\E|\wh C_{1,>L}(n,t)|^2
        \lesssim \la n\ra^{6-10\alpha+\eps}.
\end{equation}
For an explicit uniform tail, put \(\la n\ra\sim N\).  If \(N\le L\), then \eqref{C:internal-tail-fixed} gives
\[
        \sup_t\E|\wh C_{1,>L}(n,t)|^2
        \lesssim N^{-2\alpha}L^{2(3-4\alpha+\eps)}.
\]
Let \(s=5\alpha-9/2-\kappa\).  The dyadic weighted contribution from \(N\le L\) is bounded, up to harmless losses, by
\[
        L^{6-8\alpha+}
        \sum_{N\le L}N^{2s+3-2\alpha+}
        \lesssim L^{-2\kappa+},
\]
which tends to zero after choosing the losses smaller than \(\kappa\).  If \(N>L\), the full majorant gives the tail
\[
        \sum_{N>L}N^{2s+3}N^{6-10\alpha+}\to0
\]
for the same range of \(s\); more precisely, it is $O(L^{-2\kappa+})$.  Time-increment tails are identical, with the harmless additional factor \(M^{2\alpha\eta}\) and \(\eta\) chosen small compared with $\kappa$.

\emph{External remaining-color tail.}  The finite cutoff also removes the remaining color-2 Brownian mode when \(|n|\) is above the output cutoff.  This is an ordinary output-frequency tail.  The bounds \eqref{C:C1-fixed-bound} and \eqref{C:C1-incr-bound} give an $O(L^{-2\kappa+})$ weighted tail when $s=5\alpha-9/2-\kappa$.

Combining the internal and output tails and applying the cutoff part of Lemma~\ref{lem:C-coeff-path} proves the almost sure Cauchy property for the full countable cutoff family.
\end{proof}

\begin{proof}[Proof of Theorem~\ref{thm:C-cubic}]
Lemmas~\ref{lem:C-third-coeff} and \ref{lem:C-third-tail}, together with Lemma~\ref{lem:C-coeff-path}, give
\[
        \Gamma_1^{(3)}\in C([0,T];\C^s)
        \quad\text{for all }s<\frac{9\alpha}{2}-\frac92.
\]
Indeed, the coefficient exponent is \(p=6-9\alpha+\eps\), hence
\[
        -\frac32-\frac p2
        =\frac{9\alpha}{2}-\frac92-\frac\eps2.
\]

For the first-chaos contraction, Proposition~\ref{prop:C-first-chaos-fixed}, Lemma~\ref{lem:C-first-incr}, and Lemma~\ref{lem:C-first-tail} give
\[
        C_1\in C([0,T];\C^s)
        \quad\text{for all }s<5\alpha-\frac92.
\]
Here the coefficient exponent is \(p=6-10\alpha+\eps\), so
\[
        -\frac32-\frac p2
        =5\alpha-\frac92-\frac\eps2.
\]

The construction of \(\Gamma_2=V_1\circ\Psi_2\) is identical after exchanging the colors and speeds.  Since the H\"older exponents are strictly above the source indices used in Section~\ref{sec:deterministic-closure}, the cutoff convergence just proved gives
\[
        \Gamma_{i,\Lambda}\to\Gamma_i
        \quad\text{in}\quad
        L_T^1H^{s_2-\alpha}\cap L_T^1B_{2,\infty}^{\sigma-\alpha}
\]
for every admissible choice of $(s_2,\sigma)$ in Definition~\ref{def:admissible}.  This supplies the cubic component of the enhanced-data distance.

\end{proof}

\section{Weak-covariance stochastic symbols}\label{app:weak-covariance-symbols}

Assume the Fourier-diagonal covariance of Definition~\ref{def:wc-covariance} and the gain \(\delta_{\mathsf R}>0\) from Definition~\ref{def:wc-gain}.  We construct the covariance-split quadratic and cubic symbols and prove their cutoff convergence.  The zero-mode quadratic branch has dyadic size \(N^{3-3\alpha-\kappa}\).

\subsection{Real Fourier covariance algebra}

Fix once and for all a set \(\mathbb Z^3_+\) containing exactly one representative of each nonzero pair \(\{n,-n\}\).  For \(n\in\mathbb Z^3_+\), write the Hermitian color covariance as
\[
        \mathsf R(n)=A(n)+iB(n),
        \qquad A(n)^\mathsf T=A(n),\qquad B(n)^\mathsf T=-B(n).
\]
The corresponding real covariance matrix on the vector
\((X_1(n),X_2(n),Y_1(n),Y_2(n))\in\mathbb R^4\) is
\begin{equation}\label{eq:real-fourier-covariance-matrix}
        \mathsf C(n)=\frac12
        \begin{pmatrix}
        A(n)&-B(n)\\
        B(n)& A(n)
        \end{pmatrix}.
\end{equation}
At the zero mode one uses the real two-dimensional covariance \(\mathsf R(0)\); the compatibility condition \(\mathsf R(-n)=\overline{\mathsf R(n)}\) makes \(\mathsf R(0)\) real symmetric.  Positivity of Definition~\ref{def:wc-covariance} is equivalently positivity of \(\mathsf C(n)\) for every involution class.  The local Gaussian Hilbert space is the quotient
\begin{equation}\label{eq:real-fourier-quotient-space}
        \mathfrak h_n^{\mathsf R}:=\mathbb R^4/\ker\mathsf C(n),
        \qquad
        \langle v,w\rangle_{\mathfrak h_n^{\mathsf R}}=v^\mathsf T\mathsf C(n)w,
\end{equation}
with the evident two-dimensional version when \(n=0\).  If
\(Z_a(n)=X_a(n)+iY_a(n)\) and \(Z_a(-n)=\overline{Z_a(n)}\), then the local covariance identities are
\begin{equation}\label{eq:real-fourier-complex-covariances}
        \mathbb E[Z_a(n)Z_b(-n)]=\mathsf R_{ab}(n),
        \qquad
        \mathbb E[Z_a(n)Z_b(n)]=0
        \quad(n\ne0).
\end{equation}
Tensoring \(\mathfrak h_n^{\mathsf R}\) with the finite-dimensional time-kernel span generated by the stochastic convolutions gives the finite Gaussian space used below.  All Wick products are the ordinary homogeneous-chaos projections in the direct sum of these quotient Hilbert spaces.

\begin{lemma}[Local R-Wick identities in one involution class]
\label{lem:wc-real-fourier-local-wick}
Let \(F_a(n,t)\) be any finite stochastic-convolution coefficient carried by the real Fourier class \(\{n,-n\}\).  Then, at finite cutoff,
\begin{equation}\label{eq:wc-local-wick-identity}
        F_a(n,s)F_b(m,t)
        =\bigl(:F_a(n,s)F_b(m,t):\bigr)_{\mathsf R}
        +\mathbf 1_{n+m=0}\,\Sigma_{ab}^{\mathsf R}(n;s,t).
\end{equation}
If the local matrix \(\mathsf C(n)\) is degenerate, the same identity holds in the quotient space \eqref{eq:real-fourier-quotient-space}.  The formula uses covariance pairings and homogeneous-chaos projections.
\end{lemma}

\begin{proof}
For nonzero \(n\), the real covariance matrix \eqref{eq:real-fourier-covariance-matrix} gives \eqref{eq:real-fourier-complex-covariances}.  The time kernels enter by tensoring with their deterministic \(L^2\)-inner products.  Hence the covariance of \(F_a(n,s)\) and \(F_b(m,t)\) is zero unless \(m=-n\), and is \(\Sigma_{ab}^{\mathsf R}(n;s,t)\) when \(m=-n\).  The second-order Wick identity in a finite Gaussian Hilbert space is \(XY=:XY:+\mathbb E[XY]\), which gives \eqref{eq:wc-local-wick-identity}.  Passing to the quotient by the null space leaves the covariance form and the homogeneous-chaos projection unchanged.
\end{proof}

\begin{lemma}[Finite R-Wick convention with degenerate covariance]\label{lem:wc-degenerate-finite-wick}
Fix a Galerkin cutoff and finitely many time kernels.  The Gaussian variables generated by the two colors and by the finitely many Fourier labels form a finite real Gaussian vector with positive semidefinite covariance matrix.  Its Wick products are defined by first realizing this vector as a linear image of a standard Euclidean Gaussian vector, equivalently by working in the finite Gaussian Hilbert space obtained after quotienting out the zero-norm directions.  This definition is independent of the chosen square root of the covariance matrix.  In particular, for first-chaos variables one has the exact finite identity
\[
        X_aX_b=:X_aX_b:_{\mathsf R}+\mathbf E[X_aX_b],
\]
and the analogous higher-degree product formulas are the ordinary finite Gaussian Wick formulas in that quotient Hilbert space.

If the covariance is Fourier-diagonal in the sense of Definition~\ref{def:wc-covariance}, the only non-conjugated pairing between labels \((a,m)\) and \((b,m')\) is supported on \(m+m'=0\) and has coefficient \(\mathsf R_{ab}(m)\) times the corresponding time-kernel pairing.  Thus Wick contraction changes the color-chaos order but never mixes different spatial Fourier labels.
\end{lemma}

\begin{proof}
Every finite positive semidefinite covariance matrix \(\Sigma\) can be written as \(AA^*\).  If \(G\) is a standard Euclidean Gaussian vector, then \(AG\) realizes the required finite Gaussian family.  Two choices of \(A\) give the same Gaussian law, and Wick products are characterized by orthogonal projection onto homogeneous chaoses, hence are independent of the choice; this is the finite-dimensional Gaussian Hilbert-space construction, see \cite{Janson,Nualart}.  Equivalently, one may define the Gaussian Hilbert space as the span of the variables modulo the null space of the covariance form.  This quotient is taken before any cutoff limit and only within a fixed finite collection of Fourier labels and time kernels.  If a sequence of covariance matrices loses rank, the quotient simply removes directions whose variances tend to zero; Wick polynomials are still given by the same finite contraction formulas because those formulas depend only on covariance entries.  The product identities are then the standard finite Gaussian identities.  The final support statement follows directly from the factor \(\mathbf 1_{m+m'=0}\) in Definition~\ref{def:wc-covariance}.
\end{proof}

The preceding lemma fixes the finite \(\mathsf R\)-Wick convention used below.  On each finite dyadic set, Wick polynomials depend continuously on the finite covariance matrix and on the deterministic kernels.  Orthogonalization is performed inside the real Fourier involution class \(\{m,-m\}\); after quotienting zero-norm directions, it is a local coordinate choice.  Therefore the spatial labels carried by the coefficient tensor remain unchanged, and Fourier-diagonal support preserves both
\[
        n=q+\ell+r \quad\text{for centered coefficient tensors},
        \qquad
        n=q \quad\text{for lower-chaos covariance branches}.
\]
The dyadic tails are controlled by the summable envelopes proved below and in Subsection~\ref{subsec:weak-covariance-extension}.

\begin{proposition}[Covariance approximation and rank-deficient limits]\label{prop:wc-covariance-approximation}
Let \(\mathsf R^{(m)}\) and \(\mathsf R\) satisfy Definition~\ref{def:wc-covariance} with a common exponent \(\kappa\), with
\[
        d_\kappa(\mathsf R^{(m)},\mathsf R)\to0,
\]
where \(d_\kappa\) is defined in \eqref{eq:wc-covariance-distance}.  Fix a finite Galerkin cutoff and finitely many dyadic blocks and time parameters.  Then every finite \(\mathsf R^{(m)}\)-Wick component appearing in the quadratic lift, the mixed centered blocks, and the cubic contractions converges in every finite-dimensional \(L^p\) norm to the corresponding \(\mathsf R\)-Wick component.  The statement remains valid if the ranks of the local covariance matrices change in the limit.

After removing the finite cutoff, the convergence holds in the covariance-augmented enhanced-data distance \(d_T^{\rm wc}\) on every interval on which the data have a common bounded size.  The construction uses local covariance square roots or quotient Hilbert spaces and avoids inverse covariance matrices or spectral gap assumptions on \(\mathsf R(n)\).
\end{proposition}

\begin{proof}
At fixed cutoff use the common Brownian coordinates in \eqref{eq:canonical-covariance-coupling}.  On the resulting finite Euclidean space,
\[
 \mathsf C_{\mathsf R^{(m)}}(n)^{1/2}B_{[n]}
 \longrightarrow
 \mathsf C_{\mathsf R}(n)^{1/2}B_{[n]}
\]
almost surely and in every finite $L^p$ norm, because the positive semidefinite square-root map is continuous.  Wick products are finite polynomials in these Gaussian coordinates and the covariance pairings, so the same convergence holds for every quadratic, mixed, and cubic finite-cutoff component.  This argument remains valid when the limiting covariance loses rank and uses no inverse covariance matrix.

For the infinite enhanced-data convergence, choose a finite dyadic set.  The preceding paragraph gives convergence on that set.  The complement is controlled uniformly in \(m\) by the summable estimates for the zero-mode covariance branch in Proposition~\ref{prop:wc-theta-branch}, the mixed covariance diagonals in Proposition~\ref{prop:wc-deterministic-diagonals}, the centered covariance remainders in Proposition~\ref{prop:wc-centered-component}, and the cubic covariance branches in Proposition~\ref{prop:wc-cubic-symbols}.  The standard finite-set/tail argument gives convergence in \(d_T^{\rm wc}\).
\end{proof}

\subsection{The covariance-split quadratic symbol}
We use standard finite-dimensional Gaussian/Wiener-chaos algebra throughout this appendix; see \cite{Nualart,PeccatiTaqqu} for background.  At finite cutoff define the \(\mathsf R\)-Wick product by
\begin{equation}\label{eq:wc-theta-finite-split}
        \Psi_{1,\Lambda}\Psi_{2,\Lambda}
        =:\Psi_{1,\Lambda}\Psi_{2,\Lambda}:_{\mathsf R}
        +C_{\Theta,\Lambda}^{\mathsf R},
\end{equation}
where the deterministic branch is the spatial zero mode
\begin{equation}\label{eq:wc-theta-deterministic-branch}
        \widehat C_{\Theta,\Lambda}^{\mathsf R}(n,t)
        =\mathbf 1_{n=0}\sum_\ell c_\Lambda(\ell)c_\Lambda(-\ell)
        \Sigma_{12}^{\mathsf R}(\ell;t,t).
\end{equation}

\begin{proposition}[Zero-mode covariance branch]\label{prop:wc-theta-branch}
Let \(C_{\Theta,\Lambda}^{\mathsf R}\) be defined by \eqref{eq:wc-theta-deterministic-branch}.  Decompose the high loop dyadically and write
\begin{equation}
        c_{\Theta,N}^{\mathsf R}(t)
        :=\sum_{|\ell|\sim N}\rho_N(\ell)
        \Sigma_{12}^{\mathsf R}(\ell;t,t),
        \qquad
        C_{\Theta}^{\mathsf R}(t,x)=\sum_N c_{\Theta,N}^{\mathsf R}(t).
\end{equation}
For every \(\epsilon_0>0\) small enough and every dyadic \(N\),
\begin{equation}\label{eq:wc-theta-shell}
        \|c_{\Theta,N}^{\mathsf R}\|_{L_T^\infty}
        \lesssim_{T,\epsilon_0}
        \|\mathsf R\|_\kappa N^{3-3\alpha-\kappa+\epsilon_0}.
\end{equation}
Moreover, if \(0<\eta<\min\{1/2,(\kappa-(3-3\alpha))/(4\alpha)\}\), then
\begin{equation}\label{eq:wc-theta-increment}
        |c_{\Theta,N}^{\mathsf R}(t)-c_{\Theta,N}^{\mathsf R}(t')|
        \lesssim_{T,\eta,\epsilon_0}
        |t-t'|^\eta
        \|\mathsf R\|_\kappa N^{3-3\alpha-\kappa+\alpha\eta+\epsilon_0}.
\end{equation}
Consequently \(C_{\Theta,\Lambda}^{\mathsf R}\) converges in \(C_T C^r\), for every \(r\in\mathbb R\), to a deterministic zero-mode function \(C_\Theta^{\mathsf R}\), and the cutoff tail satisfies
\begin{equation}\label{eq:wc-theta-tail}
        \sup_t|C_{\Theta}^{\mathsf R}(t)-C_{\Theta,\Lambda}^{\mathsf R}(t)|
        \lesssim
        \|\mathsf R\|_\kappa \Lambda^{-(\kappa-(3-3\alpha))+\epsilon_0}
\end{equation}
for every sufficiently small \(\epsilon_0>0\).  In particular,
\begin{equation}\label{eq:wc-theta-branch-bound}
        \|C_{\Theta}^{\mathsf R}\|_{C_T C^r}
        \lesssim_{T,r,\delta_{\mathsf R}}
        \|\mathsf R\|_\kappa.
\end{equation}
\end{proposition}

\begin{proof}
Lemma~\ref{lem:wc-offdiag-covariance} gives \(|\Sigma_{12}^{\mathsf R}(\ell;t,t)|\lesssim_T \|\mathsf R\|_\kappa |\ell|^{-3\alpha-\kappa}\) on \(|\ell|\sim N\).  The number of lattice points in the shell is \(O(N^{3+\epsilon_0})\), with \(\epsilon_0\) used only as a uniform dyadic loss.  This proves \eqref{eq:wc-theta-shell}.  The increment bound \eqref{eq:wc-theta-increment} follows from \eqref{eq:wc-offdiag-cov-increment} in the same way.  The exponent in \eqref{eq:wc-theta-increment} is negative for the displayed range of \(\eta\), so the dyadic series is uniformly summable in a small Hölder time norm.  The zero spatial frequency is the only output mode in \eqref{eq:wc-theta-deterministic-branch}; therefore every spatial \(C^r\) norm is equivalent to the absolute value of the scalar coefficient.  Summing \eqref{eq:wc-theta-shell} over high-loop shells gives the limit and the tail bound \eqref{eq:wc-theta-tail}.
\end{proof}

\begin{lemma}[Picard phase sum for the weak covariance branch]\label{lem:wc-first-picard-phase-sum}
Let
\[
        \widetilde V_{i,\Lambda}^{\mathsf R}:=
        I_i\bigl(:\Psi_{1,\Lambda}\Psi_{2,\Lambda}:_{\mathsf R}\bigr)
\]
be the centered part of the weak-covariance first Picard object.  For \(\nu=0,1\), every finite cutoff \(\Lambda\), and every sufficiently small \(\eta>0\), the coefficient bounds
\begin{align}
        \sup_{t\le T}\E\bigl|\partial_t^\nu\widehat{\widetilde V}_{i,\Lambda}^{\mathsf R}(n,t)\bigr|^2
        &\lesssim_{T,\epsilon_0}\|\mathsf R\|_\kappa^{O(1)}
        \langle n\rangle^{3-(7-2\nu)\alpha+\epsilon_0},\label{eq:wc-picard-centered-fixed}\\
        \E\bigl|\partial_t^\nu\widehat{\widetilde V}_{i,\Lambda}^{\mathsf R}(n,t)-
        \partial_t^\nu\widehat{\widetilde V}_{i,\Lambda}^{\mathsf R}(n,t')\bigr|^2
        &\lesssim_{T,\epsilon_0,\eta}\|\mathsf R\|_\kappa^{O(1)}
        |t-t'|^{2\eta}
        \langle n\rangle^{3-(7-2\nu)\alpha+2\alpha\eta+\epsilon_0}
        \label{eq:wc-picard-centered-increment}
\end{align}
hold uniformly in \(\Lambda\).  The corresponding cutoff-tail estimates are controlled by the same dyadic majorants as in Lemma~\ref{lem:B-tail-sums}.  Consequently
\[
        \widetilde V_i^{\mathsf R}\in C_T\mathcal C^{\frac72\alpha-3-},
        \qquad
        \partial_t\widetilde V_i^{\mathsf R}\in C_T\mathcal C^{\frac52\alpha-3-},
\]
and the same convergence holds in the dyadic Besov bounds
\[
        \widetilde V_{i,\Lambda}^{\mathsf R}\to \widetilde V_i^{\mathsf R}
        \quad\text{in}\quad
        L_T^\infty B_{2,\infty}^{\rho_V},
        \qquad
        \partial_t\widetilde V_{i,\Lambda}^{\mathsf R}\to \partial_t\widetilde V_i^{\mathsf R}
        \quad\text{in}\quad
        L_T^\infty B_{2,\infty}^{\rho_V-\alpha}
\]
for every strict choice of \(\rho_V<7\alpha/2-3\).
\end{lemma}

\begin{proof}
At finite cutoff, the stochastic Fubini representation in Appendix~\ref{app:first-picard-phase} is used with the product of the two stochastic convolution legs read as a \(\mathsf R\)-Wick product in the finite joint Gaussian space.  By Lemma~\ref{lem:wc-degenerate-finite-wick}, this Wick product is defined in the quotient Hilbert space if the local covariance is degenerate.  Orthogonalizing the finitely many color-time kernels is performed separately inside each real Fourier involution class; hence it may change only bounded auxiliary color-time coordinates while the external Fourier labels \(k\) and \(n-k\) remain attached to the deterministic kernel.

Therefore the deterministic kernel entering the second-chaos isometry has exactly the same three-frequency phase
\[
        \sigma_0\omega_i(n)+\sigma_1\omega_1(k)+\sigma_2\omega_2(n-k)
\]
as in Proposition~\ref{prop:B-phase-sum}.  The positive semidefinite covariance matrix has normalized diagonal entries and is uniformly bounded at every label, so the finite Hilbert norm of the \(\mathsf R\)-centered tensor is bounded, up to a polynomial in \(\|\mathsf R\|_\kappa\), by the same kernel norm as in the independent-color proof.  The fixed-time estimate is therefore bounded by the deterministic phase sum \(S_i^{(\nu)}(n)\) from \eqref{B:phase-sum-def}; the time-increment estimate is bounded by \(S_i^{(\nu,\eta)}(n)\).  Proposition~\ref{prop:B-phase-sum} gives precisely \eqref{eq:wc-picard-centered-fixed} and \eqref{eq:wc-picard-centered-increment}.

For cutoff differences, one inserts the cutoff difference before the same class-wise orthonormalization.  At least one stochastic leg is then restricted to the removed Fourier tail, and the deterministic tail phase sums are those of Lemma~\ref{lem:B-tail-sums}.  The coefficient-to-path criterion Lemma~\ref{lem:B-coeff-besov} gives the stated almost sure convergence in the H\"older and Besov topologies.  All covariance dependence is carried by the bounded local Hilbert metric.
\end{proof}

\begin{proposition}[Weak-covariance quadratic lift]\label{prop:wc-quadratic-lift}
The centered quadratic symbol
\[
        \Theta^{\mathsf R}:=:\Psi_1\Psi_2:_{\mathsf R}
\]
exists almost surely in \(C_T\mathcal C^{2\alpha-3-}\), with cutoff convergence in that topology.  The full finite product is
\[
        \Psi_1\Psi_2=\Theta^{\mathsf R}+C_\Theta^{\mathsf R}.
\]
Defining
\begin{equation}\label{eq:wc-first-picard-def}
        V_i^{\mathsf R}:=I_i(\Theta^{\mathsf R}+C_\Theta^{\mathsf R}),
\end{equation}
one has both the direct and the oscillatory first-Picard estimates of Theorems~\ref{thm:stochastic-lift} and~\ref{thm:first-picard-main}.  More explicitly,
\[
        V_i^{\mathsf R}\in C_T\mathcal C^{\frac72\alpha-3-}
        \cap L_T^\infty B_{2,\infty}^{\rho_V},
        \qquad
        \partial_tV_i^{\mathsf R}\in C_T\mathcal C^{\frac52\alpha-3-}
        \cap L_T^\infty B_{2,\infty}^{\rho_V-\alpha},
\]
with cutoff convergence in both the H\"older multiplier topology and the dyadic Besov topology.  The deterministic contribution \(I_i(C_\Theta^{\mathsf R})\) is a spatial zero mode and is absorbed in these bounds.
\end{proposition}

\begin{proof}
The finite split \eqref{eq:wc-theta-finite-split} is an identity in the cutoff Gaussian space.  We first estimate the centered branch.  For a fixed Fourier coefficient,
\[
        \widehat\Theta_{\Lambda}^{\mathsf R}(n,t)
        =\sum_{k+\ell=n}c_\Lambda(k)c_\Lambda(\ell)
        :\widehat\Psi_1(k,t)\widehat\Psi_2(\ell,t):_{\mathsf R}.
\]
By the finite Wick isometry in the joint Gaussian Hilbert space, the second moment of this coefficient is a finite sum of products of two covariance entries.  Since the local two-color covariance matrices are positive semidefinite with normalized diagonal entries, Cauchy--Schwarz gives
\[
        |\E[\widehat\Psi_a(m,t)\widehat\Psi_b(m',t)]|
        \le
        \bigl(\E|\widehat\Psi_a(m,t)|^2\bigr)^{1/2}
        \bigl(\E|\widehat\Psi_b(m',t)|^2\bigr)^{1/2}.
\]
Thus every Wick-pairing contribution is dominated by the diagonal variance product used in the independent proof.  Consequently
\[
        \E|\widehat\Theta^{\mathsf R}(n,t)|^2
        \lesssim \sum_{k+\ell=n}\langle k\rangle^{-2\alpha}\langle \ell\rangle^{-2\alpha}
        \lesssim \langle n\rangle^{3-4\alpha+}.
\]
The same argument with one stochastic convolution replaced by its time increment gives the corresponding coefficient increment bound.  Lemma~\ref{lem:coeff-chaos-criterion} yields the path regularity.  For cutoff differences, the tail-indicator argument following \eqref{eq:Theta-increment} gives a polynomial weighted second-moment gain after an arbitrarily small spatial loss.  Hypercontractivity and Borel--Cantelli over the countable cutoff family therefore give almost-sure convergence of the centered quadratic lift in $C_T\mathcal C^{2\alpha-3-}$.

The lower-chaos branch is the deterministic zero mode $C_\Theta^{\mathsf R}$.  Proposition~\ref{prop:wc-theta-branch} gives an absolutely summable shell envelope $N^{3-3\alpha-\kappa+}$ and cutoff tails.  Since this branch is spatially constant, all spatial H\"older norms are equivalent to the absolute value of its scalar coefficient.

Applying $I_i$ to the centered branch gives the direct first-Picard estimates.  Lemma~\ref{lem:wc-first-picard-phase-sum} gives the phase-summed bounds: its proof uses the phase-sum theorem of Appendix~\ref{app:first-picard-phase}, because the finite $\mathsf R$-Wick projection changes only the local color Hilbert coordinates and preserves the external Fourier labels entering the three-frequency phase.  The deterministic zero-mode branch $I_i(C_\Theta^{\mathsf R})$ has only the zero spatial frequency and is therefore smoother than the stochastic branch in both the H\"older multiplier topology and the dyadic $B_{2,\infty}$ bound.  Combining these estimates proves the displayed bounds and cutoff convergence for $V_i^{\mathsf R}$ and $\partial_tV_i^{\mathsf R}$.
\end{proof}

\subsection{Cubic symbols with off-diagonal contractions}
We display \(\Gamma_1^{\mathsf R}=V_2^{\mathsf R}\circ\Psi_1\); the other symbol is obtained by exchanging the two colors and the two speeds.  The point of this subsection is that the finite \(\mathsf R\)-Wick expansion is performed before any infinite contraction kernel is introduced.  Write
\[
        V_2^{\mathsf R}=I_2(\Theta^{\mathsf R})+I_2(C_\Theta^{\mathsf R}),
        \qquad \Theta^{\mathsf R}=:\Psi_1\Psi_2:_{\mathsf R}.
\]
The first summand produces a third-chaos term and two first-chaos contractions; the second summand is a regular zero-mode contribution.  In the exact finite formulas below, \(\Sigma_{ab}^{\mathsf R}(m;s,t)\) is indexed by the Fourier label \(m\) of the first member of the contracted pair.  Hence, when the outer resonant leg is denoted by \(r\) and the paired inner leg is \(-r\), the covariance factor is \(\Sigma_{ab}^{\mathsf R}(-r;s,t)\).  The later estimates are invariant under \(r\mapsto-r\), but the sign is kept in the algebraic formulas.  Thus
\begin{equation}\label{eq:wc-cubic-decomp}
        \Gamma_i^{\mathsf R}
        =\Gamma_{i,\mathsf R}^{(3)}+C_{i,\mathrm{diag}}^{\mathsf R}
        +C_{i,\mathrm{off}}^{\mathsf R}+R_{i,\Theta}^{\mathsf R}.
\end{equation}
Here \(C_{i,\mathrm{diag}}^{\mathsf R}\) is the same-color contraction already present in the independent model, \(C_{i,\mathrm{off}}^{\mathsf R}\) is the additional off-diagonal first-chaos branch generated by \(\mathsf R_{12}\) or \(\mathsf R_{21}\), and
\[
        R_{i,\Theta}^{\mathsf R}
        := I_{3-i}(C_\Theta^{\mathsf R})\circ\Psi_i.
\]
The last formula records the index convention \(\Gamma_i=V_{3-i}\circ\Psi_i\): the zero-mode forcing is transported by \(I_{3-i}\).

\begin{lemma}[Finite R-Wick cubic product]\label{lem:wc-cubic-finite-wick}
Fix a Galerkin cutoff.  Let \(h_1,h_2,g\) be deterministic first-chaos kernels carrying colors \(1,2,1\), respectively, in the finite joint Gaussian space with covariance \(\mathsf R\).  Then
\begin{align}\label{eq:wc-cubic-wick-product}
        I^{:2:}_{\mathsf R}(h_1,h_2)\,W_1(g)
        &= I^{:3:}_{\mathsf R}(h_1,h_2,g)
        + \langle h_1,g\rangle_{\mathsf R} W_2(h_2)
        + \langle h_2,g\rangle_{\mathsf R} W_1(h_1).
\end{align}
The two contractions are supported, respectively, on the Fourier involutions of the color pairs \((1,1)\) and \((2,1)\).  In the Fourier variables of \(\Gamma_1^{\mathsf R}\), the first contraction imposes \(a+r=0\) and leaves the color-2 mode \(b=n\); the second imposes \(b+r=0\) and leaves the color-1 mode \(a=n\).  The internal \((1,2)\) contraction is the quadratic branch already subtracted in \(\Theta^{\mathsf R}=:\Psi_1\Psi_2:_{\mathsf R}\), so the displayed two branches are the lower-chaos terms in \(\Gamma_1^{\mathsf R}\).
\end{lemma}

\begin{proof}
This is the degree-three product formula in the finite \(\mathsf R\)-Wick convention of Lemma~\ref{lem:wc-degenerate-finite-wick}.  One may realize the finite vector by a covariance square root in an auxiliary Euclidean Gaussian space, apply the ordinary Wick product identity there, and then use that Wick products are determined by the resulting Gaussian law.  This also covers rank-deficient local covariance matrices.  The Fourier support of each pairing follows from the Fourier-diagonal covariance \(\mathbf 1_{m+m'=0}\mathsf R_{ab}(m)\).  Hence \(a+r=0\) in the same-color branch and \(b+r=0\) in the off-diagonal branch.  The centered quadratic \(\Theta^{\mathsf R}\) is the projection onto the corresponding centered chaos.
\end{proof}

For \(\Gamma_1^{\mathsf R}\), the off-diagonal contraction has the finite-cutoff form
\begin{equation}\label{eq:wc-cubic-offdiag-formula}
\widehat C_{1,\mathrm{off},\Lambda}^{\mathsf R}(n,t)
=\chi_\Lambda(n)\int_0^t M_{1,\mathrm{off},\Lambda}^{\mathsf R}(n;t,s)
        \widehat\Psi_{1,\Lambda}(n,s)\,ds,
\end{equation}
where
\begin{equation}\label{eq:wc-cubic-offdiag-multiplier}
M_{1,\mathrm{off},\Lambda}^{\mathsf R}(n;t,s)
=\sum_r c_\Lambda(n,r)\chi^{\rm res}(n-r,r)
        K_2(t-s,n-r)\Sigma_{21}^{\mathsf R}(-r;s,t).
\end{equation}
Here \(c_\Lambda(n,r)\) is a bounded finite-cutoff multiplier converging pointwise to one on fixed \((n,r)\).  The formula for \(C_{2,\mathrm{off},\Lambda}^{\mathsf R}\) is obtained by exchanging \(1\) and \(2\): the remaining first-chaos mode is then \(\Psi_2(n,s)\), and the Duhamel kernel is \(K_1(t-s,n-r)\).

\begin{lemma}[Integrated off-diagonal cubic kernel]\label{lem:wc-offdiag-cubic-integrated-kernel}
After inserting the remaining stochastic convolution in \eqref{eq:wc-cubic-offdiag-formula}, one has the finite identity
\begin{equation}\label{eq:wc-offdiag-cubic-first-chaos-integrated}
        \widehat C_{1,\mathrm{off},\Lambda}^{\mathsf R}(n,t)
        =\int_0^t K_{1,\mathrm{off},\Lambda}^{\mathsf R}(n;t,u)\,d\beta_n^1(u),
\end{equation}
where
\begin{equation}\label{eq:wc-offdiag-cubic-integrated-kernel}
K_{1,\mathrm{off},\Lambda}^{\mathsf R}(n;t,u)
=\int_u^t M_{1,\mathrm{off},\Lambda}^{\mathsf R}(n;t,s)
        K_1(s-u,n)\,ds.
\end{equation}
For the dyadic piece in which \(|r|\sim |n-r|\sim N\), denoted by \(K_{1,\mathrm{off},N}^{\mathsf R}\), the following bounds hold for every small \(\epsilon_0>0\):
\begin{equation}\label{eq:wc-offdiag-cubic-shell-kernel}
        \sup_{0\le u\le t\le T}
        |K_{1,\mathrm{off},N}^{\mathsf R}(n;t,u)|
        \lesssim_{T,\epsilon_0}
        \|\mathsf R\|_\kappa\,
        \langle n\rangle^{-\alpha}N^{3-4\alpha-\kappa+\epsilon_0},
        \qquad N\gtrsim \langle n\rangle.
\end{equation}
Moreover, for \(0<\eta\le1/2\),
\begin{equation}\label{eq:wc-offdiag-cubic-shell-increment}
        \sup_u
        |K_{1,\mathrm{off},N}^{\mathsf R}(n;t,u)-K_{1,\mathrm{off},N}^{\mathsf R}(n;t',u)|
        \lesssim_{T,\eta,\epsilon_0}
        |t-t'|^\eta\|\mathsf R\|_\kappa
        \langle n\rangle^{-\alpha}N^{3-4\alpha-\kappa+\alpha\eta+\epsilon_0}.
\end{equation}
The same estimates hold for the color-exchanged kernel \(K_{2,\mathrm{off},N}^{\mathsf R}\).
\end{lemma}

\begin{proof}
The identity \eqref{eq:wc-offdiag-cubic-first-chaos-integrated} is a finite-dimensional stochastic Fubini identity in the remaining color-one Gaussian Hilbert space.  Indeed
\[
        \widehat\Psi_{1,\Lambda}(n,s)
        =\chi_\Lambda(n)\int_0^s K_1(s-u,n)\,d\beta_n^1(u),
\]
and all sums are finite before the cutoff is removed.

On the shell \(|r|\sim |n-r|\sim N\), the resonant cutoff forces \(N\gtrsim\langle n\rangle\).  The Duhamel factor satisfies \(|K_2(t-s,n-r)|\lesssim N^{-\alpha}\), Lemma~\ref{lem:wc-offdiag-covariance} gives
\(|\Sigma_{21}^{\mathsf R}(-r;s,t)|\lesssim\|\mathsf R\|_\kappa N^{-3\alpha-\kappa}\), and the remaining convolution factor satisfies \(|K_1(s-u,n)|\lesssim \langle n\rangle^{-\alpha}\).  The number of lattice points in the shell is \(O(N^{3+\epsilon_0})\).  This gives \eqref{eq:wc-offdiag-cubic-shell-kernel}.  The time-increment estimate follows by replacing one of the three time-dependent factors by its increment: Lemma~\ref{lem:kernel-cov-bounds} gives the increment bound for the two Duhamel kernels, and \eqref{eq:wc-offdiag-cov-increment} gives it for the covariance factor.  In each case the fixed-time shell bound acquires only \(|t-t'|^\eta N^{\alpha\eta}\), which proves \eqref{eq:wc-offdiag-cubic-shell-increment}.  The color-exchanged proof is identical.
\end{proof}

\begin{lemma}[Off-diagonal cubic kernel tails]\label{lem:wc-offdiag-cubic-kernel-tails}
Let \(K_{i,\mathrm{off},\Lambda}^{\mathsf R}\) be the integrated kernel of Lemma~\ref{lem:wc-offdiag-cubic-integrated-kernel}.  For every internal loop cutoff \(L\ge1\),
\begin{equation}\label{eq:wc-offdiag-cubic-internal-tail}
\sup_{u,t}|K_{i,\mathrm{off},>L}^{\mathsf R}(n;t,u)|
\lesssim_{T,\epsilon_0}
\|\mathsf R\|_\kappa\langle n\rangle^{-\alpha}
\sum_{\substack{N\gtrsim\langle n\rangle\\ N>L}}
N^{3-4\alpha-\kappa+\epsilon_0}.
\end{equation}
The analogous time-increment tail is obtained by replacing the power by
\(3-4\alpha-\kappa+\alpha\eta+\epsilon_0\).  Consequently, for every
\(s<5\alpha+\kappa-9/2\) and sufficiently small \(\eta>0\), the cutoff kernels form a Cauchy family in the integrated-kernel seminorm \(\mathfrak K_T^{s,\eta}\) of Subsection~\ref{subsec:C-first-chaos-integrated}.  The external tail, where the remaining first-chaos output frequency \(n\) is removed by the Galerkin cutoff, is summable by the coefficient bound \eqref{eq:wc-cubic-offdiag-coefficient} below.
\end{lemma}

\begin{proof}
The internal tail estimate is the dyadic sum of \eqref{eq:wc-offdiag-cubic-shell-kernel}; the increment version is the dyadic sum of \eqref{eq:wc-offdiag-cubic-shell-increment}.  Choose \(\eta\) and \(\epsilon_0\) so that
\[
        3-4\alpha-\kappa+\alpha\eta+\epsilon_0<0.
\]
Then fixed Fourier labels have vanishing internal tails.  In the weighted kernel seminorm, let \(N_0=\langle n\rangle\).  If \(N_0\le L\), the contribution of \eqref{eq:wc-offdiag-cubic-internal-tail} is bounded, up to an arbitrarily small dyadic loss, by
\[
        L^{6-8\alpha-2\kappa+}
        \sum_{N_0\le L} N_0^{2s+3-2\alpha+}
        \lesssim L^{-2\iota+}
\]
when \(s=5\alpha+\kappa-9/2-\iota\).  If \(N_0>L\), the full bound obtained by summing all internal shells gives
\[
        \sum_{N_0>L}N_0^{2s+3}N_0^{6-10\alpha-2\kappa+}\to0.
\]
The increment part has the same proof with the extra factor \(N^{2\alpha\eta}\), absorbed by reducing \(\eta\).  The external tail is exactly the tail of the output dyadic series controlled by the same full coefficient and increment majorants.
\end{proof}

\begin{lemma}[Off-diagonal cubic first-chaos branch]\label{lem:wc-cubic-offdiag}
The branches \(C_{i,\mathrm{off}}^{\mathsf R}\) converge almost surely in the same source-compatible topologies as the same-color first-chaos cubic contractions.  More precisely,
\begin{equation}\label{eq:wc-cubic-offdiag-bound}
        C_{i,\mathrm{off}}^{\mathsf R}
        \in C_T\mathcal C^{5\alpha-\frac92-}
        \cap L_T^1B_{2,\infty}^{\sigma-\alpha},
\end{equation}
with cutoff convergence in these topologies.  In fact the off-diagonal branch has the stronger coefficient estimate
\begin{equation}\label{eq:wc-cubic-offdiag-coefficient}
        \sup_{t\le T}\E|\widehat C_{i,\mathrm{off}}^{\mathsf R}(n,t)|^2
        \lesssim_{T,\epsilon_0,\delta_{\mathsf R}}
        \|\mathsf R\|_\kappa^2\,
        \langle n\rangle^{6-10\alpha-2\kappa+\epsilon_0},
\end{equation}
and the corresponding time-increment estimate has the additional factor \(|t-t'|^{2\eta}\langle n\rangle^{2\alpha\eta}\).
\end{lemma}

\begin{proof}
Choose \(\epsilon_0\) and \(\eta\) so that
\[
        3-4\alpha-\kappa+\alpha\eta+\epsilon_0<0.
\]
This is possible in the present range; in fact \(\alpha>3/4\) and \(\kappa\ge0\) already suffice after taking \(\eta,\epsilon_0\) small, while the theorem keeps the stronger global condition \(\kappa>3-3\alpha+10\eps\) for the quadratic zero-mode branch.  Summing \eqref{eq:wc-offdiag-cubic-shell-kernel} over dyadic \(N\gtrsim\langle n\rangle\) yields
\begin{equation}\label{eq:wc-offdiag-cubic-kernel-summed}
        \sup_{u,t}|K_{i,\mathrm{off}}^{\mathsf R}(n;t,u)|
        \lesssim
        \|\mathsf R\|_\kappa\langle n\rangle^{3-5\alpha-\kappa+\epsilon_0}.
\end{equation}
By Itô isometry in the remaining first-chaos variable,
\[
        \E|\widehat C_{i,\mathrm{off}}^{\mathsf R}(n,t)|^2
        =\int_0^t |K_{i,\mathrm{off}}^{\mathsf R}(n;t,u)|^2\,du,
\]
which gives \eqref{eq:wc-cubic-offdiag-coefficient}.  The increment estimate follows from \eqref{eq:wc-offdiag-cubic-shell-increment}, the same dyadic summation, and the boundary interval estimate in the Itô isometry.  The coefficient-to-path criterion of Lemma~\ref{lem:C-coeff-path} then gives
\[
        C_{i,\mathrm{off}}^{\mathsf R}\in C_T\mathcal C^{5\alpha+\kappa-9/2-},
\]
which is stronger than the stated \(C_T\mathcal C^{5\alpha-9/2-}\).  The source-space inclusion follows from the admissible inequalities exactly as for the same-color first-chaos cubic term.  Cutoff convergence is Lemma~\ref{lem:wc-offdiag-cubic-kernel-tails} combined with Lemma~\ref{lem:C-K21-conv}.  Fixed dyadic shells are finite sums and converge entrywise, while the internal loop tail and the external output tail are summable in the integrated-kernel seminorm.
\end{proof}

\begin{lemma}[Centered third chaos under weak covariance]\label{lem:wc-cubic-third-chaos}
The centered component \(\Gamma_{i,\mathsf R}^{(3)}\) converges almost surely in
\[
        C_T\mathcal C^{\frac92\alpha-\frac92-}
        \cap L_T^1H^{s_2-\alpha}\cap L_T^1B_{2,\infty}^{\sigma-\alpha}
\]
for every admissible \((s_2,\sigma)\).
\end{lemma}

\begin{proof}
At finite cutoff, Lemma~\ref{lem:wc-cubic-finite-wick} writes the third component as a third homogeneous Wick chaos in the finite joint two-color Gaussian space.  By the class-wise coordinate argument of Lemma~\ref{lem:wc-classwise-covariance-coordinates}, the covariance factorization may mix colors and time-kernel coordinates only inside a fixed real Fourier involution class.  It preserves the external Fourier relation \(a+b+r=n\) and the resonant support \(|a+b|\sim |r|\).  The local covariance matrices are uniformly bounded because their diagonal entries are one in complex coordinates, so the Hilbert-kernel norm of the third-chaos coefficient is bounded, up to a polynomial in \(\|\mathsf R\|_\kappa\), by the same expression as in Lemma~\ref{lem:C-third-kernel-norm}.  Using the first-Picard coefficient estimate from Appendix~\ref{app:first-picard-phase} and the stochastic-convolution estimate gives the same dyadic shell
\[
        M^3\cdot M^{3-7\alpha+}\cdot M^{-2\alpha}=M^{6-9\alpha+}.
\]
The time-increment estimate is the same with an additional factor \(M^{2\alpha\eta}\).  Lemma~\ref{lem:C-coeff-path} and the cutoff-tail argument of Lemma~\ref{lem:C-third-tail} then give the asserted topologies.  Degenerate covariance blocks are handled by working in the quotient Hilbert space of zero-norm directions.
\end{proof}

\begin{lemma}[Regular zero-mode cubic contribution]\label{lem:wc-regular-zero-mode-cubic}
The term
\[
        R_{i,\Theta}^{\mathsf R}=I_{3-i}(C_\Theta^{\mathsf R})\circ\Psi_i
\]
is controlled in every topology required of \(\Gamma_i\) in Definition~\ref{def:enhanced-data-norm}.  The cutoff approximations converge in the same source topologies.
\end{lemma}

\begin{proof}
The function \(C_\Theta^{\mathsf R}\) is a deterministic spatial zero mode by Proposition~\ref{prop:wc-theta-branch}.  Hence \(A_i(t):=I_{3-i}(C_\Theta^{\mathsf R})(t)\) is also a deterministic spatial zero mode with finite \(C^1_T\) norm controlled by \(\|\mathsf R\|_\kappa\) and \(\delta_{\mathsf R}^{-1}\).  In the fixed inhomogeneous Bony convention, a zero-mode factor has only finitely many low Littlewood--Paley blocks, so the resonant product \(A_i\circ\Psi_i\) involves only finitely many low blocks of \(\Psi_i\).  Consequently it is a smooth spatial random function.  If one rewrites this contribution as an ordinary multiplier term, the same conclusion follows from boundedness of the smooth multiplier \(A_i(t)\).  The cutoff convergence is inherited from Proposition~\ref{prop:wc-theta-branch} and from the low-frequency convergence of \(\Psi_{i,\Lambda}\).  Thus \(R_{i,\Theta}^{\mathsf R}\) is better than the cubic source threshold and belongs to both source spaces in Definition~\ref{def:enhanced-data-norm}.
\end{proof}

\begin{proposition}[Weak-covariance cubic symbol estimates]\label{prop:wc-cubic-symbols}
Under Definition~\ref{def:wc-covariance} and \(\kappa>3-3\alpha+10\eps\), the cubic symbols \(\Gamma_i^{\mathsf R}\) have the decomposition \eqref{eq:wc-cubic-decomp}.  All components exist as cutoff limits and satisfy the same source bounds as the independent-color cubic symbols in Proposition~\ref{prop:cubic-resonance}.  More precisely,
\[
        \Gamma_{i,\mathsf R}^{(3)}\in C_T\mathcal C^{\frac92\alpha-\frac92-},
        \qquad
        C_{i,\mathrm{diag}}^{\mathsf R},\ C_{i,\mathrm{off}}^{\mathsf R},\ R_{i,\Theta}^{\mathsf R}
        \in L_T^1H^{s_2-\alpha}\cap L_T^1B_{2,\infty}^{\sigma-\alpha},
\]
with the stated stronger H\"older regularity for the first-chaos branches.  Hence the enhanced-data norm of Definition~\ref{def:enhanced-data-norm} can be augmented by \(C_{i,\mathrm{off}}^{\mathsf R}\) and \(R_{i,\Theta}^{\mathsf R}\) without changing the deterministic fixed point.
\end{proposition}

\begin{proof}
The finite algebra is the degree-three $\mathsf R$-Wick identity in Lemma~\ref{lem:wc-cubic-finite-wick}.  Applied to the finite expansion of $V_{3-i}^{\mathsf R}\circ\Psi_i$, it gives exactly the pieces in \eqref{eq:wc-cubic-decomp}: the centered third chaos, the same-color first-chaos contraction, the off-diagonal first-chaos contraction, and the zero-mode branch generated by $C_\Theta^{\mathsf R}$.  This decomposition is made before any cutoff is removed.

The centered third-chaos term is estimated in Lemma~\ref{lem:wc-cubic-third-chaos}.  Its proof is the independent third-chaos coefficient estimate with the local covariance Hilbert metric inserted.  The positive semidefinite covariance matrices have normalized diagonal entries, so every finite Wick pairing is bounded by the corresponding diagonal variance product.  The same lattice phase sums and time-increment estimates as in Appendix~\ref{app:cubic-resonance} therefore give the regularity $C_T\mathcal C^{\frac92\alpha-\frac92-}$ and the cutoff tails.

The same-color first-chaos branch $C_{i,\mathrm{diag}}^{\mathsf R}$ is the branch already constructed in Subsection~\ref{subsec:C-first-chaos-integrated}.  The diagonal entries of the covariance are one, so its integrated Volterra kernel is unchanged, except for harmless local coordinate choices in the finite $\mathsf R$-Wick representation.

The off-diagonal first-chaos branch is controlled by Lemmas~\ref{lem:wc-offdiag-cubic-integrated-kernel}--\ref{lem:wc-offdiag-cubic-kernel-tails} and Lemma~\ref{lem:wc-cubic-offdiag}.  The key shell bound is
\[
        \langle n\rangle^{-\alpha}N^{3-4\alpha-\kappa+},
        \qquad N\gtrsim\langle n\rangle,
\]
for the integrated kernel after the remaining stochastic convolution is inserted.  The internal high scale is summable already for $\alpha>3/4$ and $\kappa\ge0$; the covariance decay improves the bound.  The stronger standing condition $\kappa>3-3\alpha+10\eps$ is inherited from the quadratic zero-mode branch and therefore also suffices here.

Finally, $R_{i,\Theta}^{\mathsf R}=I_{3-i}(C_\Theta^{\mathsf R})\circ\Psi_i$ is controlled by Lemma~\ref{lem:wc-regular-zero-mode-cubic}: the factor $I_{3-i}(C_\Theta^{\mathsf R})$ is a deterministic spatial zero mode with a convergent cutoff limit, so the resonant product involves only finitely many low Littlewood--Paley blocks of the stochastic convolution.  It is therefore smoother than the cubic source threshold.

Each component is obtained as a finite-cutoff identity followed by a Cauchy argument in the displayed topology.  Fixed dyadic blocks converge by finite-dimensional convergence of the cutoff kernels and covariance Wick polynomials; the tails are controlled by the summable envelopes in the cited lemmas.  The displayed components therefore belong to the stated path/source spaces and may be added to the enhanced-data norm without changing the deterministic fixed point.
\end{proof}

\subsection{Covariance-augmented enhanced-data set}

\begin{theorem}[Covariance-augmented stochastic enhanced data]\label{thm:wc-enhanced-data}
Assume Definition~\ref{def:wc-covariance} and \(\delta_{\mathsf R}>0\).  Fix the strict admissible parameters and \(T_0\le1\).  There exists an event \(\Omega_{\rm wc}(T_0)\) with \(\Prob(\Omega_{\rm wc}(T_0))=1\) on which the covariance-augmented enhanced datum is well defined:
\begin{equation}\label{eq:wc-enhanced-data-list}
        (\Psi_i,\Theta^{\mathsf R},C_\Theta^{\mathsf R},V_i^{\mathsf R},\partial_tV_i^{\mathsf R},
        \Gamma_{i,\mathsf R}^{(3)},C_{i,\mathrm{diag}}^{\mathsf R},C_{i,\mathrm{off}}^{\mathsf R},R_{i,\Theta}^{\mathsf R},
        \mathcal D_{\mathsf R}^{i;j,k},\mathcal B_{\mathsf R}^{i;j,k}).
\end{equation}
All components have cutoff limits in the same distribution, source, and operator topologies used by the deterministic map, and the covariance-augmented cutoff data satisfy
\[
        d_T^{\rm wc}(\Xi_\Lambda^{\mathsf R}(\omega),\Xi^{\mathsf R}(\omega))\to0,
        \qquad 0<T\le T_0,
        \quad \omega\in\Omega_{\rm wc}(T_0).
\]
The additional constants are controlled by \(\|\mathsf R\|_\kappa\) and by the gap \(\delta_{\mathsf R}\).  If \(\mathsf R_{12}=\mathsf R_{21}=0\), the data reduce to the independent-color enhanced datum.
\end{theorem}

\begin{proof}
The linear stochastic convolutions are defined from the same finite Fourier Brownian coordinates, now with covariance \(\mathsf R\).  Their one-color estimates are unchanged because \(\mathsf R_{ii}=1\).  The cross quadratic product is split at finite cutoff by \eqref{eq:wc-theta-finite-split}.  Proposition~\ref{prop:wc-theta-branch} gives the deterministic zero-mode branch with shell profile \(N^{3-3\alpha-\kappa+}\), and the condition \(\delta_{\mathsf R}>0\) gives summability and cutoff tails.  The centered part \(:\Psi_1\Psi_2:_{\mathsf R}\) is estimated by Lemma~\ref{lem:wc-first-picard-phase-sum}; the proof is the same deterministic phase-sum theorem as Appendix~\ref{app:first-picard-phase}, because the finite \(\mathsf R\)-Wick projection changes only local color-time covariance coordinates and keeps the external Fourier labels.  Proposition~\ref{prop:wc-quadratic-lift} therefore constructs \(\Theta^{\mathsf R}\), \(C_\Theta^{\mathsf R}\), \(V_i^{\mathsf R}\), and \(\partial_tV_i^{\mathsf R}\) both in the H\"older multiplier topology and in the dyadic \(B_{2,\infty}\) topology.

The mixed-operator part is Theorem~\ref{thm:wc-mixed-operators}.  Its lower-chaos branches are the covariance-weighted Volterra diagonals of Proposition~\ref{prop:wc-deterministic-diagonals}; for off-diagonal covariance blocks Lemma~\ref{lem:wc-offdiag-covariance} gives the extra factor \(N^{-\kappa}\).  The centered covariance remainders are reduced by Lemma~\ref{lem:wc-classwise-covariance-coordinates} and Proposition~\ref{prop:wc-covariance-local-tensor} to the same tensor incidence \(n=q+\ell+r\) as in the independent case, and are then estimated in operator topology by Proposition~\ref{prop:wc-centered-component}.

The cubic symbols are Proposition~\ref{prop:wc-cubic-symbols}.  The finite degree-three \(\mathsf R\)-Wick formula separates the centered third chaos, the same-color first-chaos branch, the off-diagonal first-chaos branch with covariance factor \(N^{-\kappa}\), and the regular zero-mode branch generated by \(C_\Theta^{\mathsf R}\).  The first-chaos terms are defined in the integrated-kernel topology of Proposition~\ref{prop:C-integrated-first-chaos-contraction}; the off-diagonal branch has the same integrated Volterra estimate with the additional covariance decay.  The zero-mode branch is smoother because it is a deterministic scalar covariance source acted on by the Klein--Gordon propagator.

Finally, Propositions~\ref{prop:wc-covariance-approximation} and~\ref{prop:wc-enhanced-distance} give continuity under covariance approximation, including rank-deficient limits.  The Wick identities hold at finite cutoff, and the uniform dyadic tails give convergence in the displayed path, source, and operator topologies.  When \(\mathsf R_{12}=\mathsf R_{21}=0\), \(C_\Theta^{\mathsf R}\), all off-diagonal covariance diagonals, and all off-diagonal cubic branches vanish, leaving the independent-color data.
\end{proof}

\begin{definition}[Covariance-augmented bounded size and local smallness]\label{def:wc-bounded-size}
For a weak-covariance enhanced datum \(\Xi^{\mathsf R}\), define \(M_T^{\rm wc}(\Xi^{\mathsf R})\) by starting from the independent enhanced-data norm \(M_T\) in Definition~\ref{def:enhanced-data-norm}, replacing
\[
        \Theta,V_i,\partial_tV_i,\Gamma_i,\mathcal D^{i;j},\mathcal B^{i;j,k}
\]
by the covariance-augmented components
\[
        \Theta^{\mathsf R},C_\Theta^{\mathsf R},V_i^{\mathsf R},\partial_tV_i^{\mathsf R},
        \Gamma_{i,\mathsf R}^{(3)},C_{i,\mathrm{diag}}^{\mathsf R},C_{i,\mathrm{off}}^{\mathsf R},R_{i,\Theta}^{\mathsf R},
        \mathcal D_{\mathsf R}^{i;j,k},\mathcal B_{\mathsf R}^{i;j,k}.
\]
The added components are measured in the distribution, source, and operator topologies stated in Theorem~\ref{thm:wc-enhanced-data}: \(C_\Theta^{\mathsf R}\) is measured in its deterministic zero-mode topology, the cubic pieces are measured in \(L_T^1H^{s_2-\alpha}\cap L_T^1B_{2,\infty}^{\sigma-\alpha}\), and the covariance diagonals and centered remainders are measured in the same operator norms \(\mathcal L_H\) and \(\mathcal L_B\) as in Definition~\ref{def:enhanced-data-distance}.  With deterministic Cauchy data \(\mathbf y\), set
\begin{equation}\label{eq:wc-bd-size}
        M_T^{\rm bd,wc}(\Xi^{\mathsf R},\mathbf y)
        :=1+M_T^{\rm wc}(\Xi^{\mathsf R})+\|\mathbf y\|_{\mathcal H_T^{s_2,\sigma}}.
\end{equation}
The weak localized-smallness gauge \(\mu_T^{\rm wc}(\Xi^{\mathsf R},\mathbf y)\) is the quantity \(\mu_T\) from Definition~\ref{def:localized-mixed-size}, with two replacements.  First, \(M_T^{\rm bd}\) is replaced by \(M_T^{\rm bd,wc}\).  Second, the diagonal contribution is the full covariance family
\[
        \sum_{(i;j,k)\in\mathfrak M}\mathcal D_{\mathsf R}^{i;j,k}
\]
in place of the independent same-color subfamily.  The centered contribution is
\[
        \sum_{(i;j,k)\in\mathfrak M}\mathcal B_{\mathsf R}^{i;j,k}
\]
in the direct \(L_T^1B_{2,\infty}^{\sigma-\alpha}\) source operator norm.

A family \(\mathcal K^{\rm wc}\) of covariance-augmented data and deterministic data is weak localized-small on \([0,T_0]\) if
\[
        \sup_{(\Xi^{\mathsf R},\mathbf y)\in\mathcal K^{\rm wc}}
        M_{T_0}^{\rm bd,wc}(\Xi^{\mathsf R},\mathbf y)<\infty,
        \qquad
        \nu_{\mathcal K}^{\rm wc}(T)
        :=\sup_{(\Xi^{\mathsf R},\mathbf y)\in\mathcal K^{\rm wc}}
        \mu_T^{\rm wc}(\Xi^{\mathsf R}|_{[0,T]},\mathbf y)\to0
\]
as \(T\downarrow0\).  This is Definition~\ref{def:localized-smallness-class} applied to the enlarged weak-covariance enhanced datum.  Proposition~\ref{prop:wc-diagonal-local-smallness} gives the smallness of the enlarged diagonal family, and Proposition~\ref{prop:centered-besov-local-smallness} gives the centered Besov-source smallness.
\end{definition}

\paragraph{Covariance-augmented deterministic norm.}
After the components in Theorem~\ref{thm:wc-enhanced-data} have been constructed, the deterministic estimates in Section~\ref{sec:deterministic-closure} apply with the same algebraic source decomposition and the same solution spaces.  The covariance matrix enters only through the norms of the augmented enhanced components: \(M_T^{\rm bd}\) is replaced by \(M_T^{\rm bd,wc}\), and the diagonal subfamily in \(\mu_T^{\rm op}\) is enlarged from the same-color branches to the full family \(\mathcal D_{\mathsf R}^{i;j,k}\).

\subsection{Enhanced-data distance and Galerkin bridge}
The augmented weak-covariance enhanced-data distance is obtained from Definition~\ref{def:enhanced-data-distance} by replacing $\Theta,V_i,\partial_tV_i,$ and $\Gamma_i$ with their covariance-split counterparts and adding the difference norms of
\begin{equation}\label{eq:wc-distance-added-components}
        C_\Theta^{\mathsf R},\qquad
        C_{i,\mathrm{off}}^{\mathsf R},\qquad
        R_{i,\Theta}^{\mathsf R},\qquad
        \mathcal D_{\mathsf R}^{i;j,k},\qquad
        \mathcal B_{\mathsf R}^{i;j,k}.
\end{equation}
For the diagonal part of the localized-smallness norm one uses all covariance-weighted Volterra diagonals \(\mathcal D_{\mathsf R}^{i;j,k}\), in place of the independent same-color subfamily.  This is the only change in the deterministic control norm, and Proposition~\ref{prop:wc-diagonal-local-smallness} supplies the required small-time factor.

\begin{proposition}[Weak-covariance enhanced-data distance]\label{prop:wc-enhanced-distance}
Let \(\mathbb E^{\mathsf R}\) and \(\widetilde{\mathbb E}^{\widetilde{\mathsf R}}\) be two covariance-augmented enhanced data satisfying Definition~\ref{def:wc-covariance} with the same exponent \(\kappa\).  Define \(d_T^{\rm wc}\) by adding the component differences in \eqref{eq:wc-distance-added-components} to the distance \(d_T\) of Definition~\ref{def:enhanced-data-distance}, and by measuring the covariance-diagonal and centered-operator differences in the same operator norms \(\mathcal L_H\) and \(\mathcal L_B\).  The corresponding bounded-size control is \(M_T^{\rm bd,wc}\) from Definition~\ref{def:wc-bounded-size}.  If \(d_\kappa(\mathsf R,\widetilde{\mathsf R})\to0\) and the finite-cutoff enhanced objects converge in the component topologies above, then \(d_T^{\rm wc}\to0\) on every fixed interval on which the corresponding data have bounded size and a common localized-smallness modulus.
\end{proposition}

\begin{proof}
The scalar and first-chaos components use the shell and tail estimates in Propositions~\ref{prop:wc-theta-branch}, \ref{prop:wc-quadratic-lift}, and~\ref{prop:wc-cubic-symbols}.  The mixed-operator components use Proposition~\ref{prop:wc-covariance-continuity}.  The finite-dimensional covariance approximation, including rank-deficient limits, is Proposition~\ref{prop:wc-covariance-approximation}.  The deterministic data component is unchanged.  Thus each component of \(d_T^{\rm wc}\) either converges by a finite-dimensional covariance argument on finitely many dyadic blocks or is controlled by a uniform summable tail.
\end{proof}

\begin{proposition}[Finite weak-covariance Galerkin bridge]\label{prop:wc-galerkin-bridge}
At finite cutoff, the correlated Galerkin equation with covariance \eqref{eq:wc-covariance-def-main} is algebraically equivalent to the finite covariance-augmented enhanced \(X\)--\(Y\) system.  The reconstruction remains
\begin{equation}
        u_{i,\Lambda}=\Psi_{i,\Lambda}+V_{i,\Lambda}^{\mathsf R}+X_{i,\Lambda}+Y_{i,\Lambda}.
\end{equation}
Every product of stochastic convolution factors is read through the finite \(\mathsf R\)-Wick identity.  In particular, \(C_{\Theta,\Lambda}^{\mathsf R}\) enters
\(V_{i,\Lambda}^{\mathsf R}=I_i(:\Psi_{1,\Lambda}\Psi_{2,\Lambda}:_{\mathsf R}+C_{\Theta,\Lambda}^{\mathsf R})\)
as a deterministic lower-chaos forcing term.  The covariance-weighted mixed diagonals and the off-diagonal cubic first-chaos branches are likewise retained in the finite enhanced data.  Summing the component equations reconstructs the correlated Galerkin PDE.
\end{proposition}

\begin{proof}
Before removing the cutoff all Fourier sums are finite and all Gaussian products live in a finite-dimensional joint two-color Gaussian space.  The identity
\begin{equation}
        \Psi_{a,\Lambda}\Psi_{b,\Lambda}
        =:\Psi_{a,\Lambda}\Psi_{b,\Lambda}:_{\mathsf R}
        +C_{ab,\Lambda}^{\mathsf R}
\end{equation}
is an exact Wick identity, where \(C_{ab,\Lambda}^{\mathsf R}\) denotes the finite covariance pairing.  It is a finite algebraic identity in the cutoff Gaussian space.

Substitute \(u_{i,\Lambda}=\Psi_{i,\Lambda}+V_{i,\Lambda}^{\mathsf R}+X_{i,\Lambda}+Y_{i,\Lambda}\) into the finite mild equation and apply the fixed Bony decomposition.  The equation for \(\Psi_{i,\Lambda}\) removes the projected correlated noise.  The equation defining \(V_{i,\Lambda}^{\mathsf R}=I_i(:\Psi_{1,\Lambda}\Psi_{2,\Lambda}:_{\mathsf R}+C_{\Theta,\Lambda}^{\mathsf R})\) removes the full finite quadratic stochastic product, including its deterministic zero-mode covariance branch.  The low--high sources define \(X_{i,\Lambda}\), and the resonant terms define the finite mixed operators after their covariance diagonals and centered remainders are separated.  Thus the resulting finite enhanced equations are exactly the covariance-augmented \(X\)--\(Y\) system.  Conversely, adding back the equations for \(\Psi_{i,\Lambda}\), \(V_{i,\Lambda}^{\mathsf R}\), \(X_{i,\Lambda}\), and \(Y_{i,\Lambda}\) cancels the finite \(\mathsf R\)-Wick expansion term by term and recovers the unexpanded finite correlated Galerkin PDE.
\end{proof}

\begin{lemma}[Weak-covariance cutoff data are localized-small]\label{lem:wc-cutoff-localized-small}
Let \(\Xi^{\mathsf R}\) be the covariance-augmented compatible enhanced datum constructed in Theorem~\ref{thm:wc-enhanced-data}, and let \(\Xi^{\mathsf R}_\Lambda\) be the corresponding Galerkin-compatible cutoff data.  Fix deterministic Cauchy data \(\mathbf y\) with finite norm \(\mathcal H_T^{s_2,\sigma}\).  After discarding finitely many cutoff levels, the family
\[
        \{(\Xi^{\mathsf R},\mathbf y)\}
        \cup \{(\Xi^{\mathsf R}_\Lambda,\mathbf y):\Lambda\ge\Lambda_0\}
\]
is weak localized-small on a sufficiently short interval in the sense of Definition~\ref{def:wc-bounded-size}.
\end{lemma}

\begin{proof}
Componentwise cutoff convergence in \(d_T^{\rm wc}\), proved by Theorem~\ref{thm:wc-enhanced-data}, Proposition~\ref{prop:wc-enhanced-distance}, and the mixed-operator convergence in Theorem~\ref{thm:wc-mixed-operators}, gives a uniform bound for \(M_T^{\rm bd,wc}\) on the limit together with all sufficiently large cutoffs.  The ordinary deterministic source factors carry the explicit \(T^\delta\) smallness of Definition~\ref{def:localized-mixed-size}.  Proposition~\ref{prop:wc-diagonal-local-smallness} gives local smallness of the covariance diagonals after \(L_T^\infty\to L_T^1\), and Proposition~\ref{prop:centered-besov-local-smallness} gives the same property for the centered covariance remainders.  Cutoff differences are controlled in these source topologies by Proposition~\ref{prop:wc-enhanced-distance}.  Hence the bounded-size and localized-smallness conditions of Definition~\ref{def:wc-bounded-size} hold on a common short interval.
\end{proof}

\begin{corollary}[Weak-covariance cutoff stability]\label{cor:wc-cutoff-stability}
On the full-probability event of Theorem~\ref{thm:wc-enhanced-data}, the covariance-augmented Galerkin enhanced data converge in \(d_T^{\rm wc}\) and form a common localized-smallness class after restricting to a sufficiently small common interval.  Hence the correlated Galerkin solutions converge pathwise to the fixed point constructed from the limiting covariance-augmented datum.
\end{corollary}

\begin{proof}
Fixed dyadic blocks converge entrywise because the cutoff expressions are finite Fourier sums.  The tails are controlled by the summable envelopes for \(C_\Theta^{\mathsf R}\), the covariance diagonals, the centered operator remainders, and the cubic branches, so the covariance-augmented data converge in \(d_T^{\rm wc}\).  Lemma~\ref{lem:wc-cutoff-localized-small} supplies the common weak localized-smallness interval after discarding finitely many cutoff levels.  Corollary~\ref{cor:deterministic-cutoff-passage} then applies to the covariance-augmented data.
\end{proof}

\begin{proof}[Proof of Theorem~\ref{thm:weak-covariance-main}]
Theorem~\ref{thm:wc-enhanced-data} constructs $\Xi^{\mathsf R}$ on a full-probability event and proves the convergence \eqref{eq:wc-main-pathwise-dT}.  Definition~\ref{def:wc-bounded-size} and Lemma~\ref{lem:wc-cutoff-localized-small} place the limiting datum and all sufficiently large cutoffs in a common localized-smallness class.  The deterministic solution theorem, applied with the covariance-augmented norms, therefore gives the fixed point $(X_i^{\mathsf R},Y_i^{\mathsf R})$ and the reconstruction \eqref{eq:wc-main-solution-formula}.  Proposition~\ref{prop:wc-galerkin-bridge} identifies the finite enhanced equations with the correlated Galerkin equations, and Corollary~\ref{cor:wc-cutoff-stability} gives pathwise convergence of the fields and sources.  The Lipschitz estimate in Proposition~\ref{prop:uniform-deterministic-bridge} gives independence of compatible enhanced coordinates and cofinal cutoff sequences on every common interval.  When the off-diagonal covariance vanishes, each additional covariance branch is zero, so the construction reduces to the independent-color case.
\end{proof}

\end{document}